\documentclass[12pt]{amsart}      
\usepackage{amssymb}
\usepackage{eucal}
\usepackage{amsmath}
\usepackage{amscd}
\usepackage{xcolor}
\usepackage{fancybox}
\usepackage{colordvi}
\usepackage{multicol}
\usepackage{colordvi}

\usepackage{ifpdf}
\ifpdf
\usepackage[colorlinks,final,backref=page,hyperindex]{hyperref}
\else
\usepackage[colorlinks,final,backref=page,hyperindex]{hyperref}
\fi
\usepackage[active]{srcltx} 

\newtheorem{theorem}{Theorem}[section]
\newtheorem{prop}[theorem]{Proposition}
\newtheorem{proposition}[theorem]{Proposition}
\newtheorem{lemma}[theorem]{Lemma}

\newtheorem{coro}[theorem]{Corollary}

\newtheorem{conjecture}[theorem]{Conjecture}

\newcommand{\proofbegin}{\begin{proof}} 
\newcommand{\proofend}{\end{proof}} 
\newcommand{\epf}{\proofend}

\newtheorem{temprmk}[theorem]{Remark}
\newtheorem{tempdef}[theorem]{Definition}
\newtheorem{tempex}[theorem]{Example}

\newenvironment{definition}{\begin{tempdef}\em}{\end{tempdef}}

\font\cyr=wncyr10
\newcommand{\nc}{\newcommand}

\renewcommand{\Bbb}{\mathbb}
\renewcommand{\frak}{\mathfrak}

\nc{\bear}{\begin{array}}
\nc{\enar}{\end{array}}
\nc{\rar}{\rightarrow}
\nc{\la}{\longrightarrow}
\nc{\dap}[1]{\downarrow \rlap{$\scriptstyle{#1}$}}
\nc{\ola}[1]{\stackrel{#1}{\la}}
\nc{\defeq}{\stackrel{\rm def}{=}}
\nc{\mrm}[1]{{\rm #1}}
\nc{\dirlim}{\displaystyle{\lim_{\longrightarrow}}\,}
\nc{\invlim}{\displaystyle{\lim_{\longleftarrow}}\,}
\nc{\sha}{\mbox{\cyr X}}    
\newcommand{\shom}{{\mathop{{\mathcal H}om}\nolimits}}
\nc{\vep}{\varepsilon}
\nc{\sig}{{^\Sigma}}
\nc{\rlangle}{{(}}
\nc{\rrangle}{{)}}
\nc{\wc}{?}

\nc{\ab}{{\mathop{\rm ab}\nolimits}}
\nc{\ad}{{\mathop{\rm ad}\nolimits}}
\nc{\alg}{{\mathop{\rm alg}\nolimits}}
\nc{\an}{{\mathop{\rm an}\nolimits}}

\DeclareMathOperator{\Aut}{{\mathop{\rm Aut}\nolimits}}
\DeclareMathOperator{\aut}{\Aut}
\nc{\bt}{{\mathop{\rm B}\nolimits}}
\nc{\BK}{{\mathop{\rm BK}\nolimits}}
\nc{\can}{{\mathop{\rm can}\nolimits}}
\newcommand{\cris}{{\mbox{\scriptsize crys}}}
\nc{\Crys}{\cris}
\nc{\crys}{\cris}
\nc{\ecrys}{{\ell\mbox{\scriptsize -crys}}}
\nc{\lcrys}{{\lambda\mbox{\scriptsize -crys}}}
\nc{\pcrys}{{p\mbox{\scriptsize -crys}}}
\nc{\fpcrys}{{\frakp\mbox{\scriptsize -crys}}}
\renewcommand{\det}{{\mathop{\rm det}\nolimits}}
\DeclareMathOperator{\Disc}{{\mathop{\rm Disc}\nolimits}}
\nc{\dr}{{\mathop{\rm dR}\nolimits}}
\DeclareMathOperator{\End}{{\mathop{\rm End}\nolimits}}
\nc{\et}{{\mathop{\rm et}\nolimits}}
\DeclareMathOperator{\Ext}{{\mathop{\rm Ext}\nolimits}}
\nc{\ext}{\Ext}
\DeclareMathOperator{\fil}{{\mathop{\rm Fil}\nolimits}}
\nc{\Fil}{\fil}
\DeclareMathOperator{\fitt}{Fitt}
\DeclareMathOperator{\Fr}{{\mathop{\rm Fr}\nolimits}}
\DeclareMathOperator{\frob}{{\mathop{\rm Frob}\nolimits}}
\nc{\Frob}{\frob}
\DeclareMathOperator{\gal}{{\mathop{\rm Gal}\nolimits}}
\nc{\Gal}{\gal}
\DeclareMathOperator{\GL}{{\mathop{\rm GL}\nolimits}}
\nc{\gr}{{\mathop{\rm gr}\nolimits}}
\renewcommand{\hom}{{\mathop{\rm Hom}\nolimits}}
\nc{\Hom}{\hom}
\nc{\id}{{\mathop{\rm id}\nolimits}}
\renewcommand{\Im}{{\mathop{\rm Im}\nolimits}}
\DeclareMathOperator{\im}{{\mathop{\rm im}\nolimits}}
\DeclareMathOperator{\Ind}{Ind}
\DeclareMathOperator{\lcm}{{\mathop{\rm lcm}\nolimits}}
\DeclareMathOperator{\length}{{\mathop{\rm length}\nolimits}}
\renewcommand{\max}{{\mathop{\rm max}\nolimits}}
\renewcommand{\min}{{\mathop{\rm min}\nolimits}}

\nc{\naive}{{\mathop{\rm nv}\nolimits}}
\DeclareMathOperator{\ord}{{\mathop{\rm ord}\nolimits}}
\nc{\pst}{{\mathop{\rm pst}\nolimits}}
\DeclareMathOperator{\rank}{{\mathop{\rm rank}\nolimits}}

\nc{\reg}{{\mathop{\rm reg}\nolimits}}

\nc{\res}{{\mathop{\rm res}\nolimits}}
\DeclareMathOperator{\SO}{{\mathop{\rm SO}\nolimits}}
\DeclareMathOperator{\Spec}{{\mathop{\rm Spec}\nolimits}}
\nc{\spec}{\Spec}
\renewcommand{\ss}{{\mathop{\rm ss}\nolimits}}
\nc{\st}{{\rm st}}
\DeclareMathOperator{\Sym}{{\mathop{\rm Sym}\nolimits}}
\nc{\sym}{\Sym}
\nc{\Symm}{\Sym}
\nc{\symm}{\Sym}
\DeclareMathOperator{\SL}{{\mathop{\rm SL}\nolimits}}
\DeclareMathOperator{\Tam}{{\mathop{\rm Tam}\nolimits}}
\nc{\tf}{{\mathop{\rm tf}\nolimits}}
\renewcommand{\th}{{\mathop{\rm th}\nolimits}}
\nc{\tor}{{\mathop{\rm tor}\nolimits}}
\DeclareMathOperator{\Tr}{{\mathop{\rm Tr}\nolimits}}
\DeclareMathOperator{\tr}{{\mathop{\rm tr}\nolimits}}
\nc{\un}{{\rm un}}
\nc{\ur}{{\rm ur}}

\nc{\f}{{\bf f}}
\nc{\ff}{\f}
\nc{\BR}{{\bf R}}
\nc{\BH}{{\bf H}}
\nc{\II}{{\bf I}}
\nc{\m}{{\bf m}}
\DeclareMathOperator{\pms}{{\bf PM}}

\nc{\spm}{{\bf SPM}}
\nc{\dfootnote}[1]{{}}          
\nc{\ffootnote}[1]{\dfootnote{#1}}
\nc{\mfootnote}[1]{\dfootnote{#1}}
\nc{\mlabel}[1]{\label{#1}}  

\renewcommand{\AA}{{\Bbb A}}
\nc{\BA}{{\Bbb A}}
\nc{\CC}{{\Bbb C}}
\nc{\DD}{{\Bbb D}}
\nc{\FF}{{\Bbb F}}
\nc{\GG}{{\Bbb G}}
\nc{\HH}{{\Bbb H}}
\nc{\NN}{{\Bbb N}}
\nc{\PP}{{\Bbb P}}
\nc{\QQ}{{\Bbb Q}}
\nc{\RR}{{\Bbb R}}
\nc{\TT}{{\Bbb T}}
\nc{\VV}{{\Bbb V}}
\nc{\ZZ}{{\Bbb Z}}

\nc{\cal}{\mathcal}
\nc{\cala}{{\mathcal A}}
\nc{\CA}{\cala}
\newcommand{\calb}{{\mathcal B}}
\newcommand{\FC}{{\mathcal C}}
\nc{\calc}{{\mathcal C}}
\nc{\cald}{{\mathcal D}}
\nc{\cale}{{\mathcal E}}
\nc{\CE}{\cale}
\nc{\calf}{{\mathcal F}}

\nc{\CF}{\calf}
\nc{\calg}{{\mathcal G}}
\nc{\CG}{\calg}
\nc{\CH}{{\mathcal H}}
\nc{\calh}{{\mathcal H}}
\nc{\call}{{\mathcal L}}
\newcommand{\CL}{{\mathcal L}}
\nc{\calm}{{\mathcal M}}
\nc{\CM}{\calm}
\newcommand{\CN}{{\mathcal N}}
\nc{\calo}{{\mathcal O}}
\nc{\CO}{\calo}
\nc{\calos}{{\calo_{S}}}
\nc{\CP}{{\mathcal P}}
\nc{\calp}{{\mathcal P}}
\newcommand{\GP}{{\mathcal P}}
\nc{\ipms}{{\mathcal PM}}
\nc{\icpms}{{\mathcal CPM}}
\nc{\tpms}{{\mathcal TPM}}
\nc{\cals}{{\mathcal S}}
\nc{\CS}{{\mathcal S}}
\nc{\CT}{{\mathcal T}}
\nc{\calt}{{\mathcal T}}
\nc{\calu}{{\mathcal U}}
\nc{\CV}{{\mathcal V}}
\nc{\calv}{{\mathcal V}}
\nc{\calw}{{\mathcal W}}
\nc{\calx}{{\mathcal X}}
\nc{\CX}{{\mathcal X}}

\nc{\frakA}{{\frak A}}
\nc{\fraka}{{\frak a}}
\nc{\ga}{\fraka}
\nc{\frakB}{{\frak B}}
\nc{\frakb}{{\frak b}}
\nc{\gb}{{\frak b}}
\nc{\gc}{{\frak c}}
\newcommand{\uhp}{{\frak H}}
\newcommand{\gm}{{\frak{m}}}
\nc{\frakm}{{\frak m}}

\nc{\frakp}{{\frak p}}
\newcommand{\gp}{{\frak{p}}}

\newcommand{\flcat}{\mathcal{F}\mathcal{R}}
\newcommand{\mfcat}{{\MF}}
\newcommand{\MF}{{\mathcal{M}\mathcal{F}}}

\newcommand{\Qbar}{\bar{\QQ}}
\newcommand{\Qpbar}{\bar{\QQ}_\pee}
\newcommand{\Qlbar}{\bar{\QQ}_\ell}

\newcommand{\barrho}{\bar{\rho}}

\newcommand{\pee}{p}

\newcommand{\dprod}{\displaystyle\prod}
\newcommand{\dsum}{\displaystyle\sum}

\newcommand{\mat}[4]
  {\left( \begin{array}{cc} {#1} & {#2} \\ {#3} & {#4} \end{array} \right)}
\renewcommand{\vec}[2]
  {\left( \begin{array}{c} {#1} \\ {#2} \end{array} \right)}
\newcommand{\smat}[4]
  {{\mbox{\scriptsize $\mat{{#1}}{{#2}}{{#3}}{{#4}}$}}}
 
\newcommand{\svec}[2]
  {{\mbox{\scriptsize $\vec{{#1}}{{#2}}$}}}

\newcommand{\rep}{{\rm Rep}_{cris}}
\newcommand{\tw}{{\rm tw}}
\newcommand{\kd}{D}

\newcommand{\mfcato}{\text{$\CO_\lambda$-$\mfcat$}}
\newcommand{\mfcatk}{\text{$K_\lambda$-MF}}
\newcommand{\mfcatkap}{\text{$\kappa$-$\mfcat$}}
\newcommand{\MFcat}{\text{MF}}

\begin{document}

\title[Adjoint motives and the Tamagawa
number conjecture]{Adjoint motives of modular forms and the Tamagawa
number conjecture}

\author{Fred Diamond}
\address{Department of Mathematics, 
King's College London, 
London, WC2R 2LS, UK}
\email{fdiamond@math.brandeis.edu}

\author{Matthias Flach}
\address{Department of Mathematics,
California Institute of Technology,
Pasadena, CA 91125, USA}
\email{flach@its.caltech.edu}

\author{Li Guo}
\address{
Department of Mathematics and Computer Science,
University of Rutgers at Newark,
Newark, NJ 07102, USA}
\email{liguo@rutgers.edu}

\date{August 2001, with minor update December 2025. This is  the longer version of the paper: The Tamagawa number conjecture of adjoint motives of modular forms, Ann. Sci. École Norm. Sup. (4) 37 (2004), no. 5, 663-727. This version contains details in the foundational part (Section 1) that were not included in the published version}



\begin{abstract} Let $f$ be a newform of weight $k\geq 2$, level $N$ with
coefficients in a number field $K$, and $A$ the adjoint
motive of the motive $M$ associated to $f$. We carefully discuss the
construction of the realisations of $M$ and $A$, as well as natural
integral structures in these realisations. We then use the method of
Taylor and Wiles to verify the $\lambda$-part of the Tamagawa number
conjecture of Bloch and Kato for $L(A,0)$ and
$L(A,1)$. Here $\lambda$ is any prime of $K$ not dividing $Nk!$, and so
that the mod $\lambda$
representation associated to $f$ is absolutely irreducible
when restricted to the Galois group over $\QQ(\sqrt{(-1)^{(\ell-1)/2}\ell})$
where $\lambda\mid \ell$.  The method also establishes modularity of
all lifts of the mod $\lambda$ representation which are crystalline
of Hodge-Tate type $(0,k-1)$.
\end{abstract}

\maketitle

\tableofcontents

\setcounter{section}{-1}

\section{Introduction}
This paper concerns the Tamagawa number conjecture of
Bloch and Kato~\cite{bloch_kato} for adjoint motives of modular
forms of weight $k\ge 2$.  The conjecture relates the value at $0$
of the associated $L$-function to arithmetic invariants of the
motive.  We prove that it holds up to powers of certain
``bad primes.''  The strategy for achieving this is essentially due
to Wiles~\cite{wiles}, as completed with Taylor in~\cite{tw}.
The Taylor-Wiles construction
yields a formula relating the size of a certain module measuring
congruences between modular forms to that of a certain Galois
cohomology group.  This was carried out in~\cite{wiles} and
\cite{tw} in the context of modular forms
of weight~$2$, where it was used to prove results in the direction
of the Fontaine-Mazur conjecture~\cite{fontaine_mazur}.
While it was no surprise that the method could be generalized to
higher weight modular forms and that the resulting formula would
be related to the Bloch-Kato conjecture, there remained many technical
details to verify in order to accomplish this.  In particular, the
very formulation of the conjecture relies on a comparison isomorphism
between the $\ell$-adic and de Rham realizations of the motive
provided by theorems of Faltings~\cite{faltings} or Tsuji~\cite{tsuji},
and verification of the conjecture requires the careful application
of such a theorem.  We also need to generalize results on
congruences between modular forms to higher weight, and to compute
certain local Tamagawa numbers.

\subsection{Some history}

Special values of $L$-functions have long played an important
role in number theory.  The underlying principle is that the
values of $L$-functions at integers reflect arithmetic properties
of the object used to define them.
A prime example of this is Dirichlet's class number formula;
another is the Birch and Swinnerton-Dyer conjecture.
The Tamagawa number conjecture of Bloch and Kato~\cite{bloch_kato},
refined by Fontaine, Kato and Perrin-Riou~\cite{Ka1,fon_pr,fon_bour}, is a vast generalization of these.
Roughly speaking, they predict the precise value
of the first non-vanishing derivative of the $L$-function
at zero (hence any integer) for every motive over $\QQ$.
This was already done up to a rational multiple by
conjectures of Deligne and Beilinson; the additional
precision of the Bloch-Kato conjecture can be thought
of as a generalized class number formula, where
ideal class groups are replaced by groups defined
using Galois cohomology.

Dirichlet's class number formula amounts
to the conjecture for the Dedekind zeta function for
a number field at $s = 0$ or $1$.  The conjecture is
also known up to a power of $2$ for
Dirichlet $L$-functions (including the Riemann zeta function)
at any integer (\cite{mazur_wiles}, \cite{bloch_kato}, \cite{huber_kings}).
It is known up to an explicit set of bad primes for the
$L$-function of a CM elliptic curve at $s=1$ if
the order of vanishing is $\le 1$ (\cite{coates_wiles}, \cite{rubin0},
\cite{kolyvagin}). There are also partial results for $L$-functions
of other modular forms at the central
critical value (\cite{G-Z}, \cite{Kappt},
\cite{kol_log}, \cite{nekovar}, \cite{zhang}) and
for values of certain Hecke $L$-functions (\cite{Hn}, \cite{Gu2},
\cite{kings}).

Here we consider the adjoint $L$-function of a modular form of weight
$k\ge 2$ at $s=0$ and $1$.
Special values of the $L$-function associated to the
adjoint of a modular form, and more generally, twists
of its symmetric square, have been studied by many
mathematicians.  A method of Rankin relates the values to
those of Petersson inner products, and this was used by
Ogg~\cite{Ogg},
Shimura~\cite{Sh2}, Sturm~\cite{Stu1,Stu2}, Coates
and Schmidt~\cite{C-S,schmidt} to obtain nonvanishing results
and rationality results along the lines of Deligne's conjecture.
Hida~\cite{Hida81} related the precise value to a number measuring
congruences between modular forms.  In the case of forms
corresponding to (modular) elliptic curves, results relating
the value to certain Galois cohomology groups (Selmer groups)
were obtained by Coates and Schmidt in the context of
Iwasawa theory, and by one of the authors, who in~\cite{Fl1} obtained
results in the direction of the Bloch-Kato conjecture.

A key point of Wiles' paper~\cite{wiles} is that for many elliptic
curves, modularity could be deduced from a formula relating congruences
and Galois cohomology~\cite{wiles}.  This formula could be regarded as
a primitive form of the Bloch-Kato conjecture for the adjoint motive of
a modular form.  His attempt to prove it using the Euler system
method introduced in~\cite{Fl1} was not successful except in the
CM case using generalizations of results in~\cite{Gu1} and~\cite{rubin1}.
Wiles, in work completed with Taylor~\cite{tw}, eventually
proved his formula using a new construction which could be
viewed as a kind of ``horizontal Iwasawa theory.''

In this paper, we refine the method of \cite{wiles} and
\cite{tw}, generalize it to higher weight modular forms and
relate the result to the Bloch-Kato conjecture.  Ultimately,
we prove the conjecture for the adjoint of an arbitrary newform
of weight $k \ge 2$ up to an explicit finite set of bad
primes.  We should stress the importance of making this
set as small and explicit as possible; indeed the refinements
in~\cite{fd_annals}, \cite{cdt} and \cite{bcdt} which
completed the proof of the Shimura-Taniyama-Weil conjecture
can be viewed as work in this direction for weight two
modular forms.  In this paper, we make use of some of the techniques
introduced in \cite{fd_annals} and \cite{cdt}, as well as the
modification of Taylor-Wiles construction in \cite{fd_twc} and
\cite{fujiwara}.
One should be able to improve our results using current technology
in the weight two case (using  \cite{cdt}, \cite{bcdt} and \cite{savitt}),
and in the ordinary case (using \cite{dickinson}, \cite{skinner_wiles});
one just has to relate the results in those papers to the
Bloch-Kato conjecture.  Finally we remark that Wiles' method
has been generalized to the setting of Hilbert modular forms
by Fujiwara \cite{fujiwara} and Skinner-Wiles \cite{skinner_wiles2},
but it seems much harder to extract results on special values from
their work.

\subsection{The framework}
The Bloch-Kato conjecture is formulated in terms
of ``motivic structures,'' a term referring to
the usual collection of cohomological data associated
to a motive.  This data consists of:
\begin{itemize}
\item vector spaces $M_\wc$, called realizations,
for $\wc=\bt$, $\dr$ and $\ell$ for each rational prime
$\ell$, each with extra structure (involution, filtration
or Galois action);
\item
comparison isomorphisms relating the realizations;
\item
a weight filtration.
\end{itemize}

Suppose that $f$ is a newform of weight $k\geq 2$ and
level $N$.  Much of the paper is devoted to the construction
of the motivic structure $A_f$ for which we prove the conjecture.
This construction is not new; it is due for the most
part to Eichler, Shimura, Deligne, Jannsen, Scholl and Faltings
(\cite{shimura}, \cite{del_bour}, \cite{jannsen},
\cite{scholl}, \cite{faltings})  We
review it however in order to collect the facts we need, provide
proofs we could not find in the literature, and set things up in a way
suited to the formulation of the Bloch-Kato conjecture.

Let us briefly recall here how the construction works.
We start with the modular curve $X_N$ parametrizing
elliptic curves with level $N$ structure.  Then one
takes the Betti, de Rham and $\ell$-adic cohomology
of $X_N$ with coefficients in a sheaf defined as the
$(k-2)$-nd symmetric power of the relative cohomology
of the universal elliptic curve over $X_N$.  These
come with the additional structure and comparison
isomorphisms needed to define a motivic structure
$M_{N,k}$, the comparison between $\ell$-adic and
de Rham cohomology being provided by a theorem of
Faltings \cite{faltings}.  The structures $M_{N,k}$
can also be defined as in \cite{scholl} using Kuga-Sato varieties;
this has the advantage of showing they arise from
``motives'' and provides the option of applying Tsuji's
comparison theorem \cite{tsuji}.  However the construction using
``coefficient sheaves'' is better suited to defining and comparing
lattices in the realizations which play a key role in the proof.

The structures $M_{N,k}$ also come with an action of the Hecke
operators and a perfect pairing.  The Hecke action is used to
``cut out'' a piece $M_f$, which corresponds to the newform $f$
and has rank two over the field generated by the coefficients
of $f$.  The pairing comes from Poincar\'e duality,
is related to the Petersson inner product and restricts
to a perfect pairing on $M_f$.  We finally take the trace
zero endomorphisms of $M_f$ to obtain the motivic structure
$A_f=\ad^0 M_f$.  The construction also yields integral
structures $\CM_f$ and $\CA_f$, consisting of lattices
in the various realizations and integral comparison
isomorphisms outside a set of bad primes.

Our presentation of the Bloch-Kato conjecture is much
influenced by its reformulation and generalization
due to Fontaine and Perrin-Riou.  Their version assumes
the existence of a category of motives with conjectural
properties.
Without assuming conjectures however, they define a
category $\spm_{\QQ}(\QQ)$
of premotivic structures whose objects
consist of realizations with additional structure
and comparison isomorphisms.  The category of mixed motives
is supposed to admit a fully faithful functor to it,
and a motivic structure is an object of the essential
image.
Their version of the Bloch-Kato
conjecture is then stated in terms of $\ext$ groups
of motivic structures, but whenever there is an explicit
``motivic'' construction of (conjecturally) all the relevant
extensions, the conjecture can be formulated entirely
in terms of premotivic structures.  This happens in our
case, for all the relevant $\ext$'s conjecturally vanish.
There will therefore be no further mention of motives
in this paper.

We make some other slight modifications to the framework of
\cite{fon_pr}.
\begin{itemize}
\item We work with premotivic structures with coefficients
 in a number field $K$.
\item We forget about the $\ell$-adic realization and
comparison isomorphisms at a finite set of ``bad''
primes $T$.
\end{itemize}
This yields a version of the conjecture which predicts
the value of $L(A_f,0)$ up to an $T$-unit in $K$.  We make
our set $T$ explicit:
Let $T_f$ be the set of finite primes $\lambda$ in $K$ such
that either:
\begin{itemize}
\item $\lambda\mid N k!$
\item the two-dimensional residual Galois representation
$\CM_{f,\lambda}/\lambda\CM_{f,\lambda}$
is not absolutely irreducible when restricted to $G_F$,
where $F= \QQ(\sqrt{(-1)^{(\ell-1)/2}\ell})$ and
$\lambda\mid \ell$.
\end{itemize}
Note that since $T_f$ includes the set of primes dividing $Nk!$,
we will only be applying Faltings' comparison theorem in the
``easy'' case of crystalline representations whose associated
Dieudonn\'e module has short filtration length.

\subsection{The main theorems}
Our main result can be stated as follows.

\begin{theorem} (=Theorem~\ref{thm:bk})
Let $f$ be a newform of weight $k\geq 2$ and level $N$ with
coefficients in $K$.
If $\lambda$ is not in $T_f$, then the $\lambda$-part of
the Bloch-Kato conjecture holds for $A_f$ and $A_f(1)$.
\mlabel{thm:bk0}
\end{theorem}
The main tool in the proof is the construction of
Taylor and Wiles, which we
axiomatize (Theorem~\ref{thm:axiomatic}),
and apply to higher weight forms to obtain the
following generalization of their class number formula.

\begin{theorem} (=Theorem~\ref{thm:selmer})
Let $f$ be a newform of weight $k\ge2$ and level $N$
with coefficients in $K$.
Suppose $\Sigma$ is a finite set of rational primes containing those
dividing $N$.  Suppose that $\lambda$ is a prime of $K$ which is not
in $T_f$ and does not divide any prime in $\Sigma$.
Then the $\CO_{K,\lambda}$-module
$$H^1_\Sigma(G_\QQ, A_{f,\lambda}/\CA_{f,\lambda})$$
has length $v_\lambda(\eta_f^\Sigma)$.
\mlabel{thm:selmer0}
\end{theorem}
Here $\eta_f^\Sigma$, defined in \S~\ref{ss:sig.int},
is a generalization of the congruence ideal of Hida
and Wiles; it can also be viewed as measuring the
failure of the pairing on $M_f$ to be perfect on $\CM_f$.

Another consequence of Theorem~\ref{thm:selmer0},
is the following result in the direction of
Fontaine-Mazur conjecture~\cite{fontaine_mazur}.

\begin{theorem} (=Theorem~\ref{thm:modular})
Suppose $\rho: G_\QQ \to \GL_2(K_\lambda)$
is a continuous geometric representation whose
restriction to $G_\ell$ is ramified and crystalline
and its associated Dieudonn\'e module has filtration
length less than $\ell -1$.
If its residual representation is modular
and absolutely irreducible restricted to
$\QQ(\sqrt{(-1)^{(\ell-1)/2}\ell})$ where $\lambda\mid \ell$,
then $\rho$ is modular.
\mlabel{thm:modular0}
\end{theorem}

\subsection{Acknowledgements}
Research on this project was carried out while the
first author worked at Cambridge, MIT and Rutgers, visited the IAS,
IHP and Paris VII, and received support from the EPSRC, NSF and an
AMS Centennial Fellowship.
The second author would like to thank the IAS for its hospitality and
acknowledge support from the NSF and the Sloan foundation.
The third author was supported in part by an NSF grant and a research
grant from University of Georgia in the early stages of this project,
and thanks the IAS for its hospitality.

\bigskip

\tableofcontents

\bigskip

\section{Generalities and examples of premotivic structures}
\mlabel{sec:mot}

\subsection{Galois representations}
\mlabel{ss:back}
\mlabel{sss:gal.rep}
For a field $F$, $\bar{F}$ will denote
an algebraic closure, and $G_F = \gal(\bar{F}/F)$.
We fix embeddings $\Qbar \to \Qpbar$ for each prime $p$,
and an embedding $\Qbar \to \CC$.
If $F$ is a number field, we let $\II_F$ denote
the set of embeddings $F \to \Qbar$, which we
identify with the set of embeddings $F \to \CC$
via our fixed one of $\Qbar$ in $\CC$.

Suppose that $F$ is a finite extension of $\QQ_p$ in $\Qpbar$.
We let $F_0$ denote the maximal unramified subextension
of $F$ and write $F^\ur$ for the maximal unramified extension
of $F$ in $\Qpbar$.
Let $\nu$ denote the natural map $G_F \to G_{\FF_p} \cong \hat{\ZZ}$
where we view $\hat{\ZZ}$ as being generated by
the geometric Frobenius element $\frob_p$, which we also
view as an automorphism of $\QQ_p^\ur$.  We write $\phi_p$
for the arithmetic Frobenius, so $\phi_p = \Frob_p^{-1}$.
Let $I_F$ denote the inertia subgroup of $G_F$,
$W_F$ the Weil subgroup, $\nu^{-1}(\ZZ)$
and $'W_F$ the Weil-Deligne group of $F$ (see \cite[\S8]{del_ant}).
Recall that if $K$ is a field of characteristic zero, then to give
a representation of $'W_F$ with coefficients in $K$ is
equivalent to giving a finite-dimensional
$K$-vector space $V$ (with the discrete topology),
a continuous homomorphism $\rho:W_F \to \aut_K V$
and a nilpotent endomorphism $N$ of $V$ satisfying
$\rho(g)N = p^{-\nu(g)}N\rho(g)$.
If $F = \QQ_p$, we write simply $G_p$, $I_p$,
$W_p$ and $'W_p$ for these groups and identify
$I_p \subset W_p \subset G_p$ with their images
in $G_\QQ$.

If $K$ is a number field, then $S_\f(K)$ denotes the
set of finite places of $K$.  Suppose that
$\lambda\in S_\f(K)$ divides $\ell\in S_\f (\QQ)$ and
$F \subset \Qlbar$ is a finite extension of $\QQ_\ell$.
Let $B_\dr=B_{\dr,\ell}$, $B_\crys=B_{\crys,\ell}$ and
$B_\st=B_{\st,\ell}$ be the rings defined by
Fontaine~\cite[\S2]{fon_ann},\cite[I.2.1]{fon_pr}.
Suppose that $V$ is a finite-dimensional
$K_\lambda$-representation (i.e., a $\lambda$-adic
representation) of $G_F$.  Then $D_\dr(V) =
(B_\dr\otimes_{\QQ_\ell} V)^{G_F}$ is filtered
free $F\otimes_{\QQ_\ell} K_\lambda$-module of
finite rank, and $V$ is called {\em de Rham} if
$\dim_{F} D_\dr(V)=\dim_{\QQ_\ell} V$.
Similarly $D_\crys(V) =(B_\crys\otimes_{\QQ_\ell} V)^{G_F}$
is a filtered free $F_0\otimes_{\QQ_\ell} K_\lambda$-module of
finite rank with a $(\phi_\ell\otimes 1)$-semilinear endomorphism
$\phi$, and $V$ is called {\em crystalline} if
$\dim_{F_0} D_\crys(V)=\dim_{\QQ_\ell} V$.
Similarly $D_\st(V) = (B_\st\otimes_{\QQ_\ell} V)^{G_F}$
is a filtered free $F_0\otimes_{\QQ_\ell} K_\lambda$-module of
finite rank with a $(\phi_\ell\otimes 1)$-semilinear endomorphism
$\phi$ and an endomorphism $N$ satisfying
$N\phi = \ell\phi N$, and $V$ is called {\em semistable} if
$\dim_{F_0} D_\st(V)=\dim_{\QQ_\ell} V$.
One defines $D_\pst(V)$ as the direct
limit of $(B_\st\otimes_{\QQ_\ell} V)^{G_{F'}}$ where $F'$
runs over finite extensions of $F$ in $\Qlbar$; this is a filtered
free $\QQ_\ell^\ur\otimes_{\QQ_\ell} K_\lambda$-module
and one checks that $(\rho,N)$ with $\rho(g) = g\phi^{\nu(g)}$
makes $D_\pst(V)$ a $K_\lambda$-rational representation of $'W_F$,
and $V$ is called {\em potentially semistable} if
$\dim_{\QQ_\ell^\ur} D_\pst(V)=\dim_{\QQ_\ell} V$.

Recall that $V$ is crystalline if and only if $V$ is semistable
and $N = 0$ on $D_\st(V)$, in which case $D_\cris(V) = D_\st(V)$;
$V$ is semistable if and only if $V$ is potentially semistable
and $I_F$ acts trivially on $D_\pst(V)$, in which case $D_\pst(V)
= \QQ_\ell^\ur\otimes_{F_0}D_\st(V)$ with the action of $g\in W_F$
given by $\phi^{\nu(g)}$; $V$ is potentially semistable
if and only if $V' = \res_{G_F}^{G_{F'}}V$ is semistable for some
finite extension $F'$ of $F$, in which case $D_\pst(V)
= \QQ_\ell^\ur\otimes_{F_0'}D_\st(V')$; and if $V$ is
potentially semistable, then $V$ is de Rham with
$F^\ur \otimes_F D_\dr(V)
= F^\ur\otimes_{\QQ_\ell^\ur}D_\pst(V)$ (with $D_\dr(V)
= F\otimes_{F_0}D_\st(V)$ if $V$ is semistable).

A $\lambda$-adic representation $V$ of $G_\QQ$ is
{\em pseudo-geometric}~\cite[II.2]{fon_pr}
(resp.\ {\em geometric}~\cite[II.3]{fon_pr})
if it is unramified outside of a finite number of
places of $\QQ$ and its restriction to
$G_\ell$ is de Rham (resp.\ potentially semistable).
The representation $V$ is said to have
{\em good reduction at $p$} if its restriction to
$G_p$ is crystalline (resp.\ unramified) if $p = \ell$
(resp.\ $p \neq \ell$).

We also recall some of the theory of Fontaine and Laffaille
\cite{fon_laf}.  Suppose that $F\subset \Qlbar$ is a finite
unramified extension of $\QQ_\ell$.  We let
$\mfcat(\CO_F)$ denote the category whose objects are
finitely generated $\CO_F$-modules equipped with
\begin{itemize}
\item a decreasing filtration such that $\fil^{a} A = A$ and
   $\fil^b A = 0$ for some $a,b\in\ZZ$,
and for each $i \in \ZZ$, $\fil^i A$ is a direct summand of $A$;
\item $\phi_\ell$-semilinear
maps $\phi^i:\fil^i A \to A$ for $i \in \ZZ$ satisfying
$\phi^i|_{\fil^{i+1}A} = \ell \phi^{i+1}$ and $A = \sum \Im\phi^i$.
\end{itemize}
It follows from \cite[1.8]{fon_laf} that $\mfcat(\CO_F)$ is an
abelian category.  Let $\mfcat^a(\CO_F)$ denote the full subcategory
of objects $A$ satisfying $\fil^{a}A=A$ and $\fil^{a+\ell}A = 0$ and
having no non-trivial quotients $A'$ such that
$\fil^{a+\ell-1} A' = A'$, and let $\mfcat^a_\tor(\CO_F)$ denote the full
subcategory of $\mfcat^a(\CO_F)$ consisting of objects of finite length.
So $\mfcat^0_\tor(\CO_F)$ is the category denoted
$\underline{MF}_\tor^{f,\ell'}$
in \cite{fon_laf}, and it follows from \cite[6.1]{fon_laf} that
$\mfcat^a(\CO_F)$ and $\mfcat^a_\tor(\CO_F)$
are abelian categories, stable under
taking subobjects, quotients, direct products and extensions
in $\mfcat(\CO_F)$.  We write simply $\mfcat$, $\mfcat^a$, etc.\
whenever $F=\QQ_\ell$.

Fontaine and Laffaille define a contravariant functor $\underline{U}_S$
from $\mfcat^0_\tor(\CO_F)$ to the category of continuous
$\ZZ_\ell[G_F]$-modules
which are finitely generated over $\ZZ_p$, and they prove it is fully
faithful \cite[6.1]{fon_laf}.  We let $\VV$ denote the functor
defined by $\VV(A) = \hom(\underline{U}_S(A),\QQ_\ell/\ZZ_\ell)$
and we extend it to a fully faithful functor on $\mfcat^0(\CO_F)$ by
setting $\VV(A) = \projlim \VV(A/\ell^nA)$.  Then $\VV$ defines
an equivalence between $\mfcat^0(\CO_F)$ and the full subcategory of
$\ZZ_\ell[G_F]$-modules whose objects are isomorphic to quotients
of the form $L_1/L_2$, where $L_2 \subset L_1$ are finitely
generated submodules of crystalline representations $V$ with
the following properties:
\begin{itemize}
\item $\fil^0 D = D$ and $\fil^\ell D = 0$, where
$D = (B_\crys\otimes_{\QQ_\ell} V)^{G_F}$;
\item if $V'$ is a nonzero quotient of $V$,
then $V'\otimes_{\QQ_\ell}\QQ_\ell(\ell-1)$ is ramified.
\end{itemize}
In particular, the essential image of $\VV$ is closed under
taking subobjects, quotients and finite direct sums.
Furthermore, one sees from \cite[8.4]{fon_laf} that the
natural transformations
\begin{equation}
 \begin{array}{l}
F \otimes_{\CO_F}A \to
    (B_\crys\otimes_{\ZZ_\ell} \VV(A))^{G_F},\\
B_\crys\otimes_{\CO_F}A \to
B_\crys \otimes_{\ZZ_\ell} \VV(A) {\rm\ and\ } \\
\Fil^0(B_\cris\otimes_{\ZZ_\ell}A)^{\phi=1}
\to  \QQ_\ell \otimes_{\ZZ_\ell}\VV(A)
\end{array}
\mlabel{eq:VV}
\end{equation}
are isomorphisms.

If $K$ is a number field and $\lambda\in S_\ff(K)$
is a prime over $\ell$, we let $\CO_\lambda=\CO_{K,\lambda}$
and let $\mfcato^a$
denote the category
of $\CO_\lambda$-modules
in $\mfcat^a$.  We can regard $\VV$
as a functor from $\CO_{\lambda}$-$\mfcat^0$ to the category
of $\CO_{\lambda}[G_\ell]$-modules.

If $A$ and $A'$ are objects of $\CO_{\lambda}$-$\mfcat^0$ such that
$A\otimes_{\CO_{\lambda}}A'$ defines an object of
$\CO_{\lambda}$-$\mfcat^0$, then there is a
canonical isomorphism
$$\VV(A\otimes_{\CO_{\lambda}}A') \cong
\VV(A)\otimes_{\CO_{\lambda}}\VV(A').$$
Analogous assertions hold for $\hom_{\CO_{\lambda}}(A,A')$.

\subsection{Premotivic structures}
\mlabel{ssec:pre.str}
We work with categories of premotivic structures based
on notions from \cite{fon_pr} and \cite{bloch_kato}.

For a number field $K$,
we let $\pms_K$ denote the category of
premotivic structures over $\QQ$ with coefficients
in $K$.  In the notation of \cite[III.2.1]{fon_pr},
this is the category $\spm_\QQ(\QQ)\otimes K$
of $K$-modules in $\spm_{\QQ}(\QQ)$.  Thus an object
$M$ of $\pms_K$ consists of the following data:
\begin{itemize}
\item a finite-dimensional $K$-vector space $M_\bt$
with an action of $G_{\RR}$,
\item a finite-dimensional $K$-vector space $M_{\dr}$
with a finite decreasing filtration $\fil^i$, called
the Hodge filtration;
\item for each $\lambda \in S_\f(K)$,
a finite-dimensional $K_\lambda$ vector space $M_\lambda$
with a continuous pseudo-geometric (see \cite[II.2.1]{fon_pr})
action of $G_\QQ$;
\item a $\CC\otimes K$-linear isomorphism
$$I^\infty: \CC\otimes M_\dr \to \CC\otimes M_\bt$$
respecting the action of $G_{\RR}$ (where $G_\RR$
acts on $\CC\otimes M_\bt$ diagonally and acts on
$\CC\otimes M_\dr$ via the first factor);
\item for each $\lambda \in S_\f(K)$,
a $K_\lambda$-linear isomorphism
$$I^\lambda_\bt: K_\lambda\otimes_K M_\bt \to M_\lambda$$
respecting the action of $G_{\RR}$ (where the action
on $M_\lambda$ is via the restriction $G_\RR \to G_\QQ$
determined by our choice of embedding $\Qbar \to \CC$);
\item for each $\lambda \in S_\f(K)$,
$B_{\dr,\ell}\otimes_{\QQ_\ell} K_\lambda$-linear isomorphism
$$I^\lambda: B_{\dr,\ell}\otimes_{\QQ_\ell} K_\lambda \otimes_K M_{\dr}
\to B_{\dr,\ell}\otimes_{\QQ_\ell} M_\lambda$$
respecting filtrations and the action of $G_{\QQ_\ell}$
(where $\ell$ is the prime which $\lambda$ divides,
$K_\lambda$ and $M_\lambda$ are given the degree-0 filtration,
$K_\lambda$ and $M_\dr$ are given the trivial $G_{\QQ_\ell}$-action
and the action on $M_\lambda$ is determined by our choice of
embedding $\Qbar \to \Qlbar$);
\item increasing weight filtrations $W^i$ on $M_\bt$, $M_{\dr}$
and each $M_\lambda$ respecting all of the above data,
and such that $\RR\otimes M_\bt$ with its Galois action and
weight filtration, together with the Hodge filtration on
$\CC\otimes M_\bt$ defined via $I_\bt$, defines a mixed Hodge
structure over $\RR$ (see \cite[III.1]{fon_pr}).
\end{itemize}

If $S\subsetneq S_\f(K)$ is a set of primes of $K$, we let $\pms_K^S$
denote the category defined in exactly the same way, but
with $S_\f(K)$ replaced by the complement of $S$.  If $S \subseteq S'$,
we use $\cdot^{S'}$ to denote the forgetful functor from
$\pms_K^S$ to $\pms_K^{S'}$.

The category $\pms_K^S$ is equipped with a tensor product,
which we denote $\otimes_K$, and an internal hom, which
we denote $\hom_K$.  There is also a unit object, which
we denote simply by $K$.  These are defined in the obvious way;
for example, $(M\otimes_K N)_\bt$ is the $K[G_\RR]$-module
$M_\bt \otimes_K N_\bt$.  If $K \subseteq K'$, we let
$S^{K'}$ denote the set of primes in $S_\f(K')$ lying over
those in $S$, then $K'\otimes_K\cdot$ defines a functor
from  $\pms_K^S$ to $\pms_{K'}^{S^{K'}}$.

\begin{definition}  An object $M$ of $\pms_K^S$
\begin{itemize}
\item has {\em good reduction} at $p$ if for each
$\lambda$ in $S$, $M_\lambda|G_p$ is crystalline at
$p$ if $\lambda|p$ and unramified at $p$ otherwise;
\item is {\em $L$-admissible} at $p$ if $M_\lambda$ is geometric
for every $\lambda \not\in S$ dividing $p$ and the representations
$D_{\pst}(M_\lambda|G_p)^{\ss}$, for $\lambda\not\in S$,
form a compatible system of $K$-rational representations
of the Deligne-Weil group $'W_p$ (see \cite[\S8]{del_ant});
\item is {\em $L$-admissible everywhere} if it is
$L$-admissible at $p$ for all primes $p$.
\end{itemize}
\end{definition}

If $M$ is $L$-admissible at $p$, then the local
factor associated to $D_{\pst}(M_\lambda|G_p)^{\ss}$
is of the form $P(p^{-s})^{-1}$ for some polynomial
$P(u) \in K[u]$ independent of $\lambda$ not in $S$.
For an embedding $\tau:K\to\CC$ we put $L_p(M,\tau,s)=\tau P(p^{-s})$
and we regard the collection $\{L_p(M,\tau,s)\}_{\tau\in\II_K}$ as a
meromorphic function on $\CC$ with values in
$\CC^{\II_K}\cong K\otimes\CC$. If $S$ is finite and $M$ is
$L$-admissible everywhere, then its $L$-function
\[L(M,s):=\prod_p L_p(M,s)\]
is a holomorphic $K\otimes\CC$-valued function in some right
half plane ${\rm Re}(s)>r$ (with components $L(M,\tau,s)
=\prod_p L_p(M,\tau,s)$).

For each $\tau \in\II_K$ the Hodge structure
$\CC\otimes_{K,\tau}M_\bt$ over $\RR$ has attached to it an L-function
$L_\infty(M,\tau,s)$ \cite{del_cor} and we denote by
$L_\infty(M,s)$ the resulting function with values in $K\otimes\CC$.
Similarly, if $M$ is $L$-admissible everywhere
we define the $K\otimes\CC$-valued $\epsilon$-factor
$\epsilon(M,s)=\epsilon(M)c(M)^{-s}$ where the components
$\epsilon(M,\tau,s)$ are defined in \cite{del_ant}. $L_\infty(M,\tau,s)$ is
essentially
a product of $\Gamma$-functions, $c(M)$ (the conductor of $M$)
is a positive integer divisible by $p$ if and only if $M$ does not
have good reduction at $p$, and $\epsilon(M)$ is independent of $s$.
Put $\Lambda(M,s):= L_\infty(M,s)L(M,s)$.   Combining Conjecture 1
of \cite{fontaine_mazur} with Deligne's conjectured functional equation
leads to the following conjecture (see \cite[III.4.3]{fon_pr}):

\begin{conjecture} If $M$ is $L$-admissible everywhere, then
$L(M,s)$ converges to a holomorphic function on some right half-plane
and extends meromorphically to $\CC$.  Moreover
$\Lambda(M,s)$ satisfies the functional equation
\[ \Lambda(M,s) =\epsilon(M,s)\Lambda(M^*,1-s)\]
where $M^*=\Hom_K(M,K)$ is the contragredient of $M$,
and it is holomorphic if $M$ does not have a direct summand isomorphic to the
unit premotivic structure $K$.
\mlabel{fe}
\end{conjecture}

We also define a category $\ipms^S_K$ of {\em
$S$-integral premotivic structures} as follows.
We let $\CO_S =\CO_{K,S}$
denote the set of $x \in K$ with $v_\lambda(x) \ge 0$ for
all $\lambda \not\in S$.  An object
$\CM$ of $\ipms^S_K$ consists of the following data:
\begin{itemize}
\item a finitely generated $\CO_K$-module $\CM_\bt$
with an action of $G_{\RR}$,
\item a finitely generated $\CO_S$-module $\CM_{\dr}$
with a finite decreasing filtration $\fil^i$, called
the Hodge filtration;
\item for each $\lambda \in S_\f(K)$, a finitely generated
$\CO_{\lambda}$-module $\CM_\lambda$ with continuous action of $G_\QQ$
inducing a pseudo-geometric action on
$\CM_\lambda\otimes_{\CO_{\lambda}}K_\lambda$;
\item for each $\lambda \not\in S$, an object $\CM_{\lcrys}$
of $\mfcato^0$;
\item an $\RR\otimes \CO_K$-linear isomorphism
$$I^\infty: \CC\otimes \CM_\dr \to \CC\otimes \CM_\bt$$
respecting the action of $G_{\RR}$;
\item for each $\lambda$ in $S_\f(K)$ an isomorphism
$$I^\lambda_\bt:\CM_\bt\otimes_{\CO_K}\CO_{\lambda}\cong
\CM_{\lambda}$$
respecting the action of $G_\RR$;
\item for each $\lambda\not\in S$ an $\CO_\lambda$-linear isomorphism
$$I^\lambda_\dr:\CM_\dr\otimes_{\CO}\CO_\lambda \cong \CM_{\lcrys}$$
respecting filtrations;
\item for each $\lambda \not\in S$, an $\CO_{\lambda}$-linear
isomorphism
$$I^\lambda:\VV(\CM_{\lcrys}) \to \CM_{\lambda}$$
respecting the action of $G_{\QQ_\ell}$,
where $\ell$ is the prime which $\lambda$ divides.
\item increasing weight filtrations $W^i$ on $\QQ\otimes\CM_\bt$,
$\QQ\otimes\CM_{\dr}$ and each $\QQ\otimes\CM_\lambda$
respecting all of the above data and giving rise to a mixed Hodge
structure.
\end{itemize}

With the evident notion of morphism this becomes an
$\CO_K$-linear category.
In particular, it is an abelian category.
Note that the forgetful functor sending $\CM$ to $\CM_B$ is
faithful.
Note also that there is a natural functor $\QQ\otimes\cdot$
from $\ipms_K^S$ to $\pms_K^S$, where we set
$(\QQ\otimes\CM)_\wc = \QQ\otimes\CM_\wc$ for $\wc = \bt$, $\dr$
and $\lambda$ for $\lambda\not\in S$,
with induced additional structure and comparison isomorphisms.
(The comparison $I^\lambda$ for $\QQ\otimes\CM$
is defined as the composite
$$B_{\dr,\ell}\otimes_{\ZZ_\ell} \CO_\lambda \otimes_\CO \CM_{\dr} \cong
B_{\dr,\ell}\otimes_{\ZZ_\ell}\CM_{\lcrys}
\to B_{\dr,\ell}\otimes_{\ZZ_\ell}\VV(\CM_{\lcrys})
\to B_{\dr,\ell}\otimes_{\ZZ_\ell}\CM_{\lambda},$$
where the maps are respectively, $I^\lambda_\dr$,
the canonical map (\ref{eq:VV}) and
$I^\lambda$, each with scalars
extended to $B_{\dr,\ell}$.)

If $S\subset S'$, we define a functor $\cdot^{S'}$ from
$\ipms_K^{S}$ to $\ipms_K^{S'}$ in the obvious way.
We say that $\CM$ is {\em $S'$-flat} if
$\CM_\dr^{S'}=\CM_\dr\otimes_{\CO_{K,S}}\CO_{K,S'}$
is flat over $\CO_{K,S'}$.  Note that if $\CM$ is
$S'$-flat, then so is any subobject of $\CM$.
Let $K'$ be a finite extension of $K$.
We also have a natural
functor $\CO_{K'}\otimes_{\CO_K}\cdot$
from $\ipms_K^S$ to $\ipms_{K'}^{S'}$ where $S'$ is the set
of primes over those in $S$.

We say that $\CM$ has {\em good reduction} at $p$, is
{\em $L$-admissible} at $p$ or is
{\em $L$-admissible everywhere} according to whether the
same is true for $\QQ\otimes\CM$.  Note that if $\CM$
is $L$-admissible at $p$ and $p$ is not invertible in
$\CO$, then $\CM$ necessarily has good reduction at $p$.

For objects $\CM$ and $\CM'$ of $\ipms_K^S$, we can
form $\CM\otimes_{\CO_K}\CM'$ in $\ipms_K^S$
provided
$\CM_{\lcrys}\otimes_{\CO_{\lambda}}\CM'_{\lcrys}$
defines an object of $\mfcat^0_\lambda$ for all $\lambda\not\in S$.
In particular this holds if there exist positive integers $a$, $a'$
such that $\fil^a\CM_\dr = 0$, $\fil^{a'}\CM'_\dr = 0$ and $a+a' < \ell$
for all primes $\ell$ not invertible in $\CO$.
If $\CN$ and  $\CN'$ are objects of $\ipms_K^S$ such
that $\CN\otimes_{\CO_K}\CN'$ is as well, and
if $\alpha:\CM\to \CN$ and $\alpha':\CM'\to \CN'$ are morphisms
in $\ipms_K^S$, then there is a well-defined
morphism $\alpha\otimes\alpha':\CM\otimes_{\CO_K}\CN
\to \CM'\otimes_{\CO_K}\CN'$ in $\ipms_K^S$.
Analogous assertions hold
for the formation and properties of $\hom_{\CO_K}(\CM,\CM')$.

Note that if $\CM$ is an object of $\ipms_K^S$, then $\End\CM$
is a finitely generated $\CO_K$-module.  If $I$ is an
$\CO_K$-submodule of $\End\CM$, then we define
an object $\CM[I]$ of $\ipms_K^S$ as the kernel of
$$(x_1,\ldots,x_r): \CM \to \CM^r$$
where $x_1,\ldots,x_r$ generate $I$.  This is independent
of the choice of generators.  This applies in particular
when $I$ is the image in $\End \CM$ of an ideal in a commutative
$\CO_K$-algebra $R$ mapping to $\End\CM$, or the augmentation
ideal in $\CO_K[G]$ where $G$ is a group acting on $\CM$.
In the latter case, we write $\CM^G$ instead of $\CM[I]$.

\subsection{Basic examples}
\mlabel{ssec:mot.ex}
\mlabel{sssec:mot.ex.dir}

The object $\QQ(-1)$ in $\pms_\QQ$ is the weight
two premotivic structure defined by $H^1(\GG_\m)$.  To give an
explicit description, let $\varepsilon$ denote the
generator of $\ZZ_\ell(1) = \invlim \mu_{\ell^n}(\Qbar)$
defined by $(e^{2\pi i/\ell^n})_n$ via our fixed embedding
$\Qbar \to \CC$.
\begin{itemize}
\item Let $\calt_B = H^1_B(\GG_\m(\CC),\ZZ) \cong (2\pi i)^{-1}\ZZ \subset \CC$
    with complex conjugation in $G_\RR$ acting by $-1$, and let $\QQ(-1)_B
= \QQ\otimes\calt_B$.
\item  Let $\calt_{\dr} = H^1_{\dr}(\GG_\m/\ZZ)$, which with its Hodge
    filtration is isomorphic to $\ZZ[-1]$ (where $[n]$ denotes
    a shift by $n$ in the filtration, so $\fil^i V[n] = \fil^{i+n}V$).
Write $\iota$ for the canonical basis of $\calt_{\dr} \cong \ZZ[-1]$ and
let $\QQ(-1)_{\dr} = \QQ\otimes\calt_{\dr}$, so
    $\fil^1\QQ(-1)_{\dr} = \QQ\iota$ and $\fil^2\QQ(-1)_{\dr} = 0$.
\item Let $\calt_\ell = H^1_{\et}(\GG_{\m,\Qbar},\ZZ_\ell)
\cong \hom_{\ZZ_\ell}(\ZZ_\ell(1),\ZZ_\ell) = \ZZ_\ell\delta$
where $\delta(\varepsilon) = 1$, and let $\QQ(-1)_\ell
= \QQ\otimes\calt_\ell$.
\item Let $\calt_\ecrys$ denote the object of $\mfcat_\ell$
    defined by $\ZZ_\ell\otimes\calt_{\dr}
    =\ZZ_\ell~\iota$ with $\phi^1(\iota) = \iota$.
\item $I^\infty:\CC\otimes_\ZZ \calt_\dr \to
    \CC\otimes_\ZZ \calt_\bt$ is defined by
    $1\otimes\iota  \mapsto 2\pi i\otimes (2\pi i)^{-1}$.
\item $I^\ell_\bt:\ZZ_\ell\otimes\calt_\bt \to \calt_\ell$ is
    defined by $1\otimes(2\pi i)^{-1} \mapsto \delta$.
\item $I^\ell_\dr:\calt_\dr \otimes_\ZZ \ZZ_\ell \to \calt_\ecrys$
    is given by $\iota\mapsto \iota$.
\item $I^\ell: B_{\dr,\ell}\otimes\QQ[-1] \to
   B_{\dr,\ell}\otimes_{\QQ_\ell} \QQ(-1)_\ell$
   is defined by $1 \otimes \iota \mapsto t \otimes \delta$ where
$t = \log[\varepsilon]$.
\item For $\ell > 2$, $\calt_\ecrys$ is an object of $\mfcat^0$
   and $I^\ell$ is induced by an isomorphism
   $\VV(\CT_\ecrys)\cong\CT_\ell$.
\end{itemize}
The above data defines objects in $\pms_\QQ$ and $\ipms^{\{2\}}_\QQ$
which we denote by $\QQ(-1)$ and $\calt$.
These could be described equivalently by $H^2(\PP^1)$, or indeed
$H^2(X)$ for any smooth, proper, geometrically connected curve
$X$ over $\QQ$.

The {\em Tate premotivic structure} $\QQ(1)$ is the object
of $\pms_\QQ$ defined by $\hom_\QQ(\QQ(-1),$ $\QQ)$.
More generally,
for any object $M$ in $\pms_K$ and integer $n$, $M(n)$ is defined
as $M\otimes_{\QQ}\QQ(1)^{\otimes n}$.  If $M$ has good reduction at $p$, then
so does $M(n)$.  If $M$ is $L$-admissible at $p$, then so is $M(n)$.
If they are $L$-admissible everywhere, then $L(M(n),s) = L(M,s+n)$,
$\Lambda(M(n),s) = \Lambda(M,s+n)$, $c(M(n)) = c(M)$, $\epsilon(M(n),s)
= \epsilon(M,s+n)$ and Conjecture~\ref{fe} for $M$ and $M(n)$ are equivalent.
In particular, the premotivic structure $\QQ(n)$ has good reduction and is
$L$-admissible at $p$ for all primes $p$, and $L(\QQ(n),s) = \zeta(s+n)$
where $\zeta$ is the Riemann $\zeta$-function.  Furthermore $L_\infty(\QQ(n),s)
= \Gamma_\RR(s+n)$ where $\Gamma_\RR(s) = L_\infty(\QQ,s) = \pi^{-s/2}\Gamma(s/2)$,
$\epsilon(\QQ(n),s) = 1$ and Conjecture~\ref{fe} for $\Lambda(\QQ(n),s)$
is the functional equation for $\zeta(s)$.

For any integer $n\ge 0$, $\calt^{S,\otimes n}$ defines an object of
$\ipms^S_\QQ$ where $S$ is any set of primes containing those
dividing $(n+1)!$; note also that
$\QQ\otimes\calt^{S,\otimes n} \cong \QQ(-n)^S$.

\medskip
For any number field $F \subset \Qbar$ let $M_F$
denote the premotivic structure $M_F$ of weight zero
defined by $H^0(\spec F)$, called the {\em Dedekind premotivic structure}
of $F$.  To give an explicit description, let $S$ denote
the set of primes dividing  $D=\Disc(F/\QQ)$.
We let
\begin{itemize}
\item $\CM_{F,\bt} = \ZZ^{\II_F}$ with the natural action of $G_\RR$,
    and  $M_{F,\bt} = \QQ\otimes\CM_{F,\bt} = \QQ^{\II_F}$.
    (Recall that we identified $\II_F$ with the set of embeddings
    $F \to \CC$ via the chosen embedding $\Qbar \to \CC$, so for
    $\alpha:\II_F \to \ZZ$ and $\sigma \in G_\RR$, we
    define $\sigma\alpha$ by $\tau \mapsto
    \alpha(\sigma^{-1}\circ\tau)$).
\item $\CM_{F,\dr} = \CO_F[1/D]$
   with $\fil^0 \CM_{F,\dr} = \CM_{F,\dr}$ and $\fil^1 \CM_{F,\dr} = 0$,
   and $M_{F,\dr} = \QQ\otimes\CM_{F,\dr} = F$.
\item $\CM_{F,\ell}= \ZZ_\ell\otimes \CM_{F,\bt} =\ZZ_\ell^{\II_F}$
    with the natural action of $G_\QQ$, and $M_{F,\ell}
    = \QQ\otimes\CM_{F,\ell} = \QQ_\ell^{\II_F}$
    (so for $\alpha:\II_F \to \ZZ_\ell$, $\sigma \in G_\QQ$
     and $\tau:F \to \Qbar$, we have $(\sigma\alpha)(\tau) =
    \alpha(\sigma^{-1}\circ\tau)$).
\item $\CM_{F,\ecrys} = \ZZ_\ell\otimes\CM_{F,\dr}=
  \ZZ_\ell\otimes \CO_F$ for $\ell\not\in S$, with the same
  filtration as $\CM_{F,\dr}$ and with $\phi^0 = \phi_\ell$.
\item Identifying $\CC\otimes M_{F,B}$ with $\CC^{\II_F}$,
    $I^\infty$ is defined by $I^\infty(1\otimes x)(\tau) = \tau(x)$.
\item $I^\ell_\bt$ is the identity map.
\item $I^\ell_\dr$ is the identity map for $\ell\not\in S$.
\item Identifying $B_{\dr,\ell}\otimes M_{F,\ell}$ with
     $B_{\dr,\ell}^{\II_F}$ and $\II_F$ with the set of embeddings
     of $F$ in $B_{\dr,\ell}$ via the chosen embedding
     $\Qbar \to \Qpbar$, $I^\ell$, is defined by
     $I^\ell(1\otimes x)(\tau) = \tau(x)$.  For $\ell\not\in S$,
     this is induced by an isomorphism
     $\VV(\CM_{F,\ecrys}) \to \CM_{F,\ell}$ which we also
     denote $I^\ell$.
\end{itemize}
The above data defines objects $\CM_F$ of $\ipms^S_\QQ$ and
$M_F$ of $\pms_\QQ$ with $\QQ\otimes\CM_F = M_F^S$.  One finds
that $M_F$ has good reduction at all primes $p$ not dividing $D$,
is $L$-admissible everywhere, and $L(M_F,s) = \zeta_F(s)$ where
$\zeta_F$ is the $\zeta$-function of $F$.  Furthermore $c(M_F)=D$,
$\epsilon(M_F) = D^{1/2}$ and $L_\infty(M_F,s) = \Gamma_\RR(s)^{r_1}\Gamma_\CC(s)^{r_2}$
where $r_1$ (resp.\ $r_2$) is the number of real (resp.\ complex)
embeddings of $F$, $\Gamma_\RR$ was defined above and $\Gamma_\CC
= 2\cdot(2\pi)^{-s}\Gamma(s)$.  Conjecture~\ref{fe} is known
for $\Lambda(M_F,s)$ (see \cite[Ch.\ XIII]{lang} for example).

Note that if $F$ is Galois over $\QQ$, then there is a natural action of
$G = \gal(F/\QQ)$ on $\CM_F$ defined as follows:  For $\alpha \in \CM_{F,B}$,
$g \in G$ and $\tau\in \II_F$, we define $g\alpha$ by $\tau\mapsto\alpha(\tau\circ g)$.
The action on $\CM_{F,\ell}$ is defined similarly. The action on
$\CM_{F,\dr} = \CO_F[1/D]$ is the obvious one, as is the action on
$\CM_{F,\ecrys}$.  One checks that the action respects
all the additional structure, and so defines an action on
$\CM_F$ as an object in $\ipms^S_\QQ$.
We similarly have an action of $G$ on $M_F\in \pms_\QQ$.

Suppose that we are given a character $\psi:\hat{\ZZ}^\times \to
K^\times$. We regard $\psi$ also as a character of
$\AA^\times$ and $G_\QQ$ via the isomorphisms $\hat{\ZZ}^\times
\cong \AA^\times/\RR_{>0}^\times\QQ^\times \cong G_\QQ^\ab$,
where the first isomorphism is induced by the natural
inclusion $\hat{\ZZ}^\times \to \AA^\times$ and the second
is given by class field theory.  (Our convention is that
a uniformizer in $\QQ_p^\times$ maps to $\frob_p$ in the
Galois group of any abelian extension of $\QQ$ unramified
at $p$.)

Suppose we are given an embedding  of a number
field $F \to \Qbar$ such that $\psi$ is trivial on the
image of $G_F$ in $G_\QQ$.  We can then regard $\psi$
as a character of $G = \gal(F/\QQ)$ and define
the {\em Dirichlet premotivic structure}
$M_{\psi}$ as $(V\otimes M_F)^G$ where
$V = K$ with $G$ acting by $\psi$.
One can check that the construction
is independent of the choice of $F$ and embedding
$F \to \Qbar$ as above.  To describe $M_{\psi}$ explicitly,
we choose $F=\QQ(e^{2\pi i/N}) \subset \CC$ where $\psi$
has conductor $N$.  We let $\tau_0: F \to \Qbar$ denote
the embedding compatible with our fixed $\Qbar \to \CC$,
and we regard $\psi$ as a Dirichlet character via the
canonical isomorphism $(\ZZ/N\ZZ)^\times \cong G$.
Let $S$ denote the set of primes in $K$ lying over those
dividing $N$ and define an object $\CM_\psi$ of
$\ipms_K^S$ by $(\CV\otimes\CM_F)^G$ where
$\CV = \CO_K^S$ with $G$ acting by $\psi$.
We then have:
\begin{itemize}
\item $\CM_{\psi,\bt}$ is the $\CO_K$-submodule of
    $\CO_K^{\II_F}$ spanned by the map
    $b_B$ defined by $\tau_0\circ g \mapsto \psi^{-1}(g)$, where $\tau_0$
    is the inclusion of $F$ in $\Qbar$;
\item $\CM_{\psi,\dr}$ is the $\CO = \CO_K[1/N]$-submodule of
    $\CO_K\otimes \CO_F[1/N]$ spanned by
    $b_{\dr} = \dsum_a \psi(a) \otimes e^{2\pi i a/N}$,
    where $a$ runs over $(\ZZ/N\ZZ)^\times$ with
    $\fil^1 \CM_{\psi,\dr}=0$ and
    $\fil^0 \CM_{\psi,\dr}=\CM_{\psi,\dr}$;
\item $\CM_{\psi,\lambda}
    = \CO_\lambda \otimes_{\CO_K} \CM_{\psi,\bt}$ with $G_\QQ$
    acting via $\psi$;
\item for $\lambda \not| N$,
    $\CM_{\psi,\lcrys}=\CO_\lambda\otimes_{\CO_K} \CM_{\psi,\dr}$
    with
    the same filtration as $\CM_{\psi,\dr}$ and $\phi^0=\psi^{-1}(\ell)$.
\end{itemize}
The  comparison isomorphisms are induced from those of $\CM_F$.
Similarly, we get the object $M_\psi$ of $\pms_K$ by setting
$M_{\psi,\wc} = \QQ\otimes\CM_{\psi,\wc}$ with
comparison isomorphisms induced from those of $M_F$.
In particular, we have $I^\infty(1\otimes b_{\dr})=
G_\psi(1\otimes b_B)$ where $G_\psi$ is the Gauss sum
$\dsum_a e^{2\pi ia/N}\otimes\psi(a)$ in $\CC\otimes K$.

We have that $\QQ\otimes\CM_\psi \cong M_\psi^S$,
$M_{\psi}$ has good reduction at all primes not dividing
the conductor of $\psi$ and is $L$-admissible everywhere,
and $L(M_{\psi},s)$ is the Dirichlet $L$-function $L(\psi^{-1},s)$.
Moreover, Conjecture~\ref{fe} is known for $\Lambda(M_{\psi},s)$ and one has
$c(M_{\psi}) = N$, $\epsilon(M_{\psi})=i^{-\eta}G_{\psi^{-1}}$ and
$L_\infty(M_{\psi},s) = \Gamma_\RR(s+\eta)$, where $\psi(-1)=(-1)^\eta$,
$\eta\in\{0,1\}$ (see \cite[Ch.\ 4]{washington} for example).

\section{Premotivic structures for modular forms}
\mlabel{sec:levelN}
In this section we review the construction of
premotivic structures associated to the space of
modular forms of weight $k$ and level $N$.
More precisely, if $k \ge 2$ and $N\ge 3$, we
construct objects of $\ipms_\QQ^S$ whose de Rham
realization contains the space of such forms, where
$S=S_N$ is the set of primes dividing $Nk!$.

\subsection{Level $N$ modular curves}  \mlabel{sssec:modcurves}
Let $k$ and $N$ be integers with $k \ge 2$ and $N \ge 3$.
Let $T= \spec \ZZ[1/Nk!]$,
and consider the functor which associates
to a $T$-scheme $T'$ the set of isomorphism classes of generalized
elliptic curves over $T'$ with level $N$ structure \cite[IV.6.6]{del_rap}.
By \cite[IV.6.7]{del_rap}, the functor is represented by a smooth,
proper curve over $T$.  We denote this curve by $\bar{X}$, and
we let $\bar{s}:\bar{E} \to \bar{X}$ denote the universal
generalized elliptic curve with level $N$ structure.
We let $X$ denote the open subscheme of $\bar{X}$ over
which $\bar{E}$ is smooth.  Then $X$ is the complement of
a reduced divisor, called {\em the cuspidal divisor}, which
we denote by $X^\infty$.  We let $E = \bar{s}^{-1}X$, $s = \bar{s}|_E$
and $E^\infty = \bar{E}\times_{\bar{X}}X^\infty$.
Using the arguments of \cite[VII.2.4]{del_rap}, one can check
that $\bar{E}$ is smooth over $T$ and $E^\infty$ is a
reduced divisor with strict normal crossings (in the sense of
\cite[1.8]{gr_mu} as well as \cite[XIII.2.1]{sga1}).

Let us also recall the standard description of
$$\bar{s}^\an: \bar{E}^\an \to \bar{X}^\an,$$
where we use ${}^\an$ to denote the associated complex
analytic space. We let
$$X_N = \coprod_{t\in(\ZZ/N\ZZ)^\times} X_{N,t},$$
where for each $t$, $X_{N,t}$ denotes a copy of
$\Gamma(N)\backslash\uhp$, the quotient of the complex upper
half-plane $\uhp$, by the principal congruence subgroup
$\Gamma(N)$ of $\SL_2(\ZZ)$.  Similarly we let
$$\bar{X}_N = \coprod_{t\in(\ZZ/N\ZZ)^\times} \bar{X}_{N,t},$$
where $\bar{X}_{N,t}$ is the compactification of $X_{N,t}$
obtained by adjoining the cusps.  We write $\bar{X}_N^\alg$
for the corresponding algebraic curve over $\CC$.

For each $t$, we define
the complex analytic surface $E_{N,t}$ to be a copy of the
quotient
$$\Gamma(N)\backslash((\uhp\times\CC)/(\ZZ\times\ZZ)),$$
where $(m,n) \in \ZZ \times \ZZ$ acts on $\uhp\times\CC$
via $(\tau,z)\mapsto (\tau,z+m\tau+n)$, and $\gamma = \smat{a}{b}{c}{d}
\in\Gamma(N)$ acts by sending the class of $(\tau,z)$ to
that of $(\gamma(\tau),(c\tau+d)^{-1}z)$.  We can regard
$E_N = \coprod E_{N,t}$ as a complex analytic family of
elliptic curves over $X_N$.  We can then extend $E_N$ to a family
$\bar{E}_N$ of generalized elliptic curves over $\bar{X}_N$
using analytic Tate curves, as in \cite[VII.4]{del_rap}.
More precisely, for $\alpha > 0$, we let $D_\alpha$
denote the disk in $D[q^{1/N}]$ defined by
$|q| < e^{-2\pi \alpha}$ and let $D_\alpha^0$ denote
the corresponding punctured disk.
Then over a punctured neighborhood around each cusp in $\bar{X}_N$,
we have an isomorphism between $E_N$ and the $N$-sided analytic Tate
curve $E_q^\an$ over $D_1^0$, and we define $\bar{E}_N$ by
gluing $N$-sided analytic Tate curves over $D_1$ to $E_N$.
For each $a$, we define a level $N$-structure on $E_{N,t}$
by the pair of sections $X_{N,t} \to E_{N,t}$ induced
by $\tau \mapsto (\tau,\tau/N)$ and $\tau \mapsto (\tau,t/N)$.
The resulting level $N$-structure on $E_N$ in fact extends
to $\bar{E}_N$.

Using for example \cite[App.\ B, 3.4]{hart}, one checks
that $\bar{E}_N$ is algebraic, so $\bar{E}_N \to \bar{X}_N$
is the analytification of a generalized elliptic
curve $\bar{E}_N^\alg \to \bar{X}_N^\alg$.  One can check
the resulting morphism $\bar{X}_N^\alg \to \bar{X}_\CC$
induces a bijection $X_N \to X(\CC)$, and since both
are smooth, proper algebraic curves over $\CC$, the
map is an isomorphism.  The analytification of the
universal generalized elliptic curve with level $N$-structure
is therefore isomorphic to $\bar{E}_N \to \bar{X}_N$ with
the level $N$-structure defined above.

\subsection{Realizations of level $N$ modular forms}
\mlabel{ssec:levelN.real}
To construct the Betti realization,
define $\CF_B$ as the locally constant sheaf $R^1s^\an_*\ZZ$ on
$X^\an$. The stalk of $\CF_B$ at a point $x = \Gamma(N)\tau \in X_{N,t}$
can be canonically identified with $H^1(E_x^\an,\ZZ) \cong
\hom(\Lambda_\tau,\ZZ)$ where $\Lambda_\tau = \ZZ\tau\oplus\ZZ$.

Let $\CF_B^k = \sym^{k-2}_{\ZZ}\CF_B$, $\CM_B = H^1(X^\an,\CF_B^k)$,
and $\CM_{c,B} = H^1_c(X^\an,\CF_B^k)$.  Note that if $\ell$ is a
prime greater than $k-2$, then multiplication by $\ell$ is injective
on $\CF_B^k$.  Since $H^0_c(X^\an,\FF_\ell \otimes \CF_B^k) = 0$,
we see that $\CM_{c,B}$ has no $\ell$-torsion if $\ell > k-2$.
If $k =2$, then the map
$$H^0(X^\an, \CF_B^k) \to H^0(X^\an,\FF_\ell \otimes \CF_B^k)$$
is surjective, and if $k > 2$ and $\ell$ does not divide $N(k-2)!$,
then this last group vanishes.  We conclude that
$\CM_{B}$ has no $\ell$-torsion if $k = 2$ or $\ell$ does
not divide $N(k-2)!$.

If $F_E$ is complex conjugation on $E^\an$ and $F_X$ is complex
conjugation on $X^\an$, then the commutative diagram
$$
\begin{CD}
E^\an @>F_E>> E^\an\\
@V s^\an VV  @VV s^\an V\\
X^\an @> F_X >> X^\an
\end{CD}
$$
gives an isomorphism
$$F_X^*\CF_B = F_X^*R^1s^\an_*\ZZ \to
    R^1s^\an_*F_E^*\ZZ = \CF_B,$$
and this induces an isomorphism $F_X^*\CF_B^k \to \CF_B^k$.
The action of complex conjugation in $G_\RR$ on $\CM_B$ is defined
as the composite
$$H^1(X^\an,\CF_B^k) \stackrel{F_X^*}{\longrightarrow}
  H^1(X^\an,F_X^*\CF_B^k) \to H^1(X^\an,\CF_B^k),$$
and the action on $\CM_{c,B}$ is defined similarly.
We let $M_B = \QQ \otimes \CM_B$ and
$M_{c,B} = \QQ \otimes\CM_{c,B}$.

For any finite prime $\ell$, we let $\CF_\ell$ denote the
$\ell$-adic sheaf $R^1s_*\ZZ_\ell$
on $X$. This is a lisse $\ell$-adic sheaf
whose stalk at a geometric point $\bar{x}:
\spec\bar{k}\to X$ can be canonically identified
with $H^1(E_{\bar{x}},\ZZ_\ell) \cong
\hom(E_{\bar{x}}[\ell^\infty],\QQ_\ell/\ZZ_\ell)$.

Let $\CF_\ell^k = \sym^{k-2}_{\ZZ_\ell}\CF_\ell$,
$\CM_\ell = H^1(X_{\bar{\QQ}},\CF_\ell^k)$ and
$\CM_{c,\ell} = H^1_c(X_{\bar{\QQ}},\CF_\ell^k)$.
Then $G_\QQ$ acts on on $\CM_\ell$ and $\CM_{c,\ell}$
by transport of structure.  We let $M_\ell = \QQ\otimes \CM_\ell$
and $M_{c,\ell} = \QQ \otimes \CM_{c,\ell}$.

The comparison theorem between Betti and $\ell$-adic
cohomology shows that $R^1 s^\an_*\ZZ_\ell$ is canonically isomorphic
to the sheaf on $X^\an$ associated to the $\ell$-adic sheaf
$R^1 s_*\ZZ_\ell$ on $X_\CC$ (see Theorem.\ I.11.6 and p.\ 130 of
\cite{frei_kiehl}).  Identifying $x \in X^\an$ with a
geometric point $\spec\CC \to X$, the isomorphism on the
corresponding stalks is the canonical one between
$\hom(E_{x}[\ell^\infty],\QQ_\ell/\ZZ_\ell)$ and
$H^1(E_x,\ZZ_\ell)$.

Taking symmetric powers and applying the
comparison theorem again gives isomorphisms
$$\ZZ_\ell\otimes  \CM_B \cong \CM_\ell\quad\mbox{and}\quad
\ZZ_\ell \otimes \CM_{c,B} \cong \CM_{c,\ell}$$
respecting the action of $G_\RR$.

The construction of the de Rham realizations is similar to the one
given in \cite{scholl2} except we use the language of
log schemes \cite{kato}.  We let $\CN_E$ denote the log structure
on $\bar{E}$ associated to $E^\infty$, and let $\CN_X$ denote the
log structure associated to $X^\infty$ \cite[1.5]{kato}.
By \cite[3.5]{kato}, the morphisms $(\bar{E},\CN_E) \to (\bar{X},\CN_X)$
and $(\bar{X},\CN_X) \to (T,\calo_T^*)$ are smooth, so by
\cite[3.12]{kato}, we have an exact sequence
of coherent locally free $\CO_{\bar{E}}$-modules
\begin{equation}
0 \to \bar{s}^*\omega^1_{\bar{X}/T}
\to \omega^1_{\bar{E}/T} \to \omega^1_{\bar{E}/\bar{X}} \to 0,
\mlabel{kato-seq}
\end{equation}
where $\omega^1$ denotes the sheaf of logarithmic relative
differentials defined in \cite[1.7]{kato}.  The sheaves
$\omega^1_{\bar{X}/T}$  and $\omega^1_{\bar{E}/\bar{X}}$ are invertible,
and can be identified, respectively, with $\Omega^1_{\bar{X}/T}(X^\infty)$
and the sheaf of regular differentials for $\bar{s}$ (denoted
$\omega_{\bar{E}/\bar{X}}$ in \cite[I.2.1]{del_rap}).

Define $\CF_{\dr}$ as the locally free sheaf
$\BR^1\bar{s}_*\omega^\bullet_{\bar{E}/\bar{X}}$
of $\CO_{\bar{X}}$-modules
on $\bar{X}$, where $\omega^\bullet_{\bar{E}/\bar{X}}$ is the
complex $d:\CO_{\bar{E}} \to \omega^1_{\bar{E}/\bar{X}}$.
This has a decreasing filtration with $\fil^2\CF_{\dr} = 0$,
$\fil^1\CF_{\dr} = \bar{s}_*\omega^1_{\bar{E}/\bar{X}}$,
and $\fil^0\CF_{\dr} = \CF_{\dr}$.  We denote $\fil^1\CF_{\dr}$
simply as $\omega$.  The sheaf $\omega$ is invertible and can
be identified with $e^*\Omega^1_{\bar{E}/\bar{X}}$ where
$e:\bar{X}\to \bar{E}$ is the zero section.  We also have a
canonical isomorphism
$\gr^0\CF_{\dr} \cong R^1\bar{s}_*\CO_{\bar{E}} \cong \omega^{-1}$
given by Grothendieck-Serre duality.  The fiber $x^*\CF_\dr$
of $\CF_\dr$ at a point $x:\spec k \to \bar{X}$ can be identified
with $\BH^1(\bar{E}_x,\omega^\bullet_{\bar{E}_x/k})$ where
$\omega^\bullet_{\bar{E}_x/k}$ is the complex
$d:\CO_{\bar{E}_x} \to \omega^1_{\bar{E}_x/k}$
and $\omega^1_{\bar{E}_x/k}$ is the sheaf of
regular differentials on $\bar{E}_x$.  The fiber of
$\fil^1\CF_\dr$ is then identified with
$H^0(\bar{E}_x, \omega^1_{\bar{E}_x/k})
\cong e_x^*\Omega^1_{\bar{E}_x/k}$ where
$e_x:\spec k \to \bar{E}_x$ is the zero section,
and that of $\gr^0\CF_\dr$ is identified with
$H^1(\bar{E}_x, \CO_{\bar{E}_x})$.

We define $\CF^k_{\dr}$ as
the filtered sheaf of $\CO_{\bar{X}}$-modules
$\sym^{k-2}_{\CO_{\bar{X}}}\CF_{\dr}$,
and we
let $\CF_{c,\dr}^k = \CF_\dr^k(-X^\infty)$.

The exact sequence (\ref{kato-seq})
gives rise to an exact sequence of complexes
\begin{equation}
0 \to \omega^\bullet_{\bar{E}/\bar{X}}\otimes_{\CO_{\bar{E}}}
\bar{s}^*\omega^1_{\bar{X}/T}
    \to \omega^\bullet_{\bar{E}/T} \to \omega^\bullet_{\bar{E}/\bar{X}}\to
0,
\mlabel{kato-seq2}
\end{equation}
where $\bar{s}^*\omega^1_{\bar{X}/T}$ has degree 1, and
$\omega^\bullet_{\bar{E}/T}$ denotes the complex
$\wedge^\bullet_{\CO_{\bar{E}}}\omega^1_{\bar{E}/T}$ \cite[1.9]{kato}.
{}From the long exact sequence for $\BR^i\bar{s}_*$ and the projection
formula, we obtain the (logarithmic) Gauss-Manin connection
$$\nabla:\CF_{\dr} \to
\CF_{\dr}\otimes_{\CO_{\bar{X}}}\omega^1_{\bar{X}/T}.$$
This induces logarithmic connections on $\CF_{\dr}^k$ and
$\CF_{c,\dr}^k$ satisfying Griffiths transversality.
We set $\CM_{\dr} = \BH^1(\bar{X},\omega^\bullet(\CF_{\dr}^k))$ and
$\CM_{c,\dr} = \BH^1(\bar{X},\omega^\bullet(\CF_{c,\dr}^k))$,
where we write $\omega^\bullet(\calg)$ for the complex
associated to the module $\calg$ with its connection.
The filtrations on $\CM_{\dr}$ and $\CM_{c,\dr}$
are defined by those on $\CF_{\dr}^k$ and $\CF_{c,\dr}^k$.
We let $M_\dr = \QQ\otimes \CM_\dr $ and
$M_{c,\dr} = \QQ\otimes \CM_{c,\dr}$.

Recall that the composite
$$\omega  \to \CF_\dr \stackrel{\nabla}{\to}
 \CF_{\dr}\otimes_{\CO_{\bar{X}}}\omega^1_{\bar{X}/T}
  \to \omega^{-1}\otimes_{\CO_{\bar{X}}}\omega^1_{\bar{X}/T}$$
is an isomorphism (the Kodaira-Spencer isomorphism, \cite[VI.4.5.2]{del_rap}).
Hence $\omega^2\cong \omega^1_{\bar{X}/T}$, and
one deduces that
$$\gr^i \CM_{\dr}  \cong
\left\{\begin{array}{ll}
H^0(\bar{X},\omega^{k-2}\otimes \omega^1_{\bar{X}/T}),&
\mbox{if $i=k-1$;}\\
H^1(\bar{X},\omega^{2-k}),&\mbox{if $i=0$;}\\
0,&\mbox{otherwise}.\end{array}\right.$$
Similarly one finds
$$\gr^i \CM_{c,\dr}  \cong
\left\{\begin{array}{ll}
H^0(\bar{X},\omega^{k-2}\otimes \Omega^1_{\bar{X}/T}),&
\mbox{if $i=k-1$;}\\
H^1(\bar{X},\omega^{2-k}(-X^\infty)),&\mbox{if $i=0$;}\\
0,&\mbox{otherwise}.\end{array}\right.$$
In particular, it follows that $\CM_{\dr}$ and $\CM_{c,\dr}$
are torsion-free.

For lack of a reference, we now sketch the construction
of the comparison isomorphism $I^\infty$.

For $\bar{Y} = \bar{X}$ or $\bar{E}$, let $h_{\bar{Y}}$ denote the
canonical morphism $\bar{Y}^\an \to \bar{Y}$ of ringed spaces, let
$\omega^1_{\bar{Y}^\an}$ denote $h_{\bar{Y}}^*\omega^1_{\bar{Y}/S}$,
and let $\omega^\bullet_{\bar{Y}^\an}$ denote the complex
$\wedge^\bullet_{\CO_{\bar{Y}^\an}}\omega^1_{\bar{Y}^\an}$.
We also consider the complex
$$\omega^\bullet_{\bar{E}^\an/\bar{X}^\an} =
\wedge^\bullet_{\CO_{\bar{E}^\an}}\omega^1_{\bar{E}^\an/\bar{X}^\an},$$
where $\omega^1_{\bar{E}^\an/\bar{X}^\an}$ denotes
$h_{\bar{E}}^*\omega^1_{\bar{E}/\bar{X}}$.

We let
$$\CF_\dr^\an =
\BR^1\bar{s}^\an_*\omega^\bullet_{\bar{E}^\an/\bar{X}^\an},$$
$\CF_\dr^{\an,k} = \sym^{k-2}_{\CO_{\bar{X}^\an}}\CF_\dr^\an$
and $\CF_{c,\dr}^{\an,k} = \CF_\dr^{\an,k}(-X^{\infty,\an})$.
These are locally free sheaves of $\CO_{\bar{X}^\an}$-modules,
and the fiber $\CC\otimes_{\CO_{\bar{X},x}^\an}\CF_{\dr,x}^\an$
of $\CF_\dr^\an$ at a point $x\in X^\an$ can be
identified with the analytic de Rham cohomology
$\BH^1(E_x^\an,\Omega^\bullet_{E_x^\an})$ of $E_x^\an$.
By (\ref{kato-seq}) we have an exact sequence
$$0 \to \bar{s}^{\an,*}\omega^1_{\bar{X}^\an}
\to \omega^1_{\bar{E}^\an} \to
\omega^1_{\bar{E}^\an/\bar{X}^\an} \to 0,$$
and a construction like the one in the de Rham realization
gives logarithmic connections
$$\CG \to \CG \otimes_{\CO_{\bar{X}^\an}} \omega^1_{\bar{X}^\an}$$
and complexes $\omega^\bullet(\CG)$
for $\CG = \CF_\dr^{\an,k}$ and $\CF_{c,\dr}^{\an,k}$.
We let $M_\dr^\an = \BH^1(\bar{X}^\an,\omega^\bullet(\CF_\dr^{\an,k}))$
and $M_{c,\dr}^\an =
\BH^1(\bar{X}^\an,\omega^\bullet(\CF_{\dr,c}^{\an,k}))$.

Now we compare $M_\dr^\an$ with $\CC\otimes M_\dr$, and
$M_{c,\dr}^\an$ with $\CC\otimes M_{c,\dr}$.  The natural map
$$h_{\bar{E}}^{-1}\omega^\bullet_{\bar{E}/\bar{X}} \to
    \omega^\bullet_{\bar{E}^\an/\bar{X}^\an}$$
is compatible with differentiation, and so gives rise to
an $h_{\bar{X}}^{-1}\CO_{\bar{X}}$-linear map
$$h_{\bar{X}}^{-1}\CF_\dr \to
\BR^1\bar{s}^\an_*h_{\bar{E}}^{-1}\omega^\bullet_{\bar{E}/\bar{X}}
\to \CF_\dr^\an.$$
This map is compatible with connections and induces an isomorphism
$$h_{\bar{X}}^*\CF_\dr \to \CF_\dr^\an$$
which on the fiber at $x\in X(\CC) = X^\an$ becomes the
isomorphism between the algebraic de Rham cohomology of
$E_x$ and the analytic de Rham cohomology of $E_x^\an$
provided by GAGA~\cite{gaga}.
{}From this one obtains $h_{\bar{X}}^{-1}\CO_{\bar{X}}$-linear maps
$$h_{\bar{X}}^{-1}\CF_\dr^k \to \CF_\dr^{k,\an},
\quad\mbox{and}\quad
h_{\bar{X}}^{-1}\CF_{c,\dr}^k \to \CF_{c,\dr}^{k,\an}$$
compatible with connections and inducing isomorphisms on
replacing $h_{\bar{X}}^{-1}$ with $h_{\bar{X}}^*$.
We thus obtain maps
$$\begin{array}{lclcl} M_\dr & \to & \BH^1(\bar{X}^\an, h_{\bar{X}}^{-1}
\omega^\bullet(\CF_\dr^k)) & \to & M_\dr^\an\quad\mbox{and}\\
M_{c,\dr} & \to & \BH^1(\bar{X}^\an, h_{\bar{X}}^{-1}
\omega^\bullet(\CF_{c,\dr}^k)) & \to & M_{c,\dr}^\an.\end{array}$$
Applying GAGA \cite{gaga} again, we find that these induce isomorphisms
$$\CC\otimes M_\dr\cong M_\dr^\an\quad\mbox{and}\quad
\CC\otimes M_{c,\dr}\cong M_{c,\dr}^\an.$$

Now we compare $M_\dr^\an$ with $\CC\otimes M_B$ and
$M_{c,\dr}^\an$ with $\CC\otimes M_{c,\dr}$.
By the Poincar\'e Lemma, we have the resolution
$$\CC \to \Omega^\bullet_{E^\an}$$
of the constant sheaf $\CC$ on $E^\an$.
We thus obtain an isomorphism
$$\CC \otimes \CF_B \cong  R^1s^\an_*\CC \cong
\BR^1s^\an_*\Omega^\bullet_{E^\an}.$$
Applying $\BR^\bullet s^\an_*$ to the exact sequence
$$0 \to \Omega^\bullet_{E^\an/X^\an} \otimes_{\CO_{E^\an}} s^{\an,*}\Omega^1_{X^\an}
\to \Omega^\bullet_{E^\an} \to
\Omega^\bullet_{E^\an/X^\an} \to 0$$
we obtain an exact sequence
$$\CO_{X^\an} \to \Omega^1_{X^\an} \to
\BR^1s^\an_*\Omega^\bullet_{E^\an}
 \to \CG
 \to \CG \otimes_{\CO_{X^\an}}\Omega^1_{X^\an},$$
where $\CG$ is the restriction of $\CF_\dr^\an$ to $X^\an$.
The first map (differentiation) is surjective, as is the last
map (the Gauss-Manin connection) by the local solvability of
linear systems of differential equations.  We thus obtain
an exact sequence
\begin{equation} 0 \to \CC\otimes \CF_B
\to \CG
\to \CG \otimes_{\CO_{X^\an}}\Omega^1_{X^\an} \to
0.\mlabel{diff_equ}
\end{equation}
Composing the first map on the stalk at $x\in X^\an$ with the
map arising from the evaluation homomorphism $\CO_{X^\an,x} \to \CC$,
we obtain the isomorphism $H^1(E_x^\an,\CC) \cong
\BH^1(E_x^\an,\Omega^\bullet_{E_x^\an})$ given by the Poincar\'e
lemma for $E_x^\an$.

Writing $\CG^k$ for the restriction of $\CF_\dr^{k,\an}$ to $X^\an$,
we get a sequence
\begin{equation} 0 \to \CC\otimes \CF_B^k
\to \CG^k
\to \CG^k \otimes_{\CO_{X^\an}}\Omega^1_{X^\an}
\to 0
\mlabel{gk-res}
\end{equation}
whose exactness follows from the local existence and uniqueness of
solutions with initial conditions.  Extending the connection on $\CG^k$
to our logarithmic connection on $\CF_{c,\dr}^{k,\an}$ and letting $j$
denote the inclusion of $X^\an$ in $\bar{X}^\an$, we claim that
$$ 0 \to \CC\otimes j_!\CF_B^k\to \CF_{c,\dr}^{k,\an}
\to \CF_{c,\dr}^{k,\an} \otimes_{\CO_{\bar{X}^\an}}\omega^1_{\bar{X}^\an}
\to 0$$
is exact.  This amounts to the following:
\begin{lemma} The connection
$$\CF_{c,\dr}^{k,\an}
\to \CF_{c,\dr}^{k,\an} \otimes_{\CO_{\bar{X}^\an}}\omega^1_{\bar{X}^\an}$$
defines an isomorphism on stalks at the cusps.
\end{lemma}
\proofbegin  Over a neighborhood of each cusp in $\bar{X}^\an$,
$\bar{E}^\an$ is isomorphic to the $N$-sided analytic Tate
curve $E_q^\an$ over $D_1$.  We claim that there is an isomorphism
between the pull-back of $\CF_{\dr}^{k,\an}$ to $D_1$
and $\CO_{D_1}^{k-1}$ such that the corresponding
connection on $\CO_{D_1}^{k-1}$ is given by
$\nabla(f) = (A(t)f + tf')\dfrac{dt}{t}$ where $A(t)\in
M_{k-1}(\CO_{D_1})$ is such that $A(0)$ is nilpotent.
Indeed this follows from the case $k=3$ which follows
from the calculation in Appendix 1.3 of \cite{katz}
for the 1-sided Tate curve.

Now observe that the desired bijectivity follows from
the local existence of a unique solution of the differential
equation $(A(t)+I)f + tf' = g$ for any $g \in \CO_{D_1}^{k-1}$,
which is guaranteed by the fact that $A(0) + I$ has no nonpositive
integer eigenvalues.
\epf

We can now conclude that
$\CC\otimes M_{c,B} \cong H^1(\bar{X}^\an,\CC\otimes j_!\CF_B^k)$
is isomorphic to $M_{c,\dr}^\an \cong \CC \otimes M_{c,\dr}$.

To define $I^\infty$ for $M_B$, note that we have a restriction map
$$\CC \otimes M_\dr\cong M_\dr^\an
\to\BH^1(X^\an,\omega^\bullet(\CG^k)) \cong \CC \otimes M_B.$$
We postpone until \S \ref{sssec:pairings}
the proof that it is an isomorphism.
It is straightforward to check that the maps respect the action
of $G_\RR$.

Before defining the crystalline realization and comparison
isomorphisms, we review Faltings' generalization
of the theory of Fontaine and Laffaille (see
\cite{faltings}).

Suppose that $\bar{Y}$ is a smooth, proper scheme over $\spec \ZZ_\ell$
with a relative divisor $D$ with strict normal crossings, and let
$Y = \bar{Y} - D$.  For each integer $a$ with $0 \le a \le \ell-2$, Faltings
defines a category $\MF_{[0,a]}^\nabla(Y)$
and a fully faithful contravariant functor $\DD$ to the category
of finite locally constant \'etale sheaves on $Y_{\QQ_\ell}$
(\cite[Thm.\ $2.6^*$]{faltings}).  We write $\VV$ for the covariant
functor defined by
$$\VV(\CA) = \DD(\CA)^* = \hom(\DD(\CA),\QQ_\ell/\ZZ_\ell).$$
In the case of  $Y = \bar{Y} = \spec \CO_F$ where $F$
is an unramified extension of $\QQ_\ell$, Faltings' category
can be identified with the full subcategory of
$\mfcat^0_\tor (\CO_F)$
whose objects $A$ satisfy $\fil^{a+1}A = 0$.
Furthermore, identifying the category of
\'etale sheaves on $\spec F$ with that of continuous
$G_F$-modules, one finds that the functor $\DD$ agrees
with $\underline{U}_S$ (see \cite[3.8]{fon_ag}).

We let $\CO_{Y,\crys}$ denote the inverse system in
$\MF_{[0,0]}^\nabla(Y)$ (indexed by integers $n \ge 1$)
defined by reduction mod $\ell^n$ of the following data:
\begin{itemize}
\item the structure sheaf $\CO_{\bar{Y}}$;
\item the logarithmic connection
  $d:\CO_{\bar{Y}} \to \omega^1_{\bar{Y}/\ZZ_\ell}$;
\item an affine covering $\{\spec A_\alpha\}$ of $\bar{Y}$
with arithmetic Frobenius-lifts $\phi_\alpha$ of $A_\alpha$
as in \cite[\S2i]{faltings}.
\end{itemize}
Then $\VV(\CO_{Y,\crys})$ is canonically identified with the constant
$\ell$-adic sheaf $\ZZ_\ell$ on $Y$.

Assuming $\ell$ does not divide $2N$
and applying \cite[Thm.\ 6.2]{faltings} to the morphism $\bar{s}$,
we obtain a canonical isomorphism
$$\VV(\CF_\ecrys) \cong \CF_\ell,$$
where $\CF_\ecrys$ is the inverse system in
$\MF^\nabla_{[0,1]}(X_{\ZZ_\ell})$
defined by reduction mod $\ell^n$ of $\CF_\dr$
with its filtration, logarithmic Gauss-Manin connection and locally
defined Frobenius maps whose precise description we shall not need.
For  $x:\spec\CO_F \to X$ with
$F$ an unramified extension of $\QQ_\ell$, we can
identify $x^*\CF_\ecrys$ with the object $H^1_\crys(E_x,\CO_{E_x,\crys})$
of $\MF_{[0,1]}(\CO_F)$.
Furthermore $\VV(x^*\CF_\ecrys)$ is canonically identified with
the stalk of $\VV(\CF_\ecrys)$ at $\bar{x}:\bar{F} \to X$ and
the isomorphism $\VV(\CF_\ecrys) \cong \CF_\ell$ on the stalk
is the one between $H^1_\crys(E_x,\CO_{E_x,\ecrys})$
and $H^1(E_{\bar{x}},\ZZ_\ell)$ provided by
\cite[Thm.\ 5.3]{faltings}.

Using the results in \cite[IIh)]{faltings}, one finds that the functor
$\VV$ respects tensor products.  So assuming further that $\ell > k-1$,
we obtain a canonical isomorphism
$$\VV(\CF_\ecrys^k) \cong \CF_\ell^k,$$
where $\CF_\ecrys^k$ is the inverse system
in $\MF^\nabla_{[0,k-2]}(X_{\ZZ_\ell})$ defined by
$\sym^{k-2}_{\CO_{\bar{Y}}}\CF_\ecrys$.
If $\ell > k$, we obtain an object
$$\CM_\ecrys = H^1_\crys(X_{\ZZ_\ell},\CF_\ecrys^k)$$
of $\mfcat^0$ whose underlying filtered module is isomorphic
to $\ZZ_\ell\otimes\CM_\dr$, and so we have the comparison $I^\ell_\dr$.
Furthermore we conclude from \cite[Thm.\ 5.3]{faltings}
that there is a canonical isomorphism
$$I^\ell:\VV(\CM_\ecrys) \cong \CM_\ell.$$
For such primes, we similarly obtain an isomorphism
$$I^\ell:\VV(\CM_{c,\ecrys}) \cong \CM_{c,\ell}$$
where $\CM_{c,\ecrys}$ has underlying filtered module
$\CM_{c,\dr}$.

\subsection{Level $N$ modular forms} \mlabel{sssec:mot.mf.rel}
\mlabel{ssec:ad.int}

We recall the relation between
$M_{\dr}$ and the space of modular forms of weight $k$ with respect
to $\Gamma(N)$.

Let $\CL$ denote the line bundle $\omega^{k-2}\otimes
\omega^1_{\bar{X}/S}$ on $\bar{X}$.  Then we have
$$\CC\otimes \fil^{k-1}\CM_{\dr} \cong
H^0(\bar{X}^\an, \CL^\an)\cong
\bigoplus_{t\in(\ZZ/N\ZZ)^\times} H^0(\bar{X}_{N,t},i_t^*\CL^\an),$$
where
$$\CL^\an = \left(e^{\an,*}\Omega^1_{\bar{E}^\an/\bar{X}^\an}
\right)^{k-2}\otimes_{\CO_{\bar{X}^\an}}
\Omega^1_{\bar{X}^\an}(-X^{\infty,\an})$$
and $i_t$ is the inclusion $\bar{X}_{N,t} \to \bar{X}$
defined in \S\ref{sssec:modcurves}.
For each $t$, we have the natural projection $\pi_t: \uhp \to
\bar{X}_{N,t}$ giving an inclusion
$$H^0(\bar{X}_{N,t},i_t^*\CL^\an) \to
    H^0(\uhp , \pi_t^*i_t^*\CL^\an)$$
whose image is precisely the space of sections of the form
$$f(\tau)(2\pi i)^{k-1}(dz)^{\otimes(k-2)}\otimes d\tau$$
with $f(\tau)$ in $M_k(\Gamma(N))$, the space of modular forms
of weight $k$ with respect to $\Gamma(N)$.  We thus obtain
an isomorphism
$$\CC\otimes \fil^{k-1}\CM_{\dr} \cong
\bigoplus_{t\in(\ZZ/N\ZZ)^\times} M_k(\Gamma(N)).$$
Note that we have chosen our normalization so that
if $f(\tau) = g(e^{2\pi i \tau/N})$ is in $M_k(\Gamma(N))$,
then the restriction of the corresponding element of
$H^0(\bar{X}_{N,t},i_t^*\CL^\an)$ to $H^0(D_1 , j_t^*i_t^*\CL^\an)$ is
$$g(q^{1/N})\omega_\can^{\otimes(k-2)}\otimes\dfrac{dq}{q},$$
where $j_t:D_1 \to \bar{X}_{N,t}$ is the inclusion of a disk
around the cusp $\infty$ used in \S\ref{sssec:modcurves}.

By the $q$-expansion principle~\cite[VII.3.9,VII.3.13]{del_rap},
the map
$$\CC\otimes \fil^{k-1}\CM_\dr \cong \bigoplus_{t\in
(\ZZ/N\ZZ)^\times} M_k(\Gamma(N)) \to \bigoplus_{t\in
(\ZZ/N\ZZ)^\times} \CC[[q^{1/N}]],$$
where the second map is given by sending $f(\tau)$ to $g(q^{1/N})$,
gives rise to a pullback diagram
\begin{equation}
\begin{CD}
\fil^{k-1}\CM_\dr @>>> \CC\otimes \fil^{k-1}\CM_\dr
    @>\cong >> \bigoplus_{t\in (\ZZ/N\ZZ)^\times} M_k(\Gamma(N))\\
@VVV @VVV @VVV\\
R[[q^{1/N}]] @>>>
 (R\otimes\CC)[[q^{1/N}]]
 @>\cong >> \bigoplus_{t\in (\ZZ/N\ZZ)^\times} \CC[[q^{1/N}]]
\end{CD}
\mlabel{eq:qexp}
\end{equation}
where $R = \ZZ[1/Nk!,\mu_N]$.  Here the bottom right isomorphism
is induced by the isomorphism $R\otimes\CC \to \oplus_t\CC$
defined by $(e^{2\pi i t/N})_t$ and commutativity of
the right square follows from~\cite[VII.4.7]{del_rap}.

Similarly one obtains an isomorphism
$$\CC\otimes \fil^{k-1}\CM_{c,\dr} \cong
\bigoplus_{t\in(\ZZ/N\ZZ)^\times} S_k(\Gamma(N)),$$
where $S_k(\Gamma(N))$ is the space of cusp forms
of weight $k$ with respect to $\Gamma(N)$.
The $q$-expansion principle identifies $\fil^{k-1}\CM_{c,\dr}$
as the subset of $\bigoplus_{t\in (\ZZ/N\ZZ)^\times}
S_k(\Gamma(N))$ whose $q$-expansion at $\infty$ has coefficients
in $R$.

\subsection{Pairings}
\mlabel{sssec:pairings}
We shall use Poincar\'e duality to construct perfect
pairings $M_{c,\wc} \otimes M_\wc \to \QQ(1-k)_\wc$
for $\wc = B$, $\dr$ and $\ell\nmid Nk!$, compatible with the
comparison isomorphisms.
These pairings arise from a {\em duality morphism}
$\CM_{c} \to \hom(\CM, \calt^{\otimes (k-1)})$.

We begin with the pairing on Betti realizations.
The cup product defines a morphism
$$\CF_B \otimes \CF_B \to R^2s^\an_*\ZZ\cong (2\pi i)^{-1}\ZZ$$
inducing an isomorphism
$$\CF_B \to \shom(\CF_B,(2\pi i)^{-1}\ZZ)$$
of sheaves on $X^\an$.
We then consider the morphism
$$\CF_B^k \otimes \CF_B^k \to (2\pi i)^{2-k}\ZZ$$
defined on sections by
\begin{equation}
x_1\otimes \cdots \otimes x_{k-2} \otimes y_1 \otimes \cdots \otimes
y_{k-2} \mapsto
\sum_{\sigma \in \Sigma_{k-2}} \prod_{i=1}^{k-2}x_i \cup y_{\sigma(i)},
\mlabel{sympa}
\end{equation}
where $\Sigma_{k-2}$ is the symmetric group on
$\{1,\ldots,k-2\}$.  The resulting map on cohomology,
preceded by the cup product and followed by the trace, yields a pairing
$$\CM_{c,B} \otimes \CM_B \to H^2_c(X^\an,\CF_B^k\otimes\CF_B^k)
\to H^2_c(X^\an,(2\pi i)^{2-k}\ZZ) \to (2\pi i)^{1-k}\ZZ = \calt_B^{\otimes(k-1)},$$
respecting the action of $G_\RR$.
After tensoring with $\ZZ[1/(k-2)!]$, we have that
$\CM_{c,B}$ is torsion-free and
$$\CF_B^k \to \shom(\CF_B^k,(2\pi i)^{2-k}\ZZ)$$
becomes an isomorphism, and it follows that
$$\delta^L_B: \CM_{c,B} \to \hom(\CM_B,\calt_B^{\otimes(k-1)})$$
also becomes an isomorphism.

The pairing on $\ell$-adic realizations is constructed
similarly.  In particular we have morphisms
$$\CF_\ell^k \to \shom(\CF_\ell^k,\ZZ_\ell(2-k))
\quad\mbox{and}\quad
\delta^L_\ell: \CM_{c,\ell} \to \hom(\CM_\ell,\calt_\ell^{\otimes(k-1)})$$
which are isomorphisms if $\ell > k-2$, or for any $\ell$
after tensoring with $\QQ_\ell$.

To define the pairing on de Rham realizations, we
begin with the cup product
$$ \BR^1\bar{s}_*\omega^\bullet_{\bar{E}/\bar{X}}
\otimes_{\CO_{\bar{X}}} \BR^1\bar{s}_*\omega^\bullet_{\bar{E}/\bar{X}}
\to \BR^2\bar{s}_*\left(\omega^\bullet_{\bar{E}/\bar{X}}
\otimes_{\bar{s}^{-1}\CO_{\bar{X}}}\omega^\bullet_{\bar{E}/\bar{X}}
\right).$$
Combining this with the morphism induced by the wedge product on
$\omega^\bullet_{\bar{E}/\bar{X}}$,
we obtain a morphism
$$\CF_\dr \otimes_{\CO_{\bar{X}}} \CF_\dr \to
\BR^2\bar{s}_*\omega^\bullet_{\bar{E}/\bar{X}}
\cong R^1\bar{s}_*\omega^\reg_{\bar{E}/\bar{X}}[-1]
\cong \CO_{\bar{X}}[-1]$$
respecting filtrations and connections.
On graded pieces, the morphism
\begin{equation}
\CF_\dr \to \shom_{\CO_{\bar{X}}}(\CF_\dr,\CO_{\bar{X}}[-1])
\mlabel{fdr-duality}
\end{equation}
coincides with the isomorphism given by Grothendieck-Serre duality.
Using the same formula on sections as for the Betti realization,
we obtain a pairing
$$\CF_\dr^k \otimes_{\CO_{\bar{X}}} \CF_\dr^k \to
\CO_{\bar{X}}[2-k]$$
respecting filtrations and connections and inducing an isomorphism
$$\CF_\dr^k \to \shom_{\CO_{\bar{X}}}(\CF_\dr^k,\CO_{\bar{X}}[2-k]).$$
The pairing on $\CF_\dr^k$ gives rise to a morphism of
filtered complexes
$$\omega^\bullet(\CF_{c,\dr}^k) \otimes
\omega^\bullet(\CF_\dr^k) \to
\omega^\bullet(\CF_{c,\dr}^k\otimes_{\CO_{\bar{X}}}\CF_\dr^k)
\to \omega^\bullet({\cal{I}}_{X^\infty})[2-k].$$
Composing the resulting map on cohomology with the cup product, we obtain
$$\CM_{c,\dr} \otimes \CM_\dr \to
\BH^2(\bar{X},\omega^\bullet(\CF_{c,\dr}^k)\otimes\omega^\bullet(\CF_\dr^k))
\to \BH^2(\bar{X},\omega^\bullet({\cal{I}}_{X^\infty})[2-k]).$$
Since differentiation induces the zero map
$H^1(\bar{X},{\cal{I}}_{X^\infty})\to H^1(\bar{X},\Omega^1_{\bar{X}/S})$,
we conclude that
$$\BH^2(\bar{X},\omega^\bullet({\cal{I}}_{X^\infty})[2-k])
\cong H^1(\bar{X},\Omega^1_{\bar{X}/S})[1-k].$$
Composing with the trace map, we obtain a pairing
$$\CM_{c,\dr} \otimes \CM_\dr \to
H^1(\bar{X},\Omega^1_{\bar{X}/T})[1-k]
\to \calt_\dr^{S,\otimes(k-1)},$$
respecting filtrations and inducing an isomorphism
$$\delta^L_\dr: \CM_{c,\dr} \to \hom(\CM_\dr,\calt_\dr^{S,\otimes(k-1)}).$$

Finally, for $\ell\not\in S$,  one checks that
(\ref{fdr-duality}), pulled back to $\bar{X}_{\ZZ_\ell}$,
arises from an isomorphism
\begin{equation}
\CF_\ecrys \to \CF_\ecrys^t = \shom_{\CO_{\bar{X}_{\ZZ_\ell}}}
(\CF_\ecrys,\CO_{\bar{X}_{\ZZ_\ell}}[-1])
\mlabel{eq:fcrys}
\end{equation}
in the category $\mfcat^\nabla_{[0,1]}(X_{\ZZ_\ell})$,
giving an isomorphism
$$ \CF_\ecrys^k \to
\shom_{\bar{X}_{\ZZ_\ell}}(\CF_\ecrys^k,\CO_{\bar{X}_{\ZZ_\ell}}[2-k])$$
in $\mfcat^\nabla_{[0,k-1]}$ and ultimately an isomorphism
$$\delta^L_\ecrys: \CM_{c,\ecrys} \to
\CM_\ecrys^t=\Hom_{\ZZ_\ell}(\CM_\ecrys,\calt_\ecrys^{\otimes (k-1)}).$$

Now one has to check that the pairings respect the
comparison isomorphisms.  This is straightforward
for the maps $I^\ell_B$, $I^\infty$ and $I^\ell_\dr$
using the compatibility
of the comparison isomorphisms with cup products and
Poincar\'e duality.  We find that the same holds for $I^\ell$ for
$\ell \nmid Nk!$,
but we must take some care with the integral
formulation since the filtration length of $M_{c,\ecrys}
\otimes M_{\ecrys}$ is $2(k-1)$, which may be greater
than $\ell-2$.  We proceed as follows.
The isomorphism (\ref{eq:fcrys}) is compatible with
$\CF_\ell \to \CF_\ell^t = \shom_{\ZZ_\ell}(\CF_\ell,\ZZ_\ell(-1))$
in the sense that the diagram
$$\begin{array}{ccc}\VV(\CF_\ecrys) & \to & \VV(\CF_\ecrys^t)\\
\downarrow&&\downarrow\\
\CF_\ell & \to & \CF_\ell^t\end{array}$$
commutes, where the vertical arrows are given by
Faltings' comparison isomorphisms (and on the right,
we are also using the canonical isomorphism between
$\VV(\CF_\ecrys^t)$ and $\VV(\CF_\ecrys)^t$ given by
the discussion in \cite[IIh)]{faltings}).
We deduce
a similar compatibility for the dualities on
$\CF_\lcrys^k$ and $\CF_\ell^k$, and eventually a commutative diagram
$$\begin{array}{ccc}\VV(\CM_{c,\ecrys}) & \to & \VV(\CM_\ecrys^t)\\
\downarrow&&\downarrow\\
\CM_{c,\ell} & \to & \CM_\ell^t,\end{array}$$
where $\CM_\ell^t = \hom_{\ZZ_\ell}(\CM_\ell,\calt_\ell^{\otimes(k-1)})$.

Recall that we have shown $\CC\otimes M_{c,B} \to
\CC\otimes M_{c,\dr}$ is an isomorphism.  The same now
follows for $\CC\otimes M_B \to \CC\otimes M_\dr$
from the perfectness of the pairings and their
compatibility with these comparisons.  We use their
inverses to define the comparison isomorphisms $I^\infty$.

Finally, we need the relation with the Petersson inner product.
Recall that the Petersson inner product for
$g \in S_k(\Gamma(N))$, $h \in M_k(\Gamma(N))$ is defined by
$$ (g,h)_{\Gamma(N)} = (-2i)^{-1}
\int_{\Gamma(N)\backslash\uhp}g(\tau)\overline{h(\tau)}(\Im
\tau)^{k-2} d\tau\wedge d\bar{\tau}.$$

For $\wc = B$, $\dr$ or $\ell$, $g\in M_{c,\wc}$ and $h\in M_\wc$,
we write $\rlangle g,h\rrangle_\wc$ for the image of $g\otimes h$
in $\QQ(1-k)_\wc$.  For $g \in \CC\otimes \fil^{k-1}M_{c,\dr}$ and
$t\in (\ZZ/N\ZZ)^\times$, write $g_t$ for the image of
$g$ in the corresponding component of $S_k(\Gamma(N))$,
and similarly for $M_\dr$ and $M_k(\Gamma(N))$.

\begin{lemma} If $g\in\CC\otimes\fil^{k-1}M_{c,\dr}$ and $h
\in \CC\otimes\fil^{k-1}M_\dr$, then we have
$$\rlangle g, (I^\infty)^{-1} (F_\infty\otimes 1)I^\infty h\rrangle_{\dr}
=(k-2)!(4\pi)^{k-1}\sum_{t\in(\ZZ/N\ZZ)^\times}
(g_t,h_t)_{\Gamma(N)}\otimes\iota^{k-1}$$
in $\CC\otimes \QQ(1-k)_\dr$.
\mlabel{petersson1}\end{lemma}
\proofbegin  Let $\bar{\omega}^\bullet_{\bar{E}^\an/\bar{X}^\an}$
denote the complex $F_{\bar{E},*}\omega^\bullet_{\bar{E}^\an/\bar{X}^\an}$,
and let $\bar{\CF}_\dr^\an$ denote the sheaf
$F_{\bar{X},*}\CF_\dr^\an\cong
\BR^1 s_*\bar{\omega}^\bullet_{E^\an/X^\an}$ of
modules over $\bar{\CO}_{\bar{X}^\an} = F_{\bar{X},*}\CO_{\bar{X}^\an}$.
We then have an antiholomorphic resolution of
$\CC\otimes\CF_B^k$ analogous to (\ref{gk-res}),
and complex conjugation of functions gives rise
to an isomorphism $\omega^\bullet_{\bar{E}^\an/\bar{X}^\an}
\cong \bar{\omega}^\bullet_{\bar{E}^\an/\bar{X}^\an}$,
which induces an isomorphism
$$\omega^\bullet(\CF_\dr^{k,\an}) \cong
\bar{\omega}^\bullet(\bar{\CF}_\dr^{k,\an})
= F_{\bar{X},*}\omega^\bullet(\CF_\dr^{k,\an})$$
compatible with $F_\infty \otimes 1$ on $\CC \otimes \CF_B$.

Recall that $g$ corresponds to the class in $\BH^1(\bar{X}^\an,
\omega^\bullet(\CF_{c,\dr}^{k,\an}))$ that arises from
the section of $H^0(\bar{X},\FC^1_{\bar{X}}\otimes_{\CO_{\bar{X}^\an}}
\CF_{c,\dr}^{k,\an})$ whose pull-back by $\pi_t:\uhp \to X_{N,t}$
is defined by
$$\tau \mapsto 2\pi i g_t(\tau)d\tau \otimes
\pi_t^*\omega_\can^{\otimes(k-2)}.$$
Here $\FC^n$ is the sheaf of smooth $\CC$-valued $n$-forms and
$\omega_\can$ is the canonical section of $\CF_{\dr}^\an$
(i.e., $dx/x$ on the Tate curve $\GG_m/ q^\ZZ$).
Similarly $(I^\infty)^{-1} (F_\infty\otimes 1)I^\infty h$ arise from
the section of $H^0(\bar{X},\FC^1_{\bar{X}}\otimes_{\bar{\CO}_{\bar{X}^\an}}
\bar{\CF}_{\dr}^{k,\an})$ whose pull-back is
$$\tau \mapsto - 2\pi i \bar{h}_t (\tau) d(\bar{\tau})
 \otimes \pi_t^*(\bar{\omega}_\can)^{\otimes(k-2)}.$$
Note that $\omega^\can$ and $\bar{\omega}^\can$ both
define sections of $\FC^0_{\bar{X}}\otimes\CF_B$, and
under the pairing $\CF_B \otimes \CF_B \to 2\pi i^{-1}\QQ$
(pulled back to $\uhp$ and tensored with $\FC^0_\uhp$), we have
$$\omega_\can \otimes \bar{\omega}_\can \mapsto (2\pi i)^2 (-2 i \Im \tau)
\otimes (2\pi i )^{-1}.$$
It follows that $\rlangle  g, (I^\infty)^{-1}(F_\infty\otimes 1)
I^\infty h\rrangle_\dr$
is given by the class in
$\BH^2(\bar{X},\omega^\bullet({\cal{I}}_{X^\infty}^\an)[2-k])$
arising from the section of
$\Gamma(\bar{X},\FC^2_{\bar{X}}[2-k])$ whose pull-back
by $\pi_t$ is defined by
$$\tau \mapsto
(k-2)!(2\pi i)^k g_t(\tau)\overline{h_t(\tau)}(-2 i \Im \tau)^{k-2}
        d\tau \wedge d\bar{\tau} \otimes \iota^{k-2},$$
where the pairing on $\Symm^{k-2}$ is defined in (\ref{sympa}).
Under the trace map to  $\CC\otimes \QQ(1-k)_\dr$, this goes to
$(k-2)!(4\pi)^{k-1}(g_t,h_t)_{\Gamma(N)}\otimes \iota^{k-1}$.
\epf

\subsection{Weight filtrations and summary}
\mlabel{ssec:levelN.wt}

There is a natural map $\CM_{c,\wc} \to \CM_\wc$ for each realization
respecting all of the data and comparison isomorphisms.
To make $\CM$ and $\CM_c$ premotivic structures,
one checks that setting
\begin{eqnarray*}W_i M_\wc &=& \left\{\begin{array}{ll}
0,& \mbox{if $i < k-1$;}\\
\im(M_{c,\wc} \to M_\wc),& \mbox{if $k-1 \le i < 2(k-1)$;}\\
 M_\wc,& \mbox{if $2(k-1) \le i$;}\end{array}\right.\\
W_i M_{c,\wc} &=& \left\{\begin{array}{ll}
0,& \mbox{if $i < 0$;}\\
\ker(M_{c,\wc} \to M_\wc),& \mbox{if $0 \le i < k-1$;}\\
 M_{c,\wc},& \mbox{if $k-1 \le i$.}\end{array}\right.
\end{eqnarray*}
defines weight filtrations.  In each case, this
amounts to the assertion that
$$\left(\fil^{k-1} V \right)\,\cap\,
\left( (I^\infty)^{-1}(F_\infty\otimes 1)I^\infty\fil^{k-1} V\right) = 0,$$
where $V =\fil^{k-1}(\CC\otimes M_{c,\dr}) \cong
\fil^{k-1}W_{k-1}(\CC\otimes M_\dr)$, and this is immediate
from Lemma \ref{petersson1}.

We can now regard $\CM$ and $\CM_c$ as objects of $\ipms_\QQ^S$,
and $M = \QQ\otimes\CM$ and $M_c=\QQ\otimes\CM_c$ as objects
of $\pms_\QQ^S$, where $S$ contains the set of primes dividing $Nk!$.
We have also defined a {\em (left) duality morphism}
$$\delta^L: \CM_c \to \hom(\CM,\calt^{S,\otimes(k-1)})$$
whose kernel and cokernel are torsion objects $\CN$ satisfying
$\CN_\ell = 0$ for $\ell > k-2$, and which therefore induces
an isomorphism
$M_c \to \hom(M,\QQ(1-k)^S)$.
Similarly we have a
{\em (right) duality morphism}
$\delta^R=\delta^R_N:\CM \to \hom(\CM_c,\calt^{S,\otimes(k-1)})$ whose kernel and
cokernel satisfy $\CN_\ell = 0$ for $\ell\nmid (k-2)!$.
Moreover the resulting diagram
$$\begin{CD}
 \CM_c @>\delta^L >> \hom(\CM,\calt^{S,\otimes(k-1)})\\
@VVV @VVV\\
\CM @>\delta^R >>  \hom(\CM_c,\calt^{S,\otimes(k-1)})
\end{CD}$$
is commutative,
where the right-hand vertical arrow is induced by
$(-1)^{k-1}$ times the map $\CM_c \to \CM$.

We let $\CM_\tf$ denote the maximal torsion-free quotient
of $\CM$ (i.e., $\CM/\CM[r]$ where $r\in\ZZ_{>0}$ is chosen
to annihilate the the torsion in $\CM_B$ and $\CM[r]$ denotes
the kernel of multiplication by $r$ on $\CM$.)
Finally, we let $\CM_!$ denote the premotivic structure
$\im (\CM_c \to \CM_\tf)$ in
$\ipms_\QQ^S$, pure of weight $k-1$.  We let $M_! = \QQ\otimes\CM_!$.
{}From the commutative diagram above we obtain morphisms
$$\begin{array}{ccc}
\delta^L_!: \CM_! &\to &\hom(\CM_!,\calt^{S,\otimes(k-1)})
{\rm\ and\ }\\
\delta^L_!: M_! &\to &\hom(M_!,\QQ(1-k)^S)\end{array}$$
of sign $(-1)^{k-1}$.  The first map is injective and
the second is an isomorphism;
moreover the following lemma shows that the cokernel
of the first map satisfies $\calc_\ell = 0$ unless $k > 2$ and $\ell|N(k-2)!$.
\begin{lemma} Let $\CN$ denote the cokernel of $\CM_c \to \CM$.
Then $\CN[\ell] = 0$ if $k = 2$ or $\ell$ does not divide $N(k-2)!$.
\mlabel{lem:torsion}
\end{lemma}
\proofbegin It suffices to prove the assertion for the Betti realization,
and for this it suffices
(by \cite[II.17(2)]{bredon}, for example) to check that
for any cusp and any $\alpha > 1$, $H^1(D_\alpha^0,j^*\CF_B^k)$ has no
$\ell$-torsion, where $j$ is the inclusion $D_\alpha^0 \to X^\an$
in a punctured neighborhood of the cusp.  This amounts to the
surjectivity of
$$H^0(D_\alpha^0, j^*\CF_B^k) \to
    H^0(D_\alpha^0,\FF_\ell\otimes j^*\CF_B^k),$$
which holds if $k = 2$ or $\ell$ does not divide $N(k-2)!$.
\epf

We now summarize main properties of the premotivic
structures $\CM_N$ and $M_N$ for level $N$ modular forms.

\begin{theorem}
Fix  integers $N \ge 3$ and $k\geq 2$.
Let $S_N = \{\ell \nmid Nk!\}$.
\begin{enumerate}
\item
Let $S \supseteq S_N$.
Together with the induced comparison isomorphisms and weight
filtrations,
$\CM_N^S,\ \CM_{N,c}^S$ and $\CM_{N,!}^S$ are in $\ipms^S_\QQ$, and
$\CM_{N,c}^S\to \CM_{N,!}^S \to \CM_{N,\tf}^S$ are morphisms
in $\ipms^S_\QQ$.
\item
Let $S\supseteq S_N$.
There are duality homomorphisms
$$\delta^L_N:\CM_{N,c}^S \to \hom_{\ZZ}(\CM_N^S, \CT^{S,\otimes (k-1)})$$
and
$$\delta^L_{N,!}: \CM_{N,!}^S \to \hom_{\ZZ}(\CM_{N,!}^S,
    \CT^{S,\otimes (k-1)})$$
in $\ipms^S_\QQ$, giving rise to
perfect pairings
$$(\ ,\ ): M_{N,c}^S \otimes M_N^S \to \QQ(1-k)^S$$
and
$$(\ , \ )_!: M_{N,!}^S \otimes M_{N,!}^S \to \QQ(1-k)^S$$
in $\pms^S_\QQ$.
\end{enumerate}
\mlabel{thm:levelN}
\end{theorem}

\section{The action of $\GL_2(\AA_\f)$}
\mlabel{sec:action}
In this section we define the adelic action on premotivic
structures associated to modular forms.

\subsection{Adelic modular forms}
\mlabel{ssec:adelic}
We recall the adelic definition of modular curves and forms.

Suppose that $U$ is an open compact subgroup of
$\GL_2(\AA_\f)$
where $\AA_\f$ denotes the finite adeles.  Let $U_\infty$ denote
the stabilizer of $i$ in $\GL_2(\RR)$, so $U_\infty = \RR^\times\SO_2(\RR)$.
The analytic modular curve $X_U$ of level $U$ is defined as the
quotient
$$GL_2(\QQ)\backslash \GL_2(\AA) / UU_\infty.$$
The analytic structure is characterized by requiring that
if $g$ is in $\GL_2(\AA_\ff)$, then the map $\uhp \to X_U$
defined by $\gamma(i) \to \GL_2(\QQ)g\gamma UU_\infty,\
\gamma\in \GL_2^+(\RR),$
is holomorphic.

A modular form of level $U$ is a function  $\phi:\GL_2(\AA) \to \CC$
such that
\begin{itemize}
\item if $\delta \in \GL_2(\QQ)$, $x \in \GL_2(\AA)$, $u \in U$
        and $v = \smat{a}{b}{c}{d} \in U_\infty$, then
        $$\phi(\delta x uv) = \det v (ci + d)^{-k}\phi(x)$$
\item if $g\in\GL_2(\AA^\infty)$, then the function $\uhp \to \CC$
        defined by
        $$ \gamma(i) \mapsto (\det \gamma)^{-1}(ci +d)^k\phi(g\gamma)$$
        for $\gamma = \smat{a}{b}{c}{d} \in \GL_2^+(\RR)$
        is a modular form of weight $k$ with respect to
        $gUg^{-1} \cap \GL_2^+(\QQ)$.
\end{itemize}
Note that the first condition ensures that the function in the
second condition is well-defined and satisfies the usual
transformation property for modular forms.  The second condition
therefore only amounts to a holomorphy requirement.
We let $M_k(U)$ denote the space of modular forms of level $U$.
The space of cusp forms of level $U$ is defined similarly and
denoted $S_k(U)$.

Suppose now that $U$ and $U'$ are open compact subgroups of
$\GL_2(\AA_\f)$, and $g$ is an element of $\GL_2(\AA_\f)$ such
that $g^{-1}U'g\subset U$.  Note that right multiplication by $g$
induces a holomorphic map $X_{U'} \to X_U$, and inclusions
$M_k(U) \to M_k(U')$ and $S_k(U) \to S_k(U')$.
We thus obtain an action of $\GL_2(\AA_\f)$ on
$$\CA_k = \lim_{\stackrel{\to}{U}} M_k(U)\quad\mbox{and}
\quad \CA_k^0 = \lim_{\stackrel{\to}{U}}S_k(U).$$

Suppose now that $U = U_N$ for some $N \ge 3$, where $U_N \subset
\GL_2(\AA_\f)$ is the kernel of the reduction map
$\GL_2(\hat{\ZZ}) \to \GL_2(\ZZ/N\ZZ)$.  For each
class $t \in (\ZZ/N\ZZ)^\times$, we choose an element
$g_t \in \GL_2(\hat{\ZZ})$ whose image in $\GL_2(\ZZ/N\ZZ)$
is $\smat{1}{0}{0}{t^{-1}}$.  We identify $X_N$ with $X_U$
via the maps $\eta_t: X_{N,t} \to X_U$ defined by
$$\Gamma(N) \cdot \gamma(i) \mapsto
\GL_2(\QQ) \cdot g_t \gamma \cdot UU_\infty$$
for $\gamma \in \GL_2^+(\RR)$.  We identify
$M_k(U)$ with $\oplus_t M_k(\Gamma(N))$ via the
isomorphism $\beta$ defined by
\begin{equation}
\beta(\phi)_t(\gamma(i)) =
         (\det \gamma)^{-1}(ci +d)^k\phi(g_t \gamma)
\mlabel{eq:adanform}
\end{equation}
for $\gamma = \smat{a}{b}{c}{d} \in \GL_2^+(\RR)$.
We have a similar identification for cusp forms.
Note that $\eta$ and $\beta$ are independent of the
choices of the $g_t$.
Composing with the isomorphisms defined in \S\ref{ssec:ad.int}.
we obtain isomorphisms
$$\CC\otimes \fil^{k-1}\CM_{\dr} \cong  M_k(U)\quad\mbox{and}\quad
$$
$$\CC\otimes \fil^{k-1}\CM_{c,\dr}
\cong  \CC\otimes \fil^{k-1}\CM_{!,\dr}
\cong  S_k(U).$$
In particular, we have
$$\CA_k \cong  \CC\otimes \lim_{\stackrel{\to}{N}}\fil^{k-1}\CM_{N,\dr}
{\rm\ and\ }$$
$$
\CA_k^0 \cong  \CC\otimes \lim_{\stackrel{\to}{N}}
  \fil^{k-1}\CM_{N,!,\dr}\cong
 \CC\otimes \lim_{\stackrel{\to}{N}}
\fil^{k-1}\CM_{N,c,\dr}.$$ Furthermore if $N, N' \ge 3$ and
$g\in\GL_2(\AA_\f)$ are such that $g^{-1}U_{N'}g \subset U_N$, then we
have a map
$$\CC\otimes\fil^{k-1}\CM_{N,\dr} \to
\CC\otimes\fil^{k-1}\CM_{N',\dr},$$
and similarly for $\CM_c$ and $\CM_!$.
We shall explain how to recover these from maps
of premotivic structures.

\subsection{Action on curves}
\mlabel{ssec:action:curve}
For $h\in \GL_2(\AA_\f)$ and integers $N, N' \ge 3$,
we call $(h,N,N')$ an {\em admissible triple} if both
$h$ and $N'N^{-1}h^{-1} \in M_2(\hat{\ZZ})$.
We leave it to the reader to verify the following
elementary facts:
\begin{lemma}
\mlabel{lem:div}
\begin{enumerate}
\item If $(h,N,N')$ is an admissible triple, then $N|N'$ and
 $h^{-1}U_{N'} h \subset U_{N}$.
\item If $(h,N,N')$ and $(h',N',N'')$ are admissible triples, then
so is $(h'h,N,N'')$.
\item $(h,N,N')$ is an admissible triple if and only if $h = \alpha
\smat{a}{0}{0}{b}\beta$ for some $\alpha, \beta \in \GL_2(\hat{\ZZ})$
and $a,b\in\ZZ$ such that $aN$ and $bN$ divide $N'$.
\mlabel{item:rh3}
\end{enumerate}
\end{lemma}

Suppose now that $(h,N,N')$ is an admissible triple.  Let
$\bar{E}/\bar{X}$ (respectively, $\bar{E}'/\bar{X}'$) denote
the universal generalized elliptic curve with level $N$
(respectively $N'$) structure.

We shall associate to $h$ a finite flat subgroup $G$ of
$E'$ and a level $N$-structure on $E'/G$.
Note that right multiplication by $N'h^{-1}\in M_2(\hat{\ZZ})$
defines an endomorphism of
$\bar{E}'[N'] = (\ZZ/N'\ZZ)_{/X'}^2$, and we define
$G$ to be its image.  One checks that right multiplication
by $N^{-1}N'h^{-1}$ defines an injective map
$$(\ZZ/N\ZZ)^2 \rightarrow
(\ZZ/N'\ZZ)^2/((\ZZ/N'\ZZ)^2(N'h^{-1})$$
and so gives rise to a level $N$-structure on $E'/G$.

Letting $(\bar{E}'/G)_{\rm cont}$ denote the contraction of
$\bar{E}'/G$ whose cuspidal fibers are $N$-gons~\cite[IV.1]{del_rap},
we obtain a level $N$ structure on
$(\bar{E}'/G)_{\rm cont}\to \bar{X}'$.
By the universal property of
$\bar{E}_{/\bar{X}}$, this defines a map
$\chi_h=\chi_{h,N,N'}:\bar{X}'\to \bar{X}$
such that there is an isomorphism
$(\bar{E}'/G)_{\rm cont} \to \bar{E}\times_{\bar{X}}\bar{X}'$
of generalized elliptic curves with level $N$-structures.
We let $\vep_h=\vep_{h,N,N'}$ denote  the composite of the natural map
$\bar{E}' \to (\bar{E}'/G)_{\rm cont}$ with this isomorphism,
so we obtain a commutative diagram
\begin{equation}
 \begin{CD}
\bar{E}' @>\vep_h >> \bar{E}\\
@V s'VV @VV sV\\
\bar{X}' @> \chi_h >> \bar{X}.
\end{CD}
\mlabel{curve}
\end{equation}

We now give the concrete description of its analytification
in terms of the models in \S~\ref{sssec:modcurves}.

Let $U = U_N$, $\Gamma = \Gamma(N)$, $U' = U_{N'}$ and $\Gamma' =
\Gamma(N')$.  Recall that we chose elements $g_t$ in $\GL_2(\hat{\ZZ})$
and defined maps $\eta_t: X_{N,t} \to X_U$ for $t\in
(\ZZ/N\ZZ)^\times$.  For $t'\in(\ZZ/N'\ZZ)^\times$ we denote
these by $g'_{t'}$ and $\eta'_{t'}:X_{N',t'} \to X_{U'}$.
By the strong approximation theorem, we have that
$$g_{\kappa(t')} \in \GL_2(\QQ) g'_{t'}h U\GL_2^+(\RR)$$
where $\kappa=\kappa_h: (\ZZ/N'\ZZ)^\times \to (\ZZ/N\ZZ)^\times,
t'\mapsto (||\det h||\det h)^{-1} \mod N$,
so we can write $g_{\kappa(t')} = \gamma_\f g'_{t'}h u$
for some $\gamma = \gamma_{h,t'} \in \GL_2^+(\QQ)$, $u \in U$.
One checks that the diagram
$$\begin{array}{ccc}
X_{N',t'}  &  \to  & X_{N,\kappa(t')}\\
\downarrow  &&        \downarrow\\
X_{U'}      &  \to  & X_U\end{array}$$
commutes, where the top arrow is defined by $\tau \mapsto
\gamma(\tau)$, the bottom one is defined by right
multiplication by $h$ and the vertical arrows are
$\eta'_{t'}$ and $\eta_{\kappa(t)}$.
Note that $\gamma^{-1} \in M_2(\ZZ)$
and $\gamma\Gamma'\gamma^{-1} \subset \Gamma$ so that the map
$$\begin{array}{rccc}\tilde{\gamma} : & E_{N',t'} &\to& E_{N,\kappa(t')},\\
& (\tau,z) &\mapsto & (\gamma(\tau),(c\tau +d)^{-1}z),\end{array}$$
where $\gamma = \smat{a}{b}{c}{d}$, is well-defined.  Note also
that $\gamma$ is uniquely determined by $(h,N,N')$ and $t'$ up
to left multiplication by an element of $\Gamma$, so $\tilde{\gamma}$
is independent of the choice of $\gamma$.

\begin{lemma}  Suppose that $(h,N,N')$ is an admissible triple
and $t' \in (\ZZ/N'\ZZ)^\times$.
The restriction of $\chi_h^\an$ to $X_{N',t'}$ and $\vep_h^\an$
to $E_{N',t'}$ are given by $\gamma: X_{N',t'} \to X_{N,\kappa(t')}$
and $\tilde{\gamma} :  E_{N',t'}  \to E_{N,\kappa(t')}$.
\mlabel{lem:adan}
\mlabel{lem:modan}
\end{lemma}
\proofbegin
For $\tau' \in \uhp$, let $E_{\tau'} = \CC/\Lambda_{\tau'}$
denote the fiber of $E^{\prime,\an}$ over $\Gamma'\tau' \in X_{N',t'}$,
where $\Lambda_{\tau'} =\ZZ^2\svec{\tau'}{1}$  and let
$G_{\tau'}\subset E_{\tau'}$ denote the fiber of $G^\an$.  Similarly
let $E_\tau = \CC/\Lambda_\tau$ denote the fiber of
$E^\an$ over $\Gamma\tau$, where $\tau = \gamma(\tau')$ and
$\Lambda_\tau = \ZZ^2\svec{\tau}{1}$.
It suffices to check the lemma fiberwise, i.e., that
$$\tilde{\gamma}_{\tau'} : E_{\tau'} \to E_{\tau}$$
induces an isomorphism $E_{\tau'}/G_{\tau'} \cong E_{\tau}$
compatible with level $N$-structures.

Recall that the level $N'$-structure on $E_{\tau'}$ is defined by
$$(m,n)g'_{t'} \bmod N'\hat{\ZZ}^2 \ \mapsto\
(m,n)(N')^{-1}\vec{\tau'}{1}
\bmod \Lambda_{\tau'}$$
for $(m,n) \in \ZZ^2$.  Since $h^{-1} = ug_{\kappa(t')}^{-1}\gamma_\f g'_{t'}$
for some $u \in U$ and multiplication by $ug_{\kappa(t')}^{-1}$ is an automorphism
of $(\ZZ/N'\ZZ)^2$, we see that $G_{\tau'}$ is the image of
$\ZZ^2(N'\gamma_\f)g'_{t'}$.  We therefore have
$$G_{\tau'} = \ZZ^2\gamma\svec{\tau'}{1}/\Lambda_{\tau'}.$$
Furthermore, since multiplication by $u$ is the identity on
$(\ZZ/N\ZZ)^2$, the level $N$-structure on $E_{\tau'}/G_{\tau'}$ is
given by
$$(m,n) g_{\kappa(t')}\bmod N \hat{\ZZ}^2 \ \mapsto\
(m,n)N^{-1}\gamma\svec{\tau}{1}
\bmod \ZZ^2\gamma\svec{\tau'}{1},$$
and the lemma follows.
\proofend

\begin{lemma}  Suppose that $(h,N,N')$ and $(h',N',N'')$
are admissible triples.  Then $\chi_{h'h} = \chi_{h}\circ\chi_{h'}$
and $\vep_{h'h'} = \vep_{h}\circ\vep_{h'}$.
\mlabel{lem:comadm}
\end{lemma}
\proofbegin Since $\bar{E}''$ and $\bar{X}''$ are flat over $\ZZ[1/N'']$, it
suffices to check that the desired equalities after extending
scalars to
$\CC$, so by GAGA it suffices to check the analytifications coincide, so
it suffices to check they coincide on $E_{N'',t''}$ and
$X_{N'',t''}$ for each $t''\in (\ZZ/N''\ZZ)^\times$, for their unions are
dense subsets.

Fix $t''\in(\ZZ/N''\ZZ)^\times$ and let $t' = (||\det h'||\det
h')^{-1}t''$.  Note that by Lemma \ref{lem:modan}, both $\vep_{h'h}^\an$
and $\vep_h^\an\vep_{h'}^\an$ send $E_{N'',t''}$ to $E_{N,t}$ where
$$t = (||\det h||\det h)^{-1}t' = (||\det(h'h)||\det(h'h))^{-1}t'',$$
and similarly for the $\chi$'s.
Furthermore, if we choose $\gamma$ and $\gamma'\in \GL_2^+(\QQ)$ so
that $g_t \in \gamma_\f g'_{t'}h U$ and $g'_{t'} \in \gamma'_\f
g''_{t''}h'U$, then
$$g_t \in \gamma_\f\gamma'_\f g''_{t''}h'h(h^{-1}U'h)U
            = \gamma_\f\gamma'_\f g''_{t''}h'hU.$$
So we can take $\gamma'' = \gamma\gamma'$ and apply Lemma \ref{lem:modan}
to conclude the desired equalities.
\proofend

\subsection{Action on the premotivic structures}
\mlabel{ssec:action.real}
\mlabel{sss:gwt}
Suppose
that $(h,N,N')$ is an admissible
triple.  Let $S$ denote the set of primes dividing $Nk!$
and $S'$ the set of primes dividing $N'k!$.  We let
$\CM$, $\CM_c$ and $\CM_!$ (respectively, $\CM'$, $\CM'_c$
and $\CM'_!$) denote the objects of $\ipms_\QQ^S$
(respectively, $\ipms_\QQ^{S'}$)
structures defined in \S \ref{ssec:levelN.wt}
associated to modular forms of weight $k$ and level $N$
(respectively, level $N'$).  We shall define compatible
morphisms $\CM_c^{S'} \to \CM_c'$ and $\CM^{S'} \to \CM'$,
giving also $\CM_!^{S'} \to  \CM_!'$.

We first describe the action on realizations. The
analytification of~(\ref{curve}) gives a natural
transformation $(\chi_h^\an)^*s^\an_* \to
(s')^\an_*(\vep_h^\an)^*$ from which we obtain a map
$$(\chi_h^\an)^*\CF_\bt =
  (\chi_h^\an)^*R^1s^\an_*\ZZ
 \to R^1(s')^\an_*\ZZ = \CF_\bt'.$$
On stalks over $x' \in X^{\prime,\an}$ this is given
by the map $H^1(E^\an_x,\ZZ) \to H^1(E^{\prime,\an}_x,\ZZ)$
induced by $E^{\prime,\an}_{x'}\to E^\an_x$
where $x = \chi_h(x')\in X^\an$.
Taking symmetric products gives
$$(\chi_h^\an)^*\CF_\bt^k \cong
\sym^{k-2}(\chi_h^\an)^*\CF_\bt \to
\sym^{k-2}\CF_\bt',$$
and taking cohomology yields a homomorphism
$$[h]_\bt: \CM_\bt = H^1(X^\an,\CF_\bt^k)
\to  H^1(X^{\prime,\an},(\chi_h^\an)^*\CF_\bt^k)
\to  H^1(X^{\prime,\an},\sym^{k-2}\CF_\bt') =  \CM_\bt'.$$
Since $X^{\prime,\an} \to X^\an$ is proper, the same
applies to give $[h]_{c,\bt}:\CM_{c,\bt} \to \CM'_{c,\bt}$.
It is straightforward to check that these maps respect
complex conjugation.

Similarly we have a morphism of lisse $\ell$-adic sheaves
$\chi_h^*\CF_\ell \to \CF_\ell'$, which on stalks
$\bar{x}':\spec\bar{k} \to X'$ is the natural map
$$H^1(E'_{\bar{x}'},\ZZ_\ell) \to H^1(E_{\bar{x}},\ZZ_\ell)$$
induced by $E'_{\bar{x}'} \to E_{\bar{x}}$ where
$\bar{x} = \chi_h\circ\bar{x}'$.  Applying $\sym^{k-2}$
and taking cohomology yields  $G_\QQ$-linear
maps $[h]_\ell: \CM_\ell \to \CM_\ell'$
and (since $\chi_h$ is proper)
$[h]_{c,\ell}:\CM_{c, \ell} \to \CM'_{c,\ell}$.

{}From the commutative diagram of log schemes
$$\begin{array}{ccccc}
(\bar{E}',\CN_{E'}) & \to & (\bar{X}',\CN_{X'}) & \to & (T',\CO_{T'}^*)\\
\downarrow &&\downarrow &&\downarrow \\
(\bar{E},\CN_{E}) & \to & (\bar{X},\CN_{X}) & \to & (T,\CO_{T}^*),\end{array}$$
we obtain a commutative diagram of complexes of sheaves of
$(\chi_hs')^{-1}\CO_{\bar{X}}$- modules on $\bar{E}'$:
$$\begin{array}{ccccccccc}
0 & \to & \vep_h^{-1}\omega^\bullet_{\bar{E}/\bar{X}}
\otimes_{\vep_h^{-1}\CO_{\bar{E}}}
\vep_h^{-1}\bar{s}^*\omega^1_{\bar{X}/T}&
    \to& \vep_h^{-1}\omega^\bullet_{\bar{E}/T} & \to&
\vep_h^{-1}\omega^\bullet_{\bar{E}/\bar{X}}&\to&
0\\
&&\downarrow &&\downarrow &&\downarrow&& \\
0 &\to &\omega^\bullet_{\bar{E}'/\bar{X}'}\otimes_{\CO_{\bar{E}'}}
\bar{s}^{\prime, *}\omega^1_{\bar{X}'/T'}
  &  \to &
\omega^\bullet_{\bar{E}'/T'} & \to &
\omega^\bullet_{\bar{E}'/\bar{X}'} &\to &
0,\end{array}$$
where the top row is obtained by applying $\vep_h^{-1}$ to~(\ref{kato-seq2}),
the bottom row is~(\ref{kato-seq2}) with $N'$ instead of $N$, and the vertical
arrows are given by the canonical $\vep_h^{-1}\CO_{\bar{E}}$-linear maps
in each degree.  Applying $\BR^\bullet\bar{s}'_*$ and using the natural
transformation $\chi_h^{-1}\bar{s}_* \to \bar{s}'_*\vep_h^{-1}$ gives
a commutative diagram
{\scriptsize
$$\begin{array}{ccccc}
\chi_h^{-1}\BR^1\bar{s}_*\omega^\bullet_{\bar{E}/\bar{X}}&\to&
\BR^1\bar{s}'_*(\vep_h^{-1}\omega^\bullet_{\bar{E}/\bar{X}})&\to&
\BR^1\bar{s}'_*\omega^\bullet_{\bar{E}'/\bar{X}'})\\
\downarrow &&\downarrow &&\downarrow \\
\hspace{-0.15cm}
\chi_h^{-1}\BR^2\bar{s}_*(\omega^\bullet_{\bar{E}/\bar{X}}
\otimes_{\CO_{\bar{E}}}\bar{s}^*\omega^1_{\bar{X}/T})
\hspace{-0.3cm}& \to& \hspace{-0.35cm}
\BR^2s'_*(\vep_h^{-1}(\omega^\bullet_{\bar{E}/\bar{X}}
\otimes_{\CO_{\bar{E}}}\bar{s}^*\omega^1_{\bar{X}/T}))
\hspace{-0.3cm}&\to& \hspace{-0.35cm}
\BR^2s'_*(\omega^\bullet_{\bar{E}'/\bar{X}'}\otimes_{\CO_{\bar{E}'}}
\bar{s}^{\prime, *}\omega^1_{\bar{X}'/T'}).\end{array}$$
}
The top row gives a $\chi_h^{-1}\CO_{\bar{X}}$-linear
morphism $\chi_h^{-1}\CF_\dr \to \CF'_\dr$ compatible
with filtrations, becoming the map $\BH^1_\dr(E_x/k)
\to \BH^1_\dr(E'_{x'}/k)$ on a fiber $x':\spec k \to X'$ with
$x = \chi_h\circ x'$.
Moreover the commutativity of the diagram
together with the compatibility with projection formulas gives
compatibility with Gauss-Manin connections, in the sense that
the diagram
$$\begin{array}{ccccccc}
\chi_h^{-1}\nabla &:&\chi_h^{-1}\CF_\dr &\to&
\chi_h^{-1}(\CF_\dr\otimes_{\CO_{\bar{X}}}\omega^1_{\bar{X}/T})
&\cong&\chi_h^{-1}\CF_\dr
\otimes_{\chi_h^{-1}\CO_{\bar{X}}}\chi_h^{-1}\omega^1_{\bar{X}/T}\\
&&\downarrow&&&&\downarrow\\
\nabla' &:&\CF'_\dr &&\to&&\CF'_\dr\otimes_{\CO_{\bar{X}'}}\omega^1_{\bar{X'}/T'}
\end{array}$$
commutes.
This in turn gives morphisms $\chi_h^{-1}\CF_\dr^k \to \CF_\dr^{\prime, k}$
and $\chi_h^{-1}\CF_{c,\dr}^k \to \CF_{c,\dr}^{\prime, k}$ compatible
with filtrations and connections, and taking cohomology yields
the maps $[h]_\dr:\CM_\dr \to \CM'_\dr$ and
$[h]_{c,\dr}:\CM_{c,\dr} \to \CM'_{c,\dr}$
preserving filtrations.

Suppose now that $\ell$ is a prime not dividing $N'k!$.
{}From the discussion of functoriality in \cite[4c]{faltings}
and the proof of \cite[Thm.\ 6.2]{faltings}, we see that
the above construction with $T$ and $T'$ replaced by $\ZZ_\ell$
yields a morphism $\chi_{h,\ecrys}^*\CF_\ecrys \to \CF'_\ecrys$
in the category $\mfcat^\nabla_{[0,1]}(X'_{\ZZ_\ell})$.
For $x':\spec\CO_F \to X'$ and $x = \chi_h\circ x'$ with
$F$ unramified over $\QQ_\ell$, the pull-back to
$\mfcat_0(\CO_F)$ is the natural map
$H^1_\crys(E_x,\CO_{E_x,\crys}) \to
H^1_\crys(E'_{x'},\CO_{E'_{x'},\crys})$.
Taking symmetric powers and cohomology then yields the
maps $[h]_\ecrys:\CM_\ecrys \to \CM'_\ecrys$ and
$[h]_{c,\ecrys}:\CM_\ecrys \to \CM'_\ecrys$ in the
category $\mfcat_0(\ZZ_\ell)$.

Now we sketch the proof that the homomorphisms $[h]_*$
and $[h]_{c,*}$
for $\wc=\bt$, $\dr$, $\ell$ and $\ecrys$ (for $\ell \not\in S$)
are compatible with the comparison isomorphisms.

In the case of $I_\bt^\ell$, this follows from the functoriality
properties of the Betti-\'etale comparison isomorphisms.
In particular, the natural transformations
$(f_*\cdot)_\an \to f^\an_*\cdot_\an$ (\cite[I.11]{frei_kiehl})
are well-behaved under composition of maps $f:Y \to Z$
of schemes locally of finite type over $\CC$.  One
deduces that the diagram
$$\begin{array}{ccccc}
(\chi_h^*\CF_\ell/\ell^n)_\an  &\cong &
(\chi^\an_h)^*(\CF_\ell/\ell^n)_\an &\cong&
(\chi^\an_h)^*(\CF_\bt/\ell^n)\\
\downarrow&&&&\downarrow\\
(\CF_\ell'/\ell^n)_\an  && \cong &&
(\CF_\bt'/\ell^n)\end{array}$$
commutes.  (Though not necessary for the proof,
we remark that on stalks at $x'\in X^{\prime,\an}$,
the diagram becomes
$$\begin{array}{ccc}
H^1_\et(E_x,\ZZ/\ell^n\ZZ)
&\cong & H^1(E_x^\an,\ZZ/\ell^n\ZZ) \\
\downarrow&&\downarrow\\
H^1_\et(E'_{x'},\ZZ/\ell^n\ZZ)
&\cong & H^1(E_{x'}^{\prime,\an},\ZZ/\ell^n\ZZ)
\end{array}$$
where $x = \chi_h^\an(x')$ and the vertical
maps are induced by $E'_{x'} \to E_x$.)
Applying $\sym^{k-2}_{\ZZ_\ell}$
gives a similar commutative diagram for the sheaves
$\CF_\wc$ replaced by $\CF_\wc^k$, and this
gives commutativity of the right-hand square
in the diagram
$$\begin{array}{ccccc}
\CM_\ell &\to& H^1(X'_{\bar{\QQ}},\chi_h^*\CF^k_\ell)
&\to&\CM'_\ell\\
\downarrow&&\downarrow&&\downarrow\\
\ZZ_\ell\otimes\CM_\bt&\to&
\ZZ_\ell\otimes H^1(X^{\prime,\an},\chi_h^{\an,*}\CF^k_\bt)
&\to&\ZZ_\ell\otimes\CM'_\bt
\end{array}$$
where the vertical maps are comparison isomorphisms and
the horizontal maps define $[h]_\ell$ and $[h]_\bt$.
The commutativity of the left-hand square is a similar, but
easier, application of the functoriality of the
comparison.  This gives the compatibility of $[h]_\wc$
with $I_\bt^\ell$, and $[h]_{c,\wc}$ is treated similarly.

For the compatibility of $[h]_\dr$ and $[h]_\bt$ with
$I^\infty$, one first defines $[h]_\dr^\an:\CM_\dr^\an \to
\CM_\dr^{\prime,\an}$ analogously to $[h]_\dr$ and checks
that the diagram
$$\begin{array}{ccc} \CC\otimes\CM_\dr & \cong & \CM_\dr^\an \\
\downarrow&&\downarrow\\
\CC\otimes\CM_\dr' & \cong & \CM_\dr^{\prime,\an}\end{array}$$
commutes.  Furthermore the diagram
$$\begin{array}{ccc}
(\chi_h^\an)^{-1}(\CC\otimes\CF_\bt) & \to &
(\chi_h^\an)^{-1}\CG\\
\downarrow&&\downarrow\\
\CC\otimes\CF_\bt' & \to &
\CG'\end{array}$$
commutes, where the horizontal maps are those
occuring in the definitions of $[h]_\bt$ and
$[h]_\dr^\an$ and the rows arise from~(\ref{diff_equ}).
(The commutativity follows for example from that of
$$\begin{array}{ccc}
H^1(E_x^\an,\CC) & \cong & \BH^1(E_x^\an,\Omega^\bullet_{E_x^\an})\\
\downarrow&&\downarrow\\
H^1(E_{x'}^{\prime,\an},\CC) & \to &
\BH^1(E_{x'}^{\prime,\an},
\Omega^\bullet_{E_{x'}^{\prime,\an}})\end{array}$$
at each $x'\in X^{\prime,\an}$, $x = \chi_h^\an(x')$.)
One deduces from this the desired compatibility for $[h]_\wc$
and $[h]_{c,\wc}$.

The compatibility of $[h]_\dr$ and $[h]_\ecrys$ with
$I^\ell_\dr$ for $\ell \not\in S$
follows from the construction of $[h]_\ecrys$
based on~\cite{faltings}, and the same goes for $[h]_{c,\dr}$
and $[h]_{c,\ecrys}$.

Finally, using the compatibility with pull-back of the natural
transformations denoted $\alpha$ in~\cite[Va)]{faltings}, one gets
commutativity of the diagram
$$\begin{array}{ccccc}
\VV(\chi_{h,\ecrys}^*\CF_\ecrys)  &\cong &
\chi_h^*\VV(\CF_\ecrys) &\cong&
\chi_h^*\CF_\ell\\
\downarrow&&&&\downarrow\\
\VV(\CF_\ecrys')  && \cong &&
\CF'_\ell\end{array}$$
(which on the stalk $\bar{x}':\spec\bar{F}\to X'$ arising
from $x':\spec\CO_F \to X'$ for $F$ unramified over $\QQ_\ell$
is just
$$\begin{array}{ccc}
\VV(H^1_\crys(E_x,\CO_{E_x,\crys}))
&\cong&
H^1_\et(E_{\bar{x}},\ZZ_\ell)\\
\downarrow&&\downarrow\\
\VV(H^1_\crys(E'_{x'},\CO_{E'_{x'},\crys}))
& \cong &
H^1_\et(E'_{\bar{x'}},\ZZ_\ell)\end{array}$$
where $x = \chi_h\circ x'$).  This in turn
gives the commutativity of the right-hand square in
$$\begin{array}{ccccc}
\VV(\CM_\ecrys) &\to&
\VV(H^1_\crys(X'_{\ZZ_\ell},\chi_{h,\ecrys}^*\CF^k_\ecrys))
&\to&\VV(\CM'_\ecrys)\\
\downarrow&&\downarrow&&\downarrow\\
\CM_\ell&\to&
H^1_\et(X'_{\bar{\QQ}},\chi_h^*\CF^k_\ell)
&\to&\CM'_\ell
\end{array}$$
and the commutativity of the left-hand square is also
a consequence of the compatibility of $\alpha$ with
pull-back.  This gives the compatibility of
$[h]_\ecrys$ and $[h]_\ell$ with $I^\ell$, and
the proof is similar for $[h]_{c,\wc}$.

For each realization $\wc=\bt,\dr, \ell$ and
$\ecrys$,
we have a commutative diagram
$$
\begin{CD}
\CM_{c,\wc} @>>> \CM_\wc \\
@V [h]_{c,\wc} VV @VV [h]_\wc V\\
\CM'_{c,\wc} @>>> \CM'_\wc.
\end{CD}
$$
This shows that $[h]_{\wc}$ and $[h]_{c,\wc}$
respect weight filtrations and induce maps
$[h]_{!,\wc}: \CM_{!,\wc} \to \CM'_{!,\wc}$.
Since $[h]_{\wc}$, $[h]_{\wc,c}$ (and hence
$[h]_{\wc,!}$) are compatible with the comparison
isomorphisms $I^\infty$, $I_\bt^\ell$, $I_\dr^\ell$
and $I^\ell$ when $\ell\nmid N'k!$, we conclude that
$h$ defines morphisms
$[h]:\CM^{S'} \to \CM'$, $[h]_c:\CM_c^{S'}\to \CM'_c$ and
$[h]_!:\CM_!^{S'}\to \CM'_!$ in $\ipms^{S'}_\QQ$ where
$S' = \{\,p\,\mid\, p\nmid N'k!\,\}$.

\medskip

\begin{lemma}
Suppose that $(h,N,N')$ and $(h',N',N'')$ are admissible triples.
Then
$$[h']_{N',N''} \circ [h]_{N,N'}=[h'h]_{N,N''}$$ on
$\CM^S, \CM^S_c$ and $\CM^S_!$ whenever $S\supseteq S_{N''}$.
\mlabel{lem:comreal}
\end{lemma}
\proofbegin
We only need to prove
$[h']_B \circ [h]_B =[h'h]_B$.
Using Lemma~\ref{lem:comadm}, we
get the commutative diagram
$$\begin{CD}
(\chi_{h'h}^\an)^* \CF_B @>>> (\chi_{h'}^\an)^* \CF'_B\\
\hspace{-1cm} \searrow \hspace{-3cm} @.
  @VVV \\
@. \CF_B'',
\end{CD}
\mlabel{eq:rel1}
$$
which is just
$$\begin{CD}
H^1(E_x^\an,\ZZ) @>>> H^1(E^{\prime,\an}_{x'},\ZZ)\\
\hspace{-1cm} \searrow \hspace{-3cm} @.
  @VVV \\
@. H^1(E^{\prime\prime,\an}_{x''},\ZZ)
\end{CD}
$$
on the stalk at $x''\in X^{\prime\prime,an}$
(where $x'=\chi_{h'}^\an(x'')$ and
$x=\chi_h^\an (x')$).

Taking symmetric powers and cohomology,
we get the desired commutative diagram.
\proofend

Note that if $r\in \ZZ$ and $(rI,N,N')$ is an admissible triple,
then $rN|N'$ and $(I,N,N')$ is also admissible.
\begin{lemma} If $r\in \ZZ$ and $(rI,N,N')$ is an admissible triple,
then $[rI] = r^{k-2}[I]$, $[rI]_c = r^{k-2}[I]_c$ and
$[rI]_! = r^{k-2}[I]_!$.
\mlabel{lem:gscalar}
\end{lemma}
\proofbegin
It suffices to check the lemma on Betti realizations.
Applying Lemma~\ref{lem:adan} to $h= rI$ and $h=I$, we
see that $\chi_{rI}^\an = \chi_I^\an$ on $X^{\prime,\an}$ and
$\vep_{rI}^\an = \vep_I^\an r$ on $E^{\prime,\an}$.  Since
multiplication by $r$ on the fiber $E^\an_x$
over $x\in X^\an$ induces multiplication by $r$ on
$\CF_{B,x} = H^1(E^\an_x,\ZZ)$, it follows that the
map $(\chi^\an)^*\CF_B \to \CF_B'$ arising from $(rI,N,N')$
is $r$ times the map arising from $(I,N,N')$ (writing simply
$\chi^\an$ for $\chi_{rI}^\an = \chi_I^\an$).  The same then
holds for the maps $(\chi^\an)^*\CF_B^k \to \CF_B^{\prime, k}$
with $r$ replaced by $r^{k-2}$, and the lemma follows.
\epf

Now suppose that $g\in M_2(\hat{\ZZ}) \cap \GL_2(\AA_\f)$ with
$g^{-1} U_{N'} g \subseteq U_N$.  Writing
$g = \alpha\smat{a}{0}{0}{b}\beta$ for some $\alpha, \beta \in \GL_2(\hat{\ZZ})$
and $a,b\in\ZZ$ and applying Lemma~\ref{lem:div}c) shows that $g=rh$
for some $r\in \ZZ$ such that $(h,N,N')$ is admissible.  In fact, we can
let $r_0 = \gcd(a,b)$ and let $h_0 = r_0^{-1}g_0$.
Moreover if $g = rh$ is another decomposition with $r\in\ZZ$ and
$(h,N,N')$ admissible, we see that $r|r_0$, $Nr_0|N'r$ and
and $(h_0,N,N'')$ is admissible, where $N'' = r_0^{-1}|r|N'$, so
$$
r^{k-2} [h]_{N,N'}
= r^{k-2} [r_0r^{-1}I]_{N'',N'}[h_0]_{N,N''}
= r_0^{k-2} [I]_{N'',N'}[h_0]_{N,N''}
=r_0^{k-2} [h_0]_{N,N'}
$$
by Lemmas~\ref{lem:comreal} and~\ref{lem:gscalar}.
Thus we obtain a morphism $[g]=[g]_{N,N'}:\CM^{S'} \to \CM'$ by defining
$$[g]=r^{k-2}[h]:\CM^{S'} \to \CM'
$$
which is independent of the factorization
$g=rh$ with $r\in \ZZ$ and $(h,N,N')$ admissible.
We also define $[g]_c= r^{k-2}[h]_c$ and $[g]_! = r^{k-2}[h]_!$.

Similarly if $g\in \GL_2(\AA_\f)$ with
$g^{-1} U_{N'} g \subseteq U_N$, we can write $g=rh$
for some $r\in \QQ$ so that $(h,N,N')$ is admissible
and obtain morphisms $[g]$, $[g]_c$ and $[g]_!$ in
$\pms_\QQ^{S'}$.

The properties of the action of $g$ on premotivic structures can be
summarized in the following theorem.
\begin{theorem}
Fix positive integers $k\geq 2$, and $N$, $N'\geq 3$.
Let $g$ be an element in $\GL_2(\AA_\f)$
such that $g^{-1}U_{N'} g\subseteq U_N$.
\begin{enumerate}
\item  Let $S'= S_{N'}$. Then
$[g] = [g]_{N,N'}:M^{S' }\to M'$ is a morphism
in $\pms^{S'}_\QQ$.
\item  If $N''\geq 3$  is an integer
and $g'\in\GL_2(\AA_\f)$ is such that
$g^{\prime, -1}U_{N''} g'\subseteq U_N'$, then
$[g']_{N'',N'}\circ[g]_{N',N}^{S_{N''}} = [g'g]_{N'',N}$.
\item If $g\in\GL_2(\AA_\f)\cap M_2(\hat{\ZZ})$ then
$[g]$ arises from a morphism $\CM^{S' }\to \CM'$
in $\ipms^{S'}_\QQ$.
\item Analogous assertions hold for $[g]_c$ and $[g]_!$.
\end{enumerate}
\mlabel{thm:g-action}
\mlabel{prop:comp}
\end{theorem}

\subsection{Relation with the action on modular forms}
\mlabel{ss:gmod}
Given $k\geq 2$, we have established in
\S \ref{sssec:mot.mf.rel}
an isomorphism
$$
\alpha: \CC\otimes \fil^{k-1}\CM_{\dr} \cong
\bigoplus_{t\in(\ZZ/N\ZZ)^\times} M_k(\Gamma(N)).
$$
In  (\ref{eq:adanform}) we have also described an isomorphism
$$
\beta: M_k(U_N) \to \bigoplus_{t\in (\ZZ/N\ZZ)^\times}
    M_k(\Gamma(N)).
$$
Both $\fil^{k-1} M_{N,\dr} \otimes \CC$ and $M_k(U_N)$
have an action
by elements from $M_2(\hat{\ZZ})\cap \GL_2(\AA_\f)$.

\begin{prop}
The canonical isomorphism
$$\beta^{-1}\circ \alpha: \fil^{k-1}\CM_{N,\dr}\otimes \CC
    \to M_k(U(N))$$
preserves the action of $M_2(\hat{\ZZ})\cap \GL_2(\AA_\f)$.
The same holds for the isomorphism
$$\beta^{-1}\circ \alpha: \fil^{k-1}M_{N,c,\dr}\otimes \CC
    \to S_k(U(N)).$$
\mlabel{prop:mot-ad}
\end{prop}
\proofbegin
Fix a $g\in M_2(\hat{\ZZ})\cap \GL_2(\AA_\f)$
with $g^{-1}U'g\subset U$, we will show
that the actions of $g$ on $M_k(\Gamma(N))$ induced from the
isomorphisms $\alpha$ and $\beta$ are the same.
First assume that $(h,N,N')$ is an admissible triple.

Let
\begin{align*}
&\fil^{k-1}\CM_{\dr} \otimes \CC \cong
H^0(\bar{X}^\an, \CL^\an) \\
&\cong \bigoplus_{t\in(\ZZ/N\ZZ)^\times} H^0(\bar{X}_{N,t},i_t^*\CL^\an)
\hookrightarrow \bigoplus_{t\in (\ZZ/N\ZZ)^\times}
H^0(\uhp , \pi_t^*i_t^*\CL^\an)
\notag
\end{align*}
be the injective map define in \S~\ref{sssec:mot.mf.rel}.
Identifying $\fil^{k-1}\CM_{\dr} \otimes \CC $ with its
image under this map, then the isomorphism $\alpha$ is given by
$$
f(\tau)(2\pi i)^{k-1}(dz)^{\otimes(k-2)}\otimes d\tau
    \mapsto f(\tau).$$
The restriction of $\alpha$ to $\fil^{k-1}\CM_{c,\dr} \otimes \CC$
defines the isomorphism
$$\CC\otimes \fil^{k-1}\CM_{c,\dr} \cong
\bigoplus_{t\in(\ZZ/N\ZZ)^\times} S_k(\Gamma(N)).$$

For $t'\in (\ZZ/N'\ZZ)^\times$, let
$\kappa (t')\in (\ZZ/N \ZZ)^\times$
and $\gamma=\gamma_{h,t'}=\smat{a}{b}{c}{d}\in \GL_2^+(\QQ)$
be associated to $(h,N,N')$ in Lemma~\ref{lem:adan}.
Then the action on
$$\bigoplus_{t\in (\ZZ/N\ZZ)^\times}
H^0(\uhp , \pi_t^*i_t^*\CL^\an)$$
induced by $h$ sends
$\theta\defeq ((2\pi i)^{k-1}f_{t}(\tau)
(dz)^{\otimes(k-2)}\otimes d\tau)_{t}$
to
\begin{align*}
&((2\pi i)^{k-1}f_{\kappa (t')}(\gamma_{t'}(\tau))
d((c\tau+d)^{-1}z)^{\otimes(k-2)}\otimes
d\gamma(\tau))_{t'} \\
&= ((2\pi i)^{k-1}
(c\tau+d)^{-k} \det \gamma f_{\kappa(t')}(\gamma(\tau))
    dz^{\otimes(k-2)}\otimes d\tau)_{t'}.
\end{align*}
This means
\begin{equation}
\alpha(g\theta)_{t'} (\tau) =(c\tau+d)^{-k} \det \gamma
    \alpha(\theta)_{\kappa(t')} (\gamma(\tau)).
\mlabel{eq:mod.mot-an}
\end{equation}

On the other hand, for $\phi\in M_k(U_N)$,
its image $\beta(\phi)\in \oplus_{t\in (\ZZ/N\ZZ)^\times}
M_k(\Gamma(N))$ is given by
$$\beta(\phi)_t (x(i))=(\det x)^{-1} (c'i+d')^k \phi(g_t x),
\quad
x=\smat{a'}{b'}{c'}{d'}\in \GL_2^+(\RR). $$
Similarly, the image of
$h \phi $ is given by
$$
\beta(h\phi )_{t'}(x(i))
    = (\det x)^{-1} (c'i+d')^k (h \phi )(g_{t'} x)
    = (\det x)^{-1} (c'i+d')^k \phi(g_{t'} x h).
$$
Using Lemma~\ref{lem:adan} and properties of $\phi$, we verify
that
$$
(\det x)^{-1} (c'i+d')^k \phi(g_{t'} x h)
= \beta(\phi)_{\kappa(t')}(x\gamma (i))\det\gamma
    (c x(i)+d)^{-k}.
$$
Writing $\tau = x(i)$ for the right action of $x$ on $i$,
we have
\begin{equation}
\beta(h\phi)_{t'} (\tau) =(c\tau+d)^{-k} \det \gamma
    \beta(\phi)_{\kappa(t')} (\gamma(\tau)).
\mlabel{eq:mod.ad-an}
\end{equation}
Combining (\ref{eq:mod.mot-an}) and (\ref{eq:mod.ad-an})
we see that if $(\beta^{-1}\circ \alpha)(\theta)=\phi$,
then $(\beta^{-1}\circ \alpha) (h\theta)=h\phi$.

Further let $r\in \ZZ$. Then the action of $rI$ on both spaces
are multiplication by $r^{k-2}$. This proves the proposition.
\proofend

\subsection{Compatibility with the pairings}
\mlabel{g-pair}

\begin{prop}
Suppose $S$ contains $S_N$.
Let
$$\delta^L_N:\CM_{N,c}^S \to \hom(\CM_N^S,\CT^{\otimes (k-1),S})$$
and
$$\delta^L_{N'}:\CM_{N',c}^S \to \hom(\CM_{N'}^S,\CT^{\otimes (k-1),S})$$
be the duality morphisms defined in \S \ref{sssec:pairings}.
Let $[g]^*: \hom(\CM_{N'},\CT^{\otimes (k-1),S})
    \to \hom(\CM_{N'},\CT^{\otimes (k-1),S})$ be induced from
$[g]$.
Then
$$ \delta^L_N=||\det g||^{2-k} [U_N:U_{N'}] (g^* \circ \delta^L_{N'} \circ g)$$
in $\ipms^S_\QQ$.
A similar relation holds between the duality morphism
$\CM_{N,!}$ and the duality morphism on
$\CM_{N',!}$
\mlabel{pp:gpair}
\end{prop}
\proofbegin
We only need to show
the equation on the Betti realization.
Let
$$\rlangle\ ,\ \rrangle:
\CM_{c,N} \otimes \CM_N \to \CT^{\otimes (k-1),S}_\bt$$
and
$$\rlangle\ ,\ \rrangle':
\CM_{c,N'} \otimes \CM_{N'} \to \CT^{\otimes (k-1),S}_\bt$$
be the pairings associated to
$\delta^L_N$ and $\delta^L_{N'}$.
We only need to showing that,
for $x\in \CM_{N,c,\bt}$ and $y\in \CM_{N,\bt}$,
we have
\begin{equation}
 \rlangle gx, gy \rrangle' =  ||\det g||^{2-k}
[U_N:U_{N'}]\rlangle x,y\rrangle
\mlabel{eq:gpair}
\end{equation}
in $\CT^{\otimes (k-1),S}_\bt$.

Write $g=r h$ as before with $r\in \ZZ$
and $(h,N,N')$ admissible.
Assume that the lemma holds for $h$.
Then
\begin{align*}
\rlangle gx, gy\rrangle'
&= \rlangle r^{k-2} hx,r^{k-2} hy \rrangle'
=r^{2(k-2)} \rlangle hx,hy\rrangle' \\
&= ||\det (rI)||^{2-k} \rlangle hx,hy\rrangle'
= ||\det g||^{2-k} [U_N:U_{N'}] \rlangle x,y\rrangle.
\end{align*}
Thus we only need to prove (\ref{eq:gpair}) when
$(g,N,N')$ is admissible.

Using
the functoriality of cup products~\cite[II.8.2]{bredon},
we have the commutative diagram
$$
\begin{CD}
\chi^*\CF_B \otimes \chi^*\CF_B @>>> \chi^* R^2 s_*^\an \ZZ_E
    @>>> \chi^*(2\pi i)^{-1} \ZZ_X\\
@V [h]_B \otimes [h]_B VV @V [h]_B VV @VV h_\ZZ V\\
\CF_B' \otimes \CF_B' @>\cup >>
    R^2 s_*^{\prime, \an} \ZZ_{E'}
    @>>> (2\pi i)^{-1}\ZZ_{X'}
\end{CD}
$$
where the first arrow in the top row is the composite
$$ \chi^* \CF_B \otimes \chi^* \CF_B \to
    \chi^* (\CF_B \otimes \CF_B) \ola{\cup}
    \chi^* R^2 s_*^{\an} \ZZ_E $$
and the right column is the natural map induced by $h$.
Taking the stalks over $x'\in X^{\prime,\an}$, the diagram
gives
$$
\begin{CD}
H^1(E^\an_x,\ZZ) \otimes H^1(E^\an_x,\ZZ) @>>>
H^2(E^\an_x,\ZZ)
    @>Tr >> (2\pi i)^{-1} \ZZ\\
@V [h]_x \otimes [h]_x VV @V [h]_x VV @VV h_{\ZZ,x} V\\
H^1(E^{\prime,\an}_{x'},\ZZ) \otimes H^1(E^{\prime,\an}_{x'},\ZZ) @>>>
H^2(E^{\prime,\an}_{x'},\ZZ)
    @>Tr >> (2\pi i)^{-1} \ZZ
\end{CD}
$$
where $x=\chi_h(x')$.
By Lemma~\ref{lem:adan}, for any given $x'=\Gamma' \tau'\in X^{\prime,\an}$,
the fiber $G^\an_{x'}$ of $G^\an$ is given by $G_{\tau'}$
which has cardinality $\det \gamma = ||\det h||^{-1}.$
So $h_{\ZZ,x}$ is multiplication by $||\det h||^{-1}$.
Then the same is true for $h_\ZZ$.
Taking the $(k-2)$-nd symmetric power, we have the
commutative diagram
$$
\begin{CD}
\chi^*\CF^k_B \otimes \chi^*\CF^k_B
    @>>> \chi^*(2\pi i)^{2-k} \ZZ_X\\
@V [h]_B \otimes [h]_B VV
@VV ||h||^{2-k} V\\
\CF_B^{\prime,k-2} \otimes \CF_B^{\prime,k-2}
    @>>> (2\pi i)^{2-k}\ZZ_{X'}.
\end{CD}
$$

Next note that
the map $\chi: X^{\prime,\an}\to X^\an$, identified with the map
$$ G_\QQ\backslash G_\AA / U'U_\infty
    \cong G_\QQ \backslash G_\AA / h^{-1}U' h U_\infty
    \to G_\QQ \backslash G_\AA / U U_\infty$$
is a covering map of degree $[U:h^{-1} U' h]=[U:U']$.
This gives us the commutative diagram
\begin{equation}
\begin{CD}
H^2_c(X^\an, \ZZ) @>>> H^2_c(G_\QQ \backslash G_\AA / h^{-1}U' h
U_\infty, \ZZ) @>>> H^2_c(X^{\prime,\an},\ZZ) \\
@V\Tr VV @V\Tr VV @VV \Tr V\\
(2\pi i)^{-1}\ZZ @>[U:U']>> (2\pi i)^{-1}\ZZ @=
(2\pi i)^{-1}\ZZ.
\end{CD}
\mlabel{eq:gpair1}
\end{equation}
Then by functoriality of cup products of cohomology
and
the Universal Coefficient Theorem \cite[II.15.3]{bredon}
(see also \cite[Exer. III.8.3]{hart}),
we get the commutative diagram
$$
{\scriptsize
\begin{CD}
\CM_{B,c} \otimes \CM_B & \ola{\cup} &
    H^2_c(X^\an,(2\pi i)^{2-k}\ZZ_X)
@>\cong>> (2\pi i)^{2-k}\ZZ \otimes H^2_c(X^\an,\ZZ_X)
    @> \id \otimes\Tr >> (2\pi i)^{1-k}\ZZ\\
@VV [h]_{B,c}\otimes [h]_B V @V VV
@V||\det h||^{2-k}\otimes H^2_c(\chi) VV @V||\det h||^{2-k} [U:U'] VV
    \\
\CM'_{B,c} \otimes \CM'_B & \ola{\cup} &
    H^2_c(X^{\prime,\an},(2\pi i)^{2-k}\ZZ_{X'})
@>\cong >> (2\pi i)^{2-k}\ZZ\otimes H^2_c(X^{\prime,\an},\ZZ_{X'})
@> \id \otimes \Tr >> (2\pi i)^{1-k}\ZZ
\end{CD}
}
$$
where $H^2_c(\chi)$ is the composition of the top row
in diagram~(\ref{eq:gpair1}).
This is what we want.
\proofend

\subsection{Compatibility with changing of levels}
\mlabel{ssec:g-level}
Fix $g\in M_2(\ZZ)\cap \GL_2(\BA_\f)$.
Let $N, N', M$ and $M'$ be integers such that
$N|M, N'|M', g^{-1}U_{N'}g \subset U_N$
and $g^{-1}U_{M'}g\subset U_M$.
Then we have
$[g]: \CM_N\to \CM_{N'}$ and
$[g]: \CM_M\to \CM_{M'}$.
Let $I$ be the identity matrix in $\GL_2(\BA_\ff)$. From
$N|M$ we get $I^{-1}U_{M}I\subset U_N$ and we have
$$ I:
\bar{E}_{M}{}_{/\bar{X}_{M}}\to \bar{E}_N{}_{/\bar{X}_N}$$
and
$$[I]: \CM_N \to \CM_{M}\ {\rm\ for\ }S\supset S_M.$$
In the same way, we have
$$[I]: \CM_{N'} \to \CM_{M'}\ {\rm\ for\ } S\supset S_{M'}. $$

\begin{lemma}
\begin{enumerate}
\item
For $S\supseteq S_{M'}$,
the diagram
$$\begin{CD}
\CM_N^S @> [g]_{N,N'} >> \CM_{N'}^S\\
@V [I] VV @V [I] VV\\
\CM_{M}^S @> [g]_{M,M'}>> \CM_{M'}^S
\end{CD}
$$
of morphisms in $\ipms^{S}$ is commutative.
\item
$[I]:\CM_{N,\tf}^{S_M}\to \CM_{M,\tf}^{U_N/U_{M}}$ is injective
\item
$[I]_\sharp: \CM_{N,\sharp}^{S_M} \to \CM_{M,\sharp}^{U_N/U_{M}}$
is an isomorphism for $\sharp = c$ or $!$.
\item
$[I]: M_N^{S_M} \to M_{M}^{U_N/U_{M}}$
is an isomorphism.
\end{enumerate}
\mlabel{lem:gmap}
\end{lemma}
\proofbegin
(a) follows directly from
Theorem~\ref{prop:comp}.

(b). We only need to consider $[I]_B$.
Denote $s_M^\an:E_M^\an \to X_M^\an$ for the universal elliptic
curve and denote $\CF_{M,B}$ for the locally constant sheaf
$R^1 s^\an_{M,*} \ZZ$ on $X_M^\an$.
Leray spectral sequence gives us
\begin{multline}
 0\to H^1(U_N/U_{M},H^0 (X_{M},\CF_{M,B}^k)) \to
H^1 (X_N,\CF^k_{N,B}) \\
\to H^1 (X_{M},\CF_{M,B}^k)^{U_N/U_{M}}
    \to H^2(U_N/U_{M},H^0 (X_{M},\CF_{M,B}^k)).
\notag
\end{multline}
Note that the map
$H^1 (X_N,\CF^k_{N,B})
\to H^1 (X_{M},\CF_{M,B}^k)$
is induced by the projection
$p: \bar{E}_M{}_{/\bar{X}_M} \to \bar{E}_{/\bar{X}}$.
Since $U_N/U_{M}$ is finite and $H^0(X_{M},\CF_{M,B}^k)$
is torsion free, we have
\[ H^1(U_N/U_{M},H^0(X_{M},\CF_{M,B}^k)) =
\Hom (U_N/U_{M},H^0(X_{M},\CF_{M,B}^k))=0.\]
Therefore
$\CM_{N,\bt} =H^1 (X_N,\CF^k_{N,B})
\to H^1 (X_{M},\CF_{M,B}^k)^{U_N/U_{M}}=(\CM_{M,\bt})^{U_N/U_M}$
is injective.

(c).
The proof is similar to (b).
We only need to consider the Betti realization.
We know that Leray spectral
sequence\cite[IV.9.2]{bredon} gives an exact sequence

\begin{multline}
 0\to H^1(U_N/U_{M},H^0_c(X_{M},\CF_{M,B}^k)) \to
H^1_c(X_N,\CF^k_{N,B}) \notag \\
\to H^1_c(X_{M},\CF_{M,B}^k)^{U_N/U_{M}}
    \to H^2(U_N/U_{M},H^0_c(X_{M},\CF_{M,B}^k))
\end{multline}
and $H^0_c(X_{M},\CF_{M,B}^k)=0$ since $\CF_{M,B}^k$ have no
non-trivial global sections with compact support.
Thus $\CM_{N,B,c}$ is identified with $(\CM_{M,B,c})^{U_N/U_M}$.

(d).
Since $U_N/U_M$ is finite,
$H^2(U_N/U_{M},H^0 (X_{M},\CF_{M,B}^k))$ is torsion. Thus from the
proof of (b) we see that
$H^1 (X_N,\CF^k_{N,B})\otimes \QQ
\to H^1 (X_{M},\CF_{M,B}^k)^{U_N/U_{M}}\otimes \QQ$ is an
isomorphism. This prove (d) for the Betti realization which
is all we need.
\epf

\section{The premotivic structure for forms
    of level $N$ and character $\psi$}
\mlabel{sec:levchar}

\subsection{$\sigma$-constructions}
\mlabel{ssec:sigma}
Suppose $U$ is any open compact subgroup of $\GL_2(\hat{\ZZ})$.
Let $K$ be a number field with ring of integers $\CO_K$,
Let $V$ be a finite dimensional vector space over $K$
and let $\sigma: U\to \Aut_K(V)$ be a continuous
representation of $U$. Define
$$S_\sigma=\{ \lambda\in S_\f(K) |\
\lambda | k!\mbox{\ or\ }
\GL_2(\ZZ_\ell) \not\subset \ker\sigma
\ {\rm\ where\ } \lambda | \ell\}.$$

Since $U_N$ is normal in $U$, by Theorem~\ref{prop:comp},
we have a group action of $U$ on $\CM_N$.
We regard $\CV$ as an object of $\ipms_K^{S_\sigma}$
(with trivial additional structures and tautological
comparison isomorphisms).  Note that we have an action
of $U$ on $\CV$, hence on the object
$\CM_N\otimes\CV = (\CM_N\otimes\CO_K)\otimes_{\CO_K}\CV$
of $\ipms_K^{S_\sigma}$.  We can thus define an object
$$\CM(\sigma)_N = (\CM_N \otimes \CV)^U$$
of $\ipms_K^{S_\sigma}$ (see the last paragraph of
\S \ref{ssec:pre.str}).
More explicitly, for each $\wc=B,\dr$, $\lambda\in S_\f(K)$
or $\lambda$-crys with $\lambda\in S_\sigma$,
we have
$$\CM(\sigma)_{N,\wc} = (\CM_{N,\wc} \otimes \CV)^{U}.$$
We also define objects
\[ \CM(\sigma)_{N,\sharp} =(\CM_{N,\sharp}\otimes \CV)^U\]
for $\sharp=\tf$ or $c$ and define
\[ \CM(\sigma)_{N,!}= \im (\CM(\sigma)_{N,c} \to
\CM(\sigma)_{N,\tf})
\]
where
$\CM(\sigma)_{N,c} \to \CM(\sigma)_{N}$
is the restriction of the map
$\CM_{N,c}\otimes \CV \to \CM_N\otimes \CV$.
These premotivic structures will be called the
{\em $\sigma$-constructions}.

We now show that
the $\sigma$-constructions are
{\em compatible with the change of levels}.
More precisely,

\begin{lemma}
Let $N|N'$ with $U_N\subset U$.
\begin{enumerate}
\item
For $S\supseteq S_{N'}$, we have
$\CM(\sigma)_{N,c}^S \cong \CM(\sigma)_{N',c}^S$
in $\ipms^S$.
Similarly,
$\CM(\sigma)_{N,!}^S \cong \CM(\sigma)_{N',!}^S$.
\item
For $S\supseteq S_{N'}$, we have
$M(\sigma)_N^S \cong M(\sigma)_{N'}^S$.
\end{enumerate}
\mlabel{lem:coset0}
\end{lemma}
\proofbegin
(a) We only need to prove it for $S=S_{N'}$.
{}From $N\mid N'$ we have $U_{N'}\subset U_N$.
By Lemma~\ref{lem:gmap}, we have
\[ \CM_{N,c}^S=(\CM_{N',c}^{U_N/U_{N'}})^S.\]
Then since $U_N$ acts trivially on $\CV$, we have
$$ (\CM_{N',c}^S\otimes \CV)^{U_N/U_{N'}}
    =(\CM_{N',c}^{U_N/U_{N'}})^S \otimes \CV
    = \CM_{N,c}^S\otimes \CV.$$
Hence
$$\CM(\sigma)_{N',c}^S=(\CM_{N',c}^S\otimes \CV)^U=
((\CM_{N',c}^S\otimes \CV)^{U_N/U_{N'}})^U=
(\CM_{N,c}^S\otimes \CV)^U=\CM(\sigma)_{N,c}^S.$$

Since $\CV$ is a flat $\ZZ$-module,
using Lemma~\ref{lem:gmap}, we similarly find that
$\CM(\sigma)_{N,\tf}^S \to \CM(\sigma)_{N',\tf}^S$
is injective.
Consider the commutative diagram
$$\begin{CD}
\CM(\sigma)_{N,c}^S @> \cong >> \CM(\sigma)_{N',c}^S \\
@VVV @VVV\\
\CM(\sigma)_{N,!}^S @> \subset >> \CM(\sigma)_{N',!}^S \\
@VVV @VVV\\
\CM(\sigma)_{N,\tf}^S @> \subset >> \CM(\sigma)_{N',\tf}^S.
\end{CD}
$$
By the definition of $\CM(\sigma)_{N,!}$,
the two vertical maps in the upper square are surjective
and the two vertical maps in the bottom square are injective.
Since the top row is surjective, so is the map in the middle row.
Since the bottom row is injective, the map in
the middle row is also injective, hence is bijective.

The proof for (b) is the same.
\proofend

Because of the lemma, we will suppress $N$ in the notations
by using $\CM(\sigma)_c$, $\CM(\sigma)_!$ and $M(\sigma)$
when there is no danger of confusion.

\subsection{Action of double cosets}
\mlabel{ssec:coset}
Suppose we are given two open compact subgroups
$U$ and $U'$ contained in $\GL_2(\hat{\ZZ})$.
Let $\sigma:U\to \Aut_K(V)$ and $\sigma':U'\to \Aut_K(V')$
be two representations.
Suppose we are given a $g \in \GL_2(\AA_\ff) \cap M_2(\hat{\ZZ})$
and a $K$-linear homomorphism
$\tau:V \to V'$
such that $\tau(\sigma(g^{-1}ug)v) = \sigma'(u)\tau(v)$
for all $v \in V$ and $u \in U'_1 = U' \cap gUg^{-1}$.
Let $S$ be a subset of $S_\f(K)$ containing $S_\sigma
\cup S_{\sigma'} \cup S_g^K$ where $S_g$ is the set of
$\ell$ such that $g_\ell\not\in \GL_2(\ZZ_\ell)$.
Let $\CV$ and $\CV'$ be
$\calo_K$-lattices in $V$ and $V'$ that are stable
under the actions of $U$ and $U'$.
Choose integers $N$ and $N'$ so that the following hold:
\begin{itemize}
\item $N\ge 3$,
\item $U_N \subset \ker\sigma \cap \ker\sigma'$,
\item $U_{N'} \subset gU_Ng^{-1}$ and
\item $S_{N'}^K\subset S$.
\end{itemize}
To see that this is possible, first choose $N$ divisible by $4$ so
that $S_N^K \subset S$ and the second condition holds, then choose
$N'$ so that the third condition holds, and then for each $\ell\not\in S$,
replace $N'$ by $N'/\ell^r$ where $\ell^r|| N'$.
The condition that $S_N^K \subset S$ and  $g_\ell\in
\GL_2(\ZZ_\ell)$ ensures that we still have
$g^{-1} U_{N'}g \subset U_N.$  Note that $N|N'$
so $N' \ge 3$, $U_{N'}\subset U_N \subset \ker\sigma \cap \ker\sigma'$
and $S_N^K \subset S_{N'}^K\subset S$.
Then from the last section, $g$ gives a map
$[g]:\CM_N \to \CM_{N'}$.
Further $U'$ acts on $\CM_{N'}$.
Let $U' = \coprod_i g_i U_1'$ be a coset decomposition.
We define a morphism $[U'gU]_{\tau,N,N'}$
by restricting the map
$$ \begin{array}{ccc}
(\CM_{N}\otimes \CV)^S&\to& (\CM_{N'} \otimes \CV')^S\\
x \otimes v &\mapsto& \sum_i [g_i g] x \otimes \sigma'(g_i)\tau(v)
\end{array}
$$
to $\CM(\sigma)_N^S$.
Similarly define
$$[U'gU]_{\tau,N,N',c} : \CM(\sigma)_{N,c}^S \to
(\CM_{N',c}\otimes\CV')^S$$
and
$$[U'gU]_{\tau,N,N',!} : \CM(\sigma)_{N,!} \to
(\CM_{N',!}\otimes\CV')^S.$$

\begin{lemma}
$[U'gU]_{\tau,N,N',\sharp}$ defines a morphism
$\CM(\sigma)_{N,\sharp}\to \CM(\sigma')_{N',\sharp}$
which is independent of the coset decomposition
$U' = \coprod_i g_i U_1'$, for $\sharp = \emptyset$,
$\tf$, $c$ and $!$.
\mlabel{lem:coset1}
\end{lemma}
\proofbegin It suffices to prove the lemma on Betti realizations.
We assume $\sharp = \emptyset$, the proof in the other cases
being the same.

First we prove independence of the decomposition.
Let $U'=\coprod_i g'_i U_1'$ be another coset decomposition
of $U'$. Then $g'_i=g_i h_i$ for some $h_i\in U_1'$.
Thus we can write $h_i=g u_i g^{-1}$ for some $u_i\in U$.
Denote $[U'gU]'_{\tau,N,N'}$ for the map
arising from
$\ x\otimes v \mapsto \sum_i [g'_ig]x\otimes \sigma'(g'_i) \tau(v).$
Then we have

\begin{eqnarray*}
\lefteqn{[U'gU]'_{\tau,N,N'} (x\otimes v)=
\sum_i [g'_ig]x\otimes \sigma'(g'_i) \tau(v)}\\
&=& \sum_i [g_i gu_i] x \otimes \sigma'(g'_i)\tau(v)\\
&=& \sum_i [g_i gu_i] x\otimes \sigma'(g_i)
\sigma'(gu_ig^{-1})\tau(v) \\
&=& \sum_i [g_i gu_i] x\otimes \sigma'(g_i) \tau(\sigma(u_i)v)\\
&=& \sum_i ([g_i g]\otimes \sigma'(g_i)\tau)(
u_i x\otimes \sigma(u_i)v).
\end{eqnarray*}
Since $z\in \CM(\sigma)_N$ is invariant under $U$, we have
\begin{eqnarray*}
[U'gU]_{\tau,N,N'} (z)
&=& \sum_i ([g_i g]\otimes \sigma'(g_i)\tau)\circ
(u_i \otimes \sigma(u_i))(z)\\
&=& \sum_i ([g_i g]\otimes \sigma'(g_i)\tau)(z).
\end{eqnarray*}

Now we prove that the image is contained in $\CM(\sigma')_N$.
Let $u'\in U'$.
Then
\begin{eqnarray*}
u' (\sum_i [g_i g] x \otimes \sigma'(g_i) \tau(v))
&=& \sum_i [u'g_i g] x \otimes \sigma'(u') \sigma'(g_i)\tau(v)\\
&=& \sum_i [(u'g_i) g ]x\otimes \sigma'(u'g_i) \tau(v).
\end{eqnarray*}
Since $U'= u'U'=\coprod_i u'g_i U_1'$,
$\{u'g_i\}_i$ is also a complete system
of coset representatives of $U'$ by $U_1'$,
by the independence of coset decomposition,
$$\sum_i [(u'g_i) g] x\otimes \sigma'(u'g_i) \tau(v)
    = \sum_i [g_i g] x\otimes \sigma'(g_i) \tau(v).$$
\proofend

\subsection{Compositions}
\mlabel{ssec:coset-comp}
In addition to the notations in the last section,
suppose further that $U''$ is an open compact
subgroup of $\GL_2(\hat{\ZZ})$, $\sigma'': U'' \to \aut_{K} V''$
is a continuous representation,
$g'$ is in $\GL_2(\AA_\ff)\cap M_2(\hat{\ZZ})$ and
$\tau': V' \to V''$ satisfies
$\tau'(\sigma'((g')^{-1}ug')v) = \sigma'(u)\tau'(v)$
for all $v \in V'$, $u \in U''\cap g'U'g'{}^{-1}$.
Let $S$ be a subset of of $S_\f(K)$ that contains all primes
dividing each $\ell$ such that
$$
\left \{
\begin{array}{l}
\ell > k, \\
g_\ell, g'_\ell \in \GL_2(\ZZ_\ell) {\rm\ or\ } \\
\GL_2(\ZZ_\ell) \subset \ker\sigma \cap
    \ker \sigma' \cap \ker \sigma''
\end{array}
\right .
$$
Fix such an $S$. We choose stable $\CO_K$-lattices $\CV, \CV'$ and $\CV''$
as in the last section. Also choose $N,N'$ and $N''$ such
that
\begin{itemize}
\item $U_N \subset \ker \sigma \cap \ker \sigma'$,
\item $N\ge 4$,
\item $S_{N''}^K \subset S$,
\item $g ^{-1} U_{N'}g \subset U_N$ and
$g'{}^{-1} U_{N''}g' \subset U_{N'}$.
\end{itemize}
We will give a formula for
the composition
$[U''g'U']_{\tau',N',N'',\sharp}\circ [U'g U]_{\tau,N,N',\sharp}$ for
$\sharp\in \{\phi,c,!\}$.
To simplify the notation, we will suppress $N$, $N'$ and
$\sharp$ when there is no danger of confusion.

We first recall the definition of the product formula
for double cosets \cite[\S 3.1]{shimura},\cite[\S
2.7]{miyake}.
Let $U'_1=U'\cap g U g^{-1}$
and $U''_1=U''\cap g'U'{g'}^{-1}$.
Fix coset decompositions
$U'=\coprod_i g_i U'_1$ and $U''=\coprod_j g'_j U''_1$.
We have
$$ U'gU=\coprod_i g_i g U,\
    U''g'U'=\coprod_j g'_i g' U'.$$
Let $U''g'U'gU=\coprod_n U'' g'v_n g U$ be a
double coset decomposition of $U''U'U$.
The product of the double cosets $[U''g'U']$ and
$[U'gU]$ is defined by
$$ [U''g'U']\cdot [U'g U]=\sum_n c_n [U''g' v_ng U],$$
where $c_n=\#\{(i,j) \mid g'_jg' g_igU=g' v_n g U\}$.

Let $U'_2=U' \cap g'{}^{-1}U''g'$
and let $U'=\coprod_k U'_2 u_k U'_1$ be a
double coset decomposition of $U'$.
For each $k$, let $W_k=U''\cap(g' u_kg
U (g' u_k g)^{-1})$ and
$W_k'=W_k \cap g'U'g'{}^{-1}$.

\begin{prop}
$$[U''g'U']\cdot [U'gU]=\sum_k [W_k:W_k'] [U'' g'u_k g U].$$
\mlabel{prop:prod}
\end{prop}
\proofbegin
Let $U''=\coprod_m w_{k,m} W_k$ and
$W_k=\coprod_\ell w'_{k,\ell} W'_k$ be fixed coset
decompositions.
It is easy to see that we have coset decompositions
\begin{equation}
U''g'U' = \coprod_{i,j} g'_j g'g_i U'_1
\mlabel{eq:prod1}
\end{equation}
and
\begin{equation}
U''g'U' = \coprod_{k,m,\ell}
    w_{k,m} w'_{k,\ell} g' u_k U'_1.
\mlabel{eq:prod2}
\end{equation}
Further, for each $u_k$ in the double coset decomposition
$U'=\coprod_k U'_2 u_k U'_1$ of $U'$, we have
\begin{equation}
U''g'u_kgU=\coprod_m w_{k,m}g'u_kgU.
\mlabel{eq:prod2'}
\end{equation}
This is because
any $u\in U''$ can be written in the form $w_{k,m} w'$
for some $m$ and $w'\in W_k$. Then
$ug'u_kgU=w_{k,m}w'g'u_kgU=w_{k,m}g'u_kgU.$
Thus we have the union
$$U''g'u_kgU=\bigcup_m w_{k,m}g'u_kgU.$$
Suppose $w_{k,m_1}g'u_kgU=w_{k,m_2}g'u_kgU$.
Then $w_{k,m_2}^{-1}w_{k,m_1}$ is in $g'u_kgU(g'u_kg)^{-1}$.
This means that $w_{k,m_1}^{-1}w_{k,m_2}$ is in $W_k$.
Therefore $m_1=m_2$. So the union is disjoint.

~(\ref{eq:prod1}) and (\ref{eq:prod2}) give us
two decompositions of
$U''g'U'$ into cosets of $U'_1$. So we have a one-to-one
correspondence
\begin{equation}
\{ g'_jg'g_i U'_1 | j,i\} \leftrightarrow
\{w_{k,m}w'_{k,\ell}g'u_k U'_1| k,m,\ell\}.
\mlabel{eq:prod5}
\end{equation}
Fix a $v_n$ in $U''g'U'gU=\coprod_n U''g'v_ngU$.
Because of (\ref{eq:prod5}), we have
$$c_n =\# \{(k,m,\ell)|w_{k,m}w'_{k,\ell}g'u_kgU=g'v_ngU\}.$$
Given a triple $(k,m,\ell)$ with
$w_{k,m}w'_{k,\ell}g'u_kgU=g'v_ngU$,
$u_k$ must have the property $U''g'u_kgU=U''g'v_n gU$.
Fix such a $u_k$. Since $w'_{k,\ell}$ is in $W_k$, we have
$w_{k,m}w'_{k,\ell}g'u_kgU=w_{k,m}g'u_kgU$.
By (\ref{eq:prod2'}), $g'v_ngU=w_{k,m_1}g'u_kgU$
for a unique value of $m_1$.
So $w_{k,m}w'_{k,\ell}g'u_kgU=g'v_ngU$ gives
$w_{k,m}g'u_kgU   =w_{k,m_1}g'u_kgU$, independent of $\ell$.
By (\ref{eq:prod2'}),
this holds if and only if $m=m_1$ while $\ell$ takes any
value. Therefore
\begin{eqnarray*}
c_n&=& \sum_k \#\{(k,m,\ell)|w_{k,m}w'_{k,\ell}g'u_kgU=g'v_ngU\}\\
&=& \sum_k \#\{(k,m_1,\ell)|w'_{k,\ell}g'u_kgU=g'u_kgU\}\\
&=& \sum_k [W_k:W'_k]
\end{eqnarray*}
where the sum is over $k$ such that
$U''g'u_kg U=U''g' v_ng U$. This is also the coefficient
of $[U''g'v_ng U]$ given in the proposition.
\proofend

We now describe the composition of double coset actions.
As in the proof of Proposition~\ref{prop:prod}, let
$U''=\coprod_{m,\ell} w_{k,m}w'_{k,\ell}W'_k$.
For each pair $(k,\ell)$, define
$$\tau_{k,\ell}: \CV \to \CV'',$$
$$ v \mapsto \sigma''(w'_{k,\ell})\tau' \sigma'(u_k)
    \tau \sigma((g'u_kg)^{-1} w'_{k,\ell}{}^{-1}(g'u_kg))v,
    v\in \CV.$$
We will often suppress $\sigma, \sigma'$ and $\sigma''$ from the
notation.
Define
$$\tau_k =\sum_\ell \tau_{k,\ell}.$$
It is straightforward to check that $\tau_{k,\ell}$
is independent of the choice of the coset representatives
$w'_{k,\ell}$, and hence so is $\tau_k$.
Define
$$\tilde{\tau}_k: W_k \to U, u'' \mapsto
(g'u_kg)^{-1}u''(g'u_kg).$$
For $u''\in W_k$ and $v\in \CV$,
$\{u''{}^{-1}w'_{k,\ell}\}_\ell$ is another set of
coset representatives of $W'_k$ in $W_k$.
Thus we have
\begin{eqnarray*}
\lefteqn{\tau_k(\tilde{\tau}_k(u'') v)}\\
&=& \sum_\ell w'_{k,\ell}\tau'u_k\tau (g'u_kg)^{-1}
    w'_{k,\ell}{}^{-1}(g'u_kg) (g'u_kg)^{-1}u''(g'u_kg) v\\
&=& \sum_\ell w'_{k,\ell}\tau'u_k\tau (g'u_kg)^{-1}
    w'_{k,\ell}{}^{-1}u'' (g'u_kg) v\\
&=& \sum_\ell u'' (u''{}^{-1} w'_{k,\ell})\tau'u_k\tau (g'u_kg)^{-1}
    (u''{}^{-1}w'_{k,\ell})^{-1} (g'u_kg) v\\
&=& u'' \tau_k (v).
\end{eqnarray*}
Thus $[U''g'u_k gU]_{\tau_k}$ is well-defined.

\begin{prop}
$$[U''g'U']_{\tau'}\circ [U'gU]_\tau
    = \sum_k [U''g'u_kgU]_{\tau_k}.$$
More precisely,
$$[U''g'U']_{\tau',N',N''}\circ [U'gU]_{\tau,N,N'}
    = \sum_k [U''g'u_kgU]_{\tau_k,N,N''}$$
in $\ipms^{S_{N''}}_K$.
\mlabel{prop:scomp}
\end{prop}
\proofbegin
Write $U' = \coprod_i g_i U'_1$ and $U''=\coprod_j g'_jU''_1$.
Then the operators
$$[g'_jg'g_ig]: \CM_N \to \CM_{N'}\to \CM_{N''}$$
and
$$g'_j\tau' g_i \tau: \CV\to \CV' \to \CV''$$
are both well-defined.
By definition, we have
\begin{eqnarray*}
[U''g'U']_{\tau'} \circ [U'gU]_{\tau} (x\otimes v) &
=& [U''g'U']_{\tau'} (\sum_i [g_i g] x\otimes g_i \tau(v))\\
&=&\sum_i\sum_j [g_j' g'g_ig]x\otimes g_j'\tau' (g_i \tau(v)).
\end{eqnarray*}
We first prove that, for each pair $(i,j)$,
if $g'_jg'g_iU'_1=w_{k,m}w'_{k,\ell}g'u_kU'_1$ under the
bijection (\ref{eq:prod5}),
then
$$
[g'_jg' g_ig]x\otimes g'_j \tau'g_i \tau v
    =[w_{k,m}w'_{k,\ell}g'u_k g] x\otimes
        w_{k,m}w'_{k,\ell}\tau_k v.
$$
Note that the 1-1 correspondence $g'_jg'g_igU'_1 \leftrightarrow
w_{k,m}w'_{k,\ell}g'u_k gU'_1$ can be described as follows.
Given $g'_jg'g_ig$, since $U'=\coprod_k U'_2u_kU'_1$,
$g_i$ is in $U'_2u_k'U'_1$ for a unique $k$. So
$g_i=w'_2 u_k w'_1$ with $w'_1\in U'_1$ and $w'_2\in U'_2$.
So $w'_1=gw_1g^{-1}, w'_2=g'{}^{-1}w''_2g'$ with
$w_1\in U$ and $w''_2\in U''$.
Since $g'_jw''_2$ is in $U''$, it can be uniquely
written as $w_{k,m}w'_{k,\ell}w''_3$ with $w''_3\in W'_k$.
So $w''_3=(g'u_kg)w_3(g'u_kg)^{-1}$ for some $w_3\in
U$ and $g'{}^{-1}w''_3g'=(u_kg)w_3(u_kg)^{-1}$ is in $U'$.
Based on this description, we have
\begin{eqnarray*}
g'_jg'g_ig&=& w_{k,m}w'_{k,\ell}w''_3g'u_kg w_1\\
&=&w_{k,m}w'_{k,\ell}g'u_kgw_3w_1.
\end{eqnarray*}
This gives the correspondence
$g'_jg'g_igU'_1=w_{k,m}w'_{k,\ell}g'u_kgU.$
This description, together with the transformation laws observed
by $\tau$ and $\tau'$, also gives us
\begin{eqnarray*}
g'_j\tau'g_i\tau
&=& w_{k,m}w'_{k,\ell}\tau' u_k \tau w_3 w_1.
\end{eqnarray*}
Note that if $w$ is in $U$, then
$(w \otimes w)z =z$ for $z\in \CM(\sigma)_N
=(\CM_N \otimes \CV)^U$ (see the proof of Lemma~\ref{lem:coset1}).
Therefore, for $z\in \CM(\sigma)_N$,
\begin{eqnarray*}
\lefteqn{(\sum_{i,j}[g'_jg'g_ig]\otimes g'_j\tau'g_i\tau )(z)}\\
&=& (\sum_{k,m,\ell} [w_{k,m}w'_{k,\ell}g'u_kg w_3w_1]\otimes
    w_{k,m}w'_{k,\ell}\tau' u_k\tau w_3w_1)(z)\\
&=& (\sum_{k,m,\ell} [w_{k,m}g'u_kg]\otimes
    w_{k,m}\tau_{k,\ell})(z)\\
&=& (\sum_{k,m} [w_{k,m}g'u_kg]\otimes
    w_{k,m}\tau_k)(z)\\
&=& (\sum_k [U''g'u_kgU]_{\tau_k})(z).
\end{eqnarray*}
This proves the proposition.
\proofend

Using Lemma~\ref{lem:gmap} and Proposition~\ref{prop:scomp},
we conclude that $[U'gU]$ is independent of the choice of $N$ and
$N'$. More precisely, if $[U'gU]_{\tau,M,M'}$ is also defined,
then $[U'gU]_{\tau,N,N'}^S =[U'gU]_{\tau,M,M'}^S$ when
$S\supset S_{\lcm(N',M')}$.

\subsection{Pairing on $\sigma$-constructions}
\mlabel{ssec:sigma-pair}

In this section we will construct perfect pairings
between certain $\sigma$-constructions and
compute the adjoints of double coset actions with respect
to the pairings.

For $N\ge 3$, we let $H = H_N$ denote the premotivic structure $H^0(X_N)
= H^0(X) = H^0(\bar{X})$.  More
precisely, we let $\CH_\bt = H^0(X^\an,\ZZ)$, $\CF_\dr
= \BH^0(\bar{X},\omega^\bullet_{\bar{X}/T})$,
$\CH_\ell = H^0(X_{\bar{\QQ}},\ZZ_\ell)$
and $\CH_\ecrys= H^0_\crys(X_{\ZZ_\ell},\CO_{X_{\ZZ_\ell},\crys})$
for $\ell\not\in S_N$.  These come equipped with additional structure
and comparison isomorphisms making $\CH$ an object of $\ipms_\QQ^{S_N}$,
and we let $H=\QQ\otimes\CH$.

Let $F = \QQ(\mu_N)$.  The Weil pairing on
$(\ZZ/N\ZZ)^2_X \cong E[N]$ defines an isomorphism
between $(\ZZ/N\ZZ)_X$ and $\mu_{N,X}$, hence a
morphism $X \to \spec \CO_F$, which one checks
induces a morphism $\CM_F^{S_N} \to \CH$ in
$\ipms_\QQ^{S_N}$.  One checks also that the map
$X^\an \to (\spec\CO_F)^\an$ sends the component
$X_{N,t}$ to $e^{-2\pi i t/N}$, where we identify
$(\spec\CO_F)^\an = \hom(\CO_F,\CC) = \II_F$ with the
set of primitive $N^\th$ roots in $\CC$.
In particular, it follows that $\CM_F^{S_N}
\cong \CH$ since the map $\CM_{F,\bt}^{S_N} \to \CH_\bt$
is the isomorphism
$$H^0(X^\an,\ZZ) =
\bigoplus_{t\in (\ZZ/N\ZZ)^\times}H^0(X_{N,t},\ZZ)
  = \bigoplus_{t\in (\ZZ/N\ZZ)^\times}\ZZ
   \cong  \ZZ^{\II_F}$$
induced by the bijection $t\leftrightarrow e^{-2\pi i t/N}$.

Recall that if $(h,N,N')$ is an admissible triple in the
sense of \S\ref{ssec:action:curve}, then we obtain a morphism
$\chi_h:\bar{X}' \to \bar{X}$, where $\bar{X}' = \bar{X}_{N'}$.
More generally, if $N, N' \ge 3$ and $g\in\GL_2(\AA_\f)$ is
such that $g^{-1}U_{N'}g\subset U_N$, we can define $\chi_g =
\chi_h$ for any $h\in g\QQ^\times$ such that $(h,N,N')$ is
admissible.  Then $\chi_g$ induces a morphism $\chi_g^*:
\CH^{S'} \to \CH'$ where $S' = S_{N'}$ and $\CH' = \CH_{N'}$.

\begin{lemma}
For $N$, $N'$ and $g$ as above, we have the
commutative diagram
$$\begin{CD} \CH^{S'} @>\chi_g^*>> \CH'\\
@VVV @VVV\\
\CM_F^{S'} @>>> \CM_{F'}^{S'}
\end{CD}$$
where $F'=\QQ(\mu_{N'})$ and
the morphism $\CM_F \to \CM_{F'}$
is obtained from the composite
$F \to F' \to F'$ where the first map is the canonical inclusion
and the second map is the image of $\det g^{-1}$
(i.e., of $\det g^{-1} ||\det g^{-1}||$) in $\gal(F'/\QQ)$
by class field theory.
\mlabel{lem:field}
\end{lemma}
\proofbegin
It suffices to check the commutativity of the diagram
$$\begin{array}{ccc}  X_{N'} & \to & X_N \\
    \downarrow&&    \downarrow\\
    \spec\CO_{F'}& \to &\spec\CO_F\end{array}$$
on complex points, and this is immediate from
Lemma~\ref{lem:adan}.
\proofend

For the rest of the section, we assume $U$ is an open compact
subgroup of $\GL_2(\hat{\ZZ})$ satisfying $\det U = \hat{\ZZ}^\times$.
We view $U$ as acting on $\CH = \CH_N$ for any $N \ge 3$ such
that $U_N \subset U$.

We fix a continuous character $\psi:\hat{\ZZ}^\times \to K^\times$
of finite order.  We also view $\psi$ as a character of
$\AA^\times/\QQ^\times$, and write $\psi$ for the corresponding
character of $G_\QQ$.  We view $U$ as acting on $\CO_\psi = \CO_K$
via $\psi^{-1}\circ\det$.
\begin{lemma}
Suppose $N$ is such that $N \ge 3$,  $U_N \subset U$,
the conductor of $\psi$ divides $N$.
Then the morphism $\CM_F \to \CH$ induces an isomorphism
$$\CM_{\psi}^{S_N} \to (\CO_{\psi}\otimes \CH)^U$$
in $\ipms_K^{S_N}$.
\mlabel{lem:psi}
\end{lemma}
\proofbegin
It follows from Lemma~\ref{lem:field} that the isomorphism
$$ \CO_{\psi} \otimes \CM_{F}^{S_N} \cong
    \CO_{\psi} \otimes \CH$$
respects the action of $U$, where we view $U$ acting on
$\CO_\psi\otimes\CM_F$ via $\det^{-1}:U\to\hat{\ZZ}^\times
\ola{\cong}G_\QQ^\ab$.
Taking invariants gives the lemma.
\proofend

Using the compatibility of the cup product with comparison
isomorphisms, one obtains morphisms
\begin{equation}
\CH \otimes \CM_\sharp \to \CM_\sharp
\mlabel{eqn:cup1}
\end{equation}
in the category $\ipms_\QQ^{S_N}$,
where $\CM_\sharp = \CM_{N,\sharp}$ for $\sharp =
\emptyset$, $\tf$, $c$ and $!$.  Moreover, these
morphisms are compatible with the action of $g\in\GL_2(\AA_\f)
\cap M_2(\hat{\ZZ})$ in the sense that if $g^{-1}U_{N'}g\subset
U_N$, then the resulting diagram
\begin{equation}
\begin{array}{ccccc} \CH^{S'}& \otimes&\CM_\sharp^{S'}&\to &\CM_\sharp^{S'}\\
\downarrow&&\downarrow&&\downarrow\\
\CH'& \otimes&\CM'_\sharp&\to &\CM'_\sharp,\end{array}
\mlabel{eqn:g-cup}
\end{equation}
commutes, where the morphisms $[g]_\sharp:\CM_\sharp \to \CM'_\sharp$
are those of Theorem~\ref{thm:g-action}.  (Again the commutativity
can be easily checked on Betti realizations.)  In particular, the
morphism of (\ref{eqn:cup1}) respects the action of $U$.

Next we consider the composite morphism
\begin{equation}
\CM_{\psi}^{S_N^K} \otimes \CM_\sharp \to
\CO_{\psi} \otimes \CH \otimes \CM_\sharp \to
 \CO_{\psi} \otimes \CM_\sharp
\mlabel{eq:cup}
\end{equation}
in $\ipms_K^{S_N^K}$, where the first map is the
one in Lemma~\ref{lem:psi} tensored with $\CM_\sharp$
and the second is (\ref{eqn:cup1}) tensored with
$\CO_\psi$.

\begin{lemma}
The morphism in (\ref{eq:cup}) is an isomorphism
respecting the action of $U$, where the action of
of $U$ on is trivial on $\CM_\psi$, defined in
Theorem~\ref{thm:g-action} on $\CM_\sharp$, and
via $\psi^{-1}\circ\det$ on $\CO_\psi$.
\mlabel{lem:cup}
\end{lemma}
\proofbegin  The $U$-equivariance follows from that of
(\ref{eqn:cup1}) together with Lemma~\ref{lem:psi}.

It suffices to check the maps are isomorphisms on
Betti realizations.  For this, write
$$\CM_\bt = \bigoplus_{t\in (\ZZ/N\ZZ)^\times}
  H^1(X_{N,t}, \CF_\bt^k).$$
Then for a class $x\in H^1(X_{N,t},\CF_\bt^k)$,
we have $b_\bt\otimes x\mapsto \psi^{-1}(-t) \otimes x$
under (\ref{eqn:cup1}).  It follows that the map is
surjective, hence it is injective since it is between
two isomorphic finitely generated $\CO_K$-modules.

The case of $\CM_c$ is similar, and the remaining cases
follow from those of $\CM$ and $\CM_c$.
\proofend

We now construct the pairings by defining the duality
homomorphisms.
Given a representation $\sigma: U \to \aut_{\CO_K} \CV$,
let $\hat{\sigma}\otimes(\psi^{-1}\circ\det)$ denote the
representation defined by the action of $U$ on
$\hom_{\CO_K}(\CV,\CO_{\psi})$.  For $N\ge 3$ such that
 $U_N\subset \ker \sigma$ and the conductor of $\psi$ divides $N$,
Lemma (\ref{lem:cup}) and the duality homomorphism
$\delta^L_N: \CM_{N,c}\to \hom(\CM_N,  \CT^{S_N,\otimes (k-1)})$ from
\S\ref{ssec:levelN.wt} give a morphism
\begin{equation}
\begin{array}{l}
\hom_{\CO_K}(\CV,\CO_{\psi^{-1}})\otimes \CM_{N,c}
\cong \hom_{\CO_K}(\CV,\CM_{N,c}\otimes \CM_{\psi}^{S_N^K})\\
\to  \hom_{\CO_K}(\CV,\hom_{\CO_K}(\CM_N,\CT^{S_N,\otimes (k-1)}\otimes
    \CM_{\psi}^{S_N^K}))\\
\cong  \hom_{\CO_K}(\CV\otimes \CM_N, \CM_\psi(1-k)),
\end{array}
\mlabel{eqn:bigpair}
\end{equation}
where we write $\CM_\psi(1-k)$ for the object
$\CT^{S_N,\otimes (k-1)}\otimes \CM_{\psi}^{S_N^K}$
of $\ipms_K^{S_N^K}$.  Moreover applying
Proposition~\ref{pp:gpair} with $g\in U$, we find that the morphism
(\ref{eqn:bigpair})
respects the action of $U$ where we view $U$ as acting trivially
on $\CM_\psi(1-k)$.

Taking $U$-invariants, we get a morphism
\begin{equation}
\delta^L_N: \CM_c(\hat{\sigma}\otimes(\psi^{-1}\circ\det))_N
\to \hom_{\CO_K} (\CM(\sigma)_N , \CM_{\psi}(1-k)).
\mlabel{eq:spair1a}
\end{equation}
Tensoring with $\QQ$ and normalizing by dividing by $[U:U_N]$,
we get a morphism
$$
\bar{\delta}^L_N:
M_c(\hat{\sigma}\otimes(\psi^{-1}\circ\det))_N
\to \hom_{K} (M(\sigma)_N , M_{\psi}(1-k)).
$$
Similarly, using the duality morphism
$\CM_{N,!}\to \hom(\CM_{N,!},  \CT^{S_N,\otimes (k-1)})$, we
obtain morphisms
$$
\delta^L_{N,!}:
\CM_!(\hat{\sigma}\otimes (\psi^{-1}\circ \det))_N
\to \hom_{\CO_K} (\CM_!(\sigma)_N,  \CM_{\psi} (1-k))
$$
and
$$
\bar{\delta}^L_{N,!}:
M_!(\hat{\sigma}\otimes (\psi^{-1}\circ \det))_N
\to \hom_K (M_!(\sigma)_N, M_{\psi} (1-k)).
$$

\begin{prop}
The morphisms $\bar{\delta}^L_N$ and
$\bar{\delta}^L_{N,!}$ are isomorphisms in $\pms_K^{S_N^K}$
compatible with the change of levels.
\mlabel{prop:spair}
\end{prop}
\proofbegin  It follows from Theorem~\ref{thm:levelN} b) that
(\ref{eqn:bigpair}) is an isomorphism after tensoring with
$\QQ$.  The same is true for the natural map
$$\hom_{\CO_K}(\CV\otimes \CM_N, \CM_\psi(1-k))^U
  \to  \hom_{\CO_K}(\CM(\sigma)_N, \CM_\psi(1-k)),$$
so $\bar{\delta}^L_N$ is an isomorphism.
The compatibility with change of levels follows from applying
Proposition~\ref{pp:gpair} and (\ref{eqn:g-cup})
with $g=I$.  The proof for $\bar{\delta}^L_{N,!}$ is the same.
\proofend

We will therefore omit the subscript $N$ in the notation
for the morphisms $\bar{\delta}^L$.  (Note however
that $\delta^L_N$ depends on $N$, and $\delta^L_{N,!}$
depends on $N$ only in the normalization.)

We say that $U$ is {\em sufficiently small} if
$U$ acts freely on $GL_2(\QQ)\backslash \GL_2(\AA) / U_\infty$.
In particular $U$ is sufficiently small if
$U \subset U_1(d)$ for some $d \ge 4$, where $U_1(d)$ denotes the
preimage in $\GL_2(\hat{\ZZ})$ of the subgroup of $\GL_2(\ZZ/d\ZZ)$
consisting of matrices of the form $\smat{*}{*}{0}{1}$.

\begin{prop}
Let $U$ be a sufficiently small
subgroup of $\GL_2(\hat{\ZZ})$ with
$\det U=\hat{\ZZ}^\times$.  Suppose $S \supset S_N^K$
for some $N\ge 3$ such that $U_N \subset \ker\sigma$ and
$\psi$ has conductor divisible by $N$.
Then the isomorphism $\bar{\delta}^L$ arises from an injective
morphism
\begin{equation}
\CM(\hat{\sigma}\otimes(\psi^{-1}\circ\det))_!^S
\to \hom_{\CO_K} (\CM(\sigma)_!^S, \CM_{\psi}(1-k)^S)
\mlabel{eq:spair2}
\end{equation}
in $\ipms^S_K$ whose cokernel $\calc$ satisfies
$\calc_\ell = 0$ for $\ell\not|N(k-2)!$.
\mlabel{prop:spair2}
\end{prop}
\proofbegin
We only need to prove the assertion for Betti realizations,
for then $\delta^L_N  = [U:U_N]\tilde{\delta}^L$ for
a unique morphism $\tilde{\delta}^L$ with the desired property.

To do this, we give an alternate description of
the Betti realizations and pairings which is independent of $N$.
and define a perfect pairing on the new Betti realizations.
Given a sufficiently small $U$ and
a $\sigma: U\to \Aut_{\CO_K}\CV$, define a locally
constant sheaf $\CF(\sigma)$ on
$G_\QQ\backslash G_\BA/UU_\infty$
by $\CF(\sigma)=G_\QQ\backslash (G_\BA\times \CV)/UU_\infty$.
Define
$$ \CM(\sigma)^\an=H^1(X^\an_U,\CF^k_B\otimes \CF(\sigma))$$
and
$$ \CM(\sigma)_c^\an=H^1_c(X^\an_U,\CF^k_B\otimes \CF(\sigma)).$$
Also define $\CM(\sigma)_!^\an$ to be the image
of $\CM(\sigma)_c^\an$ in $\CM(\sigma)^\an_\tf$
where $\CM(\sigma)^\an_\tf$ is the largest torsion-free
quotient of $\CM(\sigma)^\an$.
Then from the pairing $\CF_B^k\otimes \CF_B^k \to (2\pi i)^{2-k}\ZZ$
we get a pairing
$$(\CF^k_B \otimes \CF(\hat{\sigma}\otimes (\psi^{-1}\circ \det)))
\otimes (\CF^k_B \otimes \CF(\sigma))
\to (2\pi i)^{2-k} \ZZ\otimes \CF(\psi^{-1}\circ \det).$$
Here
$\CF(\psi^{-1}\circ \det)
G_\QQ\backslash (G_\BA\times \CO_{\psi})/UU_\infty$
is the sheaf on
$X^\an=G_\QQ\backslash G_\BA/UU_\infty$
defined by $\CO_{\psi}$.

We have an isomorphism of locally constant sheaves
$$\begin{CD}
\CF(\psi^{-1}\circ \det) @>\cong >>
    \Gamma_U \backslash (\uhp\times \CO_{\psi^{-1}})\\
@VVV @VVV\\
X^\an @>\cong >>
        \Gamma_U\backslash \uhp
\end{CD}
$$
where $\Gamma_U$ acts on $\CO_{\psi^{-1}}$ through
its canonical embedding into $U$ and then take the inverse.
Since $\det (\Gamma_U)=1$, the action of $\Gamma_U$ on
$\CO_{\psi}$ is trivial. Thus the sheaf on the right is constant.
So the same is true for the sheaf on the left.

Thus by Poincar\'{e} duality, we have a pairing
\begin{align*}
[\ ,\ ]^\an_c:\
& \CM(\hat{\sigma}\otimes (\psi^{-1}\circ \det))_c^\an
\otimes \CM(\sigma)^\an \notag \\
& \to
H^2_c(X^\an, (\CF^k_B \otimes \CF(\hat{\sigma}\otimes (\psi^{-1}\circ \det)))
\otimes (\CF^k_B \otimes \CF(\sigma)) ) \notag\\
& \to H^2_c(X^\an, (2\pi i)^{2-k} \ZZ\otimes \CO_{\psi^{-1}}) \notag\\
& \cong
H^2_c(X,\ZZ) \otimes (2\pi i)^{2-k} \CO_{\psi^{-1}} \notag \\
& \cong H^2_c(X,\ZZ) \otimes (2\pi i)^{2-k} \CM_{\psi,B} \notag \\
& \cong \CM_{\psi} (1-k)_B
\end{align*}
which is perfect on torsion-free quotients.
We then obtain a pairing

$$
[\ ,\ ]^\an_!:
\CM(\hat{\sigma}\otimes (\psi^{-1}\circ \det))_!^\an
\otimes \CM(\sigma)_!^\an
\to \CM_{\psi} (1-k)_B
$$
by using the commutative diagram
$$
\begin{CD}
\CM(\hat{\sigma}\otimes(\psi^{-1}\circ\det))_c^\an
\otimes_{\CO_K} \CM(\sigma)_c^\an @>>>
\CM(\hat{\sigma}\otimes(\psi^{-1}\circ\det))_c^\an
\otimes_{\CO_K} \CM(\sigma)^\an\\
@VVV @VVV\\
\CM(\hat{\sigma}\otimes(\psi^{-1}\circ\det))^\an
\otimes_{\CO_K} \CM(\sigma)_c^\an @>>>
\CM_{\psi}(1-k)_B
\end{CD}
$$
where all maps are the natural ones except that the
bottom map is the pairing
$$
 \CM(\sigma)_{c}^\an \otimes_{\CO_K}\CM(\hat{\sigma}
    \otimes(\psi^{-1}\circ\det))^\an
\to \CM_{\psi}(1-k)_B$$
proceeded by the order-reversing isomorphism and followed
by multiplication by $(-1)^{k-1}$.

The same proof as that of Lemma~\ref{lem:torsion} shows
that the cokernel of $\CM(\sigma)_c^\an \to \CM(\sigma)^\an$
has no $\ell$-torsion if $\ell \not| N(k-2)!$, and it
follows that
$$
\tilde{\delta}^\an:\CM(\hat{\sigma}\otimes(\psi^{-1}\circ\det))_!^\an
\to \hom_{\CO_K} (\CM(\sigma)_!^\an, \CM_{\psi}(1-k)_B)
$$
is injective with finite cokernel $C$ satisfying $C_\ell = 0$
for $\ell\not|N(k-2)!$.
The proposition is therefore a consequence of the following lemma.
\proofend

\begin{lemma}
There is an isomorphism
$\CM(\sigma)_!^\an \cong \CM(\sigma)_{N,!,B}$
identifying $\tilde{\delta}^\an\otimes\QQ$
with the pairing $\bar{\delta}^L_!$.
\mlabel{lem:compair1}
\end{lemma}
\proofbegin
Define $\CF'(\sigma)_N=G_\QQ\backslash (G_\BA \times \CV)/U_N U_\infty$.
The same argument as in Lemma~\ref{lem:gmap}
shows that the natural projection
$X_N^\an\to X^\an=X_U^\an$ induces a map
$$\res: H^1(X^\an, \CF_B^k\otimes \CF(\sigma))
    \to H^1(X_N^\an,\CF_{N,B}^k\otimes \CF'(\sigma)_N)^{U/U_N}$$
with finite kernel and cokernel, and an isomorphism
$$\res: H^1_c(X^\an, \CF_B^k\otimes \CF(\sigma))
    \to H^1_c(X_N^\an,\CF_{N,B}^k\otimes \CF'(\sigma)_N)^{U/U_N}.$$
The sheaf $\CF'(\sigma)_N$ on $X_N$ is constant and
the right hand sides become
$\CM(\sigma)_{N,B}$ and $\CM(\sigma)_{N,c,B}$.
This proves
$\CM(\sigma)_\tf^\an \hookrightarrow \CM(\sigma)_{N,\tf,B}$
and
$\CM(\sigma)_c^\an \cong \CM(\sigma)_{N,c,B}$.
We then get the desired isomorphism from the commutative diagram

$$\begin{CD}
\CM(\sigma)_c^\an @> \cong >> \CM(\sigma)_{N,c,B} \\
@VVV @VVV\\
\CM(\sigma)_!^\an @> \subset >> \CM(\sigma)_{N,!,B} \\
@VVV @VVV\\
\CM(\sigma)_\tf^\an @> \subset >> \CM(\sigma)_{N,\tf,B}
\end{CD}
$$
in which the two top vertical maps are surjective
and two bottom vertical maps are injective.

The compatibility of the pairings follows from the commutativity
of the diagram
$$
{\scriptsize \begin{CD}
H^1_c(X^\an,\CF^k_B \otimes
    \CF(\hat{\sigma}\otimes (\psi^{-1}\circ \det)))
@. \otimes H^1(X^\an, \CF^k_B \otimes \CF(\sigma))
&\stackrel{\tilde{\cup}}{\rightarrow} &
H^2_c(X^\an, (2\pi i)^{2-k}\ZZ\otimes
\CO_{\psi^{-1}}) &\ola{\tr} & \CM_{\psi,B}\\
@V\res VV @VV \res V
 @V\res VV @V[U:U_N] VV \\
H^1_c(X^\an_N,\CF^{k}_{N,B} \otimes
    \CF'(\hat{\sigma}\otimes (\psi^{-1}\circ \det)))
@. \otimes H^1(X^\an_N, \CF^{k}_{N,B} \otimes \CF'(\sigma))
&\stackrel{\cup}{\rightarrow} &
H^2_c(X^\an_N, (2\pi i)^{2-k}\ZZ\otimes \CO_{\psi^{-1}})
&\ola{\tr} & \CM_{\psi,B},
\end{CD}}$$
which is proved
by the same argument as in Proposition~\ref{pp:gpair} (with $g=I$).
\proofend

We now compute the adjoints of double coset operators under the
pairings.  Suppose now that $\sigma$, $\sigma'$, $g$, $\tau$
and $S$ are as in \S\ref{ssec:coset}.  Recall that we defined
morphisms
$$[U'gU]_{\tau,\sharp}: M(\sigma)_\sharp^S \to M(\sigma')_\sharp^S$$
for $\sharp = \emptyset$, $c$ and $!$.
(We omit the subscripts $N,N'$ of which they are independent.)
It is easy to check that we then also have morphisms
$$[U(||\det g||g)^{-1}U']_{\tau^t\otimes \psi(\det(g))}
: M(\sigma'\otimes (\psi^{-1}\circ \det))_\sharp^S \to
M(\sigma \otimes (\psi^{-1}\circ \det))_\sharp^S.$$
Denote this operator by $[U'gU]_{\tau,\sharp}^T$.

\begin{prop} The morphism
$[U'gU]_{\tau,c}^T$ (respectively, $[U'gU]_{\tau,!}^T$)
is the adjoint of $[U'gU]_{\tau}$
(respectively,  $[U'gU]_{\tau,!}$)
with respect to the pairing $\bar{\delta}^L$,
(respectively, $\bar{\delta}^L_!$), i.e., we
have the commutative diagram
$$\begin{CD}
M(\hat{\sigma} \otimes (\psi^{-1}\circ \det))_c^S
@>\bar{\delta}^L >>
\hom (M(\sigma)^S,  M_\psi(1-k)^S) \\
@A [U'gU]_{\tau,c}^T AA @A [U'gU]_\tau^* AA   \\
M(\hat{\sigma}'\otimes (\psi^{-1}\circ \det))_c^S
@>\bar{\delta}^L >>
\hom (M(\sigma')^S, M_\psi(1-k)^S),
\end{CD}
$$
and a similar diagram for $\bar{\delta}^L_!$.
\mlabel{pp:adj}
\end{prop}
\proofbegin We first make some simplifications.
Let $\sigma:U\to \Aut_KV$, $\sigma':U'\to \Aut_KV'$
and $\tau:V\to V'$ be as in the definition of $[U'gU]_\tau$.
Consider the diagram

$$\begin{array}{ccccccc}
&& U &\xrightarrow{\sigma} & \Aut_{K} V &&  V \\
& &  && &&\dap{\tau^{iv}}\\
U_{g^{iv}}& \subset & U^{iv}& \xrightarrow{\sigma^{iv}}
&\Aut_{K} V^{iv}
    &&  V^{iv}\\
& &  && &&\dap{\tau'''}\\
U_{g'''}& \subset & U'''& \xrightarrow{\sigma'''} &\Aut_{K} V'''
    &&  V'''\\
&  &&&&& \dap{\tau''} \\
U_{g''} & \subset & U'' & \xrightarrow{\sigma''}
&\Aut_{K} V''  &&  V'' \\
&  &&&&& \dap{\tau'} \\
U_{g'} & \subset & U' & \xrightarrow{\sigma'}
&\Aut_{K} V'  &&  V'.
\end{array}
$$
Here we have used the notations
$$U^{iv}=gUg^{-1},  V^{iv}= V, $$
$$\sigma^{iv}:U^{iv}\to \Aut_{K} V^{iv},
\sigma^{iv}(u^{iv})v^{iv}=\sigma(g^{-1}u^{iv}g)v^{iv},$$
$$g^{iv}=g, \tau^{iv}=\id:  V \to  V^{iv},$$
$$U'''=U'\cap gUg^{-1},  V'''= V,
\sigma'''=\sigma_{|U'''}:U'''\to \Aut_{K} V,$$
$$g'''=1, \tau'''=\id,$$
$$U''=U'\cap gUg^{-1},  V''= V',
\sigma''=\sigma'_{|U''}:U''\to \Aut_{K} V'',$$
$$g''=1, \tau''=\tau \quad {\rm\ and\ }$$
$$g'=1, \tau'=\id: V''= V'\to  V'.$$
Then we have morphisms
$$[U^{iv}g^{iv}U]_{\tau^{iv}}: M(\sigma)\to M(\sigma^{iv}),$$
$$[U'''g'''U^{iv}]_{\tau'''}: M(\sigma^{iv})\to M(\sigma'''),$$
$$[U''g''U''']_{\tau''}: M(\sigma''')\to M(\sigma'')
\quad {\rm\ and\ }$$
$$[U'g'U'']_{\tau'}: M(\sigma'')\to M(\sigma').$$
Also
$$\tau=\tau'\circ \tau'' \circ \tau''' \circ\tau^{iv}.$$
Applying Proposition~\ref{prop:scomp}, we have

\begin{eqnarray*}
[U'g'U'']_{\tau'}\circ [U''g''U''']_{\tau''} \circ
[U''' g''' U^{iv}]_{\tau'''} \circ [U^{iv}g^{iv}U]_{\tau^{iv}}
    &=&[U'gU]_{\tau}.
\end{eqnarray*}
We similarly verify that
\begin{eqnarray*}
&&  [U(||\det g^{iv}||g^{iv})^{-1}U^{iv}]_{(\tau^{iv})^t
    \otimes \psi(\det(g^{iv})),c}\\
&& \circ
 [U^{iv}(||\det g{'''}||g{'''})^{-1}U{'''}]_{(\tau{'''})^t
    \otimes \psi(\det(g{'''})),c} \\
&& \circ
 [U{'''}(||\det g{''}||g{''})^{-1}U{''}]_{(\tau{''})^t
    \otimes \psi(\det(g{''})),c} \\
&& \circ
 [U{''}(||\det g{'}||g{'})^{-1}U{'}]_{(\tau{'})^t
    \otimes \psi(\det(g{'})),c}\\
&=&
||\det g||^{2-k} [Ug^{-1}(gUg^{-1})]_{\id\otimes \psi(\det g),c}
\circ
[(gUg^{-1}) 1 (U'\cap (gUg^{-1}))]_{\id \otimes \psi(1),c}
\\
&& \circ
[(U'\cap (gUg^{-1})) 1 (U'\cap (gUg^{-1}))]_{\tau^t
    \otimes \psi(1),c}
\circ
[(U'\cap ((gUg^{-1}))1 U{'}]_{\id \otimes \psi(1),c}\\
&=& ||\det g||^{2-k} [ U g^{-1} U']
_{\tau^t\otimes \psi(\det g),c}\\
&=&  [ U (||\det g||g)^{-1} U']
_{\tau^t\otimes \psi(\det g),c}.
\end{eqnarray*}
Thus to prove the proposition for $\bar{\delta}^L$, we only
need to consider the four composition factors, i.e.,
to prove it in each of the following four
special cases.
\begin{enumerate}
\item[]{\bf Case 1.}
$U' = gUg^{-1}$ and $\tau=\id$,
\item[]{\bf Case 2.}
$U'\subset U$, $g = 1$ and $\tau=\id$,
\item[]{\bf Case 3.}
$U'=gUg^{-1}$ and $g=1$, and
\item[]{\bf Case 4.}
$U'\supset U$, $g=1$ and $\tau=\id$.
\end{enumerate}
Furthermore, we only need to consider the Betti realization.
In this case the desired commutativity is equivalent to that
of the diagram
$$\begin{CD}
M(\hat{\sigma}\otimes (\psi^{-1}\circ \det))_{c,B}^S
@. \otimes M(\sigma)_{B}^S @>[\ ,\ ] >> M_\psi(1-k)^S_B \\
@A [U'gU]_\tau^T AA @V [U'gU]_\tau VV  || @. \\
M(\hat{\sigma}'\otimes (\psi^{-1}\circ \det))_{c,B}^S
@. \otimes M(\sigma')_{B}^S @>[\ ,\ ]' >> M_\psi(1-k)^S_B \\
\end{CD}
$$
where $[\ ,\ ]$ and $[\ ,\ ]'$ are the pairings induced by
$\bar{\delta}^L$ for $\sigma$ and $\sigma'$.

{\bf Case 1: }
When $U'=gUg^{-1}$ and $\tau=\id$, we have $U'gU=gU$.
Let $U_N$ and $U_{N'}$ be as chosen in the definition of
$[U'gU]_\id = [U'gU]_{\id,N,N'}$.
So the action of
$[U'gU]_{\id}$ on $M_{N,B}\otimes  V$ is
$$[U'gU]_{\id}(x\otimes v)=[g]x\otimes v, x\in M_{N,B},
    v\in  V.$$
Using Proposition~\ref{pp:gpair}, the definition of
$[\ ,\ ]$ and the fact that $\tau^t=\id$,
we see that $[\ ,\ ]$ (resp. $[\ ,\ ]'$) is from the
composite of the maps in the top row (resp. bottom row)
of the following commutative diagram
$$
{\scriptsize
\begin{CD}
@. (M_{N,c,B}\otimes \hat{ V}\otimes
    K_{\psi}) \otimes (M_{N,B}\otimes  V) @>>>
(M_{N,c,B}\otimes \hat{ V}\otimes
    M_{\psi,B}) \otimes (M_{N,B}\otimes  V)  \\
@. @V (g\otimes \id \otimes
    \psi^{-1}\circ \det (g)) \otimes (g\otimes \id) VV
@V (g\otimes \id) \otimes (g\otimes \id) VV \\
@. (M_{N',c,B}\otimes \hat{ V}\otimes
    K_{\psi}) \otimes (M_{N',B}\otimes  V) @>>>
(M_{N',c,B}\otimes \hat{ V}\otimes
    M_{\psi,B}) \otimes (M_{N',B}\otimes  V)  \\
    @>>> M_{\psi}(1-k)_B @>\frac{1}{[U:U_N]}>>
    M_{\psi}(1-k)_B \\
@. @V ||\det g||^{2-k}[U_N:g^{-1}U_{N'}g] VV
@V ||\det g||^{2-k} VV\\
    @>>> M_{\psi}(1-k)_B @>\frac{1}{[U':U_{N'}]}>>
    M_{\psi}(1-k)_B.
\end{CD}
}
$$
Restricted to the suitable $U$- or $U'$-invariants, it gives
\begin{equation}
 [[U'gU]_{\id\otimes (\psi^{-1}\circ \det(g)),c}x,
    [U'gU]_{\id} y]'
= ||\deg g||^{2-k} [x,y].
\mlabel{eq:adj1}
\end{equation}
We also note that
\begin{equation}
[Ug^{-1}U']_{\id \otimes (\psi\circ \det(g)),c}
=([U'gU]_{\id\otimes (\psi^{-1}\circ \det(g)),c})^{-1}
\mlabel{eq:adj2}
\end{equation}
as maps $M(\hat{\sigma}' \otimes
    (\psi^{-1}\circ \det))_{c,B} \to
M(\hat{\sigma} \otimes
    (\psi^{-1}\circ \det))_{c,B}$.
Combining  (\ref{eq:adj1}) and  (\ref{eq:adj2})
proves the formula in this case.

{\bf Case 4. }
Assume $U'\supset U$ with $g=1$ and $\tau=\id$.
So $\sigma=\sigma'_{|U}$.
Let $U_N$ be a
sufficiently small subgroup of $U'\cap (g^{-1} U g)=U$
and let $U'=\coprod_i g_i U$ be a coset decomposition.
Then the action of $[U'gU]_\tau$ on $M(\sigma)_B$ is
the restriction of the map
$$[U'gU]_\tau (x\otimes v)=\sum_i [g_i] x\otimes g_i v, \
 x\in M_{N,B},\  v\in  V. $$
For each fixed $i$, consider the setting of the proposition
with
$U=U'=U_N$, $g=g_i$, $ V= V'$ being the given $ V$ (restricted
to $U_N$) and $\tau=g_i$. Then for any $u\in U$ and $v\in  V$,
$\tau((g_i^{-1} u g_i) v)=g_i(g_i^{-1} ug_i) v=ug_i v
=u\tau(v).$
Thus $[U_N g_i U_N]_{\tau,\sharp}:
M(\sigma|_{U_N})_{\sharp,B}\to M(\sigma|_{U_N})_{\sharp,B}$
is defined. In fact, we are in Case 1 of the proposition.
Let $[\ ,\ ]_N^N$ be the pairing
$$
M((\hat{\sigma}\otimes (\psi^{-1}\circ \det))|_{U_N})_{c,B}
\otimes M(\sigma|_{U_N})_{B}
\to M_{\psi} (1-k)_B.$$
Then by Case 1 of the proposition, we have
$$ [[U_N g_i^{-1}U_N]_{g_i^t\otimes(\psi^{-1}\circ\det(g_i))}
\alpha, \beta]_N^N=[\alpha,[U_N g_i U_N]_{g_i} \beta]_N^N $$
for
$\alpha \in M(\hat{\sigma}\otimes (\psi^{-1}\circ \det)|_{U_N})_{N,c,B}$
and $\beta \in M(\sigma|_{U_N})_{N,B}.$
In other words,
$$ [(g_i^{-1}\otimes
(g_i^t\otimes (\psi^{-1}\circ \det(g_i)))
\alpha, \beta]_N=[\alpha,(g_i\otimes g_i) \beta]_N.$$
Take
$\alpha\in M(\hat{\sigma}\otimes (\psi^{-1}\circ \det))_{N,c,B}$
and $\beta \in M(\sigma)_{N,B}$.
Then we have
$$[\alpha,\beta]_N^N=[\alpha,(g_i\otimes g_i) \beta]_N^N.$$
Going back to our original $U'\supseteq U$, $g=1$ and $\tau=\id$,
we have
$$[\alpha,[U'1U]_\id \beta]_N^N=[\alpha,\sum_i (g_i\otimes g_i) \beta]_N^N
= [U':U][\alpha,\beta]_N^N.$$
Dividing it by $[U':U_N]$, we have
\begin{eqnarray*}
\lefteqn{[\alpha,[U'1U]_\id \beta]'=
\frac{1}{[U':U_N]}[\alpha,[U'1U]_\id \beta]_N^N
= \frac{[U':U]}{[U':U_N]} [\alpha,\beta]_N^N}\\
&=&
\frac{1}{[U:U_N]}[\alpha,\beta]_N^N=[\alpha,\beta]
=[[U1U']_{\id\otimes (\psi^{-1}\circ \det(1))}\alpha,\beta].
\end{eqnarray*}
This proves the proposition in this case.

{\bf Case 2: }
This case follows by exchanging $U$ and $U'$ in Case 4.

{\bf Case 3. } In this case $\tau:  V\to  V'$
is a morphism of $U$-modules and the formula follows
directly from the definitions.
\proofend

\subsection{Premotivic structure of level $N$ and character $\psi$}
\mlabel{ssec:levchar}

Suppose that $k\ge 2$ and $N\ge 1$.  Let $\psi$
be a character $\hat{\ZZ}^\times \to K^\times$ of conductor
dividing $N$.  Let $U = U_0(N)$ denote the set of matrices
$\smat{a}{b}{c}{d} \in \GL_2(\hat{\ZZ})$
with $c \in N\hat{\ZZ}$.  Define $\sigma=\sigma(N,\psi)$
by the character
$\psi: U_0(N) \to K^\times$ sending
$\smat{a}{b}{c}{d}$ to $\psi^{-1}(a_N)$,
where $a_N$ denotes th image of $a$ in
$\prod_{p|N}\ZZ_p$.
Define $V=V(N,\psi)$ to be the vector space $K$ with an action
of $U$ by $\sigma$.
Let $\CV = \CO_K \subset V$.
Note that $S_\sigma = S_N^K$.  Choose $M\ge 3$ so that
$N|M$ and $S_{M} = S_N$ (for example, take $M = 4N$).
We let $\CM(N,\psi)_{M,\sharp}$ denote the premotivic
structure $\CM(\sigma)_{M,\sharp}$ for $\sharp = \emptyset$,
$\tf$, $c$ or $!$.

Recall that the isomorphism
$$\CC\otimes\fil^{k-1}\CM_{M,\dr} \cong \oplus_{t'\in(\ZZ/M\ZZ)^\times}
M_k(\Gamma(M)) \cong M_k(U_{M})$$
defined in \S\ref{ssec:adelic} respects the action
of $\GL_2(\hat{\ZZ})$ (Proposition~\ref{prop:mot-ad}), and so
for any embedding $K \to \CC$, it identifies
$\CC\otimes_K\fil^{k-1}\CM(N,\psi)_{M,\dr}$ with the
space of forms $f \in M_k(U_{M})$ such that
$$f\left(x\smat{a}{b}{c}{d}\right) = \psi(a_N)f(x)$$
for all $\smat{a}{b}{c}{d} \in U_0(N)$.
It is straightforward to check that this is precisely
the space $M_k(N,\psi)$ of classical modular forms of weight $k$,
level $N$ and character $\psi$, i.e., the set of forms
$f\in M_k(\Gamma(M))$ such that
$$f\left(\frac{a\tau+b}{c\tau+d}\right) = \psi(d)(c\tau+d)^kf(\tau)$$
for all $\tau \in \uhp$, $\gamma = \smat{a}{b}{c}{d} \in
\Gamma_0(N) = U_0(N) \cap \GL_2^+(\QQ)$.
We thus obtain an isomorphism
$$\CC\otimes_K\fil^{k-1}\CM(N,\psi)_{M,\dr} \cong M_k(N,\psi).$$
The same holds for $\CM(N,\psi)_{M,\tf}$, and replacing this
by $\CM(N,\psi)_{M,c}$ or $\CM(N,\psi)_{M,!}$ gives the
space of cusp forms $S_k(N,\psi)$.
Taking the product over all embeddings $K \to \CC$,
we obtain isomorphisms
\begin{equation}
\begin{array}{lllll}
\CC\otimes\fil^{k-1}\CM(N,\psi)_{M,\dr}&
  \cong&\CC\otimes\fil^{k-1}\CM(N,\psi)_{M,\tf,\dr}&
  \cong& M_k(N,\psi)^{\II_K}\\
\CC\otimes\fil^{k-1}\CM(N,\psi)_{M,c,\dr}&
  \cong&\CC\otimes\fil^{k-1}\CM(N,\psi)_{M,!,\dr}&
  \cong& S_k(N,\psi)^{\II_K}\end{array}
\mlabel{eq:Npsi}
\end{equation}
of $\CC^{\II_K} \cong \CC\otimes K$-modules.

Tensoring the $q$-expansion maps of
(\ref{eq:qexp}) with $K$, we regard the
$q$-expansion of a form in $M_k(N,\psi)^{\II_K}$
as taking values in $(R'\otimes\CC\otimes K)[[q^{1/M}]]$
where $R' = \ZZ[1/M,\mu_{M}]$.
\begin{lemma}  For any $S\supset S_N^K$, the
isomorphisms of (\ref{eq:Npsi})
identify
$$\fil^{k-1}\CM(N,\psi)_{M,\dr}^S =
\fil^{k-1}\CM(N,\psi)_{M,\tf,\dr}^S$$ (respectively,
$$\fil^{k-1}\CM(N,\psi)_{M,c,\dr}^S \cong
\fil^{k-1}\CM(N,\psi)_{M,!,\dr}^S)$$ with the
set of forms in $M_k(N,\psi)^{\II_K}$
(respectively, $S_k(N,\psi)^{\II_K}$)
having $q$-expansions in $\CO_S[[q]]$.
\mlabel{lem:forms}
\end{lemma}
\proofbegin  We may assume $S=S_N^K$.
Tensoring (\ref{eq:qexp}) with $\CO_S$ and
taking $U_0(N)$-invariants gives a pull-back diagram

$$
\begin{CD}
 \fil^{k-1}\CM(N,\psi)_{M,\dr} @>>>
    (R'\otimes \CO_S)[[q^{1/M}]] \\
@VVV @VVV \\
 M_k(N,\psi)^{\II_K}@>>>
\bigoplus_{t\in (\ZZ/M\ZZ)^\times} (\CC\otimes K)[[q^{1/M}]].
\end{CD}
$$
In particular, if $f\in M_k(N,\psi)^{\II_K}$ has
$q$-expansion in $\CO_S[[q]] \subset (R'\otimes\CO_S)[[q^{1/M}]]$,
then $f$ is in the image of  $\fil^{k-1}\CM(N,\psi)_{M,\dr}$.
Conversely, if $f$ is in the image of $\fil^{k-1}\CM(N,\psi)_{M,\dr}$,
then its $q$-expansion is in $(R'\otimes\CO_S)[[q^{1/M}]]$.
Viewing $f\in M_k(U_{M})\otimes K$, we see that $f$ is invariant
under the action of $\smat{1}{1}{0}{1} \in \GL_2(\AA_\f)$ and
under that of $g_{t'}$ for each $t'\in (\ZZ/M\ZZ)^\times$.
Unravelling the effect of $\smat{1}{1}{0}{1}$ on $q$-expansions
shows that $f$ has $q$-expansion in $(R'\otimes\CC\otimes K)[[q]]$,
hence in $(R'\otimes\CO_S)[[q]]$.
The commutativity of the diagram (\ref{eq:qexp})
together with the invariance under the $g_{t'}$ shows that the
image of $f$ in $(\CC\otimes K)[[q]]$ is independent of $t'$,
i.e., independent of the map $R'\to \CC$.
We conclude that the $q$-expansion of $f$ is invariant
under $\aut R'$, hence has coefficients in $\CO_S$.
\proofend

\section{Premotivic structure of a newform}
\mlabel{sec:newform}

\subsection{Hecke actions}
Using the notations in \S\ref{ssec:levchar} and
letting $S$ be a superset of $S_N^K$,
we now consider the action of Hecke operators on
$\CM(N,\psi)_{M,!}$ as an object in $\ipms^{S}_K$.

Fix a rational prime $p$.
By Lemma~\ref{lem:coset1},
$[U\smat{p}{0}{0}{1}_p U]_{\psi(p_p)^{-1}}$ defines an endomorphism
of $\CM(N,\psi)_{M,!}$ regarded as an object
in $\ipms^{S_{pN}}$.
This induces an endomorphism
$$[U\smat{p}{0}{0}{1}_p U]_{\psi(p_p)^{-1}}: M(N,\psi)_{M,!}\to
M(N,\psi)_{M,!}$$
in $\pms^{S_{pN}}$.

\begin{lemma}
\begin{enumerate}
\item
The action of the usual Hecke operator $T_p$ on
$M_k(N,\psi)$
corresponds to that of the double coset operator
$[U\smat{p}{0}{0}{1}_p U]_{\psi(p_p)^{-1}}$ on $M(N,\psi)_M$,
regarded as an object in $\pms^{S_{pN}}_K$.
\item
If $p\nmid N$, then
the operator $[U \smat{p}{0}{0}{p}_p U]_{\psi(p_p)^{-2}}$ acts on
$M(N,\psi)_{M,c}$ via
$\psi(p_p)^{-1}p^{k-2}$.  This corresponds to the usual Hecke
operator $S_p$ on $M_k(N,\psi)$.
\end{enumerate}
\mlabel{lem:Hop}
\end{lemma}
\proofbegin
(a)
Denote $g=\smat{p}{0}{0}{1}_p$ and
let $U=\coprod_i g_i (g Ug^{-1})\cap U$ be a coset representation.
Then we have double coset decomposition
$UgU=\coprod h_i U$ with $h_i=g_i g$.
The usual Hecke operator $T_p$ on $M_k(\Gamma_1(M))$
corresponds to that of the double coset operator
$T^\ad_p=[U \smat{p}{0}{0}{1}_p U]$ on the adelic modular
forms $M_k(U_1(M))$~\cite[\S 11.1]{diamond_im}.
By definition,
$$T_p^\ad =[U\smat{p}{0}{0}{1}_p U]=\sum_i h_i.$$
By Proposition~\ref{prop:mot-ad}, this corresponds to the operator
$T_p^{\rm mot}\defeq\sum_i [h_i]$ on $\CC\otimes \fil^{k-1}\CM_{M,!,\dr}$.
Therefore the operator
$\sum_i h_i$ on $(V\otimes M_k(U_1(M)))^U$ corresponds to the
operator
$\sum_i h_i$ on $(V\otimes \CC\otimes
\fil^{k-1}\CM_{M,!,\dr})^U=\CC\otimes \CM(N,\psi)_{M,!,\dr}$.
The second operator is defined by
$$\sum_i h_i (\sum_j x_j\otimes v_j)
=\sum_{i,j} [g_i g] x_j \otimes g_i g v_j
=\sum_{i,j} [g_i g] x_j \otimes g_i \psi(p_p)^{-1} v_j
=[U g U]_{\psi(p_p)^{-1}}.$$

On the other hand, if $f$ is in $M_k(U,\psi)$, then
$1\otimes f$ is in $(V\otimes M_k(U_1(M)))^U$ and
$f\mapsto 1\otimes f$ defines an isomorphism
$S_k(U,\psi)\cong (V\otimes M_k(U_1(M)))^U$.
Then we have
$$\sum_i h_i (1\otimes f) =\sum_i (\psi(h_i)^{-1} \otimes
    h_i f)
    = 1\otimes (\sum_i \psi(h_i)^{-1} h_i f).$$

(b)
By definition, $[U\smat{p}{0}{0}{p}U]_{\psi(p_p)^{-2}}$
is the restriction of the map
$\CM_M \otimes \CV\to \CM_M \otimes \CV,\
x\otimes v \mapsto \left [\smat{p}{0}{0}{p}_p\right ] x\otimes \psi(p_p)^{-2}v$
to $\CM(N,\psi)_M$.
Let $p=||p_p||^{-1} p'$ with $p'\in \hat{\ZZ}^\times$
and
let $\sum_i x_i\otimes v_i$ be in $\CM(N,\psi)_M$.
Since $\smat{p'}{0}{0}{p'}$ is in $U$, we have
$$ \sum_i \left [\smat{p'}{0}{0}{p'}\right ] x_i\otimes \psi^{-1}(p') v_i
= \sum_i x_i\otimes v_i.$$
Thus
$$\sum_i \left [\smat{p'}{0}{0}{p'}\right ] x_i \otimes v_i
    =\sum_i x_i \otimes \psi(p') v_i
    =\psi(p') \sum_i x_i\otimes v_i.$$
Then we have
\begin{align*}
\sum_i \left [\smat{p}{0}{0}{p}_p\right ] x_i\otimes \psi(p_p)^{-2}v_i
&=\psi(p_p)^{-1} ||p_p||^{2-k} \sum_i x_i\otimes v_i.
\end{align*}
\proofend

Because of the lemma, we can write $T_p$ (resp. $S_p$)
for the endomorphisms of $\fil^{k-1}M(N,\psi)_{M,\dr}$,
$\fil^{k-1}M(N,\psi)_{M,c,\dr}$ and $\fil^{k-1}M(N,\psi)_{M,!,\dr}$
defined by the coset operator
$[U\smat{p}{0}{0}{1}_p U]_{\psi(p_p)^{-1}}$
(resp. $[U\smat{a}{0}{0}{a} U]_{\psi(p_p)^{-2}}$).

\subsection{Compatibility of the Hecke action with
    the pairings}
Let $w$ denote $\smat{0}{-1}{N}{0}_N$ in $\dprod_{p|N}\GL_2(\ZZ_p)$.
Let $V'$ be the one dimensional representation
$$\hom(V,K)\otimes K_{\psi^{-1}\circ \det}$$ of $U$.
We abbreviate $\sigma'$ for $\hat{\sigma}_\psi\otimes (\psi^{-1}\circ
\det)$. Fix a basis $v_0$ of $\CV$.
Define $\omega: \CV\to \CV'$ by sending
$v_0$ to $\hat{v}_0\otimes 1$, where $\hat{v}_0\in
\hom(\CV,\CO)$ is such that $\hat{v}_0(v_0)=1$.

\begin{lemma}
\begin{enumerate}
\item
The operator $[UwU]_\omega$ defines an isomorphism
$$\CM(N,\psi)_M \to \CM(\sigma')_{M}$$
in $\ipms^{S_N^K}_K$ and
$[UwU]_\omega^{-1}=[U w^{-1} U]_{\omega^t\otimes
\psi(\det w)}$.
Similarly for $\CM(N,\psi)_{M,c}$ and $\CM(N,\psi)_{M,!}$.
\item
The operator
$[UwU]_\omega: \CM(\psi)_{M,!} \to
    \CM(N,\sigma')_{M,!}$
is adjoint to
$$N^{k-2}[Uw^{-1}U]_{\omega^t\otimes \psi(\det w))}:
\CM(\sigma')_{M,!} \to \CM(N,\psi)_{M,!}$$ and
coincides with $\psi(-1_\infty)N^{k-2}[Uw^{-1}U]_\omega$.
\end{enumerate}
\mlabel{lem:wiso}
\end{lemma}
\proofbegin
(a)
We first note that both $wU$ and $UwU$ consist of
matrices $\smat{a}{b}{c}{d}\in \GL_2(\AA_\f)\cap
   M_2(\hat{\ZZ})$ with
$a,d\in N \hat{\ZZ}$, $b\in \hat{\ZZ}^\times$
and $c\in N\hat{\ZZ}^\times.$
So we have $UwU=wU$.

For $u=\smat{a}{b}{c}{d}\in U$ we have
$$\sigma(w^{-1}uw)v_0=\sigma\smat{d}{-\frac{c}{N_N}}{-N_N b}{a} v_0
    =\psi^{-1}(d_N)v_0$$
and
$$\sigma'(u)(\hat{v}_0\otimes 1)=
\psi^{-1}(((\det u)^{-1} d)_N) \hat{v}_0
\otimes \psi^{-1}(\det u) \cdot 1
=\psi^{-1}(d_N)\hat{v}_0\otimes 1.
$$
since $\psi(\det u/(\det u)_N)=1$.
Therefore we have
$$\omega(\sigma(w^{-1}uw)v_0) = \sigma'(u) \omega(v_0).$$
Also $w^{-1} U_{M^2} w \subseteq U_M$.
So the operator
$$[UwU]_{\omega,M,M^2}: \CM(N,\psi)_M \to \CM(\sigma')_{M^2}$$
is well-defined and defines an morphism in
$\ipms^{S_N}_K$.
Since $\CM(\sigma')_{M^2} \cong \CM(\sigma')_M$ by
Lemma~\ref{lem:coset0} (a), we have
$$[UwU]_{\omega}: \CM(N,\psi)_M \to \CM(\sigma')_M.$$

To prove that the map is an isomorphism, consider
$-w=\smat{0}{1}{-N}{0}_N$. The same argument as above shows
that we have a homomorphism
$$[U(-w)U]_{\omega'}: \CM(\sigma')_M \to \CM(\sigma'')_M$$
in $\ipms^{S_N}_K$
where $\sigma''$ is the representation of $U$ on
$\CV''=\hom(\hat{\CV}',\CO)\otimes \CO_{\psi^{-1}\circ\det}$
which is canonically isomorphic to $\sigma_\psi$
and where
$\omega':\CV'\to \CV''$ sends $\hat{v}_0$ to
$\hat{\hat{v}}_0\otimes 1$ which is identified with $v_0$
under the above isomorphism.
Thus we have
$$[U(-w)U]_{\omega'} \circ [UwU]_\omega:
    \CM(N,\psi)_M \to \CM(N,\psi)_M.$$
Since
$[U(-w)U][UwU]=[U(-w)wU]=[U(N_N I)U],$ the composition
formula in Proposition~\ref{prop:scomp} gives us
$$[U(-w)U]_{\omega'} \circ [UwU]_\omega = [U(N I)U]_1
    =N^{k-2}.$$
Since $N$ is invertible in $\CO$,
$[UwU]_\omega$ is an isomorphism in $\ipms^{S_N^K}_K$ and
$$[UwU]_\omega^{-1}=N^{2-k}[U(-w)U]_{\omega'}
    =[Uw^{-1}U]_{\omega^t\otimes \psi(\det w)}.$$

The statements for $\CM_c$ and $\CM_!$ follow since
$[UwU]_\omega$ commutes with the weight filtration maps.

\medskip
(b)
By Lemma~\ref{pp:adj},
the operator
$[UwU]_\omega: \CM(N,\psi)_{M,!} \to \CM(\sigma')_{M,!}$
is adjoint to
$N^{k-2}[Uw^{-1}U]_{\omega'}:\CM(\sigma')_{M,!} \to
    \CM(N,\psi)_{M,!}$.

On the other hand,
$$w=\smat{0}{-N}{1}{0}_N =\smat{-1}{0}{0}{-1}_N
    \smat{N}{0}{0}{N}_N w^{-1}.$$
Thus by Proposition~\ref{prop:scomp},
\begin{eqnarray*}
 [UwU]_\omega &=& [U w^{-1} U]_\omega \circ [U \smat{N}{0}{0}{N}_N U]_1
    \circ [U\smat{-1}{0}{0}{-1}_N  U]_1 \\
    &=& \psi^{-1}(-1_\infty) N^{k-2} [U w^{-1} U]_\omega
\end{eqnarray*}
since
$\sigma(\smat{-1}{0}{0}{-1}_N)=\psi^{-1}(-1_N)
=\psi^{-1}(-1_f)=\psi(-1_\infty)$.
\proofend

Composing the operator $[UwU]_\omega$ with the duality morphism
$\bar{\delta}^L_N$
defined in (\ref{eq:spair2}), we obtain a duality morphism
\begin{eqnarray}
&\hat{\delta}^L_N:
\CM(N,\psi)_{M,!}  \to \hom_\CO (\CM(N,\psi)_{M,!}, \CM_{\psi}(1-k))
    \mlabel{eq:tpair2}
\end{eqnarray}
that becomes an isomorphism after tensoring with $\QQ$.

Let $U'=U_0(N')$ and let $\psi'$ be a character of $\AA^\times$
with conductor dividing $N'$. Then we can construct the
premotivic structure
$\CM(N',\psi')_{M',!}$ associated to $\psi'$ in the same way as we
construct $\CM(N,\psi)_{M,!}$.
Given $[U'gU]_{\tau,M,M'}: \CM(N,\psi)_{M,!}\to \CM(N',\psi')_{M',!}$
where $M'$ is chosen with the additional condition
$g^{-1} U_{M'} g \subseteq U_M$,
let $[U'gU]_\tau^t$ denote the adjoint
of $[U'gU]_\tau$ with respect to the pairing morphism
$\hat{L}$.
\begin{lemma}
\begin{enumerate}
\item
$$[U'gU]_\tau^t=||\det g||^{2-k} [Uw^{-1}U]_{\omega^t\otimes
\psi(\det w)}
    \circ [Ug^{-1}U']_{\tau^t\otimes \psi(\det(g))} \circ
    [U'w'U']_{\omega'}$$
in $\ipms^{S_{M''}}$ for any $M''$ with $g^{-1}U_{M''} g\subseteq
U_{M'}$.
\item
For any prime $p$,
 $T_p$ is self-adjoint under the duality morphism
$\hat{\delta}^L$ in $\ipms^{S_{pN}}$.
\end{enumerate}
\mlabel{lem:self-adj}
\end{lemma}
\proofbegin
(a)
We only need to verify for the Betti realization.
Let $[\ ,\ ]$, $[\ ,\ ]'$, $\langle\ ,\ \rangle$ and
$\langle\ ,\ \rangle'$ be the pairings on Betti realizations
induced by $\bar{\delta}^L_N, \bar{\delta}^L_{N'},
\hat{\delta}^L_N$ and $\hat{\delta}^L_{N'}$.
By definition we have the commutative diagram
$$
\begin{CD}
\CM(N,\psi)_{!,B} @. \otimes_\CO \hspace{-1cm}@. \CM(N,\psi)_{!,B} @> {\langle\ ,\ \rangle} >> \CM_{\psi}(1-k)_B \\
\Vert @. @. @VV [UwU]_\omega V  \Vert @. \\
\CM(N,\psi)_{!,B} @. \otimes_\CO @. \CM(\hat{\sigma}\otimes
(\psi^{-1}\circ \det))_{!,B} @> [\ ,\ ] >> \CM_{\psi}(1-k)_B \\
@V [U'gU]_\tau VV @. @VV [U'gU]_\tau^T V  \Vert @.\\
\CM(N',\psi')_{!,B} @. \otimes_\CO @. \CM(\hat{\sigma'}\otimes
(\psi'{}^{-1}\circ \det))_{!,B} @> [\ ,\ ]' >> \CM_{\psi'}(1-k)_B \\
\Vert @. @. @AA [U'wU']_{\omega'} A \Vert @. \\
\CM(N',\psi')_{!,B} @. \otimes_\CO \hspace{-1cm}@. \CM(N',\psi')_{!,B}
@> {\langle\ ,\ \rangle'} >> \CM_{\psi'}(1-k)_B.
\end{CD}
$$
Then from Proposition~\ref{pp:adj} and Lemma~\ref{lem:wiso}
\begin{align*}
\langle [U'gU]_\tau & a,b\rangle'
=[[U'gU]_\tau a,[U'w'U']_{\omega'} b ]'\\
&=[a, [U'gU]_\tau^T \circ [U'w'U']_{\omega'} b]\\
&=\langle a,[UwU]_\omega^{-1} \circ[U'gU]^T\circ [U'w'U']_{\omega'} b\rangle\\
&=\langle a,[Uw^{-1}U]_{\omega^t\otimes \psi(\det w)}
    \circ [U'gU]^T \circ [U'w'U']_{\omega'} b\rangle\\
&=\langle a,||\det g||^{2-k} [Uw^{-1}U]_{\omega^t\otimes
\psi(\det w)}
    \circ [Ug^{-1}U']_{\tau^t\otimes \psi(\det(g))} \circ
    [U'w'U']_{\omega'}b\rangle.
\end{align*}
This proves (a).

(b)
Note that
$T_p^T=[U\smat{1}{0}{0}{p_p}U]_{\psi(p_p)^{-1}}$ by
Proposition~\ref{pp:adj}. Then using Proposition~\ref{prop:scomp}
we find that $[UwU][U\smat{p_p}{0}{0}{1} U]$
and $[U\smat{1}{0}{0}{p_p} U][UwU]$ are both
$[U \smat{0}{-1}{p_p N_N}{0} U]$. It follows that
$$
[UwU]_\omega T_p=T_p^T [UwU]_\omega.
$$
Then from (a) we have
$$T_p^t=[UwU]_\omega^{-1} T_p^T [UwU]_\omega = T_p.$$
\proofend

\begin{lemma}
The duality morphism
$$\hat{\delta}^L_N: \CM(N,\psi)_{!,M} \to
    \hom_{\CO_K}(\CM_!(N,\psi)_{M},\CM_{\psi}(1-k))$$
has sign $-1$.
\end{lemma}
\proofbegin
Recall that for $\CM(N,\psi)_{!,M}$ to be non-zero, we must have
$\psi(-1_\infty) = (-1)^{k-2}$.
Fix such a $\psi$. By Lemma~\ref{lem:wiso} and the fact
that $[\ ,\ ]$ has sign $(-1)^{k-1}$, we have,
for $x,\ y\in \CM(N,\psi)_{M,!,B}$,
\begin{eqnarray*}
\langle x,y\rangle &=&  [[UwU]_\omega x, y] \\
&=& \psi(-1_\infty) [N^{k-2}[Uw^{-1}U]_{\omega} x, y]\\
&=& -[y, N^{k-2}[Uw^{-1}U]_{\omega} x]\\
&=& -[[UwU]_\omega y, x]\\
\hspace{5cm} &=& -\langle y, x\rangle. \hspace{8cm}\qedhere
\end{eqnarray*}
\proofend

\subsection{The action of Hecke rings}

Let $\tilde{\TT}$ denote the polynomial algebra
over $\CO=\CO_S$ generated by the variables $t_p$, where
$p$ runs over rational primes.
Then $\fil^{k-1}M(\psi)_{N,!,\dr}$ becomes a $\tilde{\TT}$-module
with $t_p$ acting by $[U\smat{p}{0}{0}{1}_p U]_1$.
Since the modular forms with Fourier coefficients in $\CO$
is stable under the action of the Hecke operators $T_p$,
by Lemma~\ref{lem:forms} (c), $\fil^{k-1}\CM(\psi)_{N,!,\dr}$ is
also a $\tilde{\TT}$-module.
Let $\ga$ denote the annihilator
in $\tilde{\TT}$ of $\fil^{k-1} \CM(N,\psi)_{M,!,\dr}$
and let $\TT = \tilde{\TT}/\ga$.
Similarly let $\ga' \subset \ga$ be the annihilator of
$\fil^{k-1}M(N,\psi)_{M,\dr}$ and let $\TT'=\tilde{\TT}/\ga'$.
Then $\TT\subset \End_K(\fil^{k-1}\CM(N,\psi)_{M,!,\dr})$
and $\TT'\subset \End_K(\fil^{k-1}M(N,\psi)_{M,\dr})$ .

\begin{lemma}
\begin{enumerate}
\item
The $\TT$-linear homomorphism
\begin{align*}
& \varphi: \fil^{k-1} \CM(N,\psi)_{M,!,\dr}
\to \hom (\TT,\CO), \\
& \varphi(f)(T)=a_1(T(f)),\
f\in \fil^{k-1}\CM(N,\psi)_{M,!,\dr}, T\in \TT,
\end{align*}
with $a_1(T(f))$ the first Fourier coefficient of $T(f)$,
is an isomorphism.
\item
For any $p \nmid N$, $\QQ\otimes \TT$ is generated by the
images of $T_\ell$ for $\ell \neq p$.
\item
The $\TT'$-linear homomorphism
\begin{align*}
& \varphi:  \fil^{k-1} M(N,\psi)_{M,\dr}
\to \hom (\QQ\otimes \TT',K),\\
&  \varphi(f)(T)=a_1(T(f)),\
f\in \fil^{k-1}M(N,\psi)_{M,\dr}, T\in \TT'
\end{align*}
is an isomorphism.
\item
For any $p \nmid N$, $\QQ\otimes \TT'$ is generated by the
images of $T_\ell$ for $\ell \neq p$.
\end{enumerate}
\mlabel{lem:T2}
\end{lemma}
\proofbegin
(a)
The proof is the same as the proof of Theorem 2.2 in \cite{ribet}.

(b)
Fix a $p\nmid N$.
Let $\TT_1$ be the subalgebra of $\QQ\otimes \TT$ generated by
the images of $t_\ell$ for $\ell\neq p$.
Consider the map
$$ \phi: \fil^{k-1}M(N,\psi)_{M,!,\dr}
    \cong \hom(\QQ\otimes \TT,K) \to \hom(\TT_1,K)$$
where the second map is the restriction, and thus is surjective.
We only need to prove that this map is injective.

Let $f$ be in $\fil^{k-1}M(N,\psi)_{M,!,\dr}$. By definition,
$\phi(f)(T_n)=a_1(T_n f)=a_n(f)$.
Therefore, for $f\in \ker \phi$ we have
$a_n(f)=0$ when $\ell\not| n$.
It follows from~\cite[Theorem 4.6.8 (1)]{miyake} that $f=0$.
Therefore $\TT=\TT_1$.

The proof of (c) (resp. (d)) is the same as the proof of (a)
(resp. (b)).
\proofend

\begin{prop}
\begin{enumerate}
\item
$\TT$ acts on $\CM(N,\psi)_{M,!}$ as an object in $\ipms^{S_N^K}_K$.
\item
$\TT'$ acts on $M(N,\psi)_M$ as an object in $\pms^{S_N^K}_{K}$.
\item
$\ga'$ is also the annihilator of $\CM(N,\psi)_{M,c,\dr}$
and $\TT'$ acts on $\CM(N,\psi)_{M,c}$ as an object in
$\ipms^{S_N^K}_K$.
\item
The maps $\CM(N,\psi)_{M,c} \to \CM(N,\psi)_{M,!}$
in $\ipms^{S_N^K}_K$ and $M(N,\psi)_{M,!} \to M(N,\psi)_M$ in
$\pms^{S_N^K}_K$
 are $\TT'$-linear.
\end{enumerate}
\mlabel{prop:T}
\end{prop}
\proofbegin
(a)
We first prove that $\ga$ is annihilates each
of the realizations of $M(N,\psi)_{M,!}$.
Let $c:\CC\to \CC,\ v\mapsto \bar{v}$ be the complex conjugation.
By definition,
$$\overline{\fil}^{k-1}M(N,\psi)_{M,!,\dr}
=I^\infty ((\id\otimes c)((I^\infty)^{-1}(\fil^{k-1}M(N,\psi)_{M,!,\dr}))).$$
Here
$$\id\otimes c: M(N,\psi)_{M,!,B}\otimes \CC \to
M(N,\psi)_{M,!,B} \otimes \CC, \
x\otimes v \mapsto x\otimes \bar{v}$$
is clearly $\tilde{\TT}$-linear,
as is $I^\infty$ by Lemma~\ref{lem:coset1}.
Thus $\ga$ annihilates
$\overline{\fil}^{k-1}M_!(N,\psi)_{M,\dr}$.
So $\ga$ is the annihilator of
$$ M(N,\psi)_{M,!,\dr} \otimes \CC
= (\fil^{k-1} M(N,\psi)_{M,!,\dr} \otimes \CC)
    \oplus (\overline{\fil}^{k-1} M(N,\psi)_{M,!,\dr} \otimes
    \CC).$$
Since $I^\infty$ is $\tilde{\TT}$-linear,
$\ga$ is the annihilator of $M(N,\psi)_{M,!,B} \otimes \CC$.
Since the action of $\tilde{\TT}$ on $\CC$ is trivial,
$\ga$ is the annihilator of $M(N,\psi)_{M,!,\dr}$ and
$M(N,\psi)_{M,!,B}$.
Then it is also the annihilator of
$M(N,\psi)_{M,!,\lambda} =M(N,\psi)_{M,!,B} \otimes_K
K_\lambda$ and
$M(N,\psi)_{M,!,\dr} \otimes \QQ_p$.

We next show that $\TT$ acts on $\CM(N,\psi)_{M,!}$ as an object
in $\ipms^{S_N^K}_K$.

Fix a finite prime $p$.
Recall that $\CM(N,\psi)_{M,!}$ is in $\ipms^{S_{N}^K}_K$
and, by Lemma~\ref{lem:coset1},
$T_p: \CM(N,\psi)_{M,!} \to \CM(N,\psi)_{pM,!}$ is defined in
$\ipms^{S_{pM}^K}_K= \ipms^{S_{pN}^K}_K$.
If $p| Nk!$, then $S_{pN}=S_N$. So
$T_p$ is a morphism in $\ipms^{S_N^K}_K$.

If $p\nmid Nk!$, then for each realization $\CM_\wc$ of
$\CM(N,\psi)_{M,!}\in \ipms^{S_N}$, the corresponding realization of
$\CM(N,\psi)_{pM,!}\in \ipms^{S_{pN}}$ is given by
$$
\left \{ \begin{array}{ll}
    \CM_\wc, & \wc=B, \lambda {\rm\ or\ } \lambda{\rm -crys} {\rm\ with\ }
    \lambda\nmid pNk!, \\
    \CM_\wc [1/p], &
    \wc=\dr \end{array} \right . $$
Thus by Lemma~\ref{lem:coset1}, $T_p$ acts on these
realizations and observes comparison isomorphisms among them.
So to prove the proposition, we only need to give actions of $T_p$
on $\CM_\lcrys$ with $\lambda\nmid p$ and $\CM_\dr$,
and show that these
actions observe the comparison isomorphisms $I_\dr^\lambda$ and
$I^\lcrys$. Fix a $\lambda\mid p$.
By Lemma~\ref{lem:T2}, $T_p$ is a polynomial of $T_q,\ q\neq p,$
with coefficients in $K$. Thus there is $r\in \ZZ_{>0}$ such that
$r T_p$ is a polynomial of $T_q,\ q\neq p,$ with coefficients in
$\CO_K$. Since each $T_q,\ q\neq p$, acts on $\CM_\lcrys$
by Lemma~\ref{lem:coset1},
we have an action of $r T_p$ on $\CM_\lcrys$ which we provisionally
denote by $F$.
Further, $T_p:\CM_\lambda \to \CM_\lambda$ gives
$rT_p \CM_\lambda \subseteq r\CM_\lambda$.
Since each $T_q,\ q\neq p,$ observes $I^\lcrys$, we have
$$\VV(F(\CM_\lcrys))=F(\VV(\CM_\lcrys))\cong F(\CM_p)
\subseteq r\CM_p \cong r\VV(\CM_\lcrys)=\VV(r\CM_\lcrys).$$
Thus $F(\CM_\lcrys)\subseteq r \CM_\lcrys.$
Since $\CM_\lcrys$ for $\lambda\nmid Nk!$ is torsion free,
we obtain a map
$r^{-1}F: \CM_\lcrys \to \CM_\lcrys$ that observes
$I^\lcrys: \VV(\CM_\lcrys)\to \CM_\lambda.$
Then via $I_\dr^\lambda: \CM_\dr \otimes \CO_\lambda \to
\CM_\lcrys$, we obtain
a map $r^{-1}F: \CM_\dr\otimes \CO_\lambda \to
    \CM_\dr\otimes \CO_\lambda$.
Note that $T_\ell$ acts on $M_\dr$ for all $\ell$
and, by the first part of the proof,
the action of $T_p$ on $M_\dr$ agrees
with $r^{-1}F$.
Since $\CM_\dr \otimes \ZZ_p =\bigoplus_{\lambda|p}
    \CM_\dr \otimes \CO_\lambda$,
we have the following commutative diagram
$$
\begin{CD}
\CM_\dr @>\subset >> \CM_\dr[1/p] @> T_p >> \CM_\dr[1/p] \\
\Vert @. @V \cap VV @V \cap VV\\
\CM_\dr @>\subset >> M_\dr @> T_p=r^{-1}F >>
    M_\dr \\
\Vert @.  @A \cup AA @A \cup AA\\
\CM_\dr @>\subset >> \CM_\dr\otimes \ZZ_p
    @> r^{-1}F >> \CM_\dr\otimes \ZZ_p.
\end{CD}
$$
This shows that the image of $\CM_\dr$ under $T_p$ is contained in
$\CM_\dr[1/p]\cap (\CM_\dr \otimes \ZZ_p)=\CM_\dr$.
Thus if we define $T_p:\CM_{\lcrys}\to \CM_{\lcrys}$
by $r^{-1}F$, then it is compatible with $I^\lcrys$  and
with $I_\dr^\lcrys$.
So we have seen that
$T_p:\CM_!(N,\psi)_M \to \CM_!(N,\psi)_{pM}$ in
$\ipms^{S_{pN}}$ extends to a morphism
$T_p:\CM_!(N,\psi)_M\to \CM_!(N,\psi)_M$
in $\ipms^{S_N^K}_K$.

Thus we have shown that $T_p$ is a morphism in
$\ipms^{S_N^K}_K$
for every $p$ and proven (a).

The argument for (b) is the same as for (a).

To prove (c),  note that
by (b) and the fact that $T_p$ is self-adjoint under
the pairing between $M(N,\psi)_M$ and $M(N,\psi)_{M,c}$
 defined in Lemma~\ref{lem:self-adj},
$\ga'$ is also the annihilator of $M(N,\psi)_{M,c}$.
The rest of the proof is the same as for part (a).

For (d),
we only need to prove that each $T_p$ commutes with
these maps. This follows from Lemma~\ref{lem:coset1}.
\proofend

\subsection{Premotivic structure for a newform}
Now suppose that $g$ is a (normalized) eigenform in
$\fil^{k-1}\CM(N,\psi)_{N,!,\dr}$
for the action of $\TT$. So for $T\in\TT$ we have
$T(g)=a_1(T(g))g$. Let $K=K_g$ denote the field generated by
$a_1(T(g)),\ T\in \TT$ and let
$I_g$ denote the kernel of the
map $\TT \to K,\ T\mapsto a_1(T(g))$.
Let $\CM_g \subset \CM(N,\psi)_{M,!}$ denote
the intersection of the kernels of elements of $I_g$.
By Proposition~\ref{prop:T}, $\CM_g$ is in $\CP\CM^{S}_K$
for any $S\supset S_N^K$.
Let $M_g=\CM_g \otimes K\in \pms^S_K$.

\begin{lemma}
$M_g$ is a premotivic structure of rank 2 over $K$
and $\fil^{k-1}M_{g,\dr} = Kg$.  If $g$ is a new form,
then the pairing on $M(N,\psi)_{M,!}$
restricts to a perfect alternating pairing on $M_g$,
i.e., an isomorphism
$$\wedge^2_K M_g \to M_{\psi}(1-k).$$
\mlabel{lem:g}
\end{lemma}
\proofbegin
By \cite[Proposition 12.4.14]{diamond_im},
$M(N,\psi_0)_{M,!,\dr}$, where $\psi_0$ is the trivial character
of $U_1(N)$, is free of rank two over F$\TT_N$.
It follows that $M(N,\psi)_{M,!,\dr}$ is free of rank two
over $\TT$.
Thus $M_g$ is free of rank two over
$\TT/I_g\cong K$.
By the injectivity of the $q$-expansion map, we have
$\fil^{k-1}M_{g,\dr} = Kg$.

The pairing on $M_g$ is restricted from an alternating pairing.
So it is still alternating.
To show that the pairing on $M_g$ is perfect, we only need
to show that the pairing on $M_{g,\dr}\otimes \CC$
is perfect.

By  (\ref{eq:tpair2}),
$\langle M(N,\psi)_{M,!,\dr}, M(N,\psi)_{M,!,\dr}\rangle\neq 0$.
For  $\alpha\in I_g$ and $x\in M(N,\psi)_{M,!,\dr}$, we have
$\langle M_{g,\dr},\alpha x\rangle=\langle\alpha M_{g,\dr},x\rangle=0$. Thus
$\langle M_{g,\dr}, I_g M(N,\psi)_{M,!,\dr}\rangle=0$.
By strong multiplicity one theorem,
$$\fil^{k-1}M(N,\psi)_{M,!,\dr}
= Kg + I_g\fil^{k-1}M(N,\psi)_{M,!,\dr}.$$
Then
\begin{align*}
& \CC\otimes M(N,\psi)_{M,!,\dr} \\
\cong &
(\CC\otimes \fil^{k-1} M_!(N,\psi)_{M,\dr})+
(I^\infty (\id\otimes c) (I^\infty)^{-1}) (\CC \otimes \fil^{k-1} M_!(N,\psi)_{M,\dr})\\
\cong &\CC g +\CC\otimes I_g \fil^{k-1}M_!(N,\psi)_{M,\dr} +
(I^\infty(\id\otimes c)(I^\infty)^{-1})(g\otimes \CC)
\\
&+ (I^\infty (\id\otimes
c)(I^\infty)^{-1})(\CC\otimes I_g \fil^{k-1}M_!(N,\psi)_{M,\dr}).
\end{align*}
Thus
\begin{align*}
0\neq &\langle\CC g, \CC\otimes M_!(N,\psi)_{M,\dr}\rangle
=\langle\CC g, \CC g+(I^\infty(\id\otimes c)(I^\infty)^{-1})(g\otimes \CC)\rangle\\
=&\langle\CC g, \CC\otimes M_{g,\dr}\rangle.
\end{align*}
Since $M_{g,\dr}$ is two dimensional and the pairing in alternating,
this implies that the pairing on $\CC\otimes M_{g,\dr}$ is
nondegenerate, as is needed.
\proofend

\subsection{The $L$-function}
\mlabel{ssec:eulerp}

Suppose that $f(\tau) = \sum a_n e^{2\pi i n\tau}$ is a newform of
weight $k$, conductor $N_f$ and character $\psi_f$.  Associated to
$f$ is the $L$-function with Euler product factorization:
$$L(f,s) = \sum_{n\ge 1} a_n n^{-s} = \prod_{p\nmid N_f}
(1-a_p p^{-s} + \psi_f(p)p^{k-1-2s})^{-1}\prod_{p|N_f}(1-a_pp^{-s})^{-1}.$$

Recall that $\CA_k^0$ denote the representation of $\GL_2(\AA_\f)$
defined
by
$$\lim_{\stackrel{\rightarrow}{N}} \fil^{k-1} M_{N,!,\dr}\otimes \CC.$$
Then there is an irreducible subrepresentation $\pi(f)$ of $\CA_k^0$
of central character $\psi_f||\ ||^{2-k}$
such that $f$ spans the image of
$\pi(f)^{U_1(N_f)}$ under the isomorphism
$$(\CA_k^0)^{U_1(N_f)} \cong S_k(\Gamma_1(N_f))$$
(Here we view $\psi_f$ as a character
on $\AA_f^\times \subset \AA^\times$).
Moreover, we have the decomposition
$$\CA_k^0 = \bigoplus_f \pi(f)$$
where $f$ runs over  newforms of weight $k$ of any conductor
and character.
For each $f$ we have a factorization
$$\pi(f) \cong \otimes_p \pi_p(f)$$
where $\pi_p(f)$ is an irreducible admissible representation
of $\GL_2(\QQ_p)$.  For each $p$, we let $c_p = v_p(N_f)$,
so $p^{c_p}$ is the conductor of $\pi_p(f)$.  We let
$$V_p= \left\{ \left. \smat{a}{b}{c}{d} \in \GL_2(\ZZ_p)\right|
a -1 \equiv c \equiv 0 \bmod p^{c_p}\right\}.$$
Then $\pi_p(f)^{V_p}$ is one-dimensional and we consider the
double coset operators on it defined by
$$t_p = V_p \smat{p}{0}{0}{1} V_p, \quad
t_p' = V_p \smat{1}{0}{0}{p} V_p, \quad
s_p = V_p \smat{p}{0}{0}{p} V_p.$$
Then we have $t_p = \psi_f(p_p)a_p$ and $s_p = \psi_f(p_p)p^{k-2}$
on $\pi_p(f)^{V_p}$.    We recall the following general
facts about irreducible admissible representations of $\GL_2(\QQ_p)$:
\begin{lemma}  Suppose that $c_p > 0$ and $\pi_p(f)$ has minimal
conductor among its twists.
\begin{enumerate}
\item  $\pi_p(f)$ is special if and only if $t_p t_p' = s_p$;
\item  $\pi_p(f)$ is principal series if and only if
  $t_p t_p' = p s_p$;
\item $\pi_p(f)$ is supercuspidal if and only if $t_p = 0$.
\end{enumerate}
\mlabel{lem:pi-t}
\end{lemma}

Suppose now that $g$ is a newform of weight $k$, level $N_g$
and character $\psi_g$, with coefficients in a number
field $K$.  For each embedding $\tau: K \to \CC$, we have a
newform $\tau(g)$ as above, and so a representation
$\pi_{\tau(g)}$ of $\GL_2(\AA_\f)$.  The $L$-functions
attached to $M_g$ are related to those attached to the $\tau(g)$
by the formula
$$L(M_g\otimes M_{\psi_g^{-1}},\tau,s) = L(\tau(g),s).$$
In fact, a stronger statement is true:
For any primes $p$ of $\QQ$ and $\lambda$ of $K$, the representation
$D_{\pst}(M_\lambda|G_p)^{\ss}$ of $'W_p$ is $K$-rational
and corresponds via local Langlands to $\pi_p(\tau(g))$
(where we extend scalars to $\CC$ via $\tau$ and
normalize the local Langlands correspondence as in
\cite{Ca0}).  This is due to Eichler, Shimura and Igusa
for $\lambda$ not dividing $p N_g$, by
Langlands, Deligne and Carayol~\cite{Ca0}
for  $\lambda\nmid p$, and by Scholl~\cite{scholl}
and Saito~\cite{saito} in general.  In particular $M_g$
is $L$-admissible everywhere; moreover $M_g = M^{S_N}$
for an object $M$ of $\pms_K$ which is $L$-admissible everywhere.
Using the compatibility with the local Langlands correspondence,
it also follows from the functional equations for the
$L(\tau(g),s)$ that Conjecture \ref{fe} holds for $M_g$.

Finally we recall some basic properties of the local Langlands
correspondence.  In the following, we identify characters of
$\QQ_p^\times$ with those of $'W_p$ via class field theory.
\begin{lemma} Suppose that $\pi$ is an irreducible admissible
representation of $\GL_2(\QQ_p)$ and $\rho:{} 'W_p \to \aut_\CC V$
corresponds to $||\det||^{1/2}\pi$ via local Langlands.
Then the following hold:
\begin{enumerate}
\item The conductors of $\pi$ and $\rho$ coincide.
\item The central character of $\pi$ is the determinant of $\rho$.
\item For any character $\theta$ of $\QQ_p^\times$,
 $||\det||^{1/2}(\theta\circ\det)\pi$ corresponds to $\theta\rho$.
\item The trace of $\rho(\frob_p)$ on $(V/NV)^{I_p}$
is the eigenvalue of $t_p$ on $\pi^{V_p}$.
\item $\pi$ is supercuspidal if and only if
$\rho$ is irreducible.
\end{enumerate}
\end{lemma}
These properties already determine the correspondence explicitly
if $\pi$ is principal series or special.  In particular, if
$\pi = I(\mu_1,\mu_2)$, then $\rho|_{W_p}$ is equivalent to
$\mu_1 \oplus \mu_2$ and $N \neq 0$ if and only if $\pi$ is
special.  We shall also need the explicit description of the correspondence
due to G\'erardin \cite{gerardin} in certain cases where $\pi$ is
supercuspidal.  For this we refer to sections 3.2 and 4.2 of \cite{cdt}.

\section{The adjoint premotivic structure}
\mlabel{sec:adjoint}

\subsection{Realizations of the adjoint premotivic structure}
\mlabel{ss:adj.real}

$A_f =\ad^0 M_f$ is defined to be the kernel of the trace map
$$\hom_K (M_f,M_f) \to K.$$  It is a premotivic structure
in $\pms^S_K$ for $S\supseteq S_N$.

For $\wc=B, \dr$ or $\lambda$, $A_{f,\wc}$ has an integral structure
given by
\[ \CA_{f,\wc} = \{ a\in \End(\CM_{f,\wc}) | \tr (a) =0 \}. \]
The extra structures on the realizations of $\CA_f$ are obtained
by restrictions from those of $\End(\CM_f)$.
For example the filtration on $\CA_{f,\dr}$ is given by
\begin{eqnarray*}
 &&\fil^n \CA_{f,\dr} = \{a\in \CA_\dr\subseteq \End(\CM_{f,\dr}) |
    a(\Fil^i \CM_{f,\dr}) \subseteq \Fil^{n+i}\CM_{f,\dr},
    \forall j \} \\
&=&\hspace{-.3cm} \left \{ \hspace{-.2cm} \begin{array}{ll}
    \CA_{f,\dr},& n\leq -(k-1),\\
    \{a\in \CA_{f,\dr} |
        a(\Fil^0 \CM_{f,\dr})\subseteq \Fil^0 \CM_{f,\dr}\}, &
        -(k-1)<n\leq 0, \\
    \{a\in \CA_{f,\dr} |
        a( \CM_{f,\dr})\subseteq \Fil^0 \CM_{f,\dr},\
        a(\Fil^0 \CM_{f,\dr}) =0 \}, \hspace{-.2cm} & \
        0<n\leq k-1, \\
    0, & n>k-1.
    \end{array} \right .
\end{eqnarray*}
The filtration on $\End (\CM_{f,\lcrys})$ is the same.
Since the the weight is from
$1-k$ to $k-1$, we need $p-1>2(k-1)$ instead of $p>k$
(plus a twist of the
crystalline realization) in order to get an object in $\mfcat^0$
and to make $\CA_{f,\lcrys}$ well-defined.

\begin{lemma}
There is an isomorphism $\det_{K}A_f\cong K$ in $\pms^S_K$ which
restricts to isomorphisms
$\det_{\CO_K} \CA_{f,\wc}[1/N k!] \cong \CO_K[1/N k!]$
for $\wc\in \{\bt,\dr,\lambda\}$.
\mlabel{lem:adjoint}
\end{lemma}
\proofbegin
Setting
$M_f^*=\Hom_{K}(M_f,K)$
and noting that $M_f$ has $K$-rank 2 we find
\begin{align*}\det_{K}A_f\otimes_{K}K\cong&\det_{K}(A_f\oplus K)\cong\det_{K}\Hom_{K}(M_f,M_f)\cong
\det_{K}(M_f^*\otimes_{K}M_f)\\
\cong&\det_{K}(M_f^*)^{\otimes
2}\otimes_{K}\det_{K}(M_f)^{\otimes 2}\\
\cong&\left(\det_{K}(M_f)^*\otimes_{K}\det_{K}(M_f)\right)^{\otimes
2}\\
\cong&K^{\otimes 2}\cong K.
\end{align*}
The integral isomorphisms are proved in the same way.
\proofend

We note that $A_f$ and $\hom_K(A_f,K)$ are isomorphic.
There is in fact a canonical isomorphism defined by the pairing
\begin{equation}
\alpha \otimes \beta \mapsto \tr(\alpha\circ\beta)
\mlabel{eq:dual}
\end{equation}
on each realization of $A_f$.

Note that if we replace $K$
by $K' \supset K$ and $S$ by a subset $S'$ of the primes
over those in $S$, then $A_f$ is replaced by
$(A_f\otimes_K K')^{S'}$.
If $\psi'$ is a character $\AA^\times \to K^\times$ of conductor $D'$
and $f\otimes\psi'$ denotes the eigenform (of weight $k$,
conductor dividing $N_f D'{}^2$ and character $\psi(\psi')^2$) associated
to the eigenform $\sum_{(n,D')=1}\psi'(n_n) a_n e^{2\pi i nz}$.
We will relate $\CM_f$ to the premotivic structure associated
to $f\otimes \psi'$.

\begin{lemma}
Let $S\supseteq S_f$ be such that $D'$ is invertible in $\CO_S$.
We have $\CM_{f\otimes\psi'} \cong \CM_f \otimes_{\CO} \CM_{\psi'}$
in $\ipms^S_K$.
We also have
$A_{f\otimes\psi'} \cong A_f$ in $\pms^{S}_K$.
\mlabel{lem:twist}
\end{lemma}

It is sometimes convenient
to assume that $f$ has minimal conductor among its twists.

\noindent
\proofbegin
For the given character $\psi':\hat{\ZZ}^\times \to
\CO'{}^\times$ of conductor $D'$,
let $F=\QQ(\zeta_{D'M})$, where $M\geq 3$ is as in \S\ref{ssec:levchar},
and let $K'=K(\zeta_{D'{}^2})$.
Let $N'={\rm lcm}(N,D')$ and $N''={\rm lcm} (N,DD',D'{}^2)$
where $D$ is the conductor of $\psi$.
Then we have a representation given
by
$\psi\psi'{}^2: U_0(N'') \to (\CO_{K'})^\times$.

With this setup, and $U=U_0(N), U'=U_0(N'), U''=U_0(N'')$,
we define
$$\theta: \CM (N,\psi)_{D'M} \otimes \CM_{\psi'} \to
\CM(N'',\psi\psi'{}^{2})_{(D')^3M}$$
to be the composite from the left column of
the following diagram.

$$
\begin{CD}
\bear{c} \CM(N,\psi)_{D'M}\otimes \CM_{\psi'} \\
=(\CM_{D'M}\otimes
\CV_{\psi^{-1}})^{U}\otimes \CM_{\psi'}
\enar
@> \subset >>
\bear{c}
\CM_{D'M}\otimes \CV_{\psi^{-1}} \otimes \CM_F\otimes
\CO_{\psi'} \cong \\
\CM_{D'M}\otimes \CV_{\psi^{-1}} \otimes \CH^0(X_{D'M})
\otimes \CO_{\psi'} \enar \\
@V \cup VV @V \cup VV \\
(\CM_{D'M}\otimes
\CV_{\psi^{-1}}\otimes \CO_{\psi'})^{U'}
@> \subset >>
\CM_{D'M}\otimes \CV_{\psi^{-1}} \otimes \CO_{\psi'} \\
@V [U''\smat{1}{1/D'}{0}{1} U']_1 VV
@V \sum_d [g_d g'] \otimes g_d \id VV \\
\bear{c} \CM(N'',\psi\psi'{}^2)\\
=(\CM_{D'{}^3 M} \otimes
\CV_{\psi^{-1}\psi'{}^{-2}})^{U''}
\enar
@> \subset >>
\CM_{D'{}^3M}\otimes\CV_{\psi^{-1} \psi'{}^{-2}}.
\end{CD}
$$
We will show that it restricts to an isomorphism
$\CM_f\otimes \CM_{\psi'} \to \CM_{f\otimes \psi'}$
stated in the lemma.

We first show that $\theta$ is well-defined.
In the top half of the diagram,
under the isomorphism
$$ \CM_{D'M}\otimes \CV_{\psi^{-1}} \otimes \CM_F\otimes
\CO_{\psi'{}^{-1}} \cong
\CM_{D'M}\otimes \CV_{\psi^{-1}} \otimes \CH^0(X_{D'M})
\otimes \CO_{\psi'}
\to
\CM_{D'M}\otimes \CV_{\psi^{-1}} \otimes \CO_{\psi'}$$
where the first map is from the proof of Lemma~\ref{lem:psi} and
the second map is from cup product, we have
\begin{align*}\CM(N,\psi)_{D'M} \otimes \CM_{\psi'{}^{-1}}
&\cong
(\CM_{D'M}\otimes \CV_{\psi^{-1}})^U \otimes (\CH^0(X_{D'M})
\otimes \CO_{\psi'})^{U'} \\
&\subseteq
(\CM_{D'M}\otimes \CV_{\psi^{-1}} \otimes \CH^0(X_{D'M})
\otimes \CO_{\psi'})^{U'} \\
&\to
(\CM_{D'M}\otimes \CV_{\psi^{-1}} \otimes
\CO_{\psi'})^{U'}.
\end{align*}
Here the first map is by Lemma~\ref{lem:psi}, the inclusion map is
because $U'\subseteq U$ and the third map is because the cup
product is $U'$-equivariant.

We next consider the bottom half of the diagram.
Let $\sigma'$ (resp. $\sigma''$)
be the representation of $U'$ on $\CV'=\CV_{\psi^{-1}}\otimes
\CO_{\psi'{}^{-1}}$ (resp. of $U''$ on
$\CV''=\CV_{\psi^{-1}\psi'{}^{-2}}$).
Let $\tau: \CV'\to \CV''$ be the identify map of the underlying
$\CO$-modules. Let $g'=\smat{1}{1/D'}{0}{1}$.
If $u$ is an element of $U''\cap (g'U' g'{}^{-1})$,
then $u$ is of the form
$$\smat{a+D'N'c}{b+D'{}^{-1}(d-a)-D'{}^{-2}N'c}{N'c}{-D'{}^{-1}N'c+d}$$
for some $\smat{a}{b}{N'c}{d}\in U'$.
Since $u$ is in $U''$, we have
$N''|N'c$, so $\frac{N''}{N'} |c$. Then $\frac{N}{{\rm gcd}(D',N)}
\frac{N''}{N'}|\frac{N}{{\rm gcd}(D',N)} c$.
Since
$\frac{N}{{\rm gcd}(D',N)} \frac{N''}{N'}=\frac{N''}{D'}$, we have
$D'|\frac{N}{{\rm gcd}(D',N)} c$.
{}From this and the $(1,2)$-entry of $u$ we get
$d-a \equiv 0 \mod D'$.
Thus
\begin{align*}
\tau(\sigma'(g'{}^{-1}ug') v)&=\tau(\psi^{-1}(a_N)\psi'{}^{-1}\circ
\det (g'{}^{-1}ug')v)\\
&=\tau(\psi^{-1}(a_N) \psi'(ad)^{-1}v)\\
&=\tau(\psi^{-1}(a_N)\psi'(a)^{-2} v)\\
&=\psi^{-1}(a_N)
\psi'(a)^{-2} v\\
&=\sigma''(\tau(v)).
\end{align*}
This shows that
$$[U'' \smat{1}{1/D'}{0}{1} U']_1:
(\CM_{D'M}\otimes \CV_{\psi^{-1}} \otimes \CO_{\psi'})^{U'}
\to \CM(N,\psi^{-1}\psi'{}^{-2})_{(D')^3M}$$
is well-defined. We also note that,
letting $\Lambda\subset \hat{\ZZ}^\times$ be a complete set of
coset representatives of $(\ZZ/D'\ZZ)^\times,$
then $g_d=\smat{1}{0}{0}{d}, d\in \Lambda$ gives a coset
decomposition $U'' =\coprod_d g_d (U''\cap (g'U'g'{}^{-1}))$.

Thus $\theta$ is well-defined and is the restriction of the composite
of the maps in the bottom row in the above diagram.
To prove the first statement of the lemma, we only need to show
\begin{equation}
\theta_\dr (f \otimes b_\dr)= D' \psi'(-1)f_{\psi'}.
\mlabel{eq:twist}
\end{equation}
Tracing through the maps in the above diagram, we see that
$$f\otimes b_\dr  \in \fil^{k-1}\CM(N,\psi)_{D'M,\dr}\otimes
\CM_{\psi',\dr}$$
is sent to
$$\sum_{a\mod D'} \sum_{d\mod D'}
    [g_d g] \zeta_{D'}^{-a} f\otimes g_d~ (1 \otimes \psi'(a))
\in \CM_{D'{}^3M}\otimes \CV_{\psi^{-1}\psi'{}^{-2}}.$$
Under the isomorphism
$$\CC\otimes \fil^{k-1} \CM_{D'M,\dr}\cong \oplus_{t\in (\ZZ/D'M\ZZ)^\times}
M(\Gamma(D'M)),$$ the element
$\zeta_{D'}^{-a} f$ is sent to the vector
$(\zeta_{D'}^{-t^{-1}a} f_t)_t,\ t\in (\ZZ/D'M\ZZ)^\times$.
For a fixed $d\in \Lambda$ and for $t'\in (\ZZ/D'{}^3M\ZZ)^\times$,
in the notation of Lemma~\ref{lem:adan} and
Proposition~\ref{prop:mot-ad}, we have
$\kappa(t')= (\det g_d ||\det g_d||)^{-1}t' =d^{-1} t' \mod D'M$
and
$$\gamma=\gamma_{d,t'}=
\left (\smat{1}{0}{0}{t'{}^{-1}} \smat{1}{D'{}^{-1}}{0}{d}
\smat{1}{0}{0}{d^{-1}t'}\right )^{-1} = \smat{1}{-d^{-1}t'D'{}^{-1}}{0}{1}.$$
Then $[g_dg]\otimes g_d $ sends
$(\zeta_{D'}^{-t^{-1}a} f_t)_t\otimes (1\otimes \psi'(a))$ to
$$\left (\zeta_{D'}^{-d(t')^{-1} a} f_{d^{-1}t'} |_{\smat{1}{-d^{-1}t'D'{}^{-1}}{0}{1}}
    \right)_{t'}
\otimes \psi'(a)$$ in
$\CC \otimes \fil^{k-1} \CM_{D'{}^3M}
\cong \oplus_{t'\in (\ZZ/D'{}^3M\ZZ)^\times}
M_k(\Gamma_0(D'{}^3M),\psi^{-1}\psi'{}^{-2}).$
Thus the $t'$-th entry of
$$\sum_{a\mod D'} \zeta_{D'}^{-a}
    \sum_{d\mod D'} [g_d g] f\otimes g_d \tau (1 \otimes \psi'(a))
$$
in the above diagram is
$$\sum_{a\mod D'} \sum_{d\mod D'}
\zeta_{D'}^{-d^{-1}t'a} f_{d^{-1}t'} |_{\smat{1}{-d^{-1}t'D'{}^{-1}}{0}{1}}
\otimes \psi'(a).
$$
Since $f$ is from $\fil^{k-1}\CM(N,\psi)_{D'M,\dr}$, by
Lemma~\ref{lem:forms} and its proof,
the components of $(f_t)_t$ are the same.
Then we have
\begin{align*}
\sum_{a\mod D'} & \sum_{d\mod D'}
\zeta_{D'}^{-d^{-1}t'a} f_{d^{-1}t'} |_{\smat{1}{-d^{-1}t'D'{}^{-1}}{0}{1}}
\otimes \psi'(a)\\
&=\sum_{d\mod D'} (\sum_{a\mod D'}
\zeta_{D'}^{-d^{-1}t'a} \otimes \psi'(a)) f_{d^{-1}t'} |_{\smat{1}{-d^{-1}t'D'{}^{-1}}{0}{1}}
\\
&=\sum_{d\mod D'} G_{\psi'}  \psi'(-d^{-1}t')^{-1}
f_{\kappa(t')} |_{\smat{1}{-d^{-1}t'D'{}^{-1}}{0}{1}}\\
&=D'G_{\psi'{}^{-1}}^{-1} \psi'(-1)
    \sum_{d\mod D'} \psi'(-d^{-1}t')^{-1}
f_{\kappa(t')} |_{\smat{1}{-d^{-1}t'D'{}^{-1}}{0}{1}}
\notag
\end{align*}
where the last equation comes from
$D' = \psi'(-1) G_{\psi'} G_{\psi'{}^{-1}}$.
By \cite[Lemma 4.3.10]{miyake}, the last term is the image of
$D'\psi'(-1) f_{\psi'}$ under the isomorphism
$\CC\otimes \fil^{k-1} \CM_{D'{}^3M,\dr}\cong \oplus_{t'}
M(\Gamma_0(D'{}^3M),\psi\psi'{}^2)$. Then by
Proposition~\ref{prop:mot-ad} and its proof, we have
(\ref{eq:twist}), as needed.

Since $\psi'$ is a character, we have
$$\hom(M_f\otimes M_{\psi'}, M_f\otimes M_{\psi'})
=\hom(M_f, M_f\otimes M_{\psi'}\otimes
M_{\psi'}^{\otimes (-1)})
=\hom(M_f, M_f).$$
It then follows that
$A_{f\otimes\psi'} \cong A_f$.
\proofend

\subsection{Euler factors and functional equation}
\mlabel{ssec:fe}
For each prime $p$, we let $\delta_p$ denote the dimension of
$M_{f,\lambda}^{I_p}$ for any $\lambda$ not dividing $p$, so
$$\delta_p = \left\{\begin{array}{ll}
    2,\mbox{if $p\nmid N_f$,}\\
    1,\mbox{if $p | N_f$ and $a_p\neq 0$,}\\
    0,\mbox{if $p | N_f$ and $a_p = 0.$}\end{array}\right.$$
We set $L_p^\naive(A_f,s) = L_p(A_f,s)$ if $\delta_p > 0$, and
$L_p^\naive(A_f,s) = 1$ if $\delta_p = 0$.
We let $\Sigma_e=\Sigma_e(f)$ denote the set of primes $p$
such that $\delta_p = 0$ and $L_p(A_f,s) \neq 1$.
, and set
$$L^\naive(A_f,s) = \prod_p L_p^\naive(A_f,s)
    = \prod_{p\not\in\Sigma_e(f)} L_p(A_f,s).$$
We call the primes in $\Sigma_e$ {\em exceptional} for $f$.

Recall that if $\delta_p = 2$, then writing
$$L_p(f,s) = (1-a_p p^{-s} + \psi(p)p^{k-1-2s})^{-1}
=(1-\alpha_p p^{-s})^{-1}(1-\beta_p p^{-s})^{-1},$$
we have
$$
L_p(A_f, s) = (1-\alpha_p\beta_p^{-1}p^{-s})^{-1}
(1-p^{-s})^{-1}(1-\alpha_p^{-1}\beta_pp^{-s})^{-1}.$$
If $\delta_p = 1$, then

\begin{equation}
 L_p(A_f, s) = \left\{\begin{array}{ll}
(1-p^{-1-s})^{-1}&\mbox{if $\pi_p(f)$ is special;}\\
(1-p^{-s})^{-1}&\mbox{if $\pi_p(f)$ is principal series.}
\end{array}\right.
\mlabel{eq:lp1}
\end{equation}
To see this, note that if $\pi_p(f)$ is special, then $\rho_{f,\lambda}
|_{G_p}$ has matrix $\chi\otimes \smat{\chi_\ell}{*}{0}{1}$
and is indecomposible. Here $\chi$ is unramified and $\chi_\ell$ is the
cyclotomic character. So $\ad^0\rho_{f,\lambda}$ has matrix
$\left ( \begin{array}{ccc} \chi_\ell & * & *\\ &1&*\\
    && \chi_\ell^{-1} \end{array} \right )$.
So $\ad^0 \rho_{f,\lambda} ^{I_p}$ is one dimensional acted by $G_p$
through the
character $\chi_\ell$. Since $\chi_\ell(\Frob_p)=p^{-1}$, we have
$L_p(A_f,s)=(1-p^{-1-s})^{-1}$.

If $\pi_p(f)$ is principal series, then $\rho_{f,\lambda}
|_{G_p}$ has matrix $\smat{\chi_1}{0}{0}{\chi_2}$ with
$\chi_1$ ramified and $\chi_2$ unramified.
So $\ad^0\rho_{f,\lambda}$ has matrix
$\left ( \begin{array}{ccc} \chi_1\chi_2^{-1}&  & \\ &1& \\
    && \chi_1^{-1}\chi_2 \end{array} \right )$.
So $\ad^0 \rho_{f,\lambda} ^{I_p}$ is one dimensional with trivial
$G_p$-action. So
$L_p(A_f,s)=(1-p^{-s})^{-1}$.

Shimura \cite{shim_plms} proved that $L(A_f,s)$ extends to an entire
function on the complex plane.
Recall that we regard $L(M,s)$ as taking
values in $K \otimes \CC$.  Each embedding
$\tau: K \to \CC$ gives a map $K\otimes\CC \to \CC$ and
we write $L(M,\tau,s)$ for the composite with $L(M,s)$.
Moreover, the work of Gelbart-Jacquet
\cite{G-J} and others (see \cite{schmidt}) shows that
$$\begin{array}{rcl}
\Lambda(A_f,s) & =& L(A_f,s)\Gamma_\RR(s+1)\Gamma_\CC(s+k-1)\\
 & = & 2^{2-k-s}\pi^{(1-2k-3s)/2}
L(A_f,s)\Gamma(\frac{s+1}{2})\Gamma(s+k-1)
\end{array}$$
satisfies the functional equation
\[
\Lambda(A_f,s) = \epsilon(A_f,s)\Lambda(A_f,1-s),
\]
where $\epsilon(A_f,s)$ is as defined by Deligne \cite{del_ant}.
Here we have used that $A_f$ and $\hom_K(A_f,K)$ are isomorphic
(using (\ref{eq:dual})).

\subsection{$\Sigma$-variations for non-exceptional primes}
\mlabel{ss:sig}
Let $\Sigma$ denote a finite set of rational primes.  We assume
that no prime in $S$ divides a prime in $\Sigma$ and
$\Sigma_e\subseteq \Sigma$.
The case when $\Sigma_e\not\subseteq \Sigma$ will be considered in
\S \ref{sec:ep}.
 We define a premotivic
structure $M_f\sig$ exactly as above, but we replace
$N_f$ by $N\sig = N_f \dprod_{p\in\Sigma} p^{\delta_p}$
and
$f(z)$ by the eigenform $f\sig(z) = \dsum a_n' e^{2\pi i nz}$,
where $a_n' = 0$ if $n$ is divisible by a prime in $\Sigma$,
and $a_n' = a_n$ otherwise.
We will next related the premotivic structure of $f\sig$ to
the premotivic structure of $f$.

For positive integers $m$ dividing
$N^\Sigma/N_f=\prod_{p\in \Sigma} p^{\delta_p}$, we let
$$\epsilon_m = m^{-1} \left[U^\Sigma\smat{m^{-1}}{0}{0}{1}U\right]_1
    =m^{1-k}\left[ U\sig \smat{1}{0}{0}{m} U\right]_1:
    M(N,\psi) \to M(N\sig,\psi)$$
and define the endomorphism $\phi_m$ of $M(N,\psi)$ by
\begin{itemize}
\item $\phi_1 = 1$, $\phi_p = -T_p$, $\phi_{p^2} = pS_p$;
\item $\phi_{m_1m_2} = \phi_{m_1}\phi_{m_2}$ if $(m_1,m_2) = 1$.
\end{itemize}
\begin{prop}
The morphism
$$\gamma = \sum_{m}\epsilon_m\phi_m : M(N,\psi)
    \to M(N\sig,\psi)$$
induces an isomorphism $M_f \to M_f^\Sigma$ in
$\pms^{S}$.
Furthermore, $\gamma_\dr (f)=f\sig$.
\mlabel{prop:level}
\end{prop}
\proofbegin
We first reduce the proof to the case when $\Sigma$ contains
a single prime.

Let $\Sigma'$ be a subset of $\Sigma$. Define $U^{\Sigma'}$
and $M(N^{\Sigma'},\psi)$ as above.
Let $\Sigma_0\subsetneq \Sigma_1$ be subsets of $\Sigma$ and
let $m$ be a positive integer dividing $N^{\Sigma_1}/N^{\Sigma_0}$.
Define
$$\epsilon_m^{\Sigma_0,\Sigma_1}: m^{-1} [U^{\Sigma_1}
\smat{m^{-1}}{0}{0}{1} U^{\Sigma_0}]_1
=m^{1-k} [U^{\Sigma_1} \smat{1}{0}{0}{m} U^{\Sigma_0}]_1.$$
Also define
$\phi_m^{\Sigma_0}$ to be the $\phi_m$ defined above acting on
$M(N^{\Sigma_0},\psi)$
and define
$\gamma_m^{\Sigma_0,\Sigma_1}=\epsilon_m^{\Sigma_0,\Sigma_1}\circ
\phi_m^{\Sigma_0}.$

\begin{lemma}
Let $\Sigma_0\subsetneq \Sigma_1\subsetneq \Sigma_2$
be subsets of $\Sigma$.
Let $m_1$ (resp. $m_2$) be positive integers dividing
$N^{\Sigma_1}/N^{\Sigma_0}$ (resp. $N^{\Sigma_2}/N^{\Sigma_1}$).
The following equations hold.
\begin{enumerate}
\item
$$\epsilon_{m_1m_2}^{\Sigma_0,\Sigma_2}
=\epsilon_{m_2}^{\Sigma_1,\Sigma_2}
\circ \epsilon_{m_1}^{\Sigma_0,\Sigma_1}. $$
\item
$$\phi_{m_2}^{\Sigma_1} \circ
\epsilon_{m_1}^{\Sigma_0,\Sigma_1}
=\epsilon_{m_1}^{\Sigma_0,\Sigma_1}
\circ \phi_{m_2}^{\Sigma_0}.$$
\item
$$ \gamma^{\Sigma_0,\Sigma_2}_{m_1m_2} =
\gamma^{\Sigma_1,\Sigma_2}_{m_2} \circ
    \gamma^{\Sigma_0,\Sigma_1}_{m_1}.$$
\end{enumerate}
\mlabel{lem:level}
\end{lemma}
\proofbegin
(a). We only need to prove
$$ [U^{\Sigma_2} \smat{1}{0}{0}{m_2} U^{\Sigma_1}]_1\circ
[U^{\Sigma_1} \smat{1}{0}{0}{m_1} U^{\Sigma_0}]_1
= [U^{\Sigma_2} \smat{1}{0}{0}{m_1m_2} U^{\Sigma_0}]_1.$$
In Proposition~\ref{prop:scomp} take
$U=U^{\Sigma_0},
U'=U^{\Sigma_1},
U''=U^{\Sigma_2}, 
g=\smat{1}{0}{0}{m_1},g'=\smat{1}{0}{0}{m_2}$
and $\tau=\tau'=\id$.
Then we have $U'_1=U'$. So in the double coset decomposition
$U'=\coprod_k U'_2 u_k U'_1$ there is only one $u_k$ which we
can take $I$. Further $W_1=W_1'=U''$ and $\tau_1=1$.
This gives us the desired equation.

(b).
Because of (a) and similar multiplicity property of $\phi_m$,
we only need to prove (b) in the special case $\Sigma_1\backslash
\Sigma_0=\{p\}$ and $\Sigma_2\backslash \Sigma_1=\{q\}$ where $p$
and $q$ are primes.

So we just need to prove the commutativity of the diagram
$$\begin{CD}
M(N^{\Sigma_0},\psi) @>\phi_{q^i}^{\Sigma_0} >> M(N^{\Sigma_0},\psi)\\
@V \epsilon_{p^j}^{\Sigma_0,\Sigma_1} VV @V \epsilon_{p^j}^{\Sigma_0,\Sigma_1}
VV\\
M(N^{\Sigma_1},\psi) @ >>\phi_{q^i}^{\Sigma_1} > M(N^{\Sigma_1},\psi)
\end{CD}
$$
for $0\leq i, j \leq 2$.

For this we first display the following elementary relations.
Fix a prime number $p$. For each integer $i\geq 0$, let
$V_i = \{ \smat{a}{b}{c}{d} \in \GL_2(\ZZ_p) |  c\equiv 0 \mod
p^i\}.$
Let $A_t =\smat{1}{0}{0}{p^t},t=0,1,2$.
Let $B_0=I, B_1= \smat{p}{0}{0}{1}, B_2=\smat{p}{0}{0}{p}$.
Then
\begin{equation}
[ V_2 I V_0 ] [V_0 B V_0]
 = [V_2 B V_0] + [V_2 A V_0]
\mlabel{eq:B0}
\end{equation}
and
\begin{equation}
[V_2 B V_2][V_2 A_t V_0]=\left \{ \begin{array}{ll}
    [V_2 BV_0], & t=0, \\
    p[V_2 B A_t V_0]=p^{k-1}[V_2 A_{t-1} V_0],
    & 0<t\leq 2. \end{array} \right .
\mlabel{eq:AB0}
\end{equation}
If $i > 0$, then
\begin{equation}
[ V_{i+s} B_t V_{i+s} ] [V_{i+s} I V_i] = [V_{i+s} B_t V_i]
  =  [V_{i+s} I V_i] [V_i B_t V_i], 0\leq s,t \leq 2
\mlabel{eq:B}
\end{equation}
and
\begin{equation}
[ V_{i+1} B V_{i+1} ] [V_{i+1} A V_i]
  =p [V_{i+1} pI V_i] =p^{k-1} [V_{i+1} I V_i].
\mlabel{eq:AB}
\end{equation}
These relations can be easily verified applying Proposition~\ref{prop:prod}
to double cosets in $\GL_2(\QQ_p)$ instead of to $\GL_2(\AA_\f)$.

We now continue with the proof of Lemma~\ref{lem:level} (b)
and verify
\begin{equation}
\epsilon_{p^j}^{\Sigma_0,\Sigma_1} \circ \phi_{q^i}^{\Sigma_0}
= \phi_{q^i}^{\Sigma_1} \circ \epsilon_{p^j}^{\Sigma_0,\Sigma_1}.
\mlabel{eq:lem1}
\end{equation}
Since a double coset product in $\GL_2(\AA_\f)$ is determined
by its local factors at each finite prime,
we only need to verify the equation at each prime.
The equation at primes different from $p$ and $q$ follows from
$[V_{i+s} I V_i][V_i I V_i]=[V_{i+s} I V_{i+s}][V_{i+s} I V_i]$
which is obvious.
At the prime $p$, the equation follows from
$$[ V_{i+s} A_t V_{i} ] [V_{i} I V_i] = [V_{i+s} A_t V_i]
  =  [V_{i+s} I V_{i+s}] [V_{i+s} A_t V_i], 0\leq s,t \leq 2,
$$
which is again obvious.
At the prime $q$, the equation follows from (\ref{eq:B}) by
choosing $(s,t)=(0,i)$.
This proves
(\ref{eq:lem1}) and hence (b) of Lemma~\ref{lem:level}.

(c) follows from (a) and (b):
\begin{multline}
 \gamma^{\Sigma_0,\Sigma_2}_{m_1m_2} =
    \alpha^{\Sigma_0,\Sigma_2}_{m_1m_2}
    \circ \beta^{\Sigma_0,\Sigma_2}_{m_1m_2}
= \alpha^{\Sigma_1,\Sigma_2}_{m_2} \circ
    \alpha^{\Sigma_0,\Sigma_1}_{m_1} \circ
    \beta^{\Sigma_1,\Sigma_2}_{m_2} \circ
    \beta^{\Sigma_0,\Sigma_1}_{m_1} \\
= \alpha^{\Sigma_1,\Sigma_2}_{m_2} \circ
    \beta^{\Sigma_1,\Sigma_2}_{m_2} \circ
    \alpha^{\Sigma_0,\Sigma_1}_{m_1} \circ
    \beta^{\Sigma_0,\Sigma_1}_{m_1}
=\gamma^{\Sigma_1,\Sigma_2}_{m_2} \circ
    \gamma^{\Sigma_0,\Sigma_1}_{m_1}.
\notag
\end{multline}
This completes the proof of Lemma~\ref{lem:level}.
\proofend

We continue with the proof of Proposition~\ref{prop:level}.
For the given $\Sigma$ and $m$, let $p_1,\ldots,p_r$ be the
primes in $\Sigma$. Let
$\Sigma_0=\phi$ and $\Sigma_i=\{p_1,\ldots,p_i\},\ 1\leq i\leq r$.
With this convention, we have
$\epsilon_m=\epsilon_m^{\Sigma_0,\Sigma_r}$ and
$\phi_m=\phi_m^{\Sigma_p}$.
Any positive divisor of $\prod_{i=1}^r p_i^{\delta_{p_i}}$ is of the
form $\prod_{i=1}^r p^{e_i}, 0\leq e_i \leq \delta_{p_i}$.
Recall that before the lemma, we
define
$$\gamma_{\prod_i p_i^{e_i}}^{\Sigma_0,\Sigma_r}=
\epsilon_{\prod_i p_i^{e_i}}^{\Sigma_0,\Sigma_r} \circ
\phi_{\prod_i p_i^{e_i}}^{\Sigma_0,\Sigma_r}.$$
Using Lemma~\ref{lem:level}, we obtain
$$\gamma_{\prod_i p_i^{e_i}}^{\Sigma_0,\Sigma_r}=
    \circ_{i=1}^r \gamma_{p_i^{e_i}}^{\Sigma_{i-1},\Sigma_i}.$$
We therefore have
$$
\gamma=\sum_{i=1}^r \sum_{e_i=0}^{\delta_{p_i}}
    \gamma_{\prod_i p_i^{e_i}}^{\Sigma_0,\Sigma_r}
=\sum_{i=1}^r \sum_{e_i=0}^{\delta_{p_i}}
    \circ_{i=1}^r \gamma_{p_i^{e_i}}^{\Sigma_{i-1},\Sigma_i}
=\circ_{i=0}^r (\sum_{e_i=0}^{\delta_{p_i}}
\gamma_{p_i^{e_i}}^{\Sigma_{i-1},\Sigma_i}).
$$
Thus to prove Proposition~\ref{prop:level},
we only need to prove it in the case when
$\Sigma$ contains one prime $p$.
Let $\Sigma_0=\phi$ and $\Sigma_1=\Sigma$.
We consider three cases depending on
$\delta_p=0, 1$ or $2$.

{\bf Case 1: } When $\delta_p=0$, then $N^\Sigma=N$ and $\gamma$
is just the identity map. We have also seen that $a_p=0$.
Thus $a_{pn}=0$. So $f=f^\Sigma$ and the proposition is
proved in this case.

{\bf Case 2: } When $\delta_p=1$,
$$\gamma=\epsilon_1\circ \phi_1 +\epsilon_p\circ \phi_p.
$$
By Lemma~\ref{lem:level}, for $q\neq p$,
\begin{align*}
T_q \circ \gamma &=
    -\phi_q^{\Sigma_1} \circ \epsilon_1^{\Sigma_1,\Sigma_0} \circ
\phi_1^{\Sigma_0}
-\phi_q^{\Sigma_1} \circ \epsilon_p^{\Sigma_1,\Sigma_0} \circ
\phi_0^{\Sigma_0} \\
& =    - \epsilon_1^{\Sigma_1,\Sigma_0} \circ \phi_q^{\Sigma_0} \circ
\phi_1^{\Sigma_0}
-\epsilon_p^{\Sigma_1,\Sigma_0} \circ \phi_q^{\Sigma_1} \circ
\phi_0^{\Sigma_0} \\
& =    - \epsilon_1^{\Sigma_1,\Sigma_0} \circ
\phi_1^{\Sigma_0} \circ \phi_q^{\Sigma_0}
-\epsilon_p^{\Sigma_1,\Sigma_0}  \circ
\phi_0^{\Sigma_0} \circ \phi_q^{\Sigma_1}\\
& =\gamma \circ T_q.
\end{align*}
When $q=p$, by  (\ref{eq:B}) and  (\ref{eq:AB}),
the double coset product of the local factor at primes
other than $p$ involves $I$ and is obviously zero.
The local factor at $p$ is
\begin{align*}
&[V_{i+1}B V_{i+1}]([V_{i+1} A_0 V_i][V_i B_0 V_i]
-p^{1-k}[V_{i+1} A_1 V_i][V_{i+1} B V_i])\\
&=[V_{i+1} BV_i][V_i IV_i]-p^{1-k} p [V_{i+1} pI V_i][V_i B
V_i]\\
&=[V_{i+1} IV_i][V_i B V_i] -p^{2-k} p^{k-2} [V_{i+1} IV_i][V_i
BV_i]=0.
\end{align*}
So we get
\begin{align*}
T_p \circ \gamma &= -\phi_p^{\Sigma_1} \circ
\epsilon_1^{\Sigma_0,\Sigma_1}\circ \phi_1^{\Sigma_0}
-\phi_p^{\Sigma_1} \circ \epsilon_p^{\Sigma_0,\Sigma_1} \circ
\phi_p^{\Sigma_0} \\
&= - \epsilon_1^{\Sigma_0,\Sigma_1}\circ\phi_p^{\Sigma_0} \circ
 \phi_1^{\Sigma_0}
- p^{1-k} p^{k-2} \epsilon_1^{\Sigma_0,\Sigma_1}
    \circ \phi_p^{\Sigma_0} \\
&= - p\epsilon_1^{\Sigma_0,\Sigma_1}\circ\phi_p^{\Sigma_0}
+ p^{1-k} p^{k-2} \epsilon_1^{\Sigma_0,\Sigma_1}
    \circ \phi_p^{\Sigma_0}
=0.
\end{align*}
This shows that $\gamma_\dr f$ is a multiple of $f\sig$.
So to prove the proposition, we only need to prove the last
statement there.

We consider the action of
$\gamma=\epsilon_1\circ \phi_1 +\epsilon_p\circ \phi_p$
on $f$ as a classical modular form.
First
$\epsilon_1\circ \phi_1=[U\sig I U]_1$ acts on
$\CM(N,\psi)_M\subseteq \CM_{M}\otimes \CV$ by
$I\otimes 1$, where $M\geq 3$ is chosen such that $N|M$ and
$S_M=S_N$.
The action corresponding to $I$ on the classical
forms is $I$.
Since $f$ has rational coefficients, its image
in $\oplus_{t\in (\ZZ/M\ZZ)^\times} M(M,\psi)$ is of the
form $(f_0 (z))_t$ with $f_0$ independent of $t$. Then,
in the notations of Lemma~\ref{lem:adan} and
Proposition~\ref{prop:mot-ad}, for each
$t'\in (\ZZ/pM\ZZ)^\times$ we have
$\kappa(t')= t' \mod M$
and  $\tilde{\gamma}: E_{pM,t'}\to E_{M,t}$ defined by
$\gamma=\gamma_{t'}= I.$
So $(f_0(z))_t$ is sent to $(f_0(z))_{t'}$, $t'\in
(\ZZ/pM\ZZ)^\times$ by $\epsilon_1\circ \phi_1$.
Next $U\sig=U\sig \cap (\smat{1}{0}{0}{p} U
\smat{1}{0}{0}{p}^{-1}$. So
$\epsilon_p=p^{1-k}[U\sig \smat{1}{0}{0}{p} U]_1$ acts
on $\CM_{pM} \otimes \CV$ by $\smat{1}{0}{0}{p}\otimes 1$.
The corresponding action on $\oplus_{t} M_k(\Gamma_0(pN))$
is given by taking, for each $t'\in (\ZZ/p^2M\ZZ)^\times$,
$\kappa(t') = t' \mod pM$
and
$\tilde{\gamma}_{t'}: E_{p^2M,t'} \to E_{pM,t}$
given by
$\gamma_{t'}= \smat{1}{0}{0}{p^{-1}}.$
So $\epsilon_p$ sends $(f_0)_t\in \oplus_{t} M_k(\Gamma_0(pM))$ to
$(f_0(pz))_{t'}$. Since $\phi_p=-T_p$,
$\epsilon_p \circ \phi_p$ sends $(f_0(z))_t$ to
$(-f_0{} _{|T_p})_{t'}$. So by \cite[Lemma 4.6.5]{miyake}, $\gamma$
sends $(f_0(z))_t$ to $(f_0\sig(z))_{t'}$, as is needed.

{\bf Case 3: } When $\delta_p=2$, we have $p\nmid N$
 and
$$\gamma=\epsilon_1^{\Sigma_0,\Sigma_1}\circ
\phi_1^{\Sigma_0} +\epsilon_p^{\Sigma_0,\Sigma_1}\circ
\phi_p^{\Sigma_0} +\epsilon_{p^2}^{\Sigma_0,\Sigma_1}\circ
\phi_{p^2}^{\Sigma_0}.$$
When $q\neq p$, the relation
$T_q \gamma =\gamma T_q$ is proved in the same way as in Case 2.
When $q=p$, we only need to consider the local factor of
$T_p \gamma$ at the prime $p$.
The double coset product of the local factor is
\begin{align*}
& -[V_2 B V_2][V_2 I V_0][V_0 I V_0]+p^{1-k}
[V_2 B V_2][V_2 A_1 U_0][V_0 B V_0]\\
&-p^{k-2} [V_2 B V_2][V_2 A_2 V_0][V_0 B_2 V_0].
\end{align*}
By  (\ref{eq:AB0}) with $t=0$, the first term is
$-[V_2 B V_0]$. By  (\ref{eq:AB0}) with $t=1$ and  (\ref{eq:B0}),
the second term is
$-[V_2 B V_0]-[V_2 A V_0]$.
By  (\ref{eq:AB0}) with $t=2$, the third term is
$-[V_2 A V_0].$ Thus the equation is zero, as needed.
It follows that $T_p\circ \gamma =0$.

We now show that $\gamma$ sends $f$ to $f\sig$. This will complete
the proof of Proposition~\ref{prop:level}.
As in the case when $\delta_p=1$, we see that the action of
$\epsilon_1\circ \phi_1 +\epsilon_p\circ \phi_p$ on the classical
modular forms sends
$(f_0(z))_t\in \oplus_{t\in (\ZZ/M\ZZ)^\times} M_k(\Gamma_0(pM),\psi)$ to
$(f_0(z)-(f_0{}_{|T_p}(pz))_{t'}\in
\oplus_{t'\in (\ZZ/p^2M\ZZ)^\times} M_k(\Gamma_0(p^2M),\psi)$.
Further, the action of $\epsilon_{p^2}=p^{2(1-k)}[U\sig
\smat{1}{0}{0}{p^2} U]_1$ sends
$(f_0(z))_t\in \oplus_{t\in (\ZZ/M\ZZ)^\times}
M_k(M,\psi)$ to $(f_0(p^2z))_{t'}\in \oplus_{t'\in
(\ZZ/p^2M\ZZ)^\times} M_k(p^2M,\psi)$.
Thus $\epsilon_{p^2}\circ \phi_{p^2}=\psi(p_p)^{-1}p^{k-1}
\epsilon_{p^2}$ sends
$(f_0(z))_t\in \oplus_{t\in (\ZZ/M\ZZ)^\times}
M_k(M,\psi)$ to
$$(\psi(p_p)^{-1}p^{k-1} f_0(p^2z))_{t'}\in \oplus_{t'\in
(\ZZ/p^2M\ZZ)^\times} M_k(p^2M),\psi).$$
Therefore $\gamma$ send $(f_0(z))_t$ to
$$(f_0(z)-f_0{}_{|T_p}(pz)+\psi(p_p)^{-1}p^{k-1}
f_0(p^2z))_{t'}=(f_0\sig(z))_{t'}.$$
This is what we want.
\proofend

\begin{prop}
Let $\gamma^t$ denote the adjoint of $\gamma$. Then
$$\gamma^t\circ\gamma = \phi_{N\sig/N_f}\prod_{p\in\Sigma}
    L_p^\naive(A_f,1)^{-1}$$
on $M_f$.
\mlabel{pp:level2}
\end{prop}
\proofbegin
As in the proof of Proposition~\ref{prop:level},
we can again specialize to the case when $\Sigma$ has a single
prime $p$. Then we only need to consider the place at $p$.

When $\delta_p=0$, $\gamma$ is the identity map and the
formula is clear.

When $\delta_p=1$, we have
$$\gamma=\epsilon_1 \circ \phi_1
+\epsilon_p\circ \phi_p
=[\epsilon_1, \epsilon_p]
   \left [ \begin{array}{l}\phi_1\\ \phi_p
    \end{array} \right ].
$$
By
Lemma~\ref{lem:self-adj}, we find that the adjoints are given by
$$
 \begin{array}{l}
\epsilon_1^t = [U\sig I U]_1^t = [U \smat{p}{0}{0}{1} U\sig]_1,\\
\phi_1^t =[U I U]_1^t=[U I U]_1,\\
\epsilon_p^t=p^{1-k} [Uw^{-1}U]_{\omega^t\otimes \psi(\det w)} \circ
[U \smat{1}{0}{0}{p} U\sig]_1 \circ [U\sig w' U\sig]_{\omega'}
=p^{-1} [U I U\sig]_1 {\rm\ and\ }\\
\phi_p^t=-[U \smat{p_p}{0}{0}{1} U]_{\psi(p_p)^{-1}}^t
=-[U \smat{p_p}{0}{0}{1} U]_{\psi(p_p)^{-1}}.
\end{array} $$
Then by Proposition~\ref{prop:scomp} and Lemma~\ref{lem:pi-t},
\begin{align*}
\gamma^t \gamma & =
[\phi_1^t, \phi_p^t]
\left [\begin{array}{ll} \epsilon_1^t\epsilon_1 &
    \epsilon_1^t \epsilon_p\\
    \epsilon_p^t\epsilon_1 &
    \epsilon_p^t \epsilon_p\\
\end{array} \right ]
   \left [ \begin{array}{l}\phi_1\\ \phi_p
    \end{array} \right ] \\
& =
[1, -T_p]
\left [\begin{array}{cc} T_p & 1 \\
    1 & p^{-k} [U\smat{1}{0}{0}{p} U]_1
\end{array} \right ]
   \left [ \begin{array}{l} 1\\ -T_p
    \end{array} \right ] \\
& =\left [0, 1-\left\{ \begin{array}{ll}p^{-k}s_p, &
    \mbox{if $\pi_p$ is special}\\
    p^{-k} ps_p, & \mbox{if $\pi_p$ is principal series}
    \end{array} \right \} \right ]
   \left [ \begin{array}{l} 1\\ -T_p
    \end{array} \right ] \\
& =\left \{ \begin{array}{ll}
    (1-p^{-2})\phi_p, &\mbox{if $\pi_p$ is special,}\\
    (1-p^{-1})\phi_p, &\mbox{if $\pi_p$ is principal series.}
    \end{array}\right.
\end{align*}
By (\ref{eq:lp1}), this is what we need.

When $\delta_p=2$, we have $p\nmid N$ and $N\sig=p^2N$.
Also
$$\gamma=\epsilon_1 \circ \phi_1
+\epsilon_p\circ \phi_p+\epsilon_{p^2}\circ\phi_{p^2}
=[\epsilon_1, \epsilon_p, \epsilon_{p^2}]
   \left [ \begin{array}{l}\phi_1\\ \phi_p\\ \phi_{p^2}
    \end{array} \right ].
$$
By Lemma~\ref{lem:self-adj},
\begin{align*}
\epsilon_1^t &= [U\sig I U]_1^t = [Uw^{-1}U]_{\omega^t\otimes \psi(\det w)} \circ
[U I U\sig]_1 \circ [U\sig w' U\sig]_{\omega'}\\
&=\psi(p)^2[U \smat{p^2}{0}{0}{1} U\sig]_1,\\
\phi_1^t & =[U I U]_1^t=[U I U]_1,\\
\epsilon_p^t & =p^{1-k}[U\sig \smat{1}{0}{0}{p} U]_1^t
=p^{-1} \psi(p)[U \smat{p}{0}{0}{1} U\sig]_1, \\
\phi_p^t &=-[U \smat{p_p}{0}{0}{1} U]_{\psi(p_p)^{-1}}^t
=-[U \smat{p_p}{0}{0}{1} U]_{\psi(p_p)^{-1}},\\
\epsilon_{p^2}^t & =p^{2(1-k)}[U\sig \smat{1}{0}{0}{p^2} U]_1^t
=p^{2(1-k)} [U \smat{p^2}{0}{0}{p^2} U\sig]_1
=p^{-2}[U I U\sig]_1 {\rm\ and\ } \\
\phi_{p^2}^t & =p [U \smat{p_p}{0}{0}{p_p} U]_{\psi(p_p)^{-2}}^t
=p^{k-1} \psi(p_p)^{-1} [U I U]_1^t =\psi(p) p^{k-1}[U I U]_1.
\end{align*}
So by Proposition~\ref{prop:scomp},
{\scriptsize
\begin{align*}
&\gamma^t \gamma =
[\phi_1^t, \phi_p^t, \phi_{p^2}^t]
\left [\begin{array}{lll} \epsilon_1^t\epsilon_1 &
    \epsilon_1^t \epsilon_p & \epsilon_1^t \epsilon_{p^2} \\
    \epsilon_p^t\epsilon_1 &
    \epsilon_p^t \epsilon_p & \epsilon_p^t \epsilon_{p^2} \\
    \epsilon_{p^2}^t\epsilon_1 &
    \epsilon_{p^2}^t \epsilon_p & \epsilon_{p^2}^t \epsilon_{p^2}
\end{array} \right ]
   \left [ \begin{array}{l}\phi_1\\ \phi_p\\ \phi_{p^2}
    \end{array} \right ] \\
& =
[1, -T_p, \psi(p)p^{k-1}] \hspace{-.1cm}
\left [\hspace{-.15cm}\begin{array}{ccc} T_p^2-(p+1)\psi(p)p^{k-2}
    \hspace{-.8cm}&
    T_p & p^{-1}(p+1) \\
    T_p & p^{-1}(p+1) & \psi(p)^{-1} p^{-k} T_p \\
    p^{-1}(p+1) & (\psi(p) p^k)^{-1} T_p
 & \hspace{-.2cm}  (p^{2k}\psi(p)^2)^{-1}(T_p^2-(p+1)\psi(p) p^{k-2})
\end{array} \hspace{-.15cm}\right ] \hspace{-.15cm}
   \left [\hspace{-.15cm} \begin{array}{l} 1\\ -T_p\\ \psi(p)p^{k-1}
    \end{array} \hspace{-.2cm}\right ] \\
& =[0, 0, (1-p^{-1})((1+p^{-1})^2 - \psi(p)^{-1}p^{-k}T_p^2)]
   \left [ \begin{array}{l} 1\\ -T_p\\ \psi(p)p^{k-1}
    \end{array} \right ] \\
& =(1-p^{-1})((1+p^{-1})^2 - \psi(p)^{-1}p^{-k}T_p^2)
    \psi(p)p^{k-1}.
\end{align*}
}
So the action of $\gamma^t \gamma$ on $\CM(N,\psi)$ is
\begin{align*}
(1-p^{-1})&((1+p^{-1})^2  - \psi(p)^{-1}p^{-k}a_p^2)
    \psi(p)p^{k-1} \\
&=
(1-p^{-1})((1+p^{-1})^2 - \psi(p)^{-1}p^{-k}(\alpha_p+\beta_p)^2)
    \psi(p)p^{k-1} \\
&= (1-p^{-1})((1+p^{-1})^2 - \psi(p)^{-1}p^{-k}(\alpha_p^2
    +2\psi(p)p^{k-1}+\beta_p^2))
    \psi(p)p^{k-1} \\
&=
(1-p^{-1})(1+p^{-2} - \psi(p)^{-1}p^{-k}(\alpha_p^2+\beta_p^2))
    \psi(p)p^{k-1}.
\end{align*}
On the other hand,
\begin{align*}
\phi_{p^2} p[U\smat{p_p}{0}{0}{p_p} U]& (1-\alpha_p
\beta_p^{-1} p^{-1})(1-p^{-1})(1-\alpha_p^{-1}\beta_p p^{-1})\\
&=p^{k-1}\psi(p) (1-p^{-1})(1+p^{-2}-(\alpha_p^2+\beta_p^2)
/(\alpha_p\beta_p p)) \\
&=p^{k-1}\psi(p) (1-p^{-1})(1+p^{-2}-(\alpha_p^2+\beta_p^2)
\psi(p)^{-1}p^{-k}).
\end{align*}
This proves the proposition when $\delta_p=2$.
\proofend

\subsection{Integral structure for a $\Sigma$-variation}
\mlabel{ss:sig.int}
For a newform $f$ and a finite set of primes $\Sigma$
invertible in $\CO$, we set
$$\CM_f^\Sigma = M_f^\Sigma \cap \CM(N\sig,\psi^\Sigma)_!.$$
Note that $\fil^{k-1}\CM_{f,\dr}^\Sigma = \CO \cdot f^\Sigma$
by Proposition~\ref{prop:level}.
The image of $\wedge^2_\CO\CM_{f,B}^\Sigma$ defines an
integral structure for $M_\psi(1-k)$ of the form
$\eta_f^\Sigma\otimes_\CO \CM_{\psi}(1-k)$ for some
fractional
$\CO$-ideal $\eta_f^\Sigma \subset K$.
We call $\eta_f^\Sigma$
the (naive, $\Sigma$-finite) congruence $\CO$-ideal of $f$.
Note that $\eta_f^\Sigma$ is well-behaved under extension of
scalars $\CO' \otimes_\CO\cdot$ if $K\subset K'$ and
$S'$ is contained in the set of primes dividing those in $S$.

\begin{prop}
If $\CM_{f,\lambda}/\lambda$ is an irreducible
$(\CO/\lambda)[G_\QQ]$-module for every $\lambda$
in $S$, then the map $M_f \to M_f^\Sigma$ induces
an isomorphism $\CM_f \to \CM_f^\Sigma$ in $\ipms^S$.
\mlabel{prop:int-struc}
\end{prop}
\proofbegin
Fix a $\lambda\in S$.
We will identify realizations of $\CM_f$ with their images
in $\CM_f^\Sigma$ through the isomorphism $M_f\to M_f^\Sigma$.
Since $\CM_{f,\lambda}$ is compact and $\lambda^n
\CM_{f,\lambda}^\Sigma, n\in \ZZ$ is an open covering of
$M_{f,\lambda}^\Sigma$, there is $n=n_\lambda\in \ZZ$ such that
$\CM_{f,\lambda}\subset \lambda^n \CM_{f,\lambda}^\Sigma$ and
$\CM_{f, \lambda}\not\subset \lambda^{n+1}\CM_{f,\lambda}^\Sigma$.
So
$$\lambda^{n+1}\CM_{f,\lambda}^\Sigma \subset
\lambda^{n+1}\CM_{f,\lambda}^\Sigma+\CM_{f,\lambda} \subset \lambda^n
\CM_{f,\lambda}^\Sigma.$$
Since $\CM_{f,\lambda}^\Sigma/\lambda$ and hence
$\lambda^n \CM_{f,\lambda}^\Sigma/\lambda^{n+1} \CM_{f,\lambda}$
are irreducible,
one of the inclusions must be an equality.
The left inclusion being equal means $\CM_{f,\lambda}\subset
\lambda^{n+1} \CM_{f,\lambda}^\Sigma$, a contradiction. So we must
have $\lambda^{n+1}\CM_{f,\lambda}^\Sigma+\CM_{f,\lambda} \subset \lambda^n
\CM_{f,\lambda}^\Sigma.$ Then by Nakayama Lemma, we have
$\CM_{f,\lambda}=\lambda^n \CM_{f,\lambda}^\Sigma.$

Define $\gb=\prod_{\lambda\in S} \lambda^{n_\lambda}$.
Then $\gb$ is a fractional ideal of $\CO$.
Then $\CO_\lambda \otimes (\gb \otimes_\CO \CM_{f,B}^\Sigma)
\cong \CM_{f,\lambda}^\Sigma.$
Then by the comparison isomorphisms between $\CM_{f,B}^\Sigma$ and
$\CM_{f,\lambda}^\Sigma$ for $\lambda\in S$, we have
$$\CO_\lambda \otimes (\gb \otimes_\CO \CM_{f,B}^\Sigma)
    \cong \CO_\lambda \otimes \CM_B', \lambda\in S.$$
Therefore $\CM_{f,B}\cong \gb \otimes \CM_{f,B}^\Sigma.$
Similarly, using comparison isomorphisms between
$\CM\sig_{f,\lambda}$ and $\CM\sig_{f,\lcrys}$, and between
$\CM_{f,\dr}^\Sigma$
and $\CM_{f,\ecrys}^\Sigma$, we see that
$\CM_{f,\dr}\cong \gb \otimes \CM_{f,\dr}^\Sigma$
in $\ipms^{S,S}$.
In particular,
$\fil^{k-1} \CM_{f,\dr} \cong \gb \otimes
    \fil^{k-1} \CM_{f,\dr}^\Sigma.$
Since we have seen above that the image of
$\fil^{k-1}\CM_{f,\dr}$ under the isomorphism
$M_f \to M_f^\Sigma$ is $\CO \cdot f^\Sigma$ which is also
$\fil^{k-1}\CM_{f,\dr}^\Sigma$, we conclude that $\gb=\CO$.
\proofend

Combining this with Proposition \ref{prop:level},
we obtain:
\begin{coro}  For $f$ and $S$ as in proposition
\ref{prop:int-struc}, we have
$$\eta_f^\Sigma = \eta_f^\emptyset\prod_{p\in\Sigma}
    L_p^\naive(A_f,1)^{-1}.$$
\mlabel{coro:level}
\end{coro}
\proofbegin
We only need to consider the case when $\Sigma=\{p\}$.
Consider the following commutative diagram
$$
\begin{CD}
\CM_{f,B} @. \otimes_\CO @. \CM_{f,B}  @> \langle\ ,\ \rangle >>
    M_{\psi} (1-k)_B \\
@A \gamma^t AA @. @V \gamma VV  @VVV\\
\CM\sig_{f,B} @. \otimes_\CO @. \CM\sig_{f,B}
    @> \langle\ ,\ \rangle\sig >>
    M_{\psi} (1-k)_B.
\end{CD}
$$
Since $\gamma: \CM_{f,B} \to \CM\sig_{f,B}$ is an isomorphism,
we have
$$\langle \gamma \CM_{f,B}, \gamma \CM_{f,B} \rangle\sig
=\langle \CM\sig_{f,B}, \CM\sig_{f,B} \rangle\sig
=\eta_f\sig \otimes_\CO \CM_\psi(1-k)_B.$$
On the other hand, by Proposition~\ref{pp:level2},
$$\langle \gamma \CM_{f,B}, \gamma \CM_{f,B} \rangle\sig
=\langle \gamma^t \circ \gamma \CM_{f,B}, \CM_{f,B} \rangle
=\langle \phi_{p^{\delta_p}} \CM_{f,B}, \CM_{f,B} \rangle
    L_p^\naive(A_f,1)^{-1}.$$
When $\delta_p=0$, $\phi_{p^{\delta_p}}=1$ and we have
\begin{align*}
\eta_f\sig \otimes_\CO \CM_\psi(1-k)_B
&    =\langle \CM_{f,B}, \CM_{f,B} \rangle
    L_p^\naive(A_f,1)^{-1} \\
&    = \eta_f^\phi L_p^\naive(A_f,1)^{-1} \otimes_\CO \CM_\psi(1-k)_B.
\end{align*}
When $\delta_p=1$, $\phi_{p^{\delta_p}}=-T_p$ and we have
\begin{align*} \eta_f\sig \otimes_\CO \CM_\psi(1-k)_B
    &=\psi(p)^{-1} a_p \langle \CM_{f,B}, \CM_{f,B} \rangle
    L_p^\naive(A_f,1)^{-1} \\
    &= \eta_f^\phi L_p^\naive(A_f,1)^{-1}
    \otimes_\CO \CM_\psi(1-k)_B.
\end{align*}
When $\delta_p=2$, $\phi_{p^{\delta_p}}=pS_p$ and we have
\begin{align*} \eta_f\sig \otimes_\CO \CM_\psi(1-k)_B
&    =\psi(p) p^{k-1} \langle \CM_{f,B}, \CM_{f,B} \rangle
    L_p^\naive(A_f,1)^{-1} \\
&    = \psi(p) p^{k-1} \eta_f^\phi L_p^\naive(A_f,1)^{-1}
    \otimes_\CO \CM_\psi(1-k)_B.
\end{align*}
This proves the corollary since $p$ is invertible in $\CO$.
\proofend

\subsection{Exceptional primes}
\mlabel{sec:ep}
Suppose that $f$ is a newform of weight $k$ and character $\psi$
with coefficients in $K$, and that $\psi'$ is a
Dirichlet character with values in $K$.  Let $g$
denote the newform associated to $f\otimes\psi'$, so
$f\otimes\psi' = g^\Sigma$ where $\Sigma$ is the set of primes
dividing the conductor of $\psi'$.  Combining Lemma \ref{lem:twist}
and Proposition \ref{prop:level}, we have
that $M_g \cong M_f \otimes M_{\psi'}$, so $A_f \cong A_g$
in $\pms^S_K$ for any $S\supseteq S_f$ such that Cond$\psi'$
is invertible in $\CO_S$.  So when studying $A_f$, we can assume
$f$ has minimal conductor among its twists.

Let $F$ denote the quadratic unramified extension of $\QQ_p$
in $\Qpbar$, and let $H_p$ denote the subgroup of $G_p$ fixing $F$.
\begin{lemma}
Let $f$ have minimal conductor among its twists.
If $\pi_p(f)$ is supercuspidal, then the following
are equivalent:
\begin{enumerate}
\item  $p$ is in $\Sigma_e$;
\item  $L_p(A_f, s) = (1 + p^{-s})^{-1}$;
\item  $\bar{K}_\lambda \otimes_{K_\lambda} M_{f,\lambda}|_{I_p}$
    is reducible for some (hence all) $\lambda \in S$;
\item  $\bar{K}_\lambda \otimes_{K_\lambda} M_{f,\lambda}|_{G_p}$
    is induced from a character of $H_p$ for some (hence all)
    $\lambda \in S$;
\end{enumerate}
\mlabel{lem:ep2}
\end{lemma}
\proofbegin
(a) $\Rightarrow$ (c):
If $p$ is in $\Sigma_e$ then $L_p(A_f,s)\neq 1$.
This means that $A_{f,\lambda}^{I_p}$ is not zero
which means that $\Hom(M_{f,\lambda},M_{f,\lambda})^{I_p}\supsetneq
K_\lambda$. By Schur's lemma, this shows that, as a representation
of $I_p$,  $\bar{K}_\lambda \otimes_{K_\lambda} M_{f,\lambda}$ is reducible.

(c) $\Rightarrow$ (d):
Assume that $\bar{K}_\lambda\otimes_{K_\lambda}
M_{f,\lambda}|_{I_p}$ is reducible. Since the image of $I_p$ is finite,
$\bar{K}_\lambda\otimes_{K_\lambda} M_{f,\lambda}|_{I_p}$
decomposes into $V_1\oplus V_2$ of
representations of $I_p$. Let $\Frob_p$ be a lift of the Frobenius
in $G_p/I_p$. Since $M_\lambda|_{G_p}$ is irreducible, the two
terms $V_1$ and $V_2$ are exchanged by $\Frob_p$ and preserved by
$\Frob_p^2$. So $V_1$ and $V_2$ are preserved by the subgroup
of $G_p$ generated by $I_p$ and $\Frob_p^2$, which is $H_p$.
Further $M_{f,\lambda}$ is the representation of $G_p$ induced
by the representation $V_1$ of $H_p$.

(d) $\Rightarrow$ (b):
If $\bar{K}_\lambda \otimes_{K_\lambda} M_{f,\lambda}|_{G_p}$
    is induced from a character of $H_p$ for some $\lambda \in S$.
Then the action of $D_p$ on $\bar{K}_\lambda \otimes_{K_\lambda}
    A_\lambda$ has the matrix
$\left ( \begin{array}{ccc} \vep & & \\ & \chi & \\ & & \vep^{-1}
\end{array} \right )$ where $\vep$ is a ramified character and
$\chi$ is an unramified quadratic character.
Thus
$A_\lambda^{I_p}=K_\lambda(\chi)$ and $\chi(\Frob_p)=-1$.
This shows that $L_p(A,s)=(1+p^{-s})^{-1}$.

(b) $\Rightarrow$ (a): Clear.
\proofend

Suppose now that $p$ is in $\Sigma_e$ and $\tau$ is
an embedding of $K$ in $\CC$.  Note that $c_p$
is even by the lemma (see VI.2 of \cite{serre_corloc}).
Moreover the explicit description of $\pi_p(\tau(f))$ in this
case is given in Section 3 of \cite{gerardin}.  In particular,
we have that
$$\pi_p(\tau(f))^{V_p'}|_{\GL_2(\ZZ_p)} \cong \Theta(\vep)$$
for some character
$$\vep:\CO_F^\times \to (\CO_F/p^{c_p/2})^\times \to \CC^\times,$$
where $V_p'$ denotes the set of matrices in $\GL_2(\ZZ_p)$ congruent
to $I$ mod $p^{c_p/2}$, and $\Theta(\vep)$ is defined in
\S 3.2 of \cite{cdt}.  Our hypotheses
ensure that $\vep/\vep\circ\frob_p$ has conductor $(p\CO_F)^{c_p/2}$.
Let $g_p = \smat{1}{0}{0}{p^{c_p/2}}$.
We recall the following properties of $\Theta(\vep)$:
\begin{lemma}  Let $B$ denote the set of matrices
 $\smat{a}{b}{c}{d}$ in $\GL_2(\ZZ_p)$ with $c\equiv 0 \bmod p$.
 Let $D$ denote the set of matrices in $\GL_2(\ZZ_p)$ which are
 diagonal $\bmod p^{c_p/2}$, and $\CC_\psi = \CC$ with the action
 of $D$ defined by $\smat{a}{b}{c}{d} \mapsto \psi(a)$.
\begin{enumerate}
\item $\Theta(\vep)|_B$ is irreducible.
\item $\Theta(\vep^{-1}) \cong \hat{\Theta}(\vep)
 \cong \Theta(\vep)\otimes \vep^{-1}\circ\det$.
\item $\hom_{\CC[D]}(\CC_\psi,\Theta(\vep))$ is one-dimensional.
\item If $\CV$ is a model over $\CO$ for $\Theta(\vep)$ and
 $\lambda$ is a prime of $\CO$ not dividing $p-1$, then
 $\CV/\lambda\CV$ is absolutely irreducible.
\end{enumerate}
\mlabel{lem:Theta}
\end{lemma}
\proofbegin (a), (b) and (d) are part of Lemma 3.2.1 of \cite{cdt}.
(For (d), note that $p$ is already assumed to be invertible in $\CO$.)
To prove (c), we use that $\pi_p(\tau(f))^{V_p}$ is one-dimensional
with $\smat{a}{0}{0}{1}$ acting via $\psi(a)$ for $a\in\ZZ_p^\times$.
Conjugating by $g_p$ gives the result.
\proofend

Suppose now that $\Sigma_0$ is a subset of $\Sigma_e$.
For any prime $p$, we let $\CV_p^1 = \CO$ with
$\smat{a}{b}{c}{d} \in V_p$ acting via $\psi^{-1}(a)$.
\begin{lemma} There exists a field $K_0$ such that if $K_0 \subset K$,
then the following holds for all primes $p$ in $\Sigma_0$:
The representation $\Theta(\vep^{-1})$ has a locally free
model $\CV_p^0$ over $\CO$ such that
\begin{enumerate}
\item there is an isomorphism of $\CO[\GL_2(\ZZ_p)]$-modules
$$\omega_p: \CV_p^0 \cong \hom_\CO(\CV_p^0,\CO(\psi^{-1}\circ\det)).$$
\item there is a surjection of $\CO[V_p]$-modules
$$\tau_p:\res_{g_p\GL_2(\ZZ_p)g_p^{-1}}^{V_p} g_p\CV_p^0 \to \CV_p^1.$$
\end{enumerate}
\mlabel{lem:p-scalars}
\end{lemma}
\proofbegin It suffices to prove the lemma for $\Sigma_0 = \{p\}$.
We first assume $K$ is sufficiently large that $\Theta(\vep^{-1})$
is defined over $K$.  Then it has a model $\CV'$ over $\CO$
and by part (b) of Lemma \ref{lem:Theta}, there are non-zero homomorphisms
$$\omega:\CV \to \hom_\CO(\CV,\CO(\psi^{-1}\circ\det))\quad\mbox{and}\quad
\tau:\res_{g_p\GL_2(\ZZ_p)g_p^{-1}}^{V_p} g_p\CV \to \CV_p^1.$$
These are isomorphisms after tensoring with $\CO_\lambda$
for all but finitely many $\lambda$ in $S$.  Using Lemma 3.3.1
of \cite{cdt}, we may enlarge $K$ so that for the remaining
$\lambda$, there are $\GL_2(\ZZ_p]$-stable lattices $\CV'_\lambda$
in $\CV\otimes_\CO\CO_\lambda$ and isomorphisms
$$\omega'_\lambda:\CV'_\lambda \to \hom_\CO(\CV'_\lambda
,\CO(\psi^{-1}\circ\det)).$$
Now replace $\CV$ with $\CV'$
so that $\CV'\otimes_\CO\CO_\lambda = \CV'\lambda$ for these $\lambda$
and $\CV\otimes_\CO\CO_\lambda$ otherwise.
The image of $\omega$ and $\tau$ are now of the form
$$\ga \hom_\CO(\CV,\CO(\psi^{-1}\circ\det))\quad\mbox{and\ }
\gb\CV_p^1$$
for some non-zero ideals $\ga$ and $\gb$ of $\CO$.
To complete the proof, we further enlarge $K$ so that
these become principal.
\proofend

Now let $\Sigma$ be any finite set of primes. Define
$$N_{\Sigma_0}^\Sigma = \prod_{p\not\in \Sigma_0\cup \Sigma} p^{c_p}
  \prod_{p\in\Sigma} p^{c_p + \delta_p},$$
and let $U^{\Sigma}_{\Sigma_0} = U_0(N^\Sigma_{\Sigma_0})$.
Write $U^\Sigma_{\Sigma_0} = \prod_p U_p$ and consider the
representation $\sigma^\Sigma_{\Sigma_0}$ defined by
$$\CV^\Sigma_{\Sigma_0} = \otimes_\CO \CV_{p,\Sigma_0}^\Sigma,$$
where $\CV^\Sigma_{p,\Sigma_0} = \CV_p^0$ if $p \in \Sigma_0 - \Sigma$
and $\CV_p^1$ (restricted to $U_p$) otherwise.  Here the tensor product
is over all primes $p$, but the action is trivial for all but finitely
many $p$.  Note that if $\Sigma_0 \subset \Sigma$, then
$U^\Sigma_{\Sigma_0} = U^\Sigma$ and
$\sigma^\Sigma_{\Sigma_0} = \sigma(N,\psi)$.

Suppose now that $\Sigma \subset \Sigma'$.  For positive
integers $m$ dividing $N_{\Sigma'}/N_\Sigma$, define
$$\epsilon_{\Sigma_0,m}^{\Sigma,\Sigma'} =
   m^{1-k}[U^{\Sigma'}\smat{1}{0}{0}{m'}U^\Sigma]_\tau:
   \CM(\sigma^\Sigma_{\Sigma_0})_! \to
   \CM(\sigma^{\Sigma'}_{\Sigma_0})_!,$$
where $m'=\prod_p m'_p$ with
$m'_p = p^{c_p/2}$ if $p\in\Sigma_0\cap\Sigma'-\Sigma$
and $m'_p = p^{v_p(m)}$ otherwise, and $\tau = \otimes \tau_p$
where the tensor product is over primes $p\in\Sigma_0\cap\Sigma'-\Sigma$.
We now define
$$\gamma_{\Sigma_0}^{\Sigma,\Sigma'} =
\sum_{m}\epsilon_{\Sigma_0,m}^{\Sigma,\Sigma'}\phi_m.$$
Note that if $\Sigma_0 \subset \Sigma$, then
$\epsilon_{\Sigma_0,m}^{\Sigma,\Sigma'}$ is $\epsilon_m$ and
$\gamma_{\Sigma_0}^{\Sigma,\Sigma'}$ is the same
$\gamma$ as in Proposition \ref{prop:level}.

For each prime $p$, we define the Hecke operator $T_p$
on $M_!(\sigma^\Sigma_{\Sigma_0})$ (in what category)
by the double coset operator
$$\left[U_{\Sigma_0}^\Sigma\smat{p}{0}{0}{1}_p
U_{\Sigma_0}^\Sigma\right]_{\psi^{-1}(p_p)}$$
if $p \not\in\Sigma_0$ and we set $T_p = 0$ if $p \in \Sigma_0$.
For primes $p$ not dividing $N^\Sigma$, we define $S_p$ by
$$\left[U_{\Sigma_0}^\Sigma\smat{p}{0}{0}{p}_p
U_{\Sigma_0}^\Sigma\right]_{\psi^{-2}(p_p)}.$$
Note that if $\Sigma \subset \Sigma$, then this coincides with
the previous definition.

\begin{lemma}
For any $p$ in $\Sigma_0-\Sigma$, the
morphism $\epsilon = \epsilon_{\Sigma_0,p}^{\Sigma,\Sigma\cup\{p\}}$
is injective and $T_p = 0$ on its image.  Furthermore,
its cokernel is torsion-free if $p-1$ is invertible in $\CO$.
\mlabel{lem:p-inj}
\end{lemma}
\proofbegin
Let $\CM = \CM_{Np^{c_p/2}}$, $\CM' = \CM_{Np^{c_p}}$,
$U = U_{\Sigma_0}^\Sigma$, $U' = U_{\Sigma_0}^{\Sigma\cup\{p\}}$,
$\CV = \CV_{\Sigma_0}^\Sigma$,
$\CV' = \CV_{\Sigma_0}^{\Sigma\cup\{p\}}$.

We first prove the assertions concerning the kernel and
cokernel.  Note that $\epsilon$ is the map
$$g_p \otimes \tau_p : (\CM_! \otimes_\CO \CV)^U
 \to  (\CM'_! \otimes_\CO \CV')^{U'},$$
and the latter object can be identified with
$$\Ind_{U'}^U(\CM'_! \otimes_\CO \CV')^U$$
where the inclusion $U' \to U$ is defined by
$u' \mapsto g_p^{-1}u'g_p.$
It therefore suffices to prove the map
$$ \begin{array}{ccc}\CM_! \otimes_\CO \CV &\to &
\Ind_{U'}^U(\CM'_! \otimes_\CO \CV')\\
  x \otimes y & \mapsto & (u \mapsto g_pux\otimes \tau_puy)
\end{array}$$
is injective, with torsion-free cokernel if $p-1 \in \CO^\times$.
But this map can be written as a composite
$$\CM_! \otimes_\CO \CV \to
\CM_! \otimes_\CO \Ind_{U'}^U\CV' \cong
\Ind_{U'}^U(\res_U^{U'} (\CM_! \otimes_\CO \CV'))\to
\Ind_{U'}^U(\CM'_! \otimes_\CO \CV'),$$
where the first map is defined by sending
$y$ to $(u\mapsto \tau_p u y)$ for $y\in\CV$,
and the last map is the induction of $g_p \otimes 1$.
Since $\CM_!$, $\CM'_!$, $\CV$ and $\CV'$ are torsion-free, it
suffices to prove that
$$\CV \to \Ind_{U'}^U\CV' \quad\mbox{and}\quad
  \CM \to \CM'$$
are injective, with torsion-free cokernel if $p-1 \in\CO^\times$.

For the first map, the injectivity follows from the fact
that the map is non-zero and $\CV\otimes_\CO K$ is irreducible.
The torsion-freeness then amounts to the injectivity after
tensoring with $\CO/\lambda$ for all $\lambda$ in $S$.
The map is non-zero since $\tau_p\otimes_\CO\CO/\lambda$
is surjective and factors through it, and so the injectivity
follows from Lemma \ref{lem:Theta} (d).
For the second map, the injectivity and torsion-freeness
follows for example by considering the composite
$$g_p^{-1}\circ g_p : \CM \to \CM' \to \CM_{Np^{3c_p/2},!},$$
which we have already seen is injective with torsion-free
cokernel, by Lemma \ref{lem:gmap} (c).

To complete the proof of the lemma, it suffices to prove $T_p = 0$
on the image of
$$((\fil^{k-1}\CM_\dr \otimes \CC) \otimes \CV)^U
\to  ((\fil^{k-1}\CM'_\dr \otimes \CC) \otimes \CV')^{U'},$$
or equivalently, on the image of
$$g_p\otimes \tau_p: (\pi(g) \otimes \CV)^U \to (\pi(g) \otimes \CV')^{U'}$$
for each newform $g$ of weight $k$.
But by the definition of $\Theta(\vep)$,
$(\pi(g)\otimes\CV)^U = 0$ unless $\pi_p(g)$ is supercuspidal
of conductor exactly divisible by $p^{c_p}$, so the result
follows from Lemma~\ref{lem:pi-t}.
\proofend

It follows from the lemma and the definitions that
$\gamma_{\Sigma_0}^{\Sigma,\Sigma\cup\Sigma_0}$ is injective
and commutes with the action of the Hecke operators $T_p$
for all primes $p$.  In particular, for $p \in S$, the
operator $T_p$ can be written as a polynomial in $T_q$
for $q\neq p$, and the same arguments as in Proposition~\ref{prop:T}
show that $\TT$ acts on $\CM(\sigma_{\Sigma_0}^\Sigma)_!$
as an object of $\ipms^S_K$.  Now define
$\CM_{f,\Sigma_0}^\Sigma$ and $M_{f,\Sigma_0}^\Sigma$
exactly as before, i.e., as the intersection of the
kernels of elements of $I_f$.  Note that if
$\Sigma_0 \subset \Sigma$, then
$M_{f,\Sigma_0}^\Sigma = M_f^\Sigma$.

\begin{lemma}  The operator
$\gamma = \gamma_{\Sigma_0}^{\emptyset,\Sigma}$ restricts to an isomorphism
$$M_{f,\Sigma_0}^{\emptyset} \to M_{f,\Sigma_0}^{\Sigma}.$$
Furthermore if $S$ is as in Proposition \ref{prop:int-struc}
and $p-1 \in \CO^\times$ for all $p \in \Sigma_0$, then it induces
an isomorphism
$$\CM_{f,\Sigma_0}^\emptyset \to \CM_{f,\Sigma_0}^{\Sigma}.$$
\end{lemma}
\proofbegin By Proposition \ref{prop:int-struc},
Lemma \ref{lem:p-inj} and the commutativity of the diagram
$$\begin{array}{ccc}
\CM_{f,\Sigma_0}^{\emptyset}& \to &\CM_{f,\Sigma_0}^{\Sigma}\\
\downarrow && \downarrow \\
\CM_{f,\Sigma_0}^{\Sigma_0}  &\to& \CM_{f,\Sigma_0}^{\Sigma_0\cup\Sigma},
\end{array}$$
it suffices to prove that the lemma in the case $\Sigma = \Sigma_0$.
Furthermore, by lemma \ref{lem:p-inj}, it suffices to prove that
$$\fil^{k-1}M_{f,\Sigma_0,\dr}^{\emptyset} \neq 0,$$
and this follows from the definition of $\sigma_{\Sigma_0}^\Sigma$
in terms of $\pi_p(f)$.
\proofend

Define $w = w^\Sigma_{\Sigma_0} = \smat{0}{-1}{N}{0}_N \in \GL_2(\AA_\f)$
where $N = N^\Sigma_{\Sigma_0}$.  Let $\sigma = \sigma_{\Sigma_0}^\Sigma$
and $\CV = \CV_{\Sigma_0}^\Sigma = \otimes\CV_{p,\Sigma_0}^\Sigma$, and
denote by $\sigma'$ the representation of $U = U_{\Sigma_0}^\Sigma$ defined
by
$$\CV' = \hom_\CO(\CV,\CO(\psi^{-1}\circ\det))    \cong
\otimes\hom_\CO(\CV_{p,\Sigma_0}^\Sigma,\CO(\psi^{-1}\circ\det)).$$
Define $\omega:\CV \to \CV'$ as the tensor product of the maps $\omega_p$,
where $\omega_p$ is defined in Lemma \ref{lem:p-scalars} if
$p\in \Sigma-\Sigma$, and by sending a generator $v_0$ to
the map $v_0 \mapsto 1$ otherwise.  We then have that
$$\omega(\sigma(w^{-1}uw)v) = \sigma'(u) \omega(v)$$
for all $u \in U$ and $v \in \CV$, so the operator
$[UwU]_\omega$ is well-defined and induces an isomorphism
$$\CM(\sigma)_! \to \CM(\sigma')_!.$$
Composing with the isomorphism
$$M(\sigma')_! \to \hom_K(M(\sigma)_!,M_\psi(1-k)),$$
we obtain a perfect pairing
$$M(\sigma)_! \otimes_K M(\sigma)_! \to M_\psi(1-k).$$

\begin{lemma} Suppose that $p\in \Sigma_0 - \Sigma$.
 Let $\epsilon = \epsilon_{\Sigma_0,p}^{\Sigma,\Sigma\cup\{p\}}$ and
 $\epsilon^t$ its transpose with respect to the above pairings.
 Then $\epsilon^t\epsilon$ is a scalar divisible by $p+1$.
\mlabel{lem:ete}
\end{lemma}
\proofbegin
Let $U = U_{\Sigma_0}^\Sigma$, $U' = U_0(pN_{\Sigma_0}^\Sigma)$ and
$U'' = U_{\Sigma_0}^{\Sigma\cup\{p\}}$, so
$U'' \subset g_pU'g_p^{-1}$ and $[U:U'] = p+1$.  Let
$\CV = \CV_{\Sigma_0}^\Sigma$,
$\sigma = \sigma_{\Sigma_0}^\Sigma$, $\sigma' = \sigma|_{U'}$
and $\sigma'' = \sigma_{\Sigma_0}^{\Sigma\cup\{p\}}$.
Define a pairing
$$M(\sigma')_! \otimes_K M(\sigma')_! \to M_\psi(1-k)$$
exactly as we did for $M(\sigma)_!$, using the same $w$
and $\omega$.  Letting
$$\epsilon_1 = \left[U' 1 U\right ]_1;\quad
\epsilon_2 = \left[U'' g_p U'\right ]_{\tau_p},$$
we find that $\epsilon = \epsilon_2\circ\epsilon_1$,
$\epsilon_1^{t}\epsilon_1 = p+1$ and
$\epsilon_2^{t}\epsilon_2 = [U'1U']_\delta$ for
some $U'$-linear $\delta: \CV \to \CV$.
By Lemma \ref{lem:Theta} (a), $\sigma'$ is absolutely irreducible,
so $\delta$ is a scalar in $\CO$ and
$\epsilon^t\epsilon = (p+1)\delta$.
\proofend

The same proof as that of Proposition \ref{pp:level2} now gives:
\begin{proposition}  Suppose that $\Sigma_0 \subset \Sigma_e$
and $p+1 \in \CO^\times$ for all $p \in \Sigma_e - \Sigma_0$.
Suppose that $\Sigma$ is any finite set of primes invertible in $\CO$,
and let $\gamma = \gamma_{\Sigma_0}^{\emptyset,\Sigma}$.  Then
$$\gamma^t\circ\gamma = \beta_f^\Sigma\prod_{p\in\Sigma} L_p(A_f,1)^{-1}$$
on $M_{f,\Sigma_0}^\emptyset$ for some non-zero $\beta_f^\Sigma$ in $\CO$.
\mlabel{prop:exlevel}
\end{proposition}

Note also that together with the injectivity of
$\gamma_{\Sigma_0}^{\Sigma,\Sigma\cup\Sigma_0}$, Lemma \ref{lem:ete}
implies that the pairing on $M_!(\sigma_{\Sigma_0}^\Sigma)$
is alternating and induces an isomorphism
$$\wedge^2_K M_{f,\Sigma_0}^\Sigma \to M_\psi(1-k).$$
We let $\eta_{f,\Sigma_0}^\Sigma$ be the fractional ideal in $K$
such that $\wedge^2_K \CM_{f,\Sigma_0}^\Sigma$ maps
isomorphically to  $\eta_{f,\Sigma_0}^\Sigma\CM_\psi(1-k)$.
Note that if $\Sigma_0 \subset \Sigma$, then
$\eta_{f,\Sigma_0}^\Sigma = \eta_f^\Sigma$.
The same proof as that of Corollary \ref{coro:level} now gives:

\begin{coro}  Suppose that $f$ and $S$ are as in Proposition
\ref{prop:int-struc} and $\Sigma_0 \subset \Sigma_e$
is such that $p-1 \in \CO^\times$ for all $p \in \Sigma_0$,
and $p+1 \in \CO^\times$ for all $p \in \Sigma_e - \Sigma_0$. Then
$$\eta_{f,\Sigma_0}^\Sigma \subset \eta_f^\emptyset\prod_{p\in\Sigma}
    L_p(A_f,1)^{-1}.$$
\end{coro}

\section{The Taylor-Wiles construction}
\mlabel{sec:twc}
Our method for computing the Selmer group of $A_f(1)_\lambda$
is based on that of Wiles \cite{wiles} and his work
with Taylor \cite{tw}.  We begin by recalling some definitions
and results from Galois cohomology.
Then we give an axiomatic formulation of the method of
\cite{wiles} and \cite{tw}, made possible by the simplifications
due to Faltings (\cite{tw}, appendix), Lenstra \cite{lenstra}
and one of the authors \cite{fd_twc}.  Finally we verify these
axioms are satisfied in the context of modular forms of
higher weight.

\subsection{Galois cohomology}
\mlabel{ssec:gal-coh}
In this section, we fix a prime $\ell$ and suppose that
$K$ is a finite extension of $\QQ_\ell$ with ring of
integers $\CO$, uniformizer $\lambda$ and residue field $\kappa$.
We let $\Sigma_0$ denote the set of primes different
from $\ell$.  We let $\chi_\ell:G_{\QQ} \to \ZZ_\ell^\times$
denote the $\ell$-adic cyclotomic character.

Suppose that for $i = 1,2$,
$$\rho_i: G_\QQ \to \aut_K V_i$$
is a continuous geometric representation whose
restriction to $G_\ell$ is crystalline and short,
where {\em short} means that $\fil^{\ell}D_\pst(V_i)=0$
and $\fil^0 D_\pst(V_i)=D_\pst(V_i)$ and that $V_i$ has
no nonzero quotient
$V'$ so that $V'(\ell-1)$ is unramified.
We let $V = \hom_K(V_1,V_2)$.
Suppose that $L_i$ is a $G_\QQ$-stable $\CO$-lattice in $V_i$,
and let
$$A = (K/\CO)\otimes_\CO \hom_\CO(L_1,L_2).$$
Below we shall define an $\CO$-submodule $H^1_w(G_p, A) \subset
H^1(G_p,A)$ for each prime $p$.  For each finite set of
primes $\Sigma \subset \Sigma_0$, we define
$$H^1_\Sigma(G_\QQ, A) \subset H^1(G_\QQ, A)$$
as the $\CO$-submodule of elements with restrictions
in $H^1_w(G_p,A)$ for every $p\in\Sigma$.

First we define $H^1_w(G_p, A_n)$
for each prime $p$ and integer $n \ge 0$,
where $A_n = A[\lambda^{n+1}]$.
If $p \neq\ell$, then we simply set
$$H^1_w(G_p,A_n) = H^1(G_{\FF_p},A_n^{I_p})
        = \ker (H^1 (G_p,A_n) \to H^1(I_p, A_n)).$$
To treat the case $p=\ell$, we let $\flcat$ denote the
full subcategory of finite $\ZZ_\ell G_\ell$-modules
whose objects are quotients of $G_\ell$-stable lattices
in short crystalline representations.
Thus $\flcat$ is
stable under taking direct sums, subobjects and quotients,
and the $L_i/\lambda^{n+1}L_i$ are objects of it.
Recall that there is a natural $\CO$-linear isomorphism between
$H^1(G_\ell,A_n)$ and the $\CO$-module of extensions of
$(\CO/\lambda^n)[G_\ell]$-modules
$$ 0 \to \lambda^{-n-1}L_2/L_2  \to E \to  L_1/\lambda^{n+1} L_1 \to 0.$$
We define $H^1_w(G_\ell,A_n)$ to be the set of elements of
$H^1(G_\ell,A_n)$ corresponding to extensions in $\flcat$.
Using the stability of $\flcat$ under direct sums, subobjects and
quotients, one checks that $H^1_w(G_\ell,A_n)$ is an $\CO$-submodule,
and furthermore that it is the preimage of $H^1_w(G_\ell,A_{n+1})$
under the natural map.  Now define
$H^1_w(G_p,A) = \dirlim H^1_w(G_p,A_n)$ for every prime $p$.

Recall that Bloch and Kato define a divisible $\CO$-submodule
$H^1_\f(G_p,A) \subset H^1(G_p,A)$ for each $p$ (see \cite{bloch_kato}
or \cite{fon_pr}). Their {\em Selmer group}
$$H^1_\f(G_\QQ, A) \subset H^1(G_\QQ, A)$$
is then the $\CO$-submodule of elements with restrictions
in $H^1_\f(G_p,A)$ for all primes $p$.  If $p\neq\ell$,
it is clear that
$$H^1_\f(G_p,A) \subset H^1_w(G_p,A)$$
and equality holds if $A^{I_p}$ is divisible.
\begin{proposition}  $H^1_w(G_\ell,A) = H^1_\f(G_\ell,A)$
is divisible of $\CO$-corank
$$d = \dim_K H^0(G_\ell,V) + \dim_K V - \dim_K \fil^0 D_\ecrys(V).$$
\end{proposition}
\proofbegin
The divisibility of $H^1_\f(G_\ell,A)$ follows from its definition,
and the corank is computed in \cite{bloch_kato} or \cite{fon_pr}.
One shows as in \cite{wiles} that
$$H^1_\f(G_\ell,A) \subset H^1_w(G_\ell,A),$$
so it suffices to prove that
$$ \dim_\kappa H^1_w(G_\ell,A)[\lambda] \le d.$$
Since the kernel of $H^1(G_\ell,A_0) \to
H^1(G_\ell,A)$ is contained in $H^1_w(G_\ell,A_0)$,
we have
$$\begin{array}{rl} \dim_\kappa H^1_w(G_\ell,A)[\lambda] &=
        \dim_\kappa H^1_w(G_\ell,A_0) -
        \dim_\kappa(\kappa\otimes_\CO A^{G_\ell})\\&=
        \dim_\kappa H^1_w(G_\ell,A_0) -
        \dim_\kappa A_0^{G_\ell} + \dim_K H^0(G_\ell,V).\end{array}$$
Therefore we need to prove that
$$\dim_\kappa H^1_w(G_\ell,A_0) - \dim_\kappa H^0(G_\ell,A_0)
        \le  \dim_K V - \dim_K \fil^0 D_\ecrys(V).$$

We now sketch the proof that equality holds by
the theory of Fontaine and Laffaille \cite{fon_laf}.
Recall from \S \ref{ss:back} that the category $\mfcat_{\tor}^0$
is equivalent, via the functor
\[\cald\mapsto\VV(\cald):=\hom(\underline{U}_S(\cald),\QQ_\ell/\ZZ_\ell),\]
to $\flcat$.  Let $\mfcatkap$ denote the category of
$\kappa$-modules in $\mfcat$ so that the functor $\VV$ defines an equivalence between $\mfcatkap^0$ and the category of
$\kappa$-modules in $\flcat$.  Choosing objects
$\cald_1$ and $\cald_2$ in $\mfcatkap$
so that $\VV(\cald_1) \cong L_1/\lambda L_1$ and
$\VV(\cald_2) \cong \lambda^{-1}L_2/\lambda L_2$,
we have
$$\dim_\kappa \fil^j \cald_i = \dim_K \fil^j D_\cris(V_i)$$
for $i = 1,2$ and $j = 0,\ldots,\ell-1$.  The dimension
formula then follows from the exact sequence
$$
\begin{array}{rcl}
        0 \to \hom_{\mfcatkap}(\cald_1,\cald_2)
        \to& \hom_{\kappa,\fil}(\cald_1,\cald_2)&\\
        \to& \hom_\kappa(\cald_1,\cald_2)&
        \to \ext^1_{\mfcatkap}(\cald_1,\cald_2)\to 0\end{array}
$$
whose construction we leave as an exercise (see also diagram (\ref{big}) for a
similar computation).
\epf

Now we specialize to the case where
$\rho_1 = \rho_2$, each having dimension $2$
over $K$ with
$\dim_K\fil^1 D_\cris(V_i)$ $=1$
(so also $\ell > 2$).  We assume also that
$L_1 = L_2$,  and let
$$A^0 = (K/\CO)\otimes_\CO \ad^0_\CO L_1$$
where $\ad^0_\CO L_1$ denotes the kernel
of the trace map $\End_\CO L \to \CO$.
We let $A^0_n = A^0[\lambda^{n+1}]$,
identify $A^0$ with a submodule of $A$,
and define $H^1_w(G_p, A^0_n)$,
$H^1_w(G_p, A^0)$ and $H^1_w(G_\QQ,A^0)$
as the preimage of the corresponding
group for $A$.  From the proposition above, we get
\begin{coro} Under the above hypotheses,
$H^1_f(G_\ell,A^0) = H^1_w(G_\ell,A^0)$ and
$$\dim_\kappa H^1_w(G_\ell,A_0^0) = \dim_\kappa H^0(G_\ell,A_0^0) + 1.$$
Moreover if $H^0(I_p,A^0)$ is divisible for all $p \neq \ell$, then
$$H^1_\emptyset(G_\QQ,A^0) = H^1_\f(G_\QQ,A^0).$$
 \mlabel{coro:local}
\end{coro}

\subsection{An axiomatic formulation}
Fix a continuous, odd, irreducible representation
$$\rho_0 : G_\QQ \to \aut_{\kappa}(V_0)$$
where $V_0$ is two-dimensional over $\kappa$.
We impose the following two technical hypotheses
on the representations $\rho_0$ we consider:
\begin{itemize}
\item The restriction of $\rho_0$ to $G_F$ is absolutely irreducible
      where $F$ is the quadratic subfield of $\QQ(\mu_\ell)$.
\item The Serre weight $k$ of $\rho_0$ is at most $\ell - 1$.
\end{itemize}
We let $\psi$ denote the $\CO$-valued character of $G_\QQ$
having finite order not divisible by $\ell$ so that
$\psi\chi_\ell^{1-k}$ has reduction $\det\rho_0$.
Let $K_\psi$ denote $K$ with an action of $G_\QQ$ via $\psi$,
and define $\CO_\psi$ similarly.

We consider continuous pseudo-geometric
$\ell$-adic representations
$$\rho:G_\QQ \to \aut_{K_\rho}(V_\rho)$$
where $V_\rho$ is two-dimensional over a finite extension
$K_\rho$ of $K$ contained in $\bar{K}$, $\rho$ has
determinant $\psi\chi_\ell^{1-k}$ and reduction
isomorphic to $\rho_0$ over $\bar{\kappa}$.
We let $\CO_\rho$ denote the ring of integers of $K_\rho$.
We say such a representation $\rho$ is
an {\em allowable lift} of $\rho_0$
if its restriction to $G_\ell$ is short
and crystalline.

Suppose we are given a set $\CN$ of allowable lifts.
We assume the $\bar{K}$-isomorphism classes of the
elements of $\CN$ are distinct.  For each $\rho$, we let
$\Sigma_\rho$ denote the set of primes at which $\rho$
is not minimally ramified (see \cite{fd_bu}).
For each set of primes $\Sigma$ contained
in $\Sigma_0$, we let $\CN^\Sigma$ denote the set
of $\rho$ in $\CN$ such that $\Sigma_\rho \subset \Sigma$
and write $V^\Sigma$ for the direct sum over $\CN^\Sigma$
of the $\bar{V}_\rho = \bar{K}\otimes_{K_\rho}V_\rho$.
We assume that $\CN^\Sigma$ is finite
if $\Sigma$ is finite.

A {\em trellis} for $\CN$ is an $\CO G_\QQ$-submodule $L$ of
$V^{\Sigma_0}$ such that for each finite set $\Sigma\subset \Sigma_0$,
the $\CO$-module $L^\Sigma = L \cap V^\Sigma$ is
finitely generated and the map $\bar{K} \otimes_{\CO} L^\Sigma \to
V^\Sigma$ is an isomorphism.  One checks that if $\rho$ in $\CN$ is
such that $K_\rho = K$, then $L_\rho = L \cap V_\rho$ is a
lattice in $V_\rho$.  (Write $V_\rho$ as an intersection of
kernels of endomorphisms of $V^\Sigma$ defined by elements
of $\CO G_\QQ$.)

A {\em system of perfect pairings} $\varphi$ for $L$ is an
$\CO [G_\QQ]$-isomorphism
$$\varphi^\Sigma: L^\Sigma \to \hom_{\CO} (L^\Sigma, \CO_\psi(1-k))$$
for each finite $\Sigma \subset \Sigma_0$.
Note that for each $\rho$ in $\CN^\Sigma$,
$\varphi^\Sigma$ induces a $K_\rho G_\QQ$-isomorphism
$$\wedge^2_{K_\rho} V_\rho \to K_\rho\otimes_{K} K_\psi(1-k)$$
which we denote by $\varphi_\rho^\Sigma$.

We say that a prime $q$ is {\em horizontal} if the following hold
\begin{itemize}
\item $q \equiv 1 \bmod \ell$;
\item $\rho_0$ is unramified at $q$;
\item $\rho_0(\Fr_q)$ has distinct eigenvalues.
\end{itemize}
If $Q$ is a finite set of horizontal primes,
we let $\Delta_Q$ denote the maximal quotient of
$\prod\nolimits_{q\in Q}(\ZZ/q\ZZ)^\times$
of $\ell$-power order.

\begin{theorem}
Let $\CN$ be a set of allowable lifts of $\rho_0$
(with distinct $\bar{K}$-isomorphism classes and finite
$\CN^\Sigma$ for each finite $\Sigma \subset \Sigma_0$),
$L$ a trellis for $\CN$ and $\varphi$ a system of perfect
pairings for $L$.   Suppose that
\begin{itemize}
\item  $\CN^\emptyset \neq \emptyset$;
\item  if $\Sigma \subset \Sigma^0$ is a finite set of primes
        and $\rho \in \CN^\emptyset$, then
        $$\varphi_\rho^\Sigma = \varphi_\rho^\emptyset
        \beta_\rho^\Sigma \prod_{p\in\Sigma}
        L_p(\ad^0_{K} V_\rho,1)^{-1}$$
        for some $\beta_\rho^\Sigma$ in $\CO_\rho$;
\item  if $Q$ is a finite set of horizontal primes, then
        $$ \# \CN^Q \le \# \CN^\emptyset\cdot\# \Delta_Q$$
        and $\beta_\rho^Q \in \CO$ is independent of
        $\rho \in \CN^\emptyset$.
\end{itemize}
Then every allowable lift of $\rho_0$ is isomorphic over
$\bar{K}$ to some $\rho$ in $\CN$.  Furthermore if
$K_\rho = K$, then the lengths of
$$H^1_\Sigma(G_\QQ, A_\rho^0) \quad\mbox{and}\quad
\CO_\psi(1-k)/\mu_\rho^\Sigma(\wedge_{\CO}^2 L_\rho)$$
coincide for any finite subset $\Sigma$ of $\Sigma_0$
containing $\Sigma_\rho$.
\mlabel{thm:axiomatic}
\end{theorem}
\proofbegin
One checks that to prove the theorem, we can enlarge $K$ and so
assume that $\kappa$ contains the eigenvalues of the elements
of the image of $\rho_0$.  Note also that the hypotheses ensure
the existence of an element $\rho_\min$ of $\CN^\emptyset$.
We may assume also that $K = K_{\rho_\min}$ and we write
simply $V_\min$ and $L_\min$ for $V_{\rho_\min}$ and $L_{\rho_\min}$.

We first recall the results we need from the deformation
theory of Galois representations.  See \cite{dsl_bu},
\cite{mazur_bu} and the appendix of \cite{cdt}
for more details.  We let $\FC$ denote the category
of complete local Noetherian $\CO$-algebras.

Recall that if $A$ is an object of $\FC$ with maximal
ideal $\gm$, then an $A$-deformation of $V_0$ is an isomorphism
class of free $A$-modules $M$ endowed with continuous $AG_\QQ$-action
such that $M/\gm M$ is $(A/\gm) G_\QQ$-isomorphic to
$(A/\gm)\otimes_{\kappa}V_0$.

Suppose that $\Sigma$ is a finite subset of $\Sigma_0$.
We say that $M$ is of type $\Sigma$ if the following hold:
\begin{itemize}
\item the $AG_\QQ$-module $M$ is minimally ramified outside $\Sigma$;
\item for every $n > 0$, the $\ZZ_\ell G_\ell$-module
       $M/\gm^n M$ is an object of the category $\flcat$;
\item the $AG_\QQ$-module $\wedge^2_A M$ is isomorphic to
        $A\otimes_{\CO} \CO_\psi(1-k)$.
\end{itemize}

Consider the functor on $\FC$ which associates to
$A$ the set of $A$-deformations of $\rho_0$ of type $\Sigma$.
By the results of Mazur and Ramakrishna, this functor is representable
by an object of $\FC$.  We denote this object $R^\Sigma$ and let
$M^\Sigma$ denote the universal deformation.  We recall also
that $R^\Sigma$ is topologically generated by the elements
$t_g^\Sigma$ for $g$ in $G_\QQ$, where $t_g^\Sigma$ denotes
the trace of the endomorphism $g$ of the free $R^\Sigma$-module
$M^\Sigma$.  In particular,  $R^\Sigma$ has residue field $\kappa$.

If $\Sigma_1 \subset \Sigma_2$, then $M^{\Sigma_1}$ is an
$R^{\Sigma_1}$-deformation of $V_0$ of type $\Sigma_2$ and
hence gives rise to a natural surjection
$R^{\Sigma_2} \to R^{\Sigma_1}$.

Suppose now that $\rho$ is in $\CN$ and $\Sigma_\rho
\subset \Sigma$.  Then $\CO_\rho$ is an object of $\FC$ and there is
an $\CO_\rho$-deformation $M$ of $\rho_0$ of type
$\Sigma$ so that $V_\rho$ is $K_\rho G_\QQ$-isomorphic
to $K_\rho \otimes_{\CO_\rho} M$.  We thus obtain a continuous
$\CO$-algebra homomorphism
$$\theta_\rho^\Sigma: R^\Sigma \to K_\rho$$
so that $K_\rho \otimes_{R^\Sigma} M^\Sigma$ is isomorphic to
$V_\rho$.  The maps $\theta_\rho^\Sigma$ for varying $\Sigma \supset
\Sigma_\rho$ are compatible with the natural surjections
$R^{\Sigma_2} \to R^{\Sigma_1}$ defined above.
Note also that if $K_\rho = K$, then $A = \CO$
and $\theta_\rho^\Sigma$ defines a surjection $R^\Sigma \to \CO$.
In that case we have a natural isomorphism
\begin{equation}
\hom_{\CO}(\gp_\rho^\Sigma/(\gp_\rho^\Sigma)^2,K/\CO) \cong
 H^1_\Sigma(G_\QQ, A_\rho^0)
\mlabel{tangent}
\end{equation}
of $\CO$-modules where $\gp_\rho^\Sigma$ is the kernel of
$\theta_\rho^\Sigma$.  In particular this is the case for
$\rho = \rho_\min$ and any finite $\Sigma \subset \Sigma_0$.

We regard $V^\Sigma$ as a module for $R^\Sigma$ via
\begin{equation}
R^\Sigma \to \prod_{\rho \in \CN^\Sigma} K_\rho
\mlabel{thetas}
\end{equation}
defined by the maps $\theta_\rho^\Sigma$.
Note that if $g$ is in $G_\QQ$, then $t_g^\Sigma$
acts on $V^\Sigma$ via the endomorphism
$$\tr(\rho(g)) = g + \psi(g)\chi_\ell^{1-k}(g)g^{-1}$$
which is given by an element of $\CO G_\QQ$.
It follows that $L^\Sigma$ is stable under the action
of $R^\Sigma$ and that $\phi^\Sigma$ is $R^\Sigma$-linear.

We define the finite
flat $\CO$-algebra $T^\Sigma$ to be
the image of $R^\Sigma$ in $\End_{\CO} L^\Sigma$.
The maps $\theta_\rho^\Sigma$ induce an
isomorphism of finite $\bar{K}$-algebras
$$\bar{K} \otimes_{\CO} T^\Sigma \to
\prod_{\rho\in\CN^\Sigma}\bar{K}$$
such that $t_g^\Sigma \mapsto
(\tr\rho(g))_{\rho\in\CN^\Sigma}$ for $g$
in $G_\QQ$.  In particular $T^\Sigma$ is reduced and
$$\rank_{\CO} L^\Sigma = 2\cdot \# \CN^\Sigma
        = 2 \cdot \rank_{\CO} T^\Sigma.$$

Suppose that $\rho$ is an element of $\CN^\Sigma$
such that $K_\rho = K$.  Write $\GP_\rho^\Sigma$
for the image of $\gp_\rho^\Sigma$ in $T^\Sigma$
and $I_\rho^\Sigma$ for the annihilator of
$\GP_\rho^\Sigma$ in $T^\Sigma$.  Note that
$\GP_\rho^\Sigma$ (respectively, $I_\rho^\Sigma$)
is the set of elements in $T^\Sigma$ whose image in
$\prod \bar{K}$ has trivial component at $\rho$
(respectively, at each $\rho' \neq \rho$).

Now consider the $\CO$-module
$$\Omega_\rho^\Sigma =
        L^\Sigma/(L^\Sigma[\GP_\rho^\Sigma]
        + L^\Sigma[I_\rho^\Sigma]).$$
We define $\eta_\rho^\Sigma$ as the annihilator of the
finite torsion $\CO$-module
$$\CO_\psi(1-k)/ \varphi_\rho^\Sigma(\wedge_{\CO}^2 L_\rho).$$
\begin{lemma}  The $\CO$-module $\Omega_\rho^\Sigma$ is
        isomorphic to $(\CO/\eta_\rho^\Sigma)^2$.
\mlabel{cm1}
\end{lemma}
\proofbegin Note that the kernel of the projection $L^\Sigma \to V_\rho$
coincides with that of the composite
$$L^\Sigma \to \hom_{\CO} (L^\Sigma,\CO_\psi(1-k)) \to
\hom_{\CO} (L_\rho,\CO_\psi(1-k)) $$
where the first map is
$\varphi^\Sigma$ and the second is the natural surjection.
Denoting this kernel by $L_\rho^\perp$, we have
$L_\rho \subset L^\Sigma[\GP_\rho^\Sigma]$ and
$L_\rho^\perp \subset L^\Sigma[I_\rho^\Sigma]$.
Furthermore both inclusions are equalities since
$L^\Sigma[\GP_\rho^\Sigma]$ and $L^\Sigma[I_\rho^\Sigma]$
have trivial intersection.
Therefore the $\CO$-module $\Omega_\rho^\Sigma$ is
isomorphic to the cokernel of the map
$$L_\rho \to \hom_{\CO} (L_\rho,\CO_\psi(1-k))$$
induced by $\varphi^\Sigma$, and this is isomorphic to
$$\hom_{\CO}(L_\rho,\CO_\psi(1-k))\otimes_{\CO} \CO/\eta_\rho^\Sigma$$
(in fact, canonically so as an $\CO G_\QQ$-module).
\epf

Suppose that $q$ is a horizontal element of $\Sigma$ and choose
an ordering $\alpha_1, \alpha_2$ of the eigenvalues of
$\rho_0(\Fr_q)$.  We then have a unique decomposition
$$M^\Sigma = M^\Sigma_1 \oplus M^\Sigma_2$$
as an $R^\Sigma G_q$-module such that for $i = 1, 2$,
$\kappa \otimes_{R^\Sigma} M^\sigma_1$ is unramified
with $\Fr_q$ acting via $\alpha_i$.  We let
$\xi_q^\Sigma$ denote the character
$G_q \to (R^\Sigma)^\times$ defined by the action
of $G_q$ on $M^\Sigma_1$.  The restriction of
$\xi_q^\Sigma$ to the inertia group $I_q$ factors through
$$I_q \to \ZZ_q^\times \to \Delta_{\{q\}}$$
where the first map is gotten from local class field
theory and the second is the natural projection.
We thus obtain a homomorphism $\Delta_{\{q\}}
\to (R^\Sigma)^\times$.  If $Q$ is a finite set of
horizontal primes, then we regard $R_Q$ as an
$\CO\Delta_Q$-algebra via these homomorphisms
and so regard $L^Q$ as an $\CO\Delta_Q$-module.

Now let $\gp^Q$ denote the augmentation ideal of $\CO[\Delta_Q]$,
i.e., the kernel of the map $\CO[\Delta_Q] \to \CO$ defined
by $\delta \mapsto 1$ for $\delta$ in $\Delta_Q$.  Let $S^Q$
denote the image of $\CO[\Delta_Q]$ in the ring of endomorphisms
of $L^Q$ and write $\GP^Q$ for the image of $\gp^Q$ and $I^Q$
for the annihilator of $\GP^Q$ in $S^Q$.  Now consider the
$\CO$-module
$$\Omega^Q = L^Q/(L^Q[\GP^Q]+L^Q[I^Q]).$$
\begin{lemma}  The $\CO$-rank of $L^Q[\GP^Q]$ is
        $\#\CN^\emptyset$ and the $\CO$-length of $\Omega^Q$
        is at least $\#\Delta_Q\cdot\#\CN^\emptyset$.
\mlabel{cm2}
\end{lemma}
\proofbegin  First note that if $\rho$ is in $\CN^Q$, then $\rho$
is in $\CN^\emptyset$ if and only if $\Delta_Q$ acts trivially
on $V_\rho$.  It follows that $L^\emptyset\subset L^Q[\GP^Q]$
and $(L^\emptyset)^\perp \subset L^Q[I^Q]$ where
$(L^\emptyset)^\perp$ is the kernel of the projection
$L^Q \to V^\emptyset$.  This kernel is the same as that
of the composite
$$L^Q \to \hom_{\CO} (L^Q,\CO_\psi(1-k)) \to
\hom_{\CO} (L^\emptyset,\CO_\psi(1-k)) $$
where the first map is $\varphi^Q$ and
the second is the natural surjection.  Since
$L^Q[\GP^Q]$ and $L^Q[I^Q]$ have trivial intersection and
$$L^\emptyset \oplus (L^\emptyset)^\perp$$
has finite index in $L^Q$, it follows that
$L^\emptyset = L^Q[\GP^Q]$ and $(L^\emptyset)^\perp
= L^Q[I^Q]$.  Therefore the $\CO$-module $L^Q[\GP^Q]$
has rank $\#\CN^\emptyset$.  Moreover $\Omega^Q$ is
isomorphic to the cokernel of the endomorphism
$\nu_Q$ of $L^\emptyset$ obtained by composing the inclusion
$L^\emptyset \to L^Q$ with its adjoint with respect
to $\varphi^Q$ and $\varphi^\emptyset$.
Our hypotheses on the pairings ensure that
that $\nu_Q = \#\Delta_Q\nu_Q'$ for some
$\nu_Q'$ in $\CO G_\QQ$.
\epf

The following is proved exactly as in \cite{tw} except
that we use Corollary \ref{coro:local}.
\begin{lemma}
There exists an integer $r\ge 0$ and sets
of horizontal primes $Q_n$ for each $n \ge 1$ such that
the following hold:
\begin{itemize}
\item $\#Q_n = r$;
\item $q \equiv 1 \bmod \ell^n$ for each $q \in Q_n$;
\item $R^{Q_n}$ is generated by $r$ elements as an $\CO$-algebra.
\end{itemize}
 \mlabel{gal}
\end{lemma}

We are now ready to prove that $R^\emptyset$ is a complete
intersection over which $L^\emptyset$ is free of rank two.
We let $r$ and $Q = Q_n$ for $n \ge 1$ be as in lemma \ref{gal}.
Lemma \ref{cm2} allows us to apply the implication (a) $\Rightarrow$
(c) of Theorem 2.4 of \cite{fd_twc} to conclude that
$L^Q$ is a free  $\CO[\Delta_Q]$-module of
rank $\#\CN^\emptyset$.  Since the adjoint of the natural inclusion
$L^\emptyset \to L^Q$ is surjective, it induces
an isomorphism $L^Q/\gp^Q L^Q \to L^\emptyset$ of
$R^Q$-modules.

Setting $A = \kappa[[S_1,\ldots,S_r]]$,
$B= \kappa[[X_1,\ldots,X_r]]$, $R = \kappa \otimes_{\CO} R^\emptyset$
and $H = \kappa\otimes_{\CO}L^\emptyset$, we shall define $B$-modules
$H_n$ and maps $\phi_n:A\to B$, $\psi_n: B \to R$ and
$\pi_n:H_n \to H$ satisfying the hypotheses of Theorem 1.3
of \cite{fd_twc}.  We first choose surjective $\kappa$-algebra
homomorphisms $A \to \kappa[\Delta_{Q_n}]$ and
$B \to R_n$ where $R_n = \kappa\otimes_{\CO} R^{Q_n}$.
Define $\psi_n$ as the composite $B \to R_n \to R$
and define $\phi_n:A \to B$ so the diagram
$$\begin{array}{ccc} A & \to & B \\ \downarrow && \downarrow\\
\kappa[\Delta_{Q_n}] & \to & R_n \end{array}$$
commutes.  We consider $L_n = \kappa \otimes_{\CO}L^{Q_n}$ as
a $B$-module via $B \to R_n$, and define $H_n$ as
$L_n/\gm_A^n L_n$ and $\pi_n$ as the map induced by
$L^{Q_n} \to L^\emptyset$.  We can then apply Theorem 1.3
of \cite{fd_twc} to conclude that $R$ is a complete intersection
over which $H$ is a free module.  It then follows that
$R^\emptyset$ is a complete intersection over which
$L^\emptyset$ is a free module of rank 2.

We now apply the implication (c) $\Rightarrow$ (b) of
Theorem 2.4 of \cite{fd_twc} to the $R^\emptyset$-module
$L^\emptyset$ and prime ideal $\gp_{\rho_\min}^\emptyset$.
We thus obtain the formula
$$\begin{array}{rl}
2\cdot\length_{\CO} H^1_\emptyset(G_\QQ, A_\min^0)
& = 2\cdot\length_{\CO} \gp_{\rho_\min}^\emptyset
                /(\gp_{\rho_\min}^\emptyset)^2\\
= \length_{\CO} \Omega_{\rho_\min}^\emptyset&
= 2\cdot v(\eta_{\rho_\min}^\emptyset)\end{array}$$
where the first equality follows from (\ref{tangent})
and the last from Lemma \ref{cm1}.

Suppose now that $\Sigma$ is a finite subset of $\Sigma_0$.
Applying (\ref{tangent})
and Lemma \ref{cm1} again, together with the inequality
$$\length_{\CO} H^1_\Sigma(G_\QQ,A_\min^0)
 \le \length_{\CO} H^1_\emptyset(G_\QQ,A_\min^0)
                -\sum_{p\in\Sigma}v(L_p(\ad^0_{K} V_\rho,1))$$
obtained from a Galois cohomology argument, we find that
$$\begin{array}{rl}
2\cdot\length_{\CO} \gp_{\rho_\min}^\Sigma/(\gp_{\rho_\min}\Sigma)^2
& = 2\cdot\length_{\CO} H^1_\Sigma(G_\QQ,A_\min^0)\\
        \le 2\cdot v(\eta_{\rho_\min}^\Sigma) &
        = \length_{\CO} \Omega_{\rho_\min}^\Sigma.\end{array}$$
We can then apply the implication (a) $\Rightarrow$ (c)
of Theorem 2.4 of \cite{fd_twc} to conclude
$R^\Sigma$ is a complete intersection over which $L^\Sigma$
is a free module of rank 2.

The second assertion of the theorem follows from another application
of Theorem 2.4 of \cite{fd_twc}, (\ref{tangent}) and lemma \ref{cm1}.
To deduce the first assertion of the theorem, note that the map
$$\bar{K}\otimes_{\CO}R^\Sigma\to \prod_{\rho\in\CN^\Sigma} \bar{K}$$
induced by (\ref{thetas}) is an isomorphism.  Every allowable
lift of $\rho_0$ arises, up to $\bar{K}$-isomorphism, from
a $\bar{K}$-linear map $\bar{K} \otimes_{\CO}R^\Sigma$
for some $\Sigma$.  It therefore arises from $\theta_\rho^\Sigma$ for
some $\rho$ in $\CN$.
\epf

\subsection{The construction}
Fix a number field $K_0 \subset \CC$.  Let $\lambda_0$
denote the prime over $\ell$ determined by our fixed
embeddings $\bar{\QQ} \to \CC$ and $\bar{\QQ} \to \Qlbar$,
and let $S = \{\lambda_0\}$.
Let $\CO_0$ denote the valuation ring of $\lambda_0$ in $K_0$,
$K$ the completion at $\lambda_0$ of $K_0$,
$\CO$ its ring of integers, $\lambda$ a uniformizer
and $\kappa$ its residue field.

We consider a representation
$$\rho_0 : G_\QQ \to \aut_{\kappa}(V_0)$$
as in the preceding section.  We assume also that $\rho_0$
has minimal conductor among its twists by characters
$G_\QQ \to \bar{\kappa}^\times$.
Let $k$ denote the Serre weight and $\psi$
the character associated to $\rho_0$, as in the
preceding section.

If $f$ is a newform, we let $K_{0,f}$ denote the extension
of $K_0$ generated by the Fourier coefficients of $f$,
and $\lambda_f$ the prime over $\lambda_0$ determined
by our embeddings.  Define $K_f$ as the completion, $\CO_f$
as its ring of integers and $\kappa_f$ as its residue field.

We assume that $\rho_0$ is {\em modular,} in the sense that
$$
\kappa_f \otimes_{\kappa} V_0 \cong
\kappa_f \otimes_{\CO_f} \CM_{f,\lambda_f}
$$
for some newform $f$ (which can have any weight
level and character).
We let $\CN$ denote the set of representations
$$\rho_f: G_\QQ \to \aut_{K_f} M_{f,\lambda_f}$$
such that
where $f$ is a newform of weight $k$ character $\psi$
and level $N_f$ not divisible by $\ell$.
Thus $\CN$ is in bijection with the set of modular
forms $f$ such that $\rho_f$ is an allowable lift
of $\rho_0$.  From the work of Ribet and others,
one knows that $\CN_\emptyset$ is non-empty.
(Combine the main results of \cite{dt_duke} and
\cite{fd_hk}.)

Let $\Sigma_1 \subset \Sigma_0$ denote the set of
primes at which $\rho_0$ is ramified.  We let
$\Sigma_2\subset\Sigma_1$ denote the following
set of exceptional primes for $\rho_0$:
those $p \equiv -1 \bmod \ell$ such that $\rho_0|G_p$ is absolutely
irreducible, but $\rho_0|I_p$ is not
For each such $p$, $\bar{\kappa}\otimes_\kappa \rho_0$ is
induced from a character of $H_p$, and we let $\xi_p$, $\xi_p'$
denote the pair of characters $\CO_F^\times \to \bar{\kappa}^\times$
defined as in \S\ref{sec:ep}.  We then write
$$\sigma_p^\min : \GL_2(\ZZ_p) \to \aut_{K_0} V_p^\min$$
for the representation associated to $\tilde{\xi}_p$,
$e_p^\min$ for a distinguished vector and $\mu_p^\min$ for a
normalized pairing.  Here we assume $K_0$ is sufficiently large
so that
\begin{itemize}
\item $K_0$ contains the values of $\tilde{\psi}$;
\item for each $p$ in $\Sigma_2$, $\sigma_p^\min$, $e_p^\min$ and
$\mu_p^\min$ are defined over $K_0$.
\end{itemize}

Write the Artin conductor $N(\rho_0)$ as $\prod_p p^{c_p}$.
Note that $c_p$ is even if $p$ is in $\Sigma_2$.
We let $c_p^\min = 0$ if $p$ is in $\Sigma_2$ and
$c_p^\min = c_p$ otherwise.  For any prime $p\neq\ell$,
we let $c_p^{\max} = c_p+\delta_p$
where $\delta_p = \dim_{\kappa}V_0^{I_p}$.
For each finite $\Sigma \subset \Sigma_0$, we let
$c_p^\Sigma = c_p^{\max}$ or $c_p^\min$ according to
whether or not $p$ is in $\Sigma$.  We let $U_p^*$
denote the set of matrices
$\smat{a}{b}{c}{d} \in \GL_2(\ZZ_p)$ with
$c \in p^{c_p^*}$, where $*$ can be $\min$
or ${\max}$.  We let $U_p^\Sigma = U_p^{\max}$ or $U_p^\min$
according to whether or not $p$ is in $\Sigma$.  We then set
$$U^\Sigma = \prod_p U_p^\Sigma =
\prod_{p\in\Sigma} U_p^{\max} \prod_{p\not\in\Sigma} U_p^\min
=U_0(N^\Sigma),$$
where $N^\Sigma = \prod_{p} c_p^\Sigma$.

We have already defined a representation
$$\sigma_p^\min : U_p^\min \to \aut_{K_0} V_p^\min$$
for $p\in\Sigma_2$.  If $c_p^\min > 0$, we define
$\sigma_p^\min$ by $V_p^\min = K_0$ with $\smat{a}{b}{c}{d}$
acting via $\psi_p(d)$, where
$\psi_p = \tilde{\psi}|_{\ZZ_p^\times}$.  If $p$ does not
divide $N(\rho_0)$ then $\sigma_p^\min$ is trivial.
For any $p$, we define $\sigma_p^{\max} : U_p^{\max} \to K_0^\times$
by $\psi_p(d)$.  We let $\sigma_p^\Sigma = \sigma_p^{\max}$
or $\sigma_p^\min$ according to whether or not $p\in\Sigma$.
We then define the representation $V_0^\Sigma$ of $U^\Sigma$
as the tensor product over all $p$ of $V_p^\Sigma$.
Similarly we define the $\CO_0$-lattice $\CV^\Sigma$
in $V_0^\Sigma$ as the tensor product over $\CO$ of
$\CV_p^\Sigma$.

If $\rho\in\CN^\Sigma$, then $U_f^\Sigma = U^\Sigma$
and $\sigma_f^\Sigma = K_{0,f}\otimes_{K_0}\sigma^\Sigma$.
We therefore have a natural inclusion
$V_\rho \to K_\rho\otimes_K M(\sigma^\Sigma)_{!,\lambda_0}$
induced by the isomorphism $M_f \to M_f^\Sigma$.
These give rise to an inclusion
$$V^\Sigma \to \bar{K}\otimes_K M(\sigma^\Sigma)_{!,\lambda_0}.$$
We define $L^\Sigma$ as the preimage of
$\CM(\sigma^\Sigma)_{!,\lambda_0}$ in $V^\Sigma$,
where the integral structure is defined as in \S\ref{ss:sig.int}
using $\CV^\Sigma$.

Let $\tilde{\TT}^\Sigma$ denote the polynomial algebra
over $\CO$ generated by the variables $t_p$, where
$p$ runs over rational primes $p$
such that $p \in \Sigma$ or $p \not \in \Sigma_1$.
We extend the natural action of $\CO$ on $V^\Sigma$
to one of $\tilde{\TT}^\Sigma$ by defining the
action of $t_p$ on $V_\rho$ to be $0$ or
$\tr(\rho(\frob_p))$ according to whether
or not $p\in\Sigma$.  Arguing as in the proof of
Theorem \ref{thm:axiomatic}, we see that $L^\Sigma$
is stable under the action of $\tilde{\TT}^\Sigma$.

We also regard $\CM(\sigma^\Sigma)_{!,\lambda_0}$ as a module for
$\tilde{\TT}^\Sigma$ where $t_p$ acts as the
usual double coset operator.   We let $\gm^\Sigma$
denote the maximal ideal of $\tilde{\TT}^\Sigma$
which is the kernel of the homomorphism defined
by $t_p \mapsto 0$ or $\tr(\rho_0(\frob_p))$ according
to whether or not $p$ is in $\Sigma$.  We let $\CM_!^\Sigma$
denote the localization of $\CM(\sigma^\Sigma)_{!,\lambda_0}$ at
$\gm^\Sigma$.  Thus $\CM_!^\Sigma$ can be canonically
identified with a quotient (also a direct summand)
of $\CM(\sigma^\Sigma)_{!,\lambda_0}$.

\begin{lemma}
The composites of the natural maps
$$\begin{array}{ccccc}L^\Sigma &\to& \CM(\sigma^\Sigma)_{!,\lambda_0} &\to&
\CM_!^\Sigma \\
V^\Sigma &\to& \bar{K}\otimes_\CO\CM(\sigma^\Sigma)_{!,\lambda_0} &\to&
\bar{K}\otimes_\CO\CM_!^\Sigma\end{array}$$
are isomorphisms of $\tilde{\TT}^\Sigma$-modules.
\mlabel{lem:messy}
\end{lemma}
\proofbegin  Let $\bar{\TT}^\Sigma$ denote the image of $\TT^\Sigma$ in
$\End_\CO \CM(\sigma^\Sigma)_{!,\lambda_0}$.  Let $I^\Sigma$ denote
the kernel of $\bar{\TT}^\Sigma \to \bar{\TT}^\Sigma_{\gm^\Sigma}$,
or equivalently, the annihilator of $\CM_!^\Sigma$.

As in \cite[\S4.2]{ddt} (see also \cite[Lemma 5.1.1]{cdt}),
one gets an explicit description of $\bar{K}\otimes_\CO
\bar{\TT}^\Sigma_{\gm^\Sigma}$:  There is an isomorphism
$$\bar{K}\otimes_\CO\bar{\TT}^\Sigma_{\gm^\Sigma}\to
        \prod_{\rho\in\CN^\Sigma}\bar{K}$$
defined by $(1\otimes t_p) \mapsto 0$ or $\tr(\rho(\frob_p))_\rho$
according to whether or not $p$ is in $\Sigma$.  We can then
identify $V^\Sigma$ with the $\bar{\TT}^\Sigma$-module
$$\bar{K}\otimes_\CO\CM(\sigma^\Sigma)_{!,\lambda_0}
        [\bar{K}\otimes_\CO I^\Sigma].$$
It follows that the second composite in the statement of
the lemma is an isomorphism.  The first is therefore
injective and
$$L^\Sigma = \CM(\sigma^\Sigma)_{!,\lambda_0}[I^\Sigma],$$
so it is also surjective.
\epf

Suppose now that $\Sigma \subset \Sigma' \subset \Sigma_0$.
For each exceptional prime $p$ in $\Sigma' - \Sigma$, we have
a map $\tau_p:V_p^{\max} \to V_p^\min$ restricting to
$\CV_p^{\max} \to \CV_p^\min$ (see \S\ref{sec:ep}).  We
define $\tau: \CV^\Sigma \to \CV^{\Sigma'}$
as the tensor product of the $\tau_p$.  Let
$$N' = N^{\Sigma'}/\prod_{p} p^{c_p/2}\quad\mbox{and}\quad
N = N^{\Sigma}\prod_{p} p^{c_p/2}$$
where the products are over exceptional primes in $\Sigma - \Sigma'$.
For positive integers $d$ dividing $N'/N$, we let $d' = dN/N^\Sigma$
and let $\alpha_d$ denote the matrix in $M_2(\hat{\ZZ})$ whose
component at $p$ is $\smat{p^{-v_p(d')}}{0}{0}{1}$.  We then
define $\gamma:M(\sigma^\Sigma)_! \to M(\sigma^{\Sigma'})_!$ by
$$\gamma = \sum_d \frac{1}{d'}
        \left[ U^{\Sigma'}\alpha_d U^\Sigma\right]_\tau \beta_d.$$
Then for each $\rho \in \CN^\Sigma$, $\gamma$ induces the identity
map on $V_\rho$.  Furthermore $\gamma$ restricts to a map
$\CM(\sigma^\Sigma)_! \to \CM(\sigma^{\Sigma'})_!$.  It follows
that $L^\Sigma \subset L^{\Sigma'}$.  To prove that
$L = {\displaystyle\bigcup_\Sigma} L^\Sigma$ is a trellis for $\CN$,
we just need the following:

\begin{proposition}   If
$\Sigma \subset \Sigma' \subset \Sigma^0$,
then $L^{\Sigma'}/L^\Sigma$ is torsion-free.
\mlabel{lem:tfc}
\end{proposition}
\proofbegin  It suffices to consider the case $\Sigma' = \Sigma \cup \{p\}$ for
some $p$ not in $\Sigma \cup \{\ell\}$.  We let
$c= c_p$, $\gamma = {\max}\{1,c\}$ and $\delta = \delta_p$.

Choose an integer $N \ge 3$ such that $(p\ell,N) = 1$ and
$$V_N^p = \ker (\GL_2(\hat{\ZZ}) \to \GL_2(\ZZ/N\ZZ)\times \GL_2(\ZZ_p))$$
is contained in $\ker\sigma^{\Sigma}$.
We now define subgroups $V_p$ and
$V_p'$ of $\GL_2(\ZZ_p)$ such that $V_p\subset \ker\sigma_p^\min$
is normal in $U_p^\min$, $V_p' \subset \ker\sigma_p^{\max}$
is normal in $U_p^{\max}$ and
$$ V_p' \subset \bigcap_{i=0}^{\delta} \alpha_{p^i}V_p\alpha_{p^i}^{-1}.$$
(Recall that $\alpha_{p^i} = \smat{p^{-i}}{0}{0}{1}$
unless $p\in \Sigma_2$, in which case
$\alpha_1 = \smat{p^{-c/2}}{0}{0}{1}_p$.)
If $p$ is in $\Sigma_2$, we define $V_p$ as the kernel
of $\GL_2(\ZZ_p) \to \GL_2(\ZZ/p^{c/2})$ and
$V_p' = \alpha_{1}V_p\alpha_{1}^{-1}$.
Otherwise, we define $V_p$ as the set of matrices in $\GL_2(\ZZ_p)$
congruent to $\mat{1}{*}{0}{1}$ mod $p^c$, and
$$V_p' = \left\{\left.\,\mat{w}{x}{y}{z}\,\right|\,
        w\equiv z \equiv 1\bmod p^\gamma,\,y\equiv 0\bmod p^{c+\delta}
        \,\right\}.$$
We then let $V=V_N^pV_p$, $V' = V_N^pV_p'$.
Let $\tilde{\TT}$ denote
the common subring of $\tilde{\TT}^\Sigma$ and $\tilde{\TT}^{\Sigma'}$
generated by the variables $t_p$ for $p$ not dividing $Np\ell$, and let
$\gm = \gm^\Sigma \cap \tilde{\TT} = \gm^{\Sigma'} \cap \tilde{\TT}$.

We wish to prove that the map
$$\kappa\otimes_{\CO}L^\Sigma \to \kappa\otimes_{\CO}L^{\Sigma'}$$
is injective.  The natural map $L^\Sigma \to \CM(\sigma^\Sigma)_{!,\gm}$
is injective, so it suffices to prove the map
$$\kappa\otimes_{\CO}
\CM(\sigma^\Sigma)_{!,\gm} \to \kappa\otimes_\CO
\CM(\sigma^{\Sigma'})_{!,\gm}$$
defined in the same way as $\gamma$ is injective.
Therefore it suffices to prove that
$$\bigoplus_{i=0}^{\delta}\kappa\otimes_\CO\CM(\sigma^\Sigma)_{!,\gm} \to
\kappa\otimes_\CO\CM(\sigma^{\Sigma'})_{!,\gm}$$
is injective, where the map is defined by
$$\bigoplus_{i=0}^{\delta} p^{-i}
        \left[ U^{\Sigma'}\alpha_{p^i} U^\Sigma\right]_\tau.$$
Since $\rho_0$ is irreducible, $\gm$ is not in the support
of the kernel of $\CM(V)_c \to \CM(V)_!$ or $\CM(V')_c \to \CM(V')_!$.
We can therefore identify the above map with the map
$$\bigoplus_{i=0}^{\delta}\kappa\otimes_\CO(\CV^\Sigma\otimes_{\CO_0}
\CM(V))_{!,\gm}^{U^\Sigma} \to
\kappa\otimes_\CO(\CV^{\Sigma'}\otimes_{\CO_0}
\CM(V'))^{U^{\Sigma'}}_{!,\gm}$$
induced by
$$\bigoplus_{i=0}^{\delta} p^{-i}
        \left[ V' \alpha_{p^i} V\right]_\tau.$$
We shall prove that in fact if $\delta < 2$, then
$$\gamma':
\bigoplus_{i=0}^{\delta}(\CV_0^\Sigma\otimes_{\CO_0}
\CM(V))_!^{U^\Sigma} \to
(\CV_0^{\Sigma'}\otimes_{\CO_0}
\CM(V'))_!^{U^{\Sigma'}}$$
is injective,
where $\CV_0^\Sigma = \kappa\otimes_{\CO_0}\CV^\Sigma$
and $\CV_0^{\Sigma'} = \kappa\otimes_{\CO_0}\CV^{\Sigma'}$,
and that the same holds after localization at $\gm$ if
$\delta = 2$.

First suppose $\delta =0$.  If $p$ is not in $\Sigma_2$, then
$\gamma'$ is the identity and there is nothing to prove.
So assume $p$ is in $\Sigma_2$.  By Shapiro's Lemma, we
can replace the target of $\gamma'$ with
$$\left(\Ind_{\alpha U^{\Sigma'} \alpha^{-1}}^{U^\Sigma}
\alpha^{-1}\CV_0^{\Sigma'}\otimes_{\CO_0}\CM(V)_!\right)^{U^\Sigma},$$
where $\alpha = \alpha_1$ and the map is defined so the composite
$$\CV_0^\Sigma \to \Ind_{\alpha U^{\Sigma'} \alpha^{-1}}^{U^\Sigma}
\alpha^{-1}\CV_0^{\Sigma'} \to \alpha^{-1}\CV_0^{\Sigma'}$$
is $\alpha^{-1}\tau$.  Therefore it suffices to prove that
$$\kappa\otimes_{\CO_0}\CV_p^\min \to \kappa\otimes_{\CO_0}
\Ind_{\alpha U^{\Sigma'} \alpha^{-1}}^{U^\Sigma}
\alpha^{-1}\CV_p^{{\max}}$$
is injective.  This holds because  $\kappa\otimes_{\CO_0}\CV_p^\min$
is irreducible and $\kappa\otimes_{\CO_0}\tau_p$ is non-zero
(since it is surjective by its construction).

Now suppose $\delta = 1$.  Let  $\CV_0 = \CV_0^\Sigma = \CV_0^{\Sigma'}$,
$V'' = \alpha_p^{-1}V'\alpha_p$ and
$A = \CV_0\otimes_{\CO_0}\CM_!(V'')$.  Then $V''$
is normal in both $U^\Sigma$ and $\alpha_p^{-1}U^\Sigma\alpha_p$,
and both groups act on $A$, where the action of
$\alpha_p^{-1}U^\Sigma\alpha_p$ on $\CV_0$ is defined
by first conjugating by $\alpha_p$.  Note that the actions agree on
$$\alpha_p^{-1}U^{\Sigma'}\alpha_p
=U^\Sigma \cap \alpha_p^{-1}U^\Sigma\alpha_p.$$
We claim that the map
$$A^{\alpha_p^{-1}U^\Sigma\alpha_p} \oplus A^{U^\Sigma}\to
A^{\alpha_p^{-1}U^{\Sigma'}\alpha_p}$$
is injective.  To prove this, suppose there were a nonzero
vector in $A$ fixed by both $\alpha_p^{-1}U^\Sigma\alpha_p$
and $A^{U^\Sigma}$.  Then there would be a nonzero vector
in $\CM_!(V'')$ on which $\smat{a}{b}{c}{d}$ acts via
$\psi_p^{-1}(d)$ for all matrices in $U_p^\min$ or
$\alpha_p^{-1}U_p^\min\alpha_p$.  This contradicts the
fact that $\psi_p$ has conductor $p^c$.
It follows that the map
$$B^{U^\Sigma} \oplus A^{U^\Sigma} \to B^{U^{\Sigma'}}$$
where $B = \CV_0\otimes_{\CO_0}\CM_!(V')$ and the map
$A \to B$ is defined by
$p^{-1} \left[ V' \alpha_p V'' \right].$
Since $V/V'$ and $V/V''$ have order $p$, we can
identify $\CM(V')_!^V$ and $\CM(V'')_!^V$ with $\CM(V)_!$
and conclude $\gamma'$ is injective.

Finally, suppose that $\delta = 2$.  Again we let
$\CV_0 = \CV_0^\Sigma = \CV_0^{\Sigma'}$,
$V'' = \alpha_p^{-1}V'\alpha_p$ and
$A = \CV_0\otimes_{\CO_0}\CM_!(V'')$.
We also define $\tilde{U}$ to be the subgroup of
$U^\Sigma$ generated by $U^{\Sigma'}$ and
$\alpha_p^{-1} U^\Sigma \alpha_p$.
Then we have an exact sequence
$$ 0 \to A^{U^\Sigma} \to
 A^{\alpha_p^{-1}\tilde{U}\alpha_p} \oplus A^{\tilde{U}}\to
A^{\alpha_p^{-1}U^{\Sigma'}\alpha_p}$$
with maps defined by $a\mapsto (a,-a)$ and $(a_1,a_2)\mapsto
a_1 + a_2$.  This gives an exact sequence
\begin{equation}
0 \to A^{U^\Sigma} \to B^{\tilde{U}} \oplus A^{\tilde{U}}
\to B^{U^{\Sigma'}}
\mlabel{abab}
\end{equation}
where $B = \CV_0\otimes_{\CO_0}\CM(V')_!$ and
the maps are  $a\mapsto (\beta(a),-a)$ and $(b,a)\mapsto
b + \beta(a)$, with $\beta$ defined by
$p^{-1} \left[ V' \alpha_p V''\right]_1.$
Now let $C = \CV_0\otimes_{\CO_0}\CM(V)_!$ and
$D = \CV_0\otimes_{\CO_0}\CM(\tilde{V})_!$, where
$\tilde{V} = V \cap \tilde{U}$.
The proof of Lemma 3.2 of \cite{fd_ast} shows the
injectivity of the map $C \times C \to D$ defined by
$$(a,b) \mapsto \left[\tilde{V}\cdot 1 \cdot V\right]_1 (a)
+ \left[\tilde{V} \alpha_p V\right]_1 (b).$$
So we also get the exactness of
\begin{equation}
0 \to C^{U^\Sigma} \times C^{U^\Sigma} \to D^{\tilde{U}}.
\mlabel{ihara}
\end{equation}

Now let $A_c$, $B_c$, $C_c$ and $D_c$ denote the
$\tilde{\TT}$-modules defined in the same way as
$A$, $B$, $C$ and $D$, but using $\CM_c$ instead
of $\CM_!$.  Then the natural maps $C_c \to A_c^{V}$,
$D_c \to B_c^{\tilde{V}}$  and
$D_c \to A_c^{\tilde{V}}$  are isomorphisms.
Since $\gm$ is not Eisenstein, it follows that so
are the maps $C_\gm^{U^\Sigma} \to A_\gm^U$, $D_\gm^{\tilde{U}}
\to B_\gm^{\tilde{U}}$ and $D_\gm^{\tilde{U}}
\to A_\gm^{\tilde{U}}$.  We now obtain
the desired injectivity on combining the localization at $\gm$
of the exact sequences (\ref{abab}) and (\ref{ihara})
(cf.\ the discussion after Lemma 2.5 of \cite{wiles}).
\epf

Recall that in \S \ref{sec:ep}, we defined a perfect pairing on
$M(\sigma^\Sigma)_!$, and we now let $\phi^\Sigma$ denote the paring
$$M(\sigma^\Sigma)_{!,\lambda_0} \otimes_K M(\sigma^\Sigma)_{!,\lambda_0}
         \to K.$$
\begin{lemma}
The restriction of $\phi^\Sigma$ defines a perfect
pairing $L^\Sigma \otimes_\CO L^\Sigma \to \CO$.
\mlabel{lem:pp}
\end{lemma}
\proofbegin
First we apply Lemma 3 of \cite{dt_duke} to choose a prime $r \ge 5$
at which no lift of $\rho_0$ is ramified.  It follows that $L^\Sigma$
can be identified with $L^{\Sigma \cup \{r\}}$ compatibly with the map
$$\CM(\sigma^\Sigma)_{!,\lambda_0}
        \to \CM(\sigma^{\Sigma\cup\{r\}})_{!,\lambda_0}.$$
According to Proposition \ref{prop:exlevel},
the restriction to $L^\Sigma$ of the transpose of this map
with respect to our pairings is given by an element of
$\tilde{\TT}^\Sigma - \gm^{\Sigma}$, and so by
Lemma \ref{lem:messy}, is an automorphism.  So to prove
the lemma, we can assume $r$ is in $\Sigma$.

Let $U'$ denote the set of $\mat{a}{b}{c}{d}$ in $U^\Sigma$ such
that $d_r \equiv 1\bmod r$, and similarly define $U''$ requiring
$a_r \equiv 1\bmod r$.  We let $\sigma'$ and $\sigma''$ denote
the restrictions of $\sigma^\Sigma$ to $U'$ and $U''$ respectively.
If we define an action of $\tilde{\TT}^\Sigma$ on $\CM(\sigma')_!$
and $\CM(\sigma'')_!$ in the usual way and the natural inclusions
$$\CM(\sigma)_! \to \CM(\sigma')_!\quad\mbox{and}\quad
\CM(\sigma)_! \to \CM(\sigma'')_!$$
are $\tilde{\TT}^{\Sigma}$-linear.
We record the following as a lemma for future reference:
\begin{lemma} The natural maps induce isomorphisms
$$\CM(\sigma)_{!,\gm^\Sigma} \to \CM(\sigma')_{!,\gm^\Sigma}
\quad\mbox{and}\quad
\CM(\sigma)_{!,\gm^\Sigma} \to \CM(\sigma'')_{!,\gm^\Sigma}.$$
 \mlabel{sublemma}
\end{lemma}
To prove this, note that we can replace $!$ by $c$ since
$\gm^\Sigma$ is not Eisenstein.  So we must prove the inclusion
$$\CM(\sigma')_c^{U^\Sigma} \to \CM(\sigma')_c$$
localized at $\gm^\Sigma$ is an isomorphism, and similarly
for $\sigma''$.  Let $u$ be a generator for $U^\Sigma/U'$,
then $u$ acts on $\CM(\sigma')_c$ as
$p^{2-k}S_p$ for some prime $p$ at which $\rho_0$ is
unramified and $\det\rho_0(\frob_p) = p^{k-1}$.
It follows that $u$ acts as an element of $1+\gm^\Sigma$.
Since $u$ has order $r-1$, which is not divisible by $\ell$,
it follows that $u = 1$ on $\CM(\sigma')_{c,\gm^\Sigma}$.
The same argument works for $\sigma''$ and this completes
the proof of Lemma \ref{sublemma}.

Since $U'$ is sufficiently small, we have a perfect pairing
between $\CM(\sigma')_!$ and $\CM(\hat{\sigma}'\otimes
\psi^{-1}\circ\det)_!$ given by Proposition \ref{prop:spair}.
We also have an isomorphism
$$\CM(\sigma'')_! \to \CM(\hat{\sigma}'\otimes\psi^{-1}\circ\det)_!$$
defined as in \S \ref{sec:ep}, and so we obtain a
perfect pairing between $\CM(\sigma')_!$ and $\CM(\sigma'')_!$.
It is straightforward to check that the resulting isomorphism
$$\CM(\sigma'')_{!,\lambda_0} \to
\hom_{\CO}(\CM(\sigma')_{!,\lambda_0}, \CO)$$
is $\tilde{\TT}^\Sigma$-linear, and that the composite
$$\CM(\sigma)_{!,\lambda_0} \otimes_\CO \CM(\sigma)_{!,\lambda_0}
\to \CM(\sigma')_{!,\lambda_0} \otimes_\CO \CM(\sigma'')_{!,\lambda_0}
        \to \CO$$
is $(r-1)\phi^\Sigma$.  It follows from Lemmas \ref{lem:messy} and
\ref{sublemma} that $\phi^\Sigma$ restricts to a perfect
pairing on $L^\Sigma$.
\epf

In view of Proposition \ref{prop:exlevel},
we only need to prove
$$\#\CN^Q \le \#\Delta_Q \#\CN^\emptyset$$
to enable us to apply Theorem \ref{thm:axiomatic}.
We shall in fact prove:
\begin{lemma} For any finite set of primes $\Sigma \subset \Sigma_0$
and horizontal prime $q \not\in\Sigma$, we have
$$\#\CN^{\Sigma\cup\{q\}} = \#\Delta_{\{q\}} \#\CN^\Sigma.$$
\end{lemma}
\proofbegin  First note that we may assume $\Sigma$ contains an
auxiliary prime $r$ as in the proof of Lemma \ref{lem:pp}.
Let $\tilde{\TT} = \tilde{\TT}^{\Sigma} = \tilde{\TT}^{\Sigma\cup\{q\}}$.
Choose an eigenvalue $\alpha_q$ of $\rho_0(\frob_q)$ and let
$\gm$ denote the maximal ideal which is the kernel of
the map defined by $t_q \mapsto \alpha_q$, $t_p \mapsto 0$ for
$p\in\Sigma$, and $t_p \mapsto \tr\rho_0(\frob_p)$ otherwise.

We let $U_0$ denote the set of $\smat{a}{b}{c}{d}$ in $U^\Sigma$
such that $c_q \equiv 0 \bmod q$ and $d_r \equiv 1 \bmod r$, and
let $U_1$ denote the set of $\smat{a}{b}{c}{d}$ in $U_0$
such that $c_q$ has trivial image in $\Delta_{\{q\}}$.
We identify $\Delta_{\{q\}}$ with $U_0/U_1$.
Let $\sigma_i$ denote the restriction of
$\sigma^\Sigma$ to $U_i$ for $i = 0, 1$.
Let $\bar{\TT}_i$ denote the image of $\tilde{\TT}$
in $\End_K M(\sigma_i)_{!,\lambda_0}$ for $i = 0,1$.

As in \cite[\S4.2]{ddt} (see also \S6.4 of \cite{cdt}),
one can give an explicit description of
$\bar{K}\otimes_\CO\bar{\TT}_{0,\gm}$ which shows that
its dimension over $\bar{K}$ is $\#\CN^\Sigma$.  Similarly
one finds the dimension of $\bar{K}\otimes_\CO\bar{\TT}_{1,\gm}$
is $\#\CN^{\Sigma\cup\{q\}}$.
Since $M(\sigma_i)_{!,\lambda_0}$ is free of rank two over
$K\otimes_\CO\bar{\TT}_i$ (see \S \ref{sec:ep}), we need
to prove that
$$\rank_\CO\CM(\sigma_1)_{!,\gm} =
\#\Delta_{\{q\}}\cdot\rank_\CO\CM(\sigma_0)_{!,\gm}.$$
We shall show that $\CM(\sigma_1)_{!,\gm}$ is
free over $\CO[\Delta_{\{q\}}]$ and its $\Delta_{\{q\}}$-invariants
can be identified with $\CM(\sigma_0)_{!,\gm}$.
The case $k=2$, trivial $\psi$ and $\Sigma_2 = \emptyset$
is done exactly as in \S 2 of \cite{tw}, so assume
$k > 2$, $\psi$ is not trivial or $\Sigma_2$ is not empty.
It follows from the
Hochschild-Serre spectral sequence that
$$H^i(\Delta_Q, \CM(\sigma_1)_c) = \left\{\begin{array}{ll}
\CM(\sigma_0)_c, & \mbox{if $i=0$;}\\
0, & \mbox{if $i>0$.}\end{array}\right.$$
This holds as $\tilde{\TT}$-modules, so it holds after completing
at $\gm$.  We can then replace $!$ by $c$ and deduce the freeness
of $\CM(\sigma_1)_{!,\gm}$ from cohomological vanishing.
\epf

\subsection{Consequences}
\mlabel{ss:conseq}
\begin{theorem} Suppose $\rho: G_\QQ \to \aut_K V$
is a continuous geometric representation whose
restriction to $G_\ell$ is ramified, crystalline and short.
If $\rho_0$ is modular and its  restriction to $G_F$
is absolutely irreducible, where $F$ is the quadratic
subfield of $\QQ(\mu_\ell)$, then $\rho$ is modular.
\mlabel{thm:modular}
\end{theorem}
This follows from applying Theorem \ref{thm:axiomatic} to
the set $\CN$ constructed in the preceding section for the twist
$\rho_0\otimes_k \psi'$ of minimal conductor, where $\psi'$ is
unramified at $\ell$.   Writing $\tilde{\ }$ for Teichmuller
liftings, we conclude that $\rho\otimes_K\tilde{\psi}'\psi$
is modular, where $\psi$ is a character of $\ell$-power
order such that $\chi_\ell^{1-k}\psi^2\det\rho$ has order
not divisible by $\ell$.

Suppose that $K_0$ is a number field with ring of integers $\CO_0$,
and $f = \sum a_n e^{2\pi i nz}$ is a newform of weight $k\ge 2$ and level
$N$ with coefficients $a_n$ in $K_0$.
For primes $\lambda$ of $K_0$, we write $\barrho_{f,\lambda}$
for the representation of $G_\QQ$ on
$$\bar{\kappa} \otimes_{\CO}\CM_{f,\lambda},$$
where $\CO = \CO_{0,\lambda}$ and $\bar{\kappa}$ is an algebraic closure
of $\CO_0/\lambda$.
We use $\ell$ to denote the rational prime divisible by $\lambda$.
Let $S_f$ denote the set of primes $\lambda$ of $K_0$ such that the
following hold:
\begin{itemize}
\item $\lambda$ does not divide $Nk!$;
\item $\barrho_{f,\lambda}$ has irreducible restriction to $G_F$,
 where $F$ is the quadratic subfield of $\QQ(\mu_\ell)$.
\end{itemize}
The latter condition excludes only finitely many primes by the
following two lemmas:
\begin{lemma} For all but finitely many primes $\lambda$,
$\barrho_{f,\lambda}$ is irreducible.
\end{lemma}
\proofbegin Suppose that $\lambda$ does not divide $Nk!$ and $\barrho_{f,\lambda}$
is reducible.  Its semisimplification is of the form $\psi_1 \oplus \psi_2$
where $\psi_1$ and $\psi_2$ are characters of $\gal(\QQ(\mu_{N\ell})/\QQ)$.
The representation is necessarily ordinary at $\ell$ (see \cite{edixhoven}),
so one of the characters is unramified at $\ell$ and the other has restriction
$\chi_\ell^{k-1}$ on $I_\ell$.  It follows that
$$a_p \equiv p^{k-1} + 1 \bmod \lambda$$
for all $p\equiv 1\bmod N$.  If this holds for
infinitely many $\lambda$, then we get $a_p = p^{k-1}+1$
for all such $p$, violating the Ramanujan conjecture
(a theorem of Deligne \cite{del_bour}).
\epf

\begin{lemma} Suppose that $\lambda$ does not divide $N(2k-1)(2k-3)k!$.
If $\barrho_{f,\lambda}$ is irreducible, then so is its restriction to $G_F$.
\end{lemma}
\proofbegin Consider the restriction of $\barrho_{f,\lambda}$ to $I_\ell$.
By results of Deligne and Fontaine (see \cite{edixhoven}), this restriction
has semisimplification of the form $\chi_\ell^{k-1}\oplus 1$ or
$\psi_\ell^{k-1} \oplus \psi_\ell^{\ell(k-1)}$, where $\psi_\ell$ is
a fundamental character of level 2, according to whether or not
$a_\ell$ is a unit mod $\lambda$.

Suppose that $\barrho_{f,\lambda}$ is irreducible but its restriction to
$G_F$ is not.  Then $\barrho_{f,\lambda}$ is induced from a character of
$G_F$, and $\barrho_{f,\lambda}|_{I_\ell}$ is induced from a character
of its subgroup of index 2.  It follows that the ratio of the characters
into which $\barrho_{f,\lambda}|_{I_\ell}$ decomposes is quadratic.
Since $\psi_\ell$ has order $\ell^2-1$, this forces either
$(\ell-1)|2(k-1)$ or $(\ell+1)|2(k-1)$ and we arrive at a contradiction.
\epf

\begin{theorem} Let $f$ be a newform of weight $k\ge2$ and level $N$
with coefficients in $K_0$.  Suppose that $\lambda$ is a prime in $S_f$
(as defined above), and let $\CO = \CO_{0,\lambda}$ and $K = K_{0,\lambda}$.
Suppose that $\Sigma$ is a finite set of primes not containing $\ell$
such that $M_{f,\lambda}$ is minimally ramified outside
$\Sigma$.  Then the $\CO$-module
$$H^1_\Sigma(G_\QQ, K/\CO\otimes_{\CO}\ad^0_\CO\CM_{f,\lambda})$$
has length $v_\lambda(\eta_f^\Sigma)$.
\mlabel{thm:selmer}
\end{theorem}
Applying Theorem \ref{thm:axiomatic} to
the set $\CN$ constructed in the preceding section for the twist
$\rho_0\otimes_k \psi'$ of minimal conductor, we conclude that the
theorem holds for a twist of $f$, hence for $f$ itself.  See the
discussion before Proposition \ref{prop:int-struc}
for the definition of $\eta_f^\Sigma$.

\section{The Bloch-Kato conjecture for $A_f$ and $A_f(1)$}
\mlabel{sec:bk}

In this section we explain how to deduce the $\lambda$-part
of the Bloch-Kato conjecture \cite{bloch_kato} for $A_f$ and $B_f =
A_f(1)$, where $f$ is a newform of weight $k \ge 2$, conductor $N\geq 1$, with
coefficients in the number field $K$,  and $\lambda$ is a prime in $S_f$ (the
set defined in \S\ref{ss:conseq}). By Lemma \ref{lem:twist} we can assume that $f$ has
minimal conductor among its twists and we shall do so in this section. Our formulation of the
conjecture follows Fontaine and Perrin-Riou \cite{fon_pr},
generalized to motives with coefficients in $K$. For a more systematic discussion of
the Bloch-Kato conjecture for motives with coefficients we refer to \cite{burns_flach}.
If $A_f$ is the scalar extension of a premotivic structure with coefficients in a subfield
$K'\subseteq K$ then Theorem 4.1 and Lemma 9 of \cite{burns_flach} shows that the conjecture over $K$ implies
the one over $K'$ (in the context of Deligne's conjecture this was already noted in \cite[Rem 2.10]{del_cor}).
So we need not be concerned with finding the smallest coefficient field for $A_f$.

\subsection{Order of vanishing}
\mlabel{sec:ord}
Suppose that $M$ is an $L$-admissible object of $\pms_K$ and
let $M^\kd=\hom_K(M,K(1))$. We recall the conjectured order of
vanishing of $L(M,s)$ at $s=0$ \cite[III. 4.2.2]{fon_pr}.
\begin{conjecture}
Let $\tau:K \to \CC$ be an embedding and $\lambda$ any finite
prime of $K$. Then
$$\ord_{s=0} L(M,\tau,s) =
\dim_{K_\lambda} H^1_\f(\QQ,M^\kd_\lambda) - \dim_{K_\lambda}
H^0(\QQ,M^\kd_\lambda).$$
\mlabel{conj:ord}
\end{conjecture}

\begin{theorem}  Conjecture \ref{conj:ord} holds
for both $M=A_f$ and $M=B_f$ if $\lambda$ is in $S_f$. More
precisely, we have $\ord_{s=0} L(A_f,\tau,s)=\ord_{s=0}
L(B_f,\tau,s)=0$ and
\begin{equation}H^0(\QQ,A_{f,\lambda})\cong
H^1_\f(\QQ,A_{f,\lambda})\cong
H^1_\f(\QQ,B_{f,\lambda})\cong
H^0(\QQ,B_{f,\lambda})\cong
0.\notag
\end{equation}
if $\lambda\in S_f$.
\mlabel{thm:ordl}
\end{theorem}
\proofbegin Lemma \ref{hida} below shows that
$$L(A_f,\tau,1)=L^\naive(A_f,\tau,1)
\prod_{p\in\Sigma_e(f)}L_p(A_f,\tau,1)$$
is a nonzero multiple of the Petersson inner product of
$f$ with itself and hence it follows that
$L(B_f,\tau,0) = L(A_f,\tau, 1) \neq 0$
for each $\tau$. It follows from the functional
equation (\ref{ssec:fe}) that
\begin{equation}L(A_f,\tau,0) = \frac{(k-1)\epsilon(A_f)}{2\pi^2}
L(A_f,\tau,1) \neq 0
\mlabel{funatone}\end{equation}
for each $\tau$ as well.  The absolute irreducibility of $M_{f,\lambda}$ for
each $\lambda$ implies that $\End_{K_\lambda[G_\QQ]} (M_{f,\lambda}) =
K_\lambda$,
so $H^0(\QQ,A_{f,\lambda}) = 0$, and since $M_{f,\lambda}$ is not
isomorphic to $M_{f,\lambda}(1)$, we also have
$H^0(\QQ,B_{f,\lambda})=0$. It follows then from
\cite[II.2.2.2]{fon_pr} (see also \cite[Cor. 1.5]{flachthesis}) that
$$\dim_{K_\lambda} H^1_\f(\QQ,A_{f,\lambda}) =
        \dim_{K_\lambda} H^1_\f(\QQ,B_{f,\lambda})$$
for all $\lambda$ and hence that Theorem \ref{thm:ordl}
is implied by the vanishing of $H^1_\f(\QQ,A_{f,\lambda})$. Theorem
\ref{thm:selmer} shows that
$$H^1_\f(\QQ,
A_{f,\lambda}/\ad^0_{\CO_{K,\lambda}}\CM_{f,\lambda})\subset
H^1_\Sigma(\QQ,
A_{f,\lambda}/\ad^0_{\CO_{K,\lambda}}\CM_{f,\lambda})$$ is finite
for $\lambda$ in $S_f$.  Since the kernel of
$$H^1_\f(\QQ,A_{f,\lambda}) \to
H^1_\f(\QQ, A_{f,\lambda}/\ad^0_{\CO_{K,\lambda}}\CM_{f,\lambda})$$
is finitely generated over $\CO_\lambda$ we deduce
$H^1_\f(\QQ,A_{f,\lambda})=0$ and Theorem \ref{thm:ordl} follows.
\epf

\subsection{Deligne's period}
We now recall the formulation in \cite{fon_pr} of Deligne's
conjecture \cite{del_cor} for the ``transcendental part'' of
$L(M,0)$ for $M = A_f$ or $B_f$.
The authors there actually discuss the more general conjecture of
Beilinson concerning the leading coefficient $L^*(M,0)$ for
premotivic structures arising from motives, but their formulation
relies on the conjectural existence of a category of mixed motives
with certain properties. We restrict our attention to those $M$,
such as $A_f$ and $B_f$, for which $L(M,0) \neq 0$ and which are
critical in the sense of Deligne. In that case Beilinson's
conjecture reduces (conjecturally) to Deligne's, which can be stated
without reference to the category of mixed motives.

Under these hypotheses, the {\em fundamental line} for $M$ is the
$K$-line defined by
$$\Delta_\f(M) = \hom_K(\det_K M_B^+, \det_K t_M)$$
where ${\ }^+$ indicates the subspace fixed by $F_\infty$
and $t_M = M_\dr/\fil^0 M_\dr$.  Furthermore the composite
$$\RR \otimes M_B^+
\to (\CC\otimes M_B)^+
\stackrel{(I^\infty)^{-1}}{\longrightarrow} \RR\otimes M_{\dr}
\to \RR\otimes t_M$$
is an $\RR\otimes K$-linear isomorphism.  Its determinant
over $\RR\otimes K$ defines a basis for
$\RR \otimes \Delta_\f(M)$ called the Deligne
period, denoted $c^+(M)$.
\begin{conjecture}
There exists a basis $b(M)$ for $\Delta_\f(M)$ such that
$$L(M,0)(1\otimes b(M)) = c^+(M).$$
\mlabel{conj:del}
\end{conjecture}
For example, the conjecture holds for $M = \QQ(2)$ with the
basis $b(\QQ(2))$ of
$$\Delta_\f(\QQ(2)) = \hom_\QQ(\QQ(2\pi i)^2, \QQ[2])$$
defined by $(2\pi i)^2 \mapsto -24\iota^{-2}$. There are
various rationality results for $L(A_f,0)$ and $L(B_f,0)$ in the
literature (see for example \cite[Th. 2.3]{schmidt}) although
the precise relationship with Conjecture \ref{conj:del} for $M =
A_f$ or $B_f$ is not always clear. In this section we recall the
proof of Conjecture \ref{conj:del} for $M = A_f$ and $B_f$ and give
convenient natural descriptions for $b(A_f)$ and $b(B_f)$.

We begin by observing that $A_{f,B}^+$ and $t_{A_f}$
are one-dimensional over $K$.  Furthermore, complex conjugation
$$F_\infty: M_{f,B} \to M_{f,B}$$
has trace zero and commutes with $F_\infty$, so it is
a basis for $A_{f,B}^+$.  Note also that the natural map
$$A_{f,\dr} \to
         \hom_K(\fil^{k-1}M_{f,\dr}, M_{f,\dr}/\fil^{k-1}M_{f,\dr})$$
factors through an isomorphism
$$t_{A_f} \to \hom_K(\fil^{k-1}M_{f,\dr}, M_{f,\dr}/\fil^{k-1}M_{f,\dr}).$$
The fundamental line $\Delta_\f(A_f)$ can therefore be
identified with
$$\hom_K(\fil^{k-1}M_{f,\dr}\otimes \QQ\cdot F_\infty,
                M_{f,\dr}/\fil^{k-1}M_{f,\dr}).$$
We shall describe $b(A_f)$ by specifying the image of the canonical basis
$f\otimes F_\infty$ for $\fil^{k-1}M_{f,\dr}\otimes \QQ\cdot F_\infty$ where
we view $f$ as an element of $M_{f,\dr}$ by Lemma \ref{lem:g}.
Recall that we defined in (\ref{eq:tpair2}) a perfect alternating pairing
$$\langle \cdot,\cdot \rangle\ :\ M_f \otimes_K M_f \to M_{\psi}(1-k),$$
and this induces an isomorphism
$$M_{f,\dr}/\fil^{k-1}M_{f,\dr}
 \to \hom_K(\fil^{k-1}M_{f,\dr},M_{\psi}(1-k)_\dr).$$
We shall eventually define $b(A_f)$ by specifying the
element $\langle f,b(A_f)(f\otimes F_\infty) \rangle$
of $M_{\psi}(1-k)_\dr$.

We can make a similar analysis of the fundamental line $\Delta_\f(B_f)$.
One finds that $B_{f,B}^+$ and $t_{B_f}$ are two-dimensional over $K$.
Note that $B_{f,B}^+$ can be identified with $A_{f,B}^-\otimes \QQ(1)_B$
and that the natural map
$$A_{f,B}^- \to \hom_K(M_{f,B}^+,M_{f,B}^-) \oplus
                \hom_K(M_{f,B}^-,M_{f,B}^+)$$
defined by restrictions is an isomorphism.
We therefore have an isomorphism
$$\det_K B_{f,B}^+ \to K(2)_B$$
which is canonical up to sign.  To fix the choice of
sign, we use $\alpha \wedge \alpha^{-1}$ as a basis for
$\det_K A_{f,B}^-$ where $\alpha : M_{f,B}^+ \to M_{f,B}^-$
is any $K$-linear isomorphism.  Next note that the natural map
$$B_{f,\dr} \to \hom_K(M_{f,\dr},M_f(1)_\dr)
        \to \hom_K(\fil^{k-1}M_{f,\dr}, M_f(1)_\dr)$$
factors through an isomorphism
$$t_{B_f} \to \hom_K(\fil^{k-1}M_{f,\dr}, M_f(1)_\dr).$$
Using the isomorphism
$$\det_K M_{f,\dr} \to \fil^{k-1}M_{f,\dr} \otimes_K
                (M_{f,\dr}/\fil^{k-1}M_{f,\dr})$$
(with choice of sign again indicated by the ordering),
we find that $\det_K t_{B_f}$ is naturally isomorphic to
\begin{equation}\hom_K(\fil^{k-1}M_{f,\dr},
M_{f,\dr}/\fil^{k-1}M_{f,\dr})
\otimes \QQ(2)_\dr.\mlabel{tanb}\end{equation}
We can therefore identify $\Delta_\f(B_f)$ with
$$\hom_K(\fil^{k-1}M_{f,\dr}\otimes \QQ(2)_B,
                (M_{f,\dr}/\fil^{k-1}M_{f,\dr})\otimes
\QQ(2)_\dr),$$ and we arrive at a canonical isomorphism
\begin{equation}\Delta_\f(A_f) \otimes \Delta_\f(\QQ(2)) \otimes
A_{f,B}^+\xrightarrow{\sim}\Delta_\f(B_f).
\notag\end{equation}
Fixing the basis $F_\infty$ of $A_{f,B}^+$ and the basis $\beta$
of $\Delta_\f(\QQ(2))$ which sends $(2\pi i)^2$ to $\iota^{-2}$ this
defines an isomorphism of $K$-lines
\begin{equation}
\tw:\Delta_\f(A_f)\rightarrow \Delta_\f(B_f)
\mlabel{twist}\end{equation}
so that $\tw(\phi)(x\otimes y)=\phi(x\otimes F_\infty)\otimes
\beta(y)$.

\begin{lemma} We have
$$(\RR\otimes \tw)(c^+(A_f)) = -\frac{1}{2\pi^2}c^+(B_f).$$
\mlabel{tlemma}\end{lemma}

\proofbegin Let $I^\infty_M:\CC\otimes M_{f,\dr}\cong \CC\otimes M_{f,B}$
be the comparison isomorphism for $M_f$. Via the natural
isomorphism $\CC\otimes\text{End}_K(M_f)_?\cong
\text{End}_{\CC\otimes K}(\CC\otimes M_{f,?})$
where $?=B$ or $?=\dr$, $I^\infty_M$ induces the comparison
isomorphism $I^\infty$ for both $\text{End}(M_f)$ and $A_f$:
$I^\infty(\phi)= I_M^\infty\circ\phi\circ (I_M^\infty)^{-1}$. A similar formula
holds for $c^+(A_f)$.

Suppose now that $x$ is a non-zero element of
$\RR \otimes \fil^{k-1} M_{f,\dr}$ (for example, take $x=g$),
and write $I_M^\infty(x) = y^+ + y^-$ with $y^\pm \in
\CC\otimes M_{f,B}^{\pm}$.  Then
$$c^+(A_f)(x\otimes F_\infty)
    = (I_M^\infty)^{-1}(1\otimes F_\infty)I_M^\infty(x)
        = (I^\infty_M)^{-1}(y^+ - y^-) \mod \RR\otimes\fil^{k-1}M_{f,\dr}.$$
Chasing through the above isomorphisms, one finds that $c^+(B_f)$
is characterized by
\begin{align*}
x \otimes c^+(B_f)(x \otimes (2\pi i)^2) \otimes \iota^2
& =
(2\pi i)^2 (I^\infty_M)^{-1}(y^-) \wedge (I^\infty_M)^{-1}(y^+)\\
&= x\otimes\frac{1}{2}(2\pi i)^2 c^+(A_f)(x\otimes F_\infty).
\end{align*}
\epf

Recall that $\Sigma_e(f)$ is the set of primes
$p$ such that $L_p^\naive(A_f,s) = 1$ but $L_p(A_f,s) = (1+p^{-s})^{-1}$.
We write $b_{\dr}$ for the basis of $M_{\psi,\dr}$ defined in
\S \ref{sssec:mot.ex.dir}, and set $\eta = 0$ or $1$ so $\eta \equiv k
\bmod 2$.  Note that by Proposition 5.5 of \cite{del_cor},
we have $\epsilon(M_f\otimes M_{\psi^{-1}})/\epsilon(M_{\psi^{-1}}) \in K^\times$. The same
proposition together with Lemma~\ref{lem:adjoint} gives $\epsilon(A_f)\in K^\times$.

\begin{theorem}
Let $b(A_f)\in \Delta_\f(A_f)$ be defined by the formula
$$\langle f, b(A_f)(f\otimes F_\infty) \rangle =
        \frac{i^{k-\eta}((k-2)!)^2\epsilon(M_f\otimes M_{\psi^{-1}})}{2\epsilon(M_{\psi^{-1}})\epsilon(A_f)}\prod_{p\in\Sigma_e(f)}
        (1+p^{-1})
        \cdot(b_\dr\otimes \iota^{k-1}),$$
and $b(B_f) \in \Delta_\f(B_f)$ by the formula
\begin{equation} b(B_f)= (1-k)\epsilon(A_f) \tw(b(A_f)).
\mlabel{funident2}\end{equation}
Then $L(A_f,0)(1\otimes b(A_f)) = c^+(A_f)$ and
$L(B_f,0)(1\otimes b(B_f)) = c^+(B_f)$.
\mlabel{thm:del}
\end{theorem}

\proofbegin If we show
$$
\langle f, c^+(A_f)(f\otimes F_\infty)\rangle
        = \frac{i^{k-\eta}(k-1)!(k-2)! \epsilon(M_f\otimes M_{\psi^{-1}})
L^\naive(A_f,1)}{4\pi^2
        \epsilon(M_{\psi^{-1}})}\cdot(b_\dr\otimes \iota^{k-1})
$$
in $\CC\otimes M_{\psi}(1-k)_\dr$, then the statement concerning $b(A_f)$ is an immediate consequence of
the functional equation (\ref{funatone}). The identity
$L(B_f,0)(1\otimes b(B_f)) = c^+(B_f)$ then follows by applying
$(\RR\otimes\tw)$ to the identity $L(A_f,0)(1\otimes b(A_f)) =
c^+(A_f)$ and using  (\ref{funatone}) and Lemma \ref{tlemma}.

As in \S \ref{ssec:levchar} put $U=U_0(N)$, let $\sigma:U\rightarrow K^\times$ be the representation
$\left(\begin{smallmatrix}a&b\\c&d\end{smallmatrix}\right)
\mapsto \psi^{-1}(a_N)$, and choose a multiple $N'$ of
$N$ so that $U_{N'}\subseteq U$ is sufficiently small.
Put $w=\left(\begin{smallmatrix} 0 &-1\\N &
0\end{smallmatrix}\right)
\in GL_2(\AA_\f)$ and denote by
$W=[UwU]_\omega:M(\sigma)_{N'}\rightarrow M(\hat{\sigma}\otimes
(\psi^{-1}\circ\det))_{N'}$ the isomorphism in Lemma~\ref{lem:wiso}. Note
that $ww_N^{-1}\in U$ so we can work with $w$ instead of $w_N$. For any one-dimensional $K$-representation
$\sigma$ of $U$ we shall view $M(\sigma)_{N'}$ as a sub-PM-structure of $K\otimes M_{N'}$.
With $I^\infty$ denoting the comparison isomorphism for both $M_f$ and $K\otimes M_{N'}$
we have
\begin{equation}
\begin{array}{ll}
\langle f, c^+(A_f)(f\otimes F_\infty)\rangle
&= \langle f, (I^\infty)^{-1}(1\otimes F_\infty) I^\infty f\rangle
\\
&=[ f, (I^\infty)^{-1}(1\otimes F_\infty) I^\infty Wf ]_{N'} \notag\\
&=[U:U_{N'}]^{-1} ( f, (I^\infty)^{-1}(1\otimes F_\infty) I^\infty
Wf )_{N'}.
\end{array}
\mlabel{comp1}
\end{equation}
We proceed with the computation of $Wf\in\fil^{k-1}M_!(\hat{\sigma}
\otimes (\psi^{-1}\circ\det))_{N',\dr}$. Note that the field $K_f$
generated by the Fourier coefficients of the newform $f$ is
either totally real or a CM field and hence has a well defined
automorphism $\rho$ induced by complex conjugation. It is known
that the Fourier expansion
$f^\rho(z)=\sum_{n=1}^\infty a_n^\rho e^{2\pi iz n}$ is a newform of conductor
$N$ and character $\psi^{-1}$ \cite[4.6.15(2)]{miyake},
hence represents an element of $\fil^{k-1}M_!(\hat{\sigma})_{N',\dr}$.
Recall that $b_\dr$ is the canonical
basis of $M_{\psi,\dr}\cong(H^0(X_{N'})\otimes K_{\psi^{-1}})^U_\dr$
(Lemma~\ref{lem:psi})
and that $M_{\psi}\otimes_K M_{f^\rho}$ has a natural map
into $M_!(\hat{\sigma}\otimes (\psi^{-1}\circ\det))_{N'}$ via the cup
product on $X_{N'}$ (see (\ref{eq:cup})).

\begin{lemma} We have
\begin{equation} Wf =
\psi(-1)\frac{i^{k-\eta}\epsilon(M_f\otimes M_{\psi^{-1}})}{N\epsilon(M_{\psi^{-1}})}
b_\dr\cup f^\rho.
\mlabel{wid}\end{equation}
\mlabel{atkinlehner}\end{lemma}
\proofbegin We fix an embedding $\tau:K\rightarrow \CC$ and compute the
images of both sides in $S_k(U_{N'})$ (we shall suppress $\tau$ in the notation and view all elements of $K$
as complex numbers via $\tau$). Let $\phi\in (S_k(U_{N'})\otimes_\CC\CC_{\sigma})^U$ denote
the element corresponding to $f$. We then have $\phi(xu)=\sigma^{-1}(u)\phi(x)$
for all $u\in U$ and $\beta(\phi)_t(z)=f(z)$ for all $t\in(\ZZ/N'\ZZ)^\times$.
Recall that the isomorphism  $\beta: S_k(U_{N'})\cong
\bigoplus_{t\in(\ZZ/N'\ZZ)^\times}S_k(\Gamma(N'))$ was
defined in (\ref{eq:adanform}) by .
\[ \beta(F)_t(\gamma(i)) :=  (\det \gamma)^{-1}j(\gamma,i)^k F(g_t\gamma(i))).\]
for $\gamma\in GL_2(\RR)^+$,
$j(\left(\begin{smallmatrix}a&b\\c&d\end{smallmatrix}\right),z)=
cz+d$ and $g_t\equiv\left(\begin{smallmatrix}1 & 0\\0 &t^{-1}
\end{smallmatrix}\right)\mod N' $.
Note that since $g_t\in U$ we indeed have that
\begin{align*} \beta(\phi)_t(z) =
&(\det \gamma)^{-1}j(\gamma,i)^k\phi(g_t\gamma(i))\\
= &(\det \gamma)^{-1}j(\gamma,i)^k\sigma^{-1}(g_t)\phi(\gamma(i))\\
= &(\det \gamma)^{-1}j(\gamma,i)^k\phi(\gamma(i))\\
= & f(z)
\end{align*}
is independent of $t$.

Note that the element $e^{2\pi i/N'}$ of $F=\QQ(e^{2\pi i/N'})$ maps to $e^{-2\pi i/N'}$ under the isomorphism $M_F\cong H^0(X_{N'})$ of \ref{lem:field}.
Hence $b_\dr=\sum_{a\in(\ZZ/N\ZZ)^\times}\psi(a)\otimes e^{2\pi i a/N}\in
\CC\otimes \QQ(e^{2\pi i/N'})$, when viewed as an element of $H^0_\dr(X_{N'})$, i.e. a locally constant function, is given by
\begin{align*}
t\mapsto\sum_{a\in(\ZZ/N\ZZ)^\times}\psi(a)e^{-2\pi i at/N}
& =\psi(-t)^{-1}\sum_{a\in(\ZZ/N\ZZ)^\times}\psi(a)e^{2\pi i
a/N}\\
&=\psi(-t)^{-1}G_\psi\\
&=\psi(-t)^{-1} i^\eta\epsilon(M_{\psi^{-1}},\tau).
\end{align*}
If now $\phi^\rho\in S_k(U_{N'})$ corresponds to $f^\rho$ then $\beta(\phi^\rho)_t(z)=f^\rho(z)$
is again independent of $t$ and the right hand side of (\ref{wid}) is given by
\begin{equation}
t\mapsto \frac{i^k\epsilon(M_f\otimes M_{\psi^{-1}},\tau)}{N} \psi(t)^{-1}
f^\rho.
\mlabel{comp11}\end{equation}

The perfect pairing $M_f\otimes_K M_f\rightarrow M_{\psi}(1-k)$
and the identity of Hecke eigenvalues \cite[(4.6.17)]{miyake}
induce an isomorphism $M_f^*\cong M_f\otimes_K M_{\psi^{-1}}(k-1)
\cong M_{f^\rho}(k-1)$ so that the functional equation for
$\Lambda(M_f\otimes M_{\psi^{-1}},\tau,s)$ can be written
\begin{equation}
\Lambda(M_f\otimes M_{\psi^{-1}},\tau,s)=\epsilon(M_f\otimes M_{\psi^{-1}},\tau)
N^{-s}\Lambda(M_{f^\rho}\otimes M_{\psi},\tau,k-s).
\mlabel{fe11}\end{equation}
Recall that the definition $(g\vert_k\gamma)(z)=
\det(\gamma)^{k/2}j(\gamma,z)^{-k}g(\gamma(z))$
for $\gamma\in GL_2(\RR)^+$ defines a right action of
$GL_2(\RR)^+$ on functions $g:\uhp\rightarrow\CC$. Put $W_N=
\left(\begin{smallmatrix} 0 &-1\\N & 0\end{smallmatrix}\right)$.
By \cite[Th. 4.3.6]{miyake} we have
\[\Lambda(M_f\otimes M_{\psi^{-1}},\tau,s)=\Lambda(f,s)=i^k
N^{-s+k/2}\Lambda(f\vert_kW_N,k-s)\]
which together with (\ref{fe11}) yields
\[ f^\rho=\epsilon(M_f\otimes M_{\psi^{-1}},\tau)^{-1}i^k N^{k/2} f\vert_kW_N.\]
Hence (\ref{comp11}) becomes
\begin{equation}
t\mapsto (-1)^k N^{k/2-1} \psi(t)^{-1} f\vert_kW_N.
\mlabel{comp12}\end{equation}

Turning to the left hand side of (\ref{wid}) we have
$(W\phi)(x):=\phi(xw)$ and $\phi(wh)=\phi(W_\QQ W_N^{-1}h)=\phi(W_N^{-1}h)$
where $h\in GL_2(\RR)^+$ and $W_\QQ\in GL_2(\QQ)$ is the matrix
with image $w$ (resp.
$W_N$) in $GL_2(\AA_\f)$ (resp. $GL_2(\RR)$). For
$\gamma\in GL_2(\RR)^+$ we have
\begin{align} \det(h)^{-1}j(h,i)^k\phi(\gamma h) = &\det(\gamma)
\det(\gamma h)^{-1}j(\gamma,h(i))^{-k}j(\gamma h,i)^k\phi(\gamma
h)\notag\\
 = & \det(\gamma)j(\gamma,h(i))^{-k} f(\gamma h(i))\notag\\
= & \det(\gamma)^{1-k/2} (f\vert_k\gamma)(h(i)).
\notag\end{align}
Combining these equations we find that $W\phi$ corresponds to
\begin{align*} t\mapsto \det(h)^{-1}j(h,i)^k(W\phi)(g_th) =
&\det(h)^{-1}j(h,i)^k \phi(ww^{-1}g_twh)\\=&\det(h)^{-1}j(h,i)^k
\sigma^{-1}(w^{-1}g_tw)\phi(wh)\\=&\det(h)^{-1}j(h,i)^k
\sigma^{-1}(w^{-1}g_tw) \phi(W_N^{-1}h)
\\=&\sigma^{-1}(w^{-1}g_tw)
\det(W_N^{-1})^{1-k/2}(f\vert_k W_N^{-1})(h(i)).
\end{align*}
Since $f\vert_k W_N^2=(-1)^kf$ this last expression equals
\begin{equation}\sigma^{-1}(w^{-1}g_tw)
(-1)^k N^{k/2-1}(f\vert_k W_N)(h(i)).\mlabel{comp13}\end{equation}
For $g_t\equiv\left(\begin{smallmatrix}1 & 0\\0 &t^{-1}
\end{smallmatrix}\right)\mod N $ we have $w^{-1}g_tw\equiv
\left(\begin{smallmatrix}t^{-1} & *\\0 & 1
\end{smallmatrix}\right)\mod N$ and
$\sigma^{-1}(w^{-1}g_tw)=\psi(t^{-1})=\psi(t)^{-1}$. So
(\ref{comp12}) and (\ref{comp13}) agree which finishes the proof of
the Lemma.
\epf

The definition (\ref{eq:spair1a}) of the pairing on $\sigma$-constructions shows that
$(x,\alpha\cup y)_{N'}=(x,y)\otimes_K\alpha$ where $\alpha\in M_{\psi}$ and $(x,y)$ is the
$K$-linear extension of the $\QQ(1-k)$-valued pairing on $\CM_{N'}$ defined in \S \ref{sssec:pairings}.
Combining this
with Lemma \ref{atkinlehner} the last term in (\ref{comp1}) equals
\begin{equation}[U:U_{N'}]^{-1} ( f, (I^\infty)^{-1}
    (1\otimes F_\infty)I^\infty f^\rho  )\otimes_K\alpha_{\dr}
\mlabel{comp2}\end{equation}
in $\CC\otimes K(1-r)_{\dr}\otimes_K M_{\psi,\dr}$ where
\begin{align*}
\alpha_\dr & =\psi(-1)\frac{i^{k-\eta}\epsilon(M_f\otimes M_{\psi^{-1}})}{N\epsilon(M_{\psi^{-1}})}
(I^\infty)^{-1}(1\otimes F_\infty)I^\infty b_\dr\\
& =\psi(-1)\frac{i^{k-\eta}\epsilon(M_f\otimes M_{\psi^{-1}})}{N\epsilon(M_{\psi^{-1}})}
\psi(-1)^{-1} b_\dr
\end{align*}
(with $I^\infty$ also denoting the comparison isomorphism for $M_{\psi}$).
For any premotivic structure we have $(F_\infty\otimes F_\infty)I^\infty=I^\infty(F_\infty\otimes 1)$ and we have $(F_\infty\otimes 1)(f^\rho)=f^\rho$
since $f^\rho\in K\otimes M_{N',\dr}\subset \CC\otimes K\otimes M_{N',\dr}$. Hence
\[(I^\infty)^{-1}(1\otimes F_\infty)I^\infty f^\rho=(I^\infty)^{-1}(F_\infty\otimes 1)I^\infty f^\rho.\]
Under the natural isomorphism $\CC\otimes K\otimes M_{N',B}\cong (\CC\otimes M_{N',B})^{\II_K}$ the action of
$F_\infty\otimes 1\otimes 1$ on the left hand side gets transformed into the action sending $(x_\tau)$
to $\tau\mapsto (F_{\infty,\tau}\otimes 1)(x_{\overline{\tau}})$ where $F_{\infty,\tau}$ is complex conjugation acting on
$\CC$ in the factor indexed by $\tau$. Hence the $\tau$-component of (\ref{comp2}) equals
\begin{align*}&[U:U_{N'}]^{-1} ( \tau(f), (I^\infty)^{-1}
    (F_{\infty,\tau}\otimes 1)I^\infty \overline{\tau}(f^\rho)  )\otimes_\CC\tau(\alpha_{\dr})\\
=&[U:U_{N'}]^{-1}(k-2)!(4\pi)^{k-1}\phi(N')(\tau(f),\tau(f))_{\Gamma(N')}\,\tau(\alpha_\dr)\otimes\iota^{k-1}
\end{align*}
where $\phi$ is Euler's function and we have used Lemma \ref{petersson1}. Therefore (\ref{comp2}) equals
\begin{equation}
\begin{array}{l}
[U:U_{N'}]^{-1}(k-2)!(4\pi)^{k-1}\phi(N')(f,f)_{\Gamma(N')}
\cdot\alpha_\dr\otimes\iota^{k-1}\\
=
\frac{[\bar{\Gamma}_1(N):\bar{\Gamma}(N')]}{[U:U_{N'}]}\phi(N') (k-2)!
(4\pi)^{k-1}(f,f)_{\Gamma_1(N)}
\cdot\alpha_\dr\otimes\iota^{k-1}
\end{array}
\mlabel{comp3}
\end{equation}
in $\CC\otimes M_{\psi}(1-k)_\dr$ where $[\bar{\Gamma}_1(N):\bar{\Gamma}(N')]$ is the degree of the
covering
$\Gamma(N')\backslash\uhp\rightarrow\Gamma_1(N)\backslash\uhp$.
Since the maps $\det:U\rightarrow(\ZZ/N'\ZZ)^\times$ and $\text{SL}_2(\ZZ)
\rightarrow\text{SL}_2(\ZZ/N'\ZZ)$ are surjective one finds
\begin{equation}
\begin{array}{ll}
[U:U_{N'}]&=\phi(N')[\text{SL}_2(\ZZ)\cap U:\text{SL}_2(\ZZ)\cap
U_{N'}]\\
&=\phi(N')
[\Gamma_0(N):\Gamma(N')]\\
&=\phi(N')\phi(N)[\Gamma_1(N):\Gamma(N')]\\
&=\phi(N')\phi(N)\delta(N)[\bar{\Gamma}_1(N):\bar{\Gamma}(N')]
\end{array}
\mlabel{deltaN}
\end{equation}
where $\delta(N)=1$ if $N>2$ and $\delta(N)=2$ if $N\leq 2$ (note that $-1\in\Gamma_1(N)$ iff $N\leq 2$ whereas $-1\notin\Gamma(N')$).
Combining this with Lemma \ref{hida} below we find that
(\ref{comp3}) equals
\begin{align*}&\frac{(k-2)!(4\pi)^{k-1}}{\phi(N)\delta(N)}
\cdot\frac{(k-1)!\delta
(N)N\phi(N)
L^\naive(A_f,1)}{4^k\pi^{k+1}}\cdot
\frac{i^{k-\eta}\epsilon(M_f\otimes M_{\psi^{-1}})}{N\epsilon(M_{\psi^{-1}})}
\cdot b_\dr\otimes\iota^{k-1}\\
&=\frac{i^{k-\eta}(k-2)!(k-1)!L^\naive(A_f,1)\epsilon(M_f\otimes M_{\psi^{-1}})}
{4\pi^2\epsilon(M_{\psi^{-1}})}\cdot b_\dr\otimes\iota^{k-1}.
\end{align*}
This finishes the proof of Theorem \ref{thm:del}.
\proofend

\begin{lemma} If $f$ is a newform of conductor
$N$, weight $k$ and with coefficients in the number field $K$
we have
\[(\tau(f),\tau(f))_{\Gamma_1(N)} =
\frac{(k-1)!\delta(N)N\phi(N)
L^\naive(A_f,\tau,1)}{4^k\pi^{k+1}}\]
for any embedding $\tau:K\rightarrow\CC$ and
$\delta(N)$ as in (\ref{deltaN}).
\mlabel{hida}\end{lemma}
\proofbegin We fix $\tau$ and write $f$ for $\tau(f)$ to
ease notation. By Theorem 5.1 of
\cite{Hida81} (essentially a reformulation of a Theorem
of Rankin and Shimura) we have
\[L(k,f,\bar{\psi}) =
\frac{4^k\pi^{k+1}(f,f)_{\Gamma_1(N)}}{(k-1)!\delta(N)NN_\psi\phi(N/N_\psi)}
\]
where $L(s,f,\bar{\psi})=\prod_p L_p(s,f,\bar{\psi})$,
\[L_p(s,f,\bar{\psi})^{-1}=(1-\bar{\psi}(p)\alpha_p^2p^{-s})
(1-\bar{\psi}(p)\alpha_p\beta_pp^{-s})(1-\bar{\psi}(p)\beta_p^2p^{-s})\]
and $\alpha_p,\beta_p$ are defined as in \S\ref{ssec:eulerp}
for $p\nmid N$
and $\alpha_p+\beta_p=a_p, \alpha_p\beta_p=0$ for $p\vert N$.
Denote by $M_p$ the exact power of $p$ dividing an integer $M$.
To show the Lemma it suffices to show that
\begin{equation}
L_p(k,f,\bar{\psi})\frac{\phi(N_p/N_{\psi,p})}{N_p/N_{\psi,p}}=
L_p^\naive(A_f,\tau,1)\frac{\phi(N_p)}{N_p}\mlabel{eulereq}
\end{equation}
for all primes $p$. If $p\nmid N$ this is immediate from
\S\ref{ssec:eulerp}. If $N_p=p$ and $N_{\psi,p}=1$ we have
$a_p^2=\psi(p)p^{k-2}$ by \cite[Th. 4.6.17(2)]{miyake} and $\pi_p(f)$
is special so that (\ref{eulereq}) holds true by (\ref{eq:lp1}).
The only other case in which
$a_p\neq 0$ is when $N_p=N_{\psi,p}$ \cite[Th. 4.6.17]{miyake}.
In this case $\bar{\psi}(p)=0$ and hence $L_p(k,f,\bar{\psi})=1$
whereas $\pi_p(f)$ is principal series so that $L_p^\naive(A_f,\tau,1)=(1-p^{-1})^{-1}=N_p/\phi(N_p)$ by
(\ref{eq:lp1}). Finally, if $N_p>1$ and $a_p=0$ then
$L_p(k,f,\bar{\psi})=L_p^\naive(A_f,\tau,1)=1$,
$N_p/N_{\psi,p}>1$ and both sides in (\ref{eulereq}) equal
$(1-p^{-1})$.
\epf

\noindent
{\em Remark.} In the following, we shall not need the full precision
of Theorem \ref{thm:del} but only the fact that $i^{k-\eta}((k-2)!)^2
\epsilon(M_f\otimes M_{\psi^{-1}})/2\epsilon(M_{\psi^{-1}})\epsilon(A_f)$ is a unit in
$\CO = \CO_K[(N k!)^{-1}]$. This in turn is a consequence of Lemma
\ref{lem:epsilon} below.
\bigskip

\begin{lemma} Let $M$ be an object of $\pms_K$ which is $L$-admissible
everywhere and let $\tau:K\rightarrow\CC$ be an embedding. Then $\epsilon(M,\tau)=\epsilon(M,\tau,0)$ is a unit in $\overline{\ZZ}[c(M)^{-1}]$
where $\overline{\ZZ}$ is the ring of algebraic integers.
\mlabel{lem:epsilon}
\end{lemma}
\proofbegin By definition $\epsilon(M,\tau)=\prod_p\epsilon(D_{\pst}(M_\lambda\vert G_p)\otimes_{K_\lambda,\tau'}\CC,\psi_p,dx_p)$
is a product over all places $p$ of $\QQ$ where the additive characters $\psi_p$ and the Haar measures
$dx_p$ are chosen as in \cite[5.3]{del_cor} and $\tau':K_\lambda\rightarrow\CC$ is any extension of $\tau$. The assumption
that $M$ is $L$-admissible at $p$ implies that the isomorphism class of
$D_{\pst}(M_\lambda\vert G_p)\otimes_{K_\lambda,\tau'}\CC$ is independent of $\tau'$. The definition
of $\epsilon$ in \cite[(8.12)]{del_ant} and \cite[Th. 6.5 (a),(b)]{del_ant} show that
\[\epsilon(D_{\pst}(M_\lambda\vert G_p)\otimes_{K_\lambda,\tau'}\CC,\psi_p,dx_p)=
\tau'\epsilon(D_{\pst}(M_\lambda\vert G_p),\psi_p,dx_p)\in\tau'(K_\lambda(\mu_{p^\infty})).\]
Replacing $\tau'$ by $\gamma\tau'$, $\gamma\in\text{Aut}(\CC/K(\mu_{p^\infty}))$, and using the $L$-admissibility again,
we deduce from this formula that $\epsilon(D_{\pst}(M_\lambda\vert G_p)\otimes_{K_\lambda,\tau'}\CC,\psi_p,dx_p)
\in K(\mu_{p^\infty})$. The remark after \cite[(8.12.4)]{del_ant} shows that $\epsilon$ can be
directly expressed in terms of the $\lambda$-adic representation $M_\lambda$ for $\lambda\nmid p$. Namely
\[\epsilon(D_{\pst}(M_\lambda\vert G_p),\psi_p,dx_p)=\epsilon_0((M_\lambda\vert W_p)^{ss},\psi_p,dx_p)
\det(-\Frob\vert M_\lambda^{I_p})^{-1}\]
where $\epsilon_0$ is introduced in \cite[\S 5]{del_ant} and $(M_\lambda\vert W_p)^{ss}$ is the semisimplification
of $M_\lambda$ as a representation of $W_p$. Now for any $\lambda\nmid p$
the $W_p$-representation $M_\lambda$ is the restriction of a continuous $G_p$-representation, hence
carries a $W_p$-stable $\CO_\lambda$-lattice. This implies, on the one hand, that
$\det(-\Frob\vert M_\lambda^{I_p})\in\CO_\lambda^\times$ and on the other hand, via \cite[Th. 6.5(c)]{del_ant},
that $\epsilon_0((M_\lambda\vert W_p)^{ss},\psi_p,dx_p)\in\CO_\lambda[\mu_{p^\infty}]^\times $.
Noting that with our choice of $\psi_p,dx_p$ the epsilon factor equals 1 (resp. a power of $i$) for $p\nmid c(M)$ (resp. $p=\infty$)
the Lemma follows.
\epf

\subsection{Bloch-Kato conjecture}
We now recall the formulation of the $\lambda$-part of the
Bloch-Kato conjecture. We assume that $M$ is a premotivic
structure in $\pms_K$ such that $M$ is critical, $L(M,0) \neq 0$
and Conjecture \ref{conj:del} holds. We assume that
$\lambda$ is a prime of $K$ such that
\begin{equation}H^0(\QQ,M_\lambda)\cong
H^1_\f(\QQ,M_\lambda)\cong
H^1_\f(\QQ,M_\lambda^\kd)\cong
H^0(\QQ,M_\lambda^\kd)\cong 0.\mlabel{vanishing}\end{equation}
This is conjectured to hold for all $\lambda$ under
our hypotheses on $M$ and it implies Conjecture \ref{conj:ord}.
If $M = A_f$ or $A_f(1)$  and $\lambda$ is a prime in
$S_f$ then (\ref{vanishing}) holds by Theorem \ref{thm:ordl}.

Fontaine and Perrin-Riou (\cite[II.4]{fon_pr}) define an $\CO_\lambda$-lattice
$\delta_{\f,\lambda}(M)$ in $K_\lambda\otimes_K\Delta_\f(M)$.
They assume $K = \QQ$, denote their lattice $\Delta_S(T)$ (where $S$
is a finite set of primes and $T$ is a Galois-stable lattice in
$M_\lambda$) and then prove it is independent of the choice
of $S$ and $T$.
One checks that the definition and independence argument carry over
to arbitrary $K$ by taking determinants relative to $\CO_\lambda$
and $K_\lambda$ instead of $\ZZ_\ell$ and $\QQ_\ell$.  The arguments
of \cite[II.5]{fon_pr} carry over as well, giving another description
of $\delta_{\f,\lambda}(M)$ for which we need more notation.
Choose a Galois stable lattice $\CM_\lambda \subset M_\lambda$
and a free rank one $\CO_\lambda$-module $\omega \subset
K_\lambda \otimes_K \det_{K}t_M$.  We let
$\theta(\CM_\lambda) = \det_{\CO_\lambda}\CM_\lambda^+$,
regarded as a lattice in $K_\lambda\otimes_K\det_K M_B^+$
via the comparison isomorphism $I_\lambda^B$.
We let $\CM_\lambda^\kd =
\hom_{\CO_\lambda}(\CM_\lambda,\CO_\lambda(1)) \subset M_\lambda^\kd$.
Under our hypotheses, the Tate-Shafarevich group $\sha(\CM_\lambda)$
can be identified with $H^1_\f(\QQ,M_\lambda/\CM_\lambda)$
(which we know is an $\CO_\lambda$-module of finite length, as
in the proof of Theorem \ref{thm:ordl}).  The same
holds for $\CM_\lambda^\kd$.  Furthermore, by the main result of
\cite{fl3} (also \cite[II.5.4.2]{fon_pr}), $\sha(\CM_\lambda)$ and
$\sha(\CM_\lambda^\kd)$ have the same length.  In fact,
there is an $\CO_\lambda$-linear isomorphism
\begin{equation}\sha(\CM_\lambda^\kd) \cong
 \hom_{\ZZ_\ell}(\sha(\CM_\lambda),\QQ_\ell/\ZZ_\ell).
\mlabel{casselstate}\end{equation}
Finally, the  Tamagawa ideal of $\CM_\lambda$
relative to $\omega$ is defined as
$$\Tam_\omega^0(\CM_\lambda) = \Tam_{\ell,\omega}^0(\CM_\lambda)
\cdot\Tam_{\infty}^0(\CM_\lambda)
\cdot\prod_{p\neq \ell}\Tam_p^0(\CM_\lambda),$$
where the factors are defined as in I.4.1 (and III.5.3.3)
of \cite{fon_pr}. Recall that $\Tam_p^0(\CM_\lambda)=1$ if
$\CM_\lambda$ is unramified at $p\neq \ell$ and that
$$\Tam_{\infty}^0(\CM_\lambda) =
\fitt_{\CO_\lambda}H^1(\RR,\CM_\lambda) = \CO_\lambda$$
if $\ell$ is odd. The argument of \cite[I.4.2.2]{fon_pr} shows that
if $p \neq \ell$, then
$$\Tam_p^0(\CM_\lambda) =
\fitt_{\CO_\lambda} H^1(I_p,\CM_\lambda)^{G_{\QQ_p}}_{\tor}$$
from which it is not hard to deduce that
\begin{equation}
\Tam_p^0(\CM_\lambda) = \Tam_p^0(\CM_\lambda^\kd).
\mlabel{locdual}\end{equation}

Viewing $\hom_{\CO_\lambda}(\theta(\CM_\lambda),\omega)$
as a lattice in $K_\lambda\otimes_K \Delta_\f(M)$, we have by
\cite[Th. II.5.3.6]{fon_pr}
\begin{equation}
\delta_{\f,\lambda}(M)\! =\! \frac{\fitt_{\CO_\lambda}\!\! H^0(\QQ,
M_\lambda/\CM_\lambda)\cdot
\fitt_{\CO_\lambda}\!\! H^0(\QQ,M_\lambda^\kd/\CM_\lambda^\kd)}{
\fitt_{\CO_\lambda}\sha(\CM_\lambda^\kd)\cdot
\Tam_\omega^0(\CM_\lambda)}
\hom_{\CO_\lambda}(\theta(\CM_\lambda),\omega).\mlabel{deltadef}
\end{equation}

The $\lambda$-part of the Bloch-Kato conjecture can then be
formulated as follows:

\begin{conjecture} Let $M$ in $\pms_K$ be critical,
$b(M)$ as in Conjecture \ref{conj:del}, $\lambda$ a place of $K$
such that (\ref{vanishing})
holds and $\delta_{\f,\lambda}(M)$ as in (\ref{deltadef}). Then
$$\delta_{\f,\lambda}(M) = (1\otimes b(M))\CO_\lambda.$$
\mlabel{conj:bk}\end{conjecture}

\begin{theorem} Let $f$ be a newform of weight
$k\geq 2$ and $S_f$ the set of places defined in
\S\ref{ss:conseq}.
Then Conjecture \ref{conj:bk} holds for both $M =
A_f$ and $M=A_f(1)$ if $\lambda\in S_f$.
\mlabel{thm:bk}\end{theorem}
\proofbegin We first prove a Lemma which  relates Conjecture
\ref{conj:bk} to
the value $L^\Sigma(M,0)$ for a suitable finite set of primes
$\Sigma$.

\begin{lemma} Suppose $M$, $b(M)$ and $\lambda$ are as in Conjecture
\ref{conj:bk} and $S$ is a set of places of $\QQ$ containing $\ell$,
$\infty$ and those where $M_\lambda$ is ramified. Assume
$\Sigma:=S\setminus \{\ell,\infty\}$ is nonempty and
$L_p(M,0)^{-1}\neq 0$ for all $p\in\Sigma$. Put
$b^\Sigma(M)=\prod_{p\in\Sigma}L_p(M,0)b(M)$. Then Conjecture
\ref{conj:bk} is equivalent  to
$$\frac{\fitt_{\CO_\lambda}H^0(\QQ,M_\lambda^\kd/\CM_\lambda^\kd)}
{\fitt_{\CO_\lambda}  H^1_{\Sigma}(\QQ,
M_\lambda^\kd/\CM_\lambda^\kd)
\Tam_{\ell,\omega}^0(\CM_\lambda)}
\hom_{\CO_\lambda}(\theta(\CM_\lambda),\omega)
=(1\otimes b^\Sigma(M))\CO_\lambda.
$$
where $H^1_{\Sigma}(\QQ, M_\lambda^\kd/\CM_\lambda^\kd)$ was defined
in \S \ref{ssec:gal-coh}.
\mlabel{LSigma}
\end{lemma}
\proofbegin  By \cite[Prop 1.4]{flachthesis} there is a long exact
sequence
\begin{multline} 0\rightarrow
H^1_\emptyset(\QQ,M_\lambda^\kd/\CM_\lambda^\kd) \rightarrow
H^1(G_S,M_\lambda^\kd/\CM_\lambda^\kd)
\rightarrow
\bigoplus_{p\in S}\frac{H^1(\QQ_p,M_\lambda^\kd/\CM_\lambda^\kd)}
{H^1_f(\QQ_p,M_\lambda^\kd/\CM_\lambda^\kd)}
\xrightarrow{\rho^*} \\ H^1_f(\QQ,\CM_\lambda)^*
\rightarrow H^2(G_S,M_\lambda^\kd/\CM_\lambda^\kd)\rightarrow
\bigoplus_{p\in S}H^2(\QQ_p,M_\lambda^\kd/\CM_\lambda^\kd)\rightarrow
H^0(\QQ,\CM_\lambda)^*\rightarrow 0
\notag
\end{multline} where $G_S$ is
the Galois group of the maximal extension of $\QQ$ unramified outside
$S$ and $^*$ denotes the Pontryagin dual.
{}From (\ref{vanishing}) we know that the groups
$H^1_\emptyset(\QQ,M_\lambda^\kd/\CM_\lambda^\kd)$ and
$H^1_f(\QQ,\CM_\lambda)$ are finite. The map
$\rho^*$ is Pontryagin dual to the restriction map
\[ H^1_f(\QQ,\CM_\lambda) \xrightarrow{\rho}
\bigoplus_{p\in S}H^1_f(\QQ_p,\CM_\lambda).\]
Clearly $\rho$ is injective as
$H^1_f(\QQ,\CM_\lambda)\cong H^0(\QQ,M_\lambda/\CM_\lambda)$ injects
into $H^0(\QQ_p,M_\lambda/\CM_\lambda)\cong
H^1_f(\QQ_p,\CM_\lambda)_{tors}$ for any $p$, for example some
$p\in\Sigma$. This argument also shows that
$\rho^*$ restricted to \[\bigoplus_{p\in\Sigma\cup\{\infty\}}
\frac{H^1(\QQ_p,M_\lambda^\kd/\CM_\lambda^\kd)}
{H^1_f(\QQ_p,M_\lambda^\kd/\CM_\lambda^\kd)}\cong
\bigoplus_{p\in\Sigma\cup\{\infty\}}H^1_f(\QQ_p,\CM_\lambda)^*\]
is still surjective since the dual map is still injective.
Hence we find an exact sequence
\begin{multline} 0\rightarrow
H^1_\emptyset(\QQ,M_\lambda^\kd/\CM_\lambda^\kd) \rightarrow
H^1_\Sigma(\QQ,M_\lambda^\kd/\CM_\lambda^\kd) \\
\rightarrow
\bigoplus_{p\in\Sigma\cup\{\infty\}}H^1_f(\QQ_p,\CM_\lambda)^*
\xrightarrow{\rho^*} H^0(\QQ,M_\lambda/\CM_\lambda)^*\rightarrow 0.
\mlabel{lex2}\end{multline}
By \cite[Proof of I.4.2.2]{fon_pr} we have for $p\in\Sigma$
\[\fitt_{\CO_\lambda}H^1_f(\QQ_p,\CM_\lambda)=L_p(M,0)^{-1}
\Tam_p^0(\CM_\lambda) \]
since $L_p(M,0)^{-1}\neq 0$. In particular all terms in
(\ref{lex2}) are finite
$\CO_\lambda$-modules, and the statement of Lemma \ref{LSigma}
follows easily by taking Fitting ideals in (\ref{lex2}), together
with (\ref{deltadef}).
\epf

For a prime
$\lambda\in S_f$, recall that
$\cala_{f,\lambda}=\ad^0_{\CO_\lambda} \CM_{f,\lambda}$.
Put $\calb_{f,\lambda} =
\cala_{f,\lambda}(1)$.
Identifying $\det_K t_{A_f}$ with
$\hom_K(\fil^{k-1}M_{f,\dr}, M_{f,\dr}/\fil^{k-1}M_{f,\dr})$,
we let
\begin{align*}
\omega_A = & \CO_\lambda \otimes_\CO
\hom_\CO(\fil^{k-1}\CM_{f,\dr},
\CM_{f,\dr}/\fil^{k-1}\CM_{f,\dr})\\
\cong &\hom_{\CO_\lambda}
(\fil^{k-1}\CM_{f,\dr}\otimes_\CO \CO_\lambda,
\CM_{f,\dr}\otimes_\CO \CO_\lambda/\fil^{k-1}\CM_{f,\dr}
\otimes_\CO \CO_\lambda).
\end{align*}
Similarly, identifying  $\det_K t_{B_f}$
with $\det_K t_{A_f} \otimes \QQ(2)_\dr$
we define $\omega_B$ as $\omega_A \otimes \iota^{-2}$.

\begin{proposition} We have
$\Tam_{\ell,\omega_A}^0(\cala_{f,\lambda})=
\Tam_{\ell,\omega_B}^0(\calb_{f,\lambda})=
\CO_\lambda$ for $\lambda\in S_f$.
\mlabel{tamatl}\end{proposition}
\proofbegin
With the notation in \S \ref{ss:back}, we further denote by $\MFcat$ the
category of filtered $\phi$-modules as defined in \cite[1.2.1]{fon_ast}.
We denote by $\mfcato$ (resp. $\mfcatk$)
the category of $\CO_\lambda$- (resp. $K_\lambda$-) modules in the abelian category
$\mfcat$ (resp. additive category $\MFcat$) and by $\mfcato^a$ (resp. $\mfcatk^a$) the full subcategories
with filtration restrictions as in \S\ref{ss:back}. Scalar extension $-\otimes_{\ZZ_\ell}\QQ_\ell$
induces an exact functor $\mfcato\rightarrow\mfcatk$ where the notion of exactness in $\MFcat$ is
defined in \cite[1.2.3]{fon_ast}.

Now assume that
$\cald_1,\cald_2$ are torsion free objects of $\mfcato^a$ for some
$a$ and put
$D_i=\cald_i\otimes_{\ZZ_\ell}\QQ_\ell$.
Set $\cald=\Hom_{\CO_\lambda}(\cald_1,\cald_2)$ and $D=\Hom_{K_\lambda}(D_1,D_2)$
  which are objects of $\mfcat$ and $\MFcat$ respectively. We also have
$D\cong\cal D\otimes_{\ZZ_\ell}\QQ_\ell$. An elementary computation shows that the
first two rows in the following commutative diagram are exact
\minCDarrowwidth1em
\begin{equation}
{\scriptsize
\begin{CD}
0 @>>> \Hom_{\mfcato}(\cald_1,\cald_2) @>>> \fil^0 \cald
@>{1-\phi^0}>>\cald @>\pi>>
\Ext^1_{\mfcato}(\cald_1,\cald_2) @>>>0\\
@. @V\epsilon^0 VV @VVV @VVV @V\epsilon^1 VV @.\\
0 @>>> \Hom_{\mfcatk}(D_1,D_2) @>>> \fil^0 D @>{1-\phi}>> D @>\pi>>
\Ext^1_{\mfcatk}(D_1,D_2) @>>>0\\
@. @V\theta^0 VV @V\iota VV @V\iota VV @V\theta^1 VV @.\\
0 @>>> H^0(\QQ_\ell,V) @>>> \fil^0 D(V) @>{1-\phi}>>
D(V)  @>e>> H^1_f(\QQ_\ell,V) @>>> 0 \\
@. \Vert @. @VVV @V(\text{id},0)VV \Vert @.\\
0 @>>> H^0(\QQ_\ell,V) @>>> D(V) @>>> D(V)
\oplus t_V @>>> H^1_f(\QQ_\ell,V) @>>> 0
\end{CD}
}
\mlabel{big}\end{equation}
where $\pi$ is defined as follows. For $\eta\in\cald$ define an
extension $E_\eta$ of $\cald_1$ by $\cald_2$ in $\mfcato$ with
underlying $\CO_\lambda$-module $\cald_2\oplus\cald_1$,
filtration
\[ \fil^iE_\eta:=\fil^i \cald_2\oplus\fil^i\cald_1\]
and Frobenius maps $\phi^i:\fil^iE_\eta\rightarrow E_\eta$
\begin{equation}
\phi^i(x,y)=(\phi^i(x)+\eta\phi^i(y),\phi^i(y)).
\mlabel{phiaction}\end{equation}
The same definitions for $\eta\in D$ lead to an extension in
$\mfcatk$. Then $\pi(\eta)$ is the  class of the Yoneda extension
$E_\eta$ in $\Ext^1$ (we shall identify $\Ext^1$ with the group of
Yoneda extensions throughout).

To explain the remaining part of diagram (\ref{big}) we first recall the notion of admissibility from \cite[3.6.4]{fon_ast}.
A filtered $\phi$-module $D'$ in $\MFcat$ is called admissible if
the natural map $B_{cris}\otimes_{\QQ_\ell}D'\cong B_{cris}\otimes_{\QQ_\ell}V(D')$ is an isomorphism where $V(D')$ is the $G_\ell$-representation
\[V(D')=\fil^0(D'\otimes B_{cris})^{\phi\otimes\phi=1}.\]
The functor $D'\rightarrow V(D')$ is fully faithful and exact on the category of admissible filtered $\phi$-modules, and
induces an equivalence of this category with the category $\rep(G_\ell)$ of crystalline $K_\lambda[G_l]$-representations
(see \cite[3.6.5]{fon_ast}). If $D'=\cald'\otimes_{\ZZ_l}\QQ_l$ for some object $\cald'$ of $\mfcat^0$ then $D'$ is admissible by
\cite[Th. 8.4]{fon_laf}, and for such $D'$ we have a natural isomorphism
$V(D')\cong \VV(\cald')\otimes_{\ZZ_\ell}\QQ_\ell$ by (\ref{eq:VV}). If $D'=\cald'\otimes_{\ZZ_l}\QQ_l$ for some object $\cald'$ of $\mfcat^a$ then
we can extend the definition of $\VV$ by $\VV(\cald')=\VV(\cald'[-a])(a)$ (Tate twist) and we deduce again that $D'$ is admissible.
In particular $D_1$ and $D_2$ are admissible, and then $D$ is admissible
by \cite[Prop. 3.4.3]{fon_ast}. Putting $V:=V(D)$ and $V_i:=V(D_i)$ we have an isomorphism of
$G_\ell$-representations $V=\Hom_{K_\lambda}(V_1,V_2)$ by \cite[3.6]{fon_ast}.

Coming back to diagram (\ref{big}), the map $\iota$ is just the natural isomorphism
\[D\xrightarrow{1\otimes-} B_{cris}\otimes_{\QQ_\ell}D \cong
B_{cris}\otimes_{\QQ_\ell}V
\leftarrow H^0(\QQ_\ell,B_{cris}\otimes_{\QQ_\ell}V)=:D(V)\]
and $e$ is the boundary map in Galois cohomology induced from the
short exact sequence of
$G_{\QQ_\ell}$-modules
\begin{equation} 0\rightarrow V(D)\rightarrow
\fil^0(B_{cris}\otimes_{\QQ_\ell}D)\xrightarrow{1-\phi\otimes
\phi}B_{cris}\otimes_{\QQ_\ell}D\rightarrow 0
\mlabel{bk_seq}\end{equation}
as in the proof of \cite[Lemma4.5(b)]{bloch_kato}.  It is clear
that $\mfcatk^a$ is closed under extensions inside $\mfcatk$
hence we obtain a chain of isomorphisms
\begin{multline}\theta^i: \Ext^i_{\mfcatk}(D_1,D_2)\leftarrow
\Ext^i_{\mfcatk^a}(D_1,D_2)\xrightarrow{v^i}\notag \\
\Ext^i_{\rep(G_\ell)}(V_1,V_2)\xrightarrow{\Delta^i}
\Ext^i_{\rep(G_\ell)}(K_\lambda,\Hom_{K_\lambda}(V_1,V_2))
\rightarrow H^i_f(\QQ_\ell,V)
\end{multline}
for $i=0,1$. Here $\Delta^1$
sends a Yoneda extension
\[0\rightarrow V_2\rightarrow V_3\rightarrow
V_1\rightarrow 0\]
to the pull back to $K_\lambda\cdot
1_{V_1}\subseteq \Hom_{K_\lambda}(V_1,V_1)$ of the induced
extension
\[0\rightarrow \Hom_{K_\lambda}(V_1,V_2)\rightarrow
\Hom_{K_\lambda}(V_1,V_3)\rightarrow\Hom_{K_\lambda}(V_1,V_1)\rightarrow
0.\]
The maps $v^i$ (defined by applying $V$ to a Yoneda extension) are
isomorphisms because $V$ is fully faithful and exact.

The three lower rows in
(\ref{big}) with the indicated maps form a
commutative diagram, and all
these rows are exact (see
\cite[Lemma 4.5(b)]{bloch_kato} for the two lower
rows). We shall verify the
identity $\theta^1\pi=e\iota$, all the others being
straightforward. Consider the commutative diagram
\begin{equation}
{\scriptsize
\begin{CD} 0 @>>> V_2 @>>>
\fil^0(B_{cris}\otimes_{\QQ_\ell}D_2)@>{1-\phi\otimes
\phi}>>B_{cris}\otimes_{\QQ_\ell}D_2 @>>> 0\\
@.\Vert @. @AAA @AA{1\otimes\psi}A\\
0 @>>> V_2 @>>>
\fil^0(B_{cris}\otimes_{\QQ_\ell}(D_2\oplus
D_1))^{\phi=1}@>>>\fil^0(B_{cris}\otimes_{\QQ_\ell}D_1)^{\phi=1}
@>>> 0\\
\end{CD}
}
\mlabel{big2}\end{equation}
where all unnamed arrows are natural projection or
inclusion maps, the top row is (\ref{bk_seq}) with $D$
replaced by $D_2$, and the action of
$\phi$ on $D_2\oplus D_1$ is given by (\ref{phiaction}).
For $\psi\in D$, the extension $e\iota(\psi)$ is the
pullback of (\ref{bk_seq}) under $K_\lambda(1\otimes\psi)
\subset B_{cris}\otimes_{\QQ_\ell}D$. To compute
$(\Delta^1)^{-1}e\iota(\psi)$ apply the exact functor
$\Hom_{K_\lambda}(V_1,-)$ to diagram (\ref{big2}). Via
the isomorphisms
\begin{align*}
\Hom_{K_\lambda}(V_1,B_{cris}\otimes_{\QQ_\ell}D_2)
\cong &
\Hom_{B_{cris}\otimes
K_\lambda}(B_{cris}\otimes_{\QQ_\ell}V_1,B_{cris}\otimes_{\QQ_\ell}D_2)\\
\cong & \Hom_{B_{cris}\otimes
K_\lambda}(B_{cris}\otimes_{\QQ_\ell}D_1,B_{cris}\otimes_{\QQ_\ell}D_2)\\
\cong & B_{cris}\otimes_{\QQ_\ell}\Hom_{K_\lambda}(D_1,D_2),
\\
\Hom_{K_\lambda}(V_1,\fil^0(B_{cris}\otimes_{\QQ_\ell}D_2))\cong
&
\fil^0(B_{cris}\otimes_{\QQ_\ell}\Hom_{K_\lambda}(D_1,D_2))
\end{align*}
the first row becomes isomorphic to (\ref{bk_seq}) and
the image of
\[1_{V_1}\in\Hom_{K_\lambda}(V_1,V_1)=
\Hom_{K_\lambda}(V_1,\fil^0(B_{cris}\otimes_{\QQ_\ell}D_1)^{\phi=1})
\]
in $B_{cris}\otimes_{\QQ_\ell}D$
is $1\otimes\psi$. Hence $(\Delta^1)^{-1}e\iota(\psi)$
is represented by the lower row in (\ref{big2}). But
from the definition of $\pi$ it is immediate that the
lower row in (\ref{big2}) is the image of $\pi(\psi)$
under the functor $V$. This gives the identity
$\theta^1\pi=e\iota$.

Put $T_i=\VV(\cald_i)$ for $i=1,2$
and $T=\Hom_{\CO_\lambda}(T_1,T_2)$ with its natural
$G_{\ell}$-action. Then, since
$T_1$ is torsion free, $H^1(\QQ_\ell,T)$ naturally identifies with the
set of equivalence classes of extensions of
$\CO_\lambda[G_\ell]$-modules
\begin{equation}
0\rightarrow T_2\rightarrow T_3\rightarrow T_1\rightarrow 0.
\mlabel{text}\end{equation}
Since the functor $\VV$ is exact on $\mfcato^a$, and since $\mfcato^a$ is
closed under extensions in $\mfcato$, we obtain maps
$$\Theta^i: \Ext^i_{\mfcato}(\cald_1,\cald_2)\xleftarrow{\sim}
\Ext^i_{\mfcato^a}(\cald_1,\cald_2)\rightarrow
\Ext^i_{\CO_\lambda[G_\ell]}(T_1,T_2)\cong H^i(\QQ_\ell,T)
$$
analogous to the maps $\theta^i$.
The faithfulness of $\VV$ implies that $\Theta^0$ is injective and fullness
of $\VV$ implies that $\Theta^0$ is surjective and that $\Theta^1$ is injective.
The image of $\Theta^1$ lies in the subgroup
\[
H^1_f(\QQ_\ell,T) :=\{[T_3]\in H^1(\QQ_\ell,T)\vert
[V_3]:=[T_3\otimes_{\ZZ_\ell}\QQ_\ell]\in H^1_f(\QQ_\ell,V)\}\]
since $\VV(\cald_3)\otimes_{\ZZ_\ell}\QQ_\ell\cong V(\cald_3\otimes_{\ZZ_\ell}\QQ_\ell)$
is a crystalline representation. Conversely, if (\ref{text}) lies in $H^1_f(\QQ_\ell,T)$,
the $G_\ell$-module $T_3$ is a submodule of a crystalline representation $V_3$
so that $D(V_3)$ lies in $\mfcatk^a$ and hence $T_3$ lies in the essential image of the
Fontaine Laffaille functor $\VV$, $T_3=\VV(\cald_3)$, say. Since $\VV$ is full the
extension (\ref{text}) is the image of a sequence
$0\rightarrow \cald_2\rightarrow \cald_3\rightarrow \cald_1\rightarrow 0$
in $\mfcato^a$ and since $\VV$ is fully faithful and exact this sequence is exact, hence
represents an element of $\Ext^1_{\mfcato}(\cald_1,\cald_2)$. We conclude that
\begin{equation}
\Theta^1:\Ext^1_{\mfcato}(\cald_1,\cald_2)\cong H^1_f(\QQ_\ell,T)
\mlabel{h1f}\end{equation}
is an isomorphism. It is clear that $\theta^i\epsilon^i=\tilde{\epsilon}^i\Theta^i$ where
$\tilde{\epsilon}^i:H^i(\QQ_\ell,T)\rightarrow H^i(\QQ_\ell,V)$ are the natural maps.
The last row in (\ref{big}) induces an isomorphism
\begin{align*}\det_{K_\lambda}H^0(\QQ_\ell,V)\otimes_{K_\lambda}
\det_{K_\lambda}^{-1}H^1_f(\QQ_\ell,V)&\cong \det_{K_\lambda}D
\otimes_{K_\lambda}\det_{K_\lambda}^{-1}D\otimes_{K_\lambda}
\det_{K_\lambda}^{-1}t_V\\&\cong\det_{K_\lambda}^{-1}t_V\end{align*}
and the Tamagawa ideal is defined in \cite[I.4.1.1]{fon_pr} so that
\begin{equation}
\det_{\CO_\lambda}H^0(\QQ_\ell,T)\otimes_{\CO_\lambda}
\det_{\CO_\lambda}^{-1}H^1_f(\QQ_\ell,T)\cong
\Tam_{\ell,\omega}^0(T)\omega^{-1}.
\mlabel{tamno}\end{equation}
Using the fact that $\Theta^0$ is an isomorphism together with (\ref{h1f}) and (\ref{big}) one computes that
the left hand side in (\ref{tamno}) equals
$\det_{\CO_\lambda}^{-1}\cald/\fil^0\cald$ so that
$\Tam_{\ell,\omega}^0(T)=\CO_\lambda$ if $\omega$ is a basis
of $\det_{\CO_\lambda}\cald/\fil^0\cald$.

These arguments apply to $\cald_1=\cald_2
=\CM_{f,\lcrys}$ which is an object of
$\mfcato^0$ if
$\lambda\in S_f$, more specifically if $\ell\nmid N$
and $\ell > k$. We have $T_1=T_2=\CM_{f,\lambda}$
and $T=\cala_{f,\lambda}\oplus\CO_\lambda$. Our choice
of $\omega_A$ then ensures that
$\Tam_{\ell,\omega_A}^0(\cala_{f,\lambda})=\CO_\lambda$.
For $B_f=A_f(1)$ we can use the same argument as
long as both
$\cald_1=\CM_{f,\lcrys}$ and
$\cald_2=\CM_{f,\lcrys}[1]$ are objects
of $\mfcato^1$. This is the case if $\ell>k+1$ or if
$l=k+1$ and $\CM_{f,\lcrys}$ has no
nonzero quotient $A$ in $\mfcato$ with
$\fil^{k-1}A=A$.

\begin{lemma} If $a_\ell\equiv 0\mod\lambda$ then
$\CM_{f,\lcrys}$ has no
nonzero quotient $A$ in $\mfcato$ with
$\fil^{k-1}A=A$.
\end{lemma}
\proofbegin By \cite{scholl}
we know that the characteristic polynomial of
$\phi$ on
$\CM:=\CM_{f,\lcrys}$ is
$X^2-a_\ell X+\psi(\ell)\ell^{k-1}$ hence $\bar{\phi}$ has
characteristic polynomial $X^2$ on
$\bar{\CM}:=\CM\otimes_\CO
(\CO_\lambda/\lambda)$. Since
\begin{equation}
\bar{\CM}=\bar{\phi}(\bar{\CM})+\bar{\phi}^{k-1}(\fil^{k-1}\bar{\CM})
\mlabel{mfcond}\end{equation}
and
$\dim_{\CO_\lambda/\lambda}\fil^{k-1}\bar{\CM}=1$
the map $\bar{\phi}$ is nonzero, hence conjugate to
$\left(\begin{smallmatrix} 0 & 1\\ 0 & 0
\end{smallmatrix}\right)$. Since
$\bar{\phi}=l^{k-1}\bar{\phi}^{k-1}=0$ on
$\fil^{k-1}\bar{\CM}$ and because of
(\ref{mfcond}) we have
$\fil^{k-1}\bar{\CM}=\text{ker}(\bar{\phi})=
\bar{\phi}(\bar{\CM})$. It is now easy to see
that $\bar{\CM}$ is a simple object in $\mfcato$:
Any proper subobject
$N\subset\bar{\CM}$ is $\bar{\phi}$-stable,
hence contained
in $\text{ker}(\bar{\phi})=\fil^{k-1}\bar{\CM}$ and
we have $\fil^{k-1}N=N$. But again by
(\ref{mfcond}) we find
$\bar{\phi}^{k-1}(\fil^{k-1}\bar{\CM})\not\subseteq
\text{ker}(\bar{\phi})=\bar{\phi}(\bar{\CM})$
so that $N=\fil^{k-1}N=0$.  If $A$ is a
nonzero quotient of $\CM$ then
$\bar{A}$ is a nonzero quotient of $\bar{\CM}$
hence equal to $\bar{\CM}$ and we find
$\fil^{k-1}\bar{A}\neq\bar{A}$ and
$\fil^{k-1}A\neq A$.
\epf

It remains to prove Proposition \ref{tamatl} for $\calb_{f,\lambda}$
in the ordinary case $a_\ell \not\equiv
0\mod \lambda$ (and
$l=k+1$). We use the fact that
$B_f\cong A_f^*(1)$ and appeal to the following
conjecture, a slight generalisation (from $\ZZ_\ell$
to $\CO_\lambda$) of conjecture $C_{EP}(V)$ of
\cite{periou94} (we also use a similar generalisation of
\cite[Prop. C.2.6]{periou94}).

Let $V$ be a crystalline
representation of $G_{\ell}$ over $K_\lambda$ and
$T\subseteq V$ a $G_{\ell}$-stable
$\CO_\lambda$-lattice. Let
$\omega$ (resp. $\omega^*$) be a lattice of
$$\det_{K_\lambda}D(V)/\fil^0D(V) {\rm\ (resp.\ }
\det_{K_\lambda}D(V^*(1))/\fil^0D(V^*(1)))$$
so that we
obtain a lattice
$\omega\otimes\omega^{*,-1}$ of $\det_{K_\lambda}D(V)$
via the exact sequence
\begin{equation}0\rightarrow
(D(V^*(1))/\fil^0D(V^*(1)))^*\rightarrow D(V)\rightarrow
D(V)/\fil^0D(V)\rightarrow 0.\mlabel{drseq}\end{equation}
Let $\eta(T,\omega,\omega^*)\in
B_{cris}\otimes_{\QQ_\ell}K_\lambda$ be such that
$t^{-t_H(V)}\otimes \det_{\CO_\lambda}T=\eta(T,\omega,\omega^*)\omega\otimes\omega^{*,-1}$
under the comparison isomorphism
$B_{cris}\otimes_{\QQ_\ell}\det_{K_\lambda}V\cong
B_{cris}\otimes_{\QQ_\ell}\det_{K_\lambda}D(V)$. One shows
that
$\eta(T,\omega,\omega^*)\in\QQ_\ell^{ur}\otimes_{\QQ_\ell}K_\lambda$
\cite[Lemme C.2.8]{periou94} and that in fact
$\eta(T,\omega,\omega^*)\in 1\otimes K_\lambda$ up to
an element in
$(\ZZ_\ell^{ur}\otimes_{\ZZ_\ell}\CO_\lambda)^\times$.

\begin{conjecture}
For $j\in\ZZ$, put
$h_j(V)=\dim_{K_\lambda}\fil^jD(V)/\fil^{j+1}D(V)$,
and put $\Gamma^*(j)=(j-1)!$ if $j>0$ and
$\Gamma^*(j)=(-1)^j((-j)!)^{-1}$ if $j\leq 0$.
Then
\[\CO_\lambda\frac{\Tam_{\ell,\omega}^0(T)}
{\Tam_{\ell,\omega^*}^0(T^*(1))}=\CO_\lambda\prod_j\Gamma^*(-j)^{-h_j(V)}
\eta(T,\omega,\omega^*)^{-1}\]
\mlabel{ep}\end{conjecture}

\noindent
{\em Remark.} One can show that upon taking the norm
from $K_\lambda$ to $\QQ_\ell$ all quantities in this
formula transform into the corresponding quantities
obtained by viewing $V$ as a representation over
$\QQ_\ell$ rather than $K_\lambda$. Since the norm map
$K_\lambda^\times/\CO_\lambda^\times \rightarrow
\QQ_\ell^\times/\ZZ_\ell^\times$ is injective it suffices to
prove the conjecture for $K_\lambda=\QQ_\ell$.
\bigskip

We make Conjecture \ref{ep} more explicit for
$V=A_{f,\lambda}$. In this case we have
$h_j(V)=1$ for $i=1-k,0,k-1$ and $h_j(V)=0$ otherwise so
that $\prod_j\Gamma^*(-j)^{-h_j(V)}=(-1)^{k-1}(k-1)$. For $\lambda\in
S_f$ Lemma~\ref{lem:adjoint} shows that the isomorphism
$B_{cris}\otimes_{\QQ_\ell}\det_{K_\lambda}V\cong
B_{cris}\otimes_{\QQ_\ell}\det_{K_\lambda}D(V)$ is induced
by the functor $\VV$  for the unit object in $\ipms^S_K$, hence sends
$\det_{\CO_\lambda}\CA_{f,\lambda}$ to
$\det_{\CO_\lambda}\CA_{f,\dr}\otimes_\CO \CO_\lambda$.
The computation of $t_{B_f}$ before (\ref{tanb}) works
with $M_{f,\dr}$ replaced by $\CM_{f,\dr}\otimes_\CO \CO_\lambda$
and the pairing (\ref{eq:dual}) on $A_{f}$
gives a perfect pairing
$(\CA_{f,\dr}\otimes_\CO \CO_\lambda) \otimes
(\CA_{f,\dr}\otimes_\CO \CO_\lambda) \rightarrow
\CO_\lambda$. Hence we find that
$\omega_A\otimes\omega_B^{-1}$ is a basis of
$\det_{\CO_\lambda}\CA_{f,\dr}\otimes_\CO \CO_\lambda$ via the exact
sequence (\ref{drseq}). We conclude that
$\eta(\CA_{f,\lambda},\omega_A,\omega_B)=1$ and that
Conjecture \ref{ep} reduces to the assertion
\[\CO_\lambda\Tam_{\ell,\omega_A}^0(\CA_{f,\lambda})=
\CO_\lambda\Tam_{\ell,\omega_B}^0(\CA_{f,\lambda}(1)).\]
Moreover, we know from the first part of the proof
that the left hand side equals $\CO_\lambda$ if
$\lambda\in S_f$. Now Conjecture \ref{ep} is shown in
\cite{periouinv}
for $K_\lambda=\QQ_\ell$ and $V$ an
ordinary representation of $G_{\QQ_\ell}$ (combine
Proposition 4.2.5, Theorem 3.5.4 of loc. cit.) under the
assumption of another conjecture Rec(V) which has
meanwhile been proved in \cite{colmez}. Ordinarity of
$A_{f,\lambda}$ is implied by ordinarity of
$M_{f,\lambda}$ which in turn is implied by $a_\ell
\not\equiv 0\mod\lambda$.
This finishes
the proof of Proposition~\ref{tamatl}.
\epf

We shall first prove Theorem \ref{thm:bk} for $M=B_f=A_f(1)$ in which
case Lemma \ref{LSigma} applies. Fix a prime $\lambda\in S_f$ and let $\Sigma$ be the set of
primes dividing $N$ if $N>1$ or put $\Sigma=\{p\}$ for some prime $\lambda\nmid p$ if $N=1$. The isomorphism $\gamma:M_f\rightarrow M_f^\Sigma$
of Proposition \ref{prop:level} satisfies
\[\gamma^t=\gamma^{-1}\phi\prod_{p\in\Sigma}L_p^{\naive}(B_f,0)^{-1}\]
by Proposition \ref{pp:level2} where $\phi=\prod_{\delta_p=1}(-a_p)\prod_{\delta_p=2}\psi(p)p^{k-1}\in\CO_\lambda^\times$ (it is well
known that $a_p^2=\psi(p)p^{k-1}$ or $a_p^2=\psi(p)p^{k-2}$ if $\delta_p=1$ \cite[4.6.17]{miyake}). Moreover, $\gamma$
induces an isomorphism \[B_f=\Hom_K(M_f,M_f(1))\rightarrow
B_f^\Sigma:=\Hom_K(M_f^\Sigma,M_f^\Sigma(1))\] and an isomorphism
$\gamma:\Delta_\f(B_f)\rightarrow\Delta_\f(B_f^\Sigma)$ so that
\[\gamma(b)(x\otimes (2\pi i)^2)\otimes\iota^2=\gamma_{\dr}b(\gamma_{\dr}^{-1}(x)\otimes (2\pi i)^2)\otimes\iota^2\]
for $b\in\Delta_\f(B_f)$ and $x\in M_{f,\dr}^\Sigma$. For such $b$ and $x$ we have
\begin{align}\langle x,
\gamma(b)(x\otimes (2\pi i)^2)\otimes\iota^2
\rangle^\Sigma = &\langle \gamma_{\dr}^{-1}(x),\gamma_{\dr}^t
\gamma(b)(x\otimes (2\pi i)^2)\otimes\iota^2
\rangle\notag\\
= &\langle \gamma_{\dr}^{-1}(x),\gamma_{\dr}^{-1}
\gamma(b)(x\otimes (2\pi i)^2)\otimes\iota^2
\rangle\phi\prod_{p\in\Sigma}L_p^{\naive}(B_f,0)^{-1}\notag\\
= &\langle \gamma_{\dr}^{-1}(x),b(\gamma_{\dr}^{-1}(x)\otimes (2\pi i)^2)\otimes\iota^2
\rangle\phi\prod_{p\in\Sigma}L_p^{\naive}(B_f,0)^{-1}.
\label{Bfsigma}\end{align}
Recall that $\fil^{k-1}\CM_{f,\dr}^\Sigma=\CO\cdot f^\Sigma=\CO\cdot\gamma(f)$ by Propositions
\ref{prop:int-struc} and \ref{prop:level} where $\CO=\bigcap_{\lambda\in S_f}K\cap\CO_\lambda$. Note that if $b'$ is an $\CO_\lambda$-basis for
\begin{align*}
&\hom_{\CO_\lambda}(\theta(\calb_{f,\lambda}),\omega_B) \\
&\cong
\hom_{\CO_\lambda}(\fil^{k-1}\CM_{f,\dr}^\Sigma\otimes_\CO
\theta(\calb_{f,\lambda}),(\CM_{f,\dr}^\Sigma/\fil^{k-1}\CM_{f,\dr}^\Sigma)
\otimes_\CO\CO_\lambda\otimes\iota^{-2})
\end{align*}
then
$$\CO_\lambda\cdot b'(f^\Sigma\otimes(2\pi
i)^2)\otimes\iota^2 =
(\CM_{f,\dr}^\Sigma/\fil^{k-1}\CM_{f,\dr}^\Sigma)
\otimes_\CO\CO_\lambda$$
and hence
$$\CO_\lambda\cdot\langle
f^\Sigma,b'(f^\Sigma\otimes(2\pi
i)^2)\otimes\iota^2\rangle^\Sigma
= \CO_\lambda\eta^\Sigma_f\CM_{\psi}(1-k)_{\dr}
$$
where $\eta^\Sigma_f$ was defined before Proposition \ref{prop:int-struc}.
On the other hand by (\ref{Bfsigma}),
Theorem \ref{thm:del} and the remark
after the proof of Theorem \ref{thm:del} we have for $\lambda\in S_f$
\begin{align*}
&\CO_\lambda\cdot\langle f^\Sigma,
\gamma b^\Sigma(B_f)(f^\Sigma\otimes (2\pi i)^2)\otimes\iota^2
\rangle^\Sigma \\
&=\CO_\lambda\cdot\langle f,
b^\Sigma(B_f)(f\otimes (2\pi i)^2)\otimes\iota^2
\rangle \prod_{p\in\Sigma} L_p^{\naive}(B_f,0)^{-1} \\
&= \CO_\lambda\cdot\langle f,
b(B_f)(f\otimes (2\pi i)^2)\otimes\iota^2
\rangle \prod_{p\in\Sigma_e(f)} L_p(B_f,0) \\
&= \CO_\lambda(b_\dr \otimes \iota^{k-1}) =
\CO_\lambda\CM_{\psi}(1-k)_\dr.
\end{align*}
Applying Lemma \ref{LSigma} we must therefore show
$$\frac{\fitt_{\CO_\lambda}H^0(\QQ,B_{f,\lambda}^\kd/\calb_{f,\lambda}^\kd)}
{\fitt_{\CO_\lambda}
H^1_{\Sigma}(\QQ, B_{f,\lambda}^\kd/\calb_{f,\lambda}^\kd)
\Tam_{\ell,\omega_B}^0(\calb_{f,\lambda})}\,
\eta_f^\Sigma=\CO_\lambda.
$$
Using Proposition \ref{tamatl}, the fact that
$\cala_f=\calb_f^\kd$ and the vanishing of
$H^0(\QQ, A_{f,\lambda}/\cala_{f,\lambda})$ for $\lambda\in S_f$
this reduces to the identity
$$
 \fitt_{\CO_\lambda} H^1_\Sigma(\QQ,
A_{f,\lambda}/\cala_{f,\lambda}) = \CO_\lambda\eta^\Sigma_f
$$
which is Theorem \ref{thm:selmer}.

By (\ref{casselstate}), (\ref{locdual}) and Lemma \ref{tamatl} the
factor in front of $\hom(\theta(\CM_\lambda),\omega)$ in
(\ref{deltadef}) is the same  for $M=A_f$ and $M=B_f$. The
isomorphism $\tw$ defined in (\ref{twist}) maps
$\hom_{\CO_\lambda}(\theta(\cala_{f,\lambda}),\omega_A)$ to
$\hom_{\CO_\lambda}(\theta(\calb_{f,\lambda}),\omega_B)$, hence
$\delta_{\f,\lambda}(A_f)$ to
$\delta_{\f,\lambda}(B_f)$. Theorem \ref{thm:bk}
for $M=A_f$ therefore follows from Theorem
\ref{thm:bk} for $M=B_f$, together with Theorem \ref{thm:del} and
the fact that $(1-k)\epsilon(A_f)$ is a unit in $\CO_\lambda$.
\proofend


\end{document}